\documentclass{amsart}[10pt]
\usepackage{circuitikz}
\usepackage{graphicx}
\usepackage{tikz}
\usepackage{todonotes}
\usepackage{scalerel}
\usepackage{amssymb}
\usepackage{amsfonts}
\usepackage{color}
\usepackage{mathtools}
\usepackage{tikz} 
\usepackage{ytableau}
\usepackage{pgfplots}
\usepackage{hyperref}
\usepackage{xcolor}
\usepackage{tikz-cd}
\usepackage[utf8]{inputenc}
\usepackage[T1]{fontenc}
\usepackage{mathtools}
\usepackage{comment}

\newcommand{\Ytodo}[1]{\todo[color=teal!40,bordercolor=teal,size=\footnotesize]{\textbf{Todo: }#1}}

\theoremstyle{plain}

\theoremstyle{definition}
\newtheorem{theorem}{Theorem}[section]
\newtheorem{corollary}[theorem]{Corollary}

\newtheorem{lemma}[theorem]{Lemma}

\newtheorem{proposition}[theorem]{Proposition}
\newtheorem{example}[theorem]{Example}

\newtheorem{conjecture}[theorem]{Conjecture}

\newtheorem{remark}[theorem]{Remark}

\newtheorem{definition}[theorem]{Definition}

\newcommand{\Hilb}{\mathrm{Hilb}}

\newcommand{\monomialIdealPicture}[1]{
    \begin{tikzpicture}[x=(220:1cm), y=(-40:1cm), z=(90:0.707cm),scale=0.35]
        \foreach \m [count=\y] in {#1}{
          \foreach \n [count=\x] in \m {
          \ifnum \n>0
              \foreach \z in {1,...,\n}{
                \draw [fill=orange!30] (\x+1,\y,\z) -- (\x+1,\y+1,\z) -- (\x+1, \y+1, \z-1) -- (\x+1, \y, \z-1) -- cycle;
                \draw [fill=orange!40] (\x,\y+1,\z) -- (\x+1,\y+1,\z) -- (\x+1, \y+1, \z-1) -- (\x, \y+1, \z-1) -- cycle;
                \draw [fill=orange!10] (\x,\y,\z)   -- (\x+1,\y,\z)   -- (\x+1, \y+1, \z)   -- (\x, \y+1, \z) -- cycle;  
              }
             \fi
          }
        }
    \end{tikzpicture}
}
\newcommand{\monomialIdealPictureWrongTangent}[1]{
    \begin{tikzpicture}[x=(220:1cm), y=(-40:1cm), z=(90:0.707cm),scale=0.5]
        \foreach \m [count=\y] in {#1}{
          \foreach \n [count=\x] in \m {
          \ifnum \n>0
              \foreach \z in {1,...,\n}{
                \draw [fill=gray!30] (\x+1,\y,\z) -- (\x+1,\y+1,\z) -- (\x+1, \y+1, \z-1) -- (\x+1, \y, \z-1) -- cycle;
                \draw [fill=gray!40] (\x,\y+1,\z) -- (\x+1,\y+1,\z) -- (\x+1, \y+1, \z-1) -- (\x, \y+1, \z-1) -- cycle;
                \draw [fill=gray!10] (\x,\y,\z)   -- (\x+1,\y,\z)   -- (\x+1, \y+1, \z)   -- (\x, \y+1, \z) -- cycle;  
              }
             \fi
          }
        }
    \end{tikzpicture}
}

\definecolor{staircasegray}{RGB}{128,127,127}
\definecolor{staircasedarkgray}{RGB}{51,49,49}

\title{Singularities of ideals with minimal tangent space dimension in the Hilbert scheme of points on a threefold}
\author{Yuze Luan}
\begin{document}

\maketitle

\section{introduction}
The geometry of the Hilbert scheme of points on $\mathbb{A}^3$, $\Hilb^d(\mathbb{A}^3)$, has long been a mystery, because many key properties, such as smoothness, or dimension, components, are not yet fully understood or classified. In this paper, we use combinatorial methods to peek into the singularities of $\Hilb^d(\mathbb{A}^3)$ at the monomial ideals. The monomial ideals in $\Hilb^n(\mathbb{A}^3)$ are in one-to-one correspondence with three dimensional partitions, so many algebraic computations can be translated into combinators.

The smooth monomial ideals of $\Hilb^d(\mathbb{A}^3)$ are classified in \cite{AJR1} by Jelisiejew, Ramkumar, and Sammartano, describing the dimensions of the tangent spaces at the monomial ideals as a combinatorial count, and smoothness as combinatorial conditions on the corners of the partitions. They further prove that all the tangent spaces at the singular monomial ideals are of dimension at least $3d+6$.

This paper studies the singularities of monomial ideals whose tangent spaces have the minimal dimension $3d+6$, motivated by the following conjecture, proposed and verified for $d \leq 7$ by Hu \cite{Hu}, verified for tripod monomial ideals in \cite{AJR1}.

\begin{conjecture} \cite{Hu}
\label{Hu}
If $[S/I] \in \Hilb^d(\mathbb{A}^3)$ is smoothable with $\dim_k T(I)=3d+6$, then the singularity at $[S/I] \in \Hilb^d(\mathbb{A}^3)$ is smoothly equivalent to the vertex of a cone over the Grassmannian
$\mathrm{Gr}(2,6) \hookrightarrow \mathbb{P}^{14}$ in its Pl\"ucker embedding.
\end{conjecture}
 
\subsection{Results}
We classify into two families all singular monomial ideals with tangent space dimension $3d+6$, and we call these two families opera houses and rotating staircases, and we give their combinatorial and algebraic descriptions.

\begin{theorem}
A minimized monomial ideal $I \in \Hilb^d(\mathbb{A}^3)$ whose tangent space has dimension $3d+6$ is either an opera house or a rotating staircase. 

All the opera houses are in one linkage class linking to the ideal $(x,y,z)^2$, and the singularities at the opera houses are smoothly equivalent to the vertex of a cone over the Grassmannian
$\mathrm{Gr}(2,6) \hookrightarrow \mathbb{P}^{14}$ in its Pl\"ucker embedding.
\end{theorem} 

See Definition \ref{def_opera_house} for the combinatorial definition of opera houses, Definition \ref{alg_opera_house} for the explicit monomial generators of opera house ideals, and Example \ref{example_red_opera_house} for an example of linking more complicated opera houses to simpler ones. 

See Definition \ref{rot_stair_definition_comb} for the combinatorial definition of rotating staircases, Example \ref{rot_star_examples}) for different rotating staircases, and Example \ref{rot_star_link_class} for rotating staircases in different linkage classes.
 
\subsection{Methods.}\subsubsection{Minimization.}
For a partition corresponding to an ideal in the Hilbert scheme $\Hilb^d(\mathbb{A}^3)$, we first inductively remove the corners that have the biggest two out of the three $x,y,$ and $z$ coordinates (see Example \ref{minimized} as an example of such reduction), to obtain what we call minimized monomial ideals. Algebraically, this corresponding to finding a smaller ideal in the same linkage class with the same singularity type.

\subsubsection{Classification lemmas.}The classification part of the monomial ideals is mostly combinatorial. For an arbitrary monomial ideal $I \in \Hilb^d(\mathbb{A}^3)$ with minimal tangent space dimension $3d+6$, we develop a collection of lemmas in Section 3, 4 and 5 that each contributes as a necessary condition of restricting the tangent space dimension. Combining all the lemmas together, we explicitly construct the classification conditions. Then, we prove the converse, that those staircases satisfying all of the conditions have the desired tangent space dimension, which concludes the proof of the classification theorems.

The key lemma that we use repeatedly, Lemma \ref{bigimportantlemma} follows from the intuition that the dimension of the tangent space is combinatorially determined by the interaction between the outer corners and the inner corners of the partition, such that to have the smallest number of tangent dimension among all singular monomial ideals, we need the interaction to be minimal. so after a failed attempt at directly developing a single comprehensive classification theorem, the author finds it more efficient to develop many small lemmas, and each one would describe some control in the interaction of the corners either in terms of shapes or in terms of the numerical sizes of the partition. And all the small classification lemmas work together leading to the full classification.

\subsubsection{Controlling the middle layers by the top and bottom layers.} Three dimensional partitions can have many complicated structures than two dimensional ones, which is the difficulty of such classification lemmas, and we simplify the computation largely by realizing that we can deal first with the top and bottom layer of the three dimensional partition, and develop lemmas on just the interaction between these two layers, and that confines what the middle layers can be because the middle layers must both contain the top and be contained in the bottom layer. 

When the top layer is a rectangle (Section 4), it's easy enough to discuss case by case. When the top layer has more than two outer corners (Section 5), the bottom layer must be very similar in size as the top layer, therefore the middle layers can not have much variation, eliminating many complications. 

\subsubsection{Linkage.} To compute linkage, we test many examples with Macaulay2 and develop an algorithm to reduce a complicated opera house to another smaller and simpler opera house, and eventually the algorithm reduces any arbitrary opera house via linkage to the ideal $(x,y,z)^2$, which has Grassmannian singularity. For the rotating staircases, the singularity types are not fully classified, and we use some general linkage classification theorem to find within the rotating staircases family, many pairs of rotating staircase in linkage pair, that are each in different linkage classes.

\subsection{Acknowledgement.} This work was conducted under the supervision of Joachim Jelisiejew during the author's stay in Warsaw. The author thanks Joachim Jelisiejew for suggesting of the problem and the name ``opera houses'', and many helpful discussions.

No AI was used for this paper.

\section{background}

\subsection{Staircases.\\}
For a monomial ideal $I \subset S:=\mathbb{C}[x,y,z]$, we define its
\textbf{staircase} $E \subset \mathbb{N}^3$ 
to be the monomials of $S$ that are not in $I$, which are in bijection with a three dimensional partition via the map that sends a monomial $x^ay^bz^c$ to the box with coordinates $(a,b,c)$. 

\begin{remark}
We feel the necessity to put a warning here, that the coordinates of the staircases start from $(0,0,0)$ corresponding to the monomial $1$. So staircase consisting of a row of $a$ boxes, for example, have coordinates $(0,0,0), ... (a-1, 0,0)$. The fact that we must use the convention $a -1$ with a minus one instead of just $a$ will show up many times in the computation later, but it's unavoidable and also consistently confusing so the reader will eventually get used to it.

We also remark that, to draw pictures of staircases, sometimes we lattice points, and sometimes we use boxes.
\end{remark}

\begin{definition}
Algebraically, the socle of a commutative local ring $R$ is the annihilator of its maximal ideal $m$, i.e.
$soc(R)=(0:m)$. Combinatorially, the \textbf{outer corners} $e$ of $E$ are the maximal
elements in their degrees with respect to the usual partial order of $E$, i.e. the outer corners are $e \in E$ that satisfy that all of $xe,ye,ze \in I$. 
The socle $soc(S/I)$ is spanned by all the outer corners of $E$.

The \textbf{inner corners} $i$ of a staircase $E$ are the monomials $i \in I$ that are minimal in their degrees with respect to the usual partial order, i.e. a monomial $i \in I$ is an inner corner of $E$ if for all monomials $m$ that satisfy any one of the conditions $xm=i$, $ym=i$ or $zm=i$, we have that $m \in E$.
\end{definition}

\begin{example}
\label{staircaseex}
In the following example, the staircase $E$ of the the ideal $I=(x^3, y^5, z^4, x^2y, xy^2,y^4z,xyz^3)$ is the following partition. The outer corners are $x^2z^3,xyz^2,y^3z^3,y^4$. The inner corners are $x^2y, xy^2,y^4z,xyz^3,x^3,y^5, z^4$.

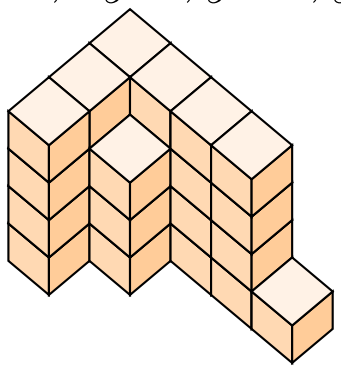
\monomialIdealPicture{{4,4,4},{4,3},{4},{4},{1}}

%\monomialIdealPicture{{4,4,4,4,4,2,2},{4,4,4,4,4,2,2},{4,4,4},{4,4},{4,4},{4},{1}}
\end{example}

\begin{definition}
Let $I \subset S$ be a co-finite monomial ideal, i.e. the quotient ring $S/I$ has finite dimension over $k$. A singularizing triple for $I$ is a triple of outer corners $\{a,b,c\}$, where $a=(a_1,a_2,a_3)$, $b=(b_1,b_2,b_3)$, $c=(c_1,c_2,c_3)$ such that all of the three conditions hold: (1) $a_1 > b_1, c_1$; (2) $ b_2 > a_2,c_2$; (3) $ c_3 > a_3,b_3$.
\end{definition}

\begin{example}
For the staircase $E'$ of the ideal $I'=(x^3, y^4,z^3,xy, xz, yz)$, there is a singularizing triple $(2,0,0), (0,3,0), (0,0,2)$.

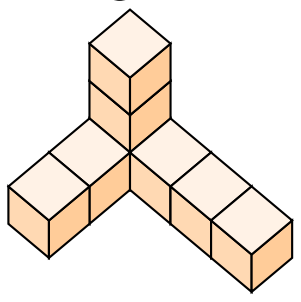
\monomialIdealPicture{{3,1,1},{1},{1},{1}}

The staircase $E$ in the Example \ref{staircaseex} above has no singularizing triple.
\end{example}

The smooth monomial ideals of $\Hilb^d(\mathbb{A}^3)$ are classified as follows.

\begin{theorem}
\cite{AJR1}
\label{smoothconditionsAJR}
Let $I \subseteq S = k[x,y,z]$ be a monomial ideal and $[S/I] \in \Hilb^d(\mathbb{A}^3)$ the corresponding point. The following conditions are equivalent
\begin{enumerate}
\item The point $[S/I]$ is a smooth point of the Hilbert scheme.
\item The ideal $I$ admits no singularizing triple.
\item The ideal $I$ is licci.
\end{enumerate}
\end{theorem}

\begin{definition}
A \textbf{dominating corner} of a cofinite monomial ideal $I$ is an outer corner $d$ where $d=(d_1,d_2,d_3)$, for which there exist two
distinct indices $i,j\in\{1,2,3\}$ such that for every other
outer corner $e=(e_1,e_2,e_3)$ of $I$, we have $d_i \geq e_i$ and $d_j \geq e_j$.
\end{definition}

\begin{example}
In the above \ref{staircaseex}, the outer corner $(2,0,3)$ is a dominating corner. 
\end{example}

\subsection{Linkage and smoothness}
\begin{definition}
Two ideals $I,J \in k[x,y,z]$ are linked if there is a regular sequence $g_1, g_2, g_3$ such that for the ideal $G$ generated as $G=(g_1, g_2, g_3)$, we have $G \subset I \cap J$, $J=G : I$ and $I=G : J$. 
\end{definition}

\begin{definition}
\cite[Definition 5.5]{AJR1}
Two pointed schemes $(X,x)$, $(Y,y)$ are \textbf{smoothly equivalent} if there exists a third pointed scheme $(Z,z)$ with smooth maps
\[
(Z,z)\to (X,x), \qquad (Z,z)\to (Y,y).
\]
\end{definition}

\begin{theorem}
\label{AJR_singularity_type}
\cite[Theorem 5.6] {AJR1}
Let $[Z]\in \operatorname{Hilb}^d(X)$ and
$[Z'']\in \operatorname{Hilb}^{d''}(X)$.
If $Z$ is linked to $Z''$, then the singularity at $[Z]$ is smoothly equivalent to the singularity at $[Z'']$.
\end{theorem}

\begin{proposition}
Every staircase with at least one dominating corner is linked to some staircase with no dominating corner. 
\end{proposition}

\begin{proof}
Suppose the staircase $E$ of the ideal $I$ has a dominating corner $d=(d_1,d_2,d_3)$, and assume without loss of generality that $d_1 \geq e_1$, and $d_2 \geq e_2$ for all outer corners $e=(e_1,e_2,e_3)$ of $E$. Take linkage of $I$ by the regular sequence $(x^{d_1+1}, y^{d_2+1},z^{h_z+1})$ where $h_z$ is the largest $z$-coordinate of the staircase $E$. The result $E'$ of $I':=(x^{d_1+1}, y^{d_2+1},z^{h_z+1}):I$ is the complement of $E$ in the cube of size $(d_1+1)\times(d_2+1)\times(h_z+1),$ with the dominating corner $d$ removed.
%\Ytodo{maybe add explanation/reference, maybe no} 

If there are more dominating corners of the staircase $E'$ after taking linkage, we inductively repeat the process above, until there are no dominating corners in the staircase.
\end{proof}

\begin{example}
\label{minimized}
The left staircase of the ideal $(x^3, y^5, z^4, x^2yz^2, xy^2z^2,y^4z^3,x^2z^3)$ has dominating corner $x^3y^5z^2$. The right staircase is taking its linkage by the complete intersection $(x^3, y^5, z^4)$, with ideal $(z^2,x^2z,x^2y,x^3,xy^3,y^4z,y^5)$, that has no dominating corner.

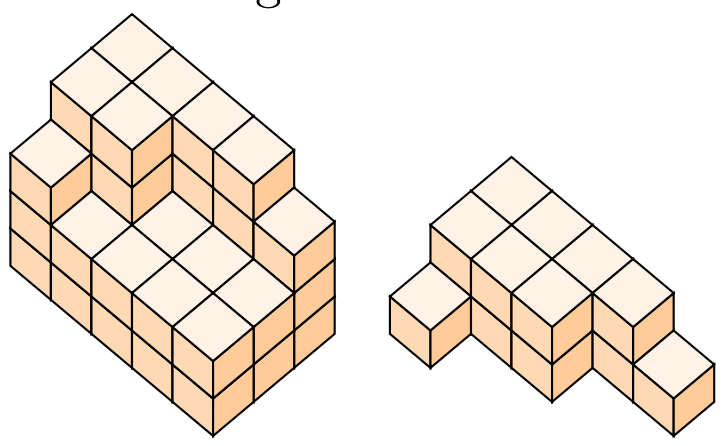
\monomialIdealPicture{{4,4,3},{4,4,2},{4,2,2},{4,2,2},{3,2,2}}
\monomialIdealPicture{{2,2,1},{2,2},{2,2},{2},{1}}
\end{example}

\begin{definition}
A staircase with no dominating corner is denoted as a \textbf{minimized staircase}. 
\end{definition}

From this point on, we always consider minimized staircases as a representative of a linkage class,  and if not otherwise specified, a staircase refers to a minimized one.

\subsection{Tangent space decomposition}
In this section we decompose the tangent space at a monomial ideal into graded pieces, and state some combinatorial theorems on the dimension of each graded tangent space.

\begin{definition}\cite[Definition 0.1.]{AR}
A \textbf{signature} is a non-constant triple on the two-element set $\{\text{p}, \text{n}\}$, where
\[
\text{p} = \text{``positive or 0''}, \quad \text{n} = \text{``negative''}.
\]
Let $\mathfrak{S}$ denote the set of signatures, and define for each $\mathfrak{s} \in \mathfrak{S}$

$$\mathbb{Z}_{\mathfrak{s}}^3 = \{(\alpha_1, \alpha_2, \alpha_3) \in \mathbb{Z}^3 : \alpha_i \geqslant 0 \text{ if } \mathfrak{s}_i = \text{p}; \alpha_i < 0 \text{ if } \mathfrak{s}_i = \text{n}\}.$$
we decompose the tangent space $T(I)$ of $I$ in the Hilbert scheme into subspaces defined in terms of $\mathbb{Z}$ graded directions, where $|T(I)|_\alpha$ denotes the graded component of $T(I)$ of degree $\alpha \in \mathbb{Z}^3$.
$$T_{\mathfrak{s}}(I) = \bigoplus_{\alpha \in \mathbb{Z}_{\mathfrak{s}}^3} |T(I)|_\alpha \subseteq T(I).$$
\end{definition}

\begin{remark}\cite[Proposition 1.9.]{AR} 
The tangent spaces
$T_{\mathrm{ppp}}(I)=T_{\mathrm{nnn}}(I)=0$,
and therefore
$T(I)=\bigoplus_{\mathfrak{s} \in \mathfrak{S}}T_{\mathfrak{s}}(I)$.    
\end{remark}

\begin{theorem}
\cite[Theorem 2.4.]{AR}
For a field $k$ and any monomial point $[I] \in \Hilb^d (\mathbb{A}^3)$, 
\begin{align*}
\dim_{{k}} T_{\text{ppn}}(I) &= \dim_{{k}} T_{\text{nnp}}(I) + d, \\
\dim_{{k}} T_{\text{pnp}}(I) &= \dim_{{k}} T_{\text{npn}}(I) + d, \\
\dim_{{k}} T_{\text{npp}}(I) &= \dim_{{k}} T_{\text{pnn}}(I) + d.
\end{align*}
\end{theorem}

\begin{corollary}
For any monomial point $I \in \Hilb^d(\mathbb{A}^3)$, 
\begin{align*}
\dim_{k} T(I)
%&= \dim_{k} T_{nnp}(I) +\dim_{k} T_{npn}(I)+\dim_{k} T_{pnn}(I)\\
%&+\dim_{k} T_{npp}(I)+\dim_{k} T_{pnp}(I)+\dim_{k} T_{ppn}(I)\\
&=2(\dim_{k} T_{nnp}(I) +\dim_{k} T_{npn}(I)+\dim_{k} T_{pnn}(I))+3d. 
\end{align*}

\end{corollary}
Therefore to compute the dimension of the tangent space of $I \in \Hilb^d(\mathbb{A}^3)$, it is enough to compute the dimension of $\dim T_{nnp}(I), \dim T_{npn}(I)$ and $\dim T_{pnn}$. We use the following combinatorial condition on the staircase $E$ to compute the dimensions of these tangent spaces.

\begin{proposition}
\label{proptanspacedecom}
\cite[Proposition 1.5.]{AR}
Let $\alpha \in \mathbb{Z}^n$ and $[I] \in \Hilb^d (\mathbb{A})^n$. The set of bounded connected components of $I \setminus (I - \alpha)$ corresponds to a basis of $|T(I)|_\alpha$.
\end{proposition}

Here, $I \setminus (I-\alpha)$ denotes the set of lattice points of $I$ that do not belong to the translate $I-\alpha$.  The next sections of the paper focus on determining the combinatorial conditions under which such bounded connected components arise.
%In the above proposition, $I \setminus (I - \alpha)$ is the collection of lattice points of $I$ and removing the set $I$ shifted by the vector $-\alpha$. Proposition \ref{proptanspacedecom} reduces the computation of dimension of tangent space of an ideal to counting bounded connected components of $I\setminus(I - \alpha)$. Therefore, throughout the paper our goal is to determine when such bounded components appear.

\subsection{Bounded connected components}

Denote $U$ as the union of lattice points in $\mathbb{Z}^3$ where at least one of the $x$ or $y$ axis is negative, and the $z$ coordinate is non-negative. Denote $I^c$ as the complement of $I$ in the lattice $\mathbb{Z}^3$.

\begin{definition}
We decompose the set $I \setminus (I-\alpha)$ into two disjoint subsets. Define the \textbf{infinite component} and \textbf{finite component} of $I \setminus (I - \alpha)$
$$
I_{\infty}:=I\cap(U-\alpha),\qquad
I_{\mathrm{fin}}:=I\cap(E-\alpha).
$$
Since $E \sqcup I \sqcup U=\mathbb{Z}^3$, we obtain the decomposition
$$
I\setminus(I-\alpha)=I_{\mathrm{fin}}\sqcup I_{\infty}.
$$
\end{definition}

\begin{remark}
Consider $\alpha=(-n_1, -n_2, p), n_1, n_2 \geq 1, p \geq 0$ with signature $nnp$, and the proof for vectors of other signatures follows similarly. The translation $U-\alpha$ shifts $U$ to the positive $x$, positive $y$, and negative $z$ direction. So the infinite component $I_\infty$ is the collection of boxes in $I$ whose $x$ and $y$ coordinates satisfy $0 \leq x < n_1$ or $0 \leq y < n_2$, and the $z$ coordinate is non-negative, which is always one connected piece that extends to infinity in the positive $x,y$ and $z$ directions. Intuitively, the infinite component $I_\infty$ has is an L-shaped wall, with thickness $n_1$ along the $x$ axis, and thickness $n_2$ along the $y$ axis (see Example \ref{example_wall} below). 
\end{remark}

\begin{remark}
\label{rmkpathconnect}
Since $E$ is bounded, every connected component of $I\cap(E-\alpha)$ is bounded. Therefore, the bounded connected components of $I\setminus(I-\alpha)$ are precisely the connected components of $I\cap(E-\alpha)$ that are not connected by a lattice path to the infinite component $I\cap(U-\alpha)$.
\end{remark}

\begin{example}
\label{example_wall}
Let $I=(x^2, y^2, z^2, xy, xz, yz)$ and $\alpha=(-1,-1,1)$. The staircase $E$ is marked by the black dots. The infinite component $I_\infty$ is marked by the \textcolor{teal!90!black}{blue} dots. We've only drawn a few dots, but the "wall" shaped component has thickness $1$, extends infinitely in the positive $x, y, z$ directions. The finite component $I_{\textit{fin}}$ is marked by the \textcolor{red!70!black}{red} dots. The finite component is not path connected to the infinite component, so there is a unique bounded connected component.
\end{example}

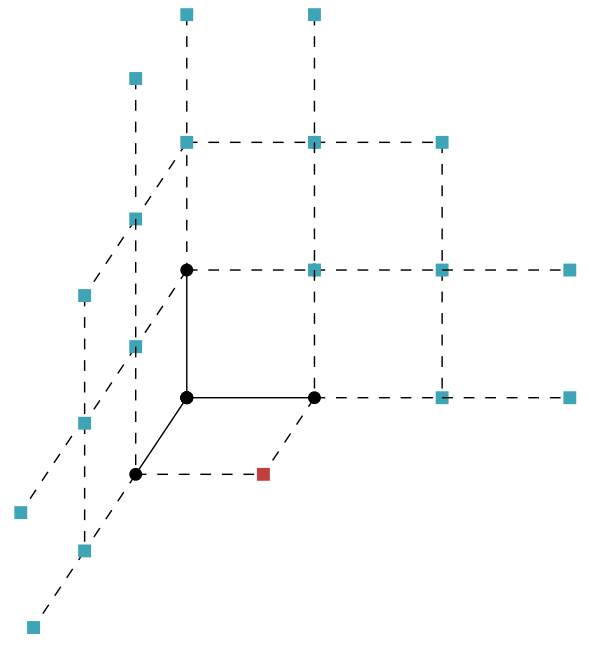
\resizebox{0.3\textwidth}{!}{%
\begin{circuitikz}
\tikzstyle{every node}=[font=\fontsize{18.2pt}{23.7pt}\selectfont]
\begin{scope}[transparency group, opacity=1]
\node at (6.25,11.125) [circ, color={rgb,255:red,0; green,0; blue,0}] {};
\end{scope}
\begin{scope}[transparency group, opacity=1]
\node at (6.25,11.125) [circ, color={rgb,255:red,0; green,0; blue,0}] {};
\end{scope}
\begin{scope}[transparency group, opacity=1]
\node at (6.25,11.125) [circ, color={rgb,255:red,0; green,0; blue,0}] {};
\end{scope}
\begin{scope}[transparency group, opacity=1]
\node at (7.5,11.125) [circ, color={rgb,255:red,0; green,0; blue,0}] {};
\end{scope}
\begin{scope}[transparency group, opacity=1]
\node at (6.25,11.125) [circ, color={rgb,255:red,0; green,0; blue,0}] {};
\end{scope}
\draw (6.25,11.125) to[short] (7.5,11.125);
\begin{scope}[transparency group, opacity=1]
\node at (6.25,12.375) [circ, color={rgb,255:red,0; green,0; blue,0}] {};
\end{scope}
\draw (6.25,12.375) to[short] (6.25,11.125);
\draw (7.5,12.375) to[short, -o] (7.5,12.375) ;

\draw [short] (6.25,11.125) -- (5.75,10.375);
\draw [dashed] (5.75,11.625) -- (5.25,10.875);
\draw [dashed] (8.75,12.375) -- (8.75,11.125);
\draw [dashed] (5.25,9.625) -- (4.75,8.875);
\draw [dashed] (6.25,12.375) -- (7.5,12.375);
\draw [dashed] (6.25,13.5) -- (6.25,12.25);
\draw [dashed] (6.25,12.375) -- (5.75,11.625);
\draw [dashed] (5.25,10.875) -- (5.25,9.625);
\draw [dashed] (5.125,10.75) -- (4.625,10);
\draw [dashed] (7.5,11.125) -- (8.75,11.125);
\begin{scope}[transparency group, opacity=1]

\end{scope}

\begin{scope}[transparency group, opacity=1]
\node at (5.75,10.375) [circ, color={rgb,255:red,0; green,0; blue,0}] {};
\end{scope}
\draw [dashed] (5.75,10.375) -- (5.25,9.625);
\draw (7.5,13.625) to[short, -o] (7.5,13.625) ;
\draw [dashed] (8.75,13.625) -- (8.75,12.375);

\draw [dashed] (7.5,12.375) -- (8.75,12.375);
\draw [dashed] (5.75,10.375) -- (7,10.375);
\draw [dashed] (7.5,11.125) -- (7,10.375);
\draw [dashed] (6.25,13.625) -- (7.5,13.625);
\draw [dashed] (7.5,13.625) -- (8.75,13.625);
\draw [dashed] (6.25,13.625) -- (5.75,12.875);
\draw [dashed] (5.75,12.875) -- (5.25,12.125);
\draw [dashed] (5.25,12.125) -- (5.25,10.875);
\begin{scope}[transparency group, opacity=1]
\node at (5.25,10.875) [squarepole, color={rgb,255:red,62; green,165; blue,182}] {};
\end{scope}
\begin{scope}[transparency group, opacity=1]
\node at (5.75,11.625) [squarepole, color={rgb,255:red,62; green,165; blue,182}] {};
\end{scope}
\begin{scope}[transparency group, opacity=1]
\node at (7.5,12.375) [squarepole, color={rgb,255:red,62; green,165; blue,182}] {};
\end{scope}
\begin{scope}[transparency group, opacity=1]
\node at (8.75,12.375) [squarepole, color={rgb,255:red,62; green,165; blue,182}] {};
\end{scope}
\begin{scope}[transparency group, opacity=1]
\node at (5.75,12.875) [squarepole, color={rgb,255:red,62; green,165; blue,182}] {};
\end{scope}
\begin{scope}[transparency group, opacity=1]
\node at (5.25,12.125) [squarepole, color={rgb,255:red,62; green,165; blue,182}] {};
\end{scope}
\begin{scope}[transparency group, opacity=1]
\node at (7.5,13.625) [squarepole, color={rgb,255:red,62; green,165; blue,182}] {};
\end{scope}
\begin{scope}[transparency group, opacity=1]
\node at (8.75,13.625) [squarepole, color={rgb,255:red,62; green,165; blue,182}] {};
\end{scope}
\begin{scope}[transparency group, opacity=1]
\node at (6.25,13.625) [squarepole, color={rgb,255:red,62; green,165; blue,182}] {};
\end{scope}
\begin{scope}[transparency group, opacity=1]
\node at (5.25,9.625) [squarepole, color={rgb,255:red,62; green,165; blue,182}] {};
\end{scope}
\begin{scope}[transparency group, opacity=1]
\node at (8.75,11.125) [squarepole, color={rgb,255:red,62; green,165; blue,182}] {};
\end{scope}
\draw [dashed] (8.75,11.125) -- (10,11.125);
\draw [dashed] (8.75,12.375) -- (10,12.375);
\draw [dashed] (6.25,14.875) -- (6.25,13.625);
\draw [dashed] (5.75,11.625) -- (5.75,10.375);
\draw [dashed] (5.75,12.875) -- (5.75,11.625);
\draw [dashed] (7.5,13.625) -- (7.5,12.375);
\draw [dashed] (7.5,12.375) -- (7.5,11.125);
\begin{scope}[transparency group, opacity=1]
\node at (4.625,10) [squarepole, color={rgb,255:red,62; green,165; blue,182}] {};
\end{scope}
\begin{scope}[transparency group, opacity=1]
\node at (4.75,8.875) [squarepole, color={rgb,255:red,62; green,165; blue,182}] {};
\end{scope}
\begin{scope}[transparency group, opacity=1]
\node at (6.25,14.875) [squarepole, color={rgb,255:red,62; green,165; blue,182}] {};
\end{scope}
\begin{scope}[transparency group, opacity=1]
\node at (10,12.375) [squarepole, color={rgb,255:red,62; green,165; blue,182}] {};
\end{scope}
\begin{scope}[transparency group, opacity=1]
\node at (10,11.125) [squarepole, color={rgb,255:red,62; green,165; blue,182}] {};
\end{scope}
\draw [dashed] (5.75,14.25) -- (5.75,13);
\draw [dashed] (7.5,14.875) -- (7.5,13.625);
\begin{scope}[transparency group, opacity=1]
\node at (5.75,14.25) [squarepole, color={rgb,255:red,62; green,165; blue,182}] {};
\end{scope}
\begin{scope}[transparency group, opacity=1]
\node at (7.5,14.875) [squarepole, color={rgb,255:red,62; green,165; blue,182}] {};
\end{scope}
\begin{scope}[transparency group, opacity=1]
\node at (7,10.375) [squarepole, color={rgb,255:red,193; green,63; blue,63}] {};
\end{scope}
\end{circuitikz}
}%
 
Assume for the rest of this section that the vector $\alpha=(-n_1, -n_2, p)$. The cases of other signatures $npn, pnn$ are analogous by symmetry.

\begin{definition}
Let $h_x$, $h_y$, and $h_z$ denote the maximal $x$, $y$, and $z$ coordinates of a staircase $E$. Unless otherwise stated, we canonically take the $z$ direction as the vertical direction and $h_z$ the height of the staircase. (Technically, the staircase has $h_z+1$ layers, but we abuse the notation because it's easier to work with than having to put plus one all the time.)

\begin{enumerate}
\item The boxes with $z$-coordinate $0$ form the \textbf{bottom layer}, and those
with $z$ coordinate $h_z$ form the \textbf{top layer}. We similarly define
the top and bottom layers in the $x$ and $y$ directions. Unless
otherwise specified, the terms top layer and bottom layer
refer to the $z$ direction.
\item For $0 \leq i\leq h_z$, let $P_i$ denote the two-dimensional
partition given by the $i$-th $z$-layer of $E$. The \textbf{inner corners of the layer} $P_i$ are the inner corners of $E$
whose $z$ coordinate is $i$, excluding those with maximal
$x$ or $y$ coordinate in that layer. Equivalently, the inner
corners of $P_i$ are the internal inner corners of the
two-dimensional partition $P_i$.
\end{enumerate}
\end{definition}

\begin{example}
For the following staircase, the inner corners of the bottom layer are $x^3y, xy^2$.  The inner corner of the middle layer is $xy$. There are no inner corners of the top layer.

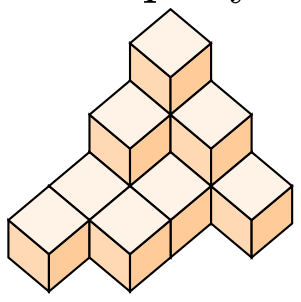
\monomialIdealPicture{{3,2,1,1},{2,1,1},{1}}
\end{example}

\begin{lemma}
\label{innercorbottomlayer}
The bottom layer of a minimized staircase must have at least one inner corner.
\end{lemma}

\begin{proof}
Suppose there is no inner corner of the bottom layer, and then all the rows of boxes of the bottom layer have the same length, and denote their coordinates as $a_1 =... =a_k$, which means that $(a_1,k,0)$ is a dominating corner, so this is a contradiction.
\end{proof}

\section{tool box}
This section develops some combinatorial lemmas that will be used repeatedly in the classification theorems in the next few sections. Each lemma provides a geometric configuration that forces the existence of one or more bounded connected components. We start with some definitions.

\begin{definition}
Let $O=(o_1,o_2,o_3)$ be an outer corner of the staircase $E$.
An \textbf{extended rectangular outer corner} of size
$(r_1,r_2,r_3)$ at $O$ is a rectangular cuboid $R^+$ %$R:=R^+\cap E,$
such that
$$R^+ =
\left\{
m=(m_1,m_2,m_3):
\begin{array}{l}
o_1-r_1\leq m_1\leq o_1+1,\\
o_2-r_2\leq m_2\leq o_2+1,\\
o_3-r_3\leq m_3\leq o_3+1
\end{array}
\right\},$$
such that every lattice point $m\in R^+$ satisfies,
\begin{enumerate}
    \item $m\in I$ if $m_i=o_i+1$ for any $i =1,2,3$;
    \item $m\in E$ if $o_i-r_i+1 \leq m_i \leq o_i$ for all $i=1,2,3$.
\end{enumerate}
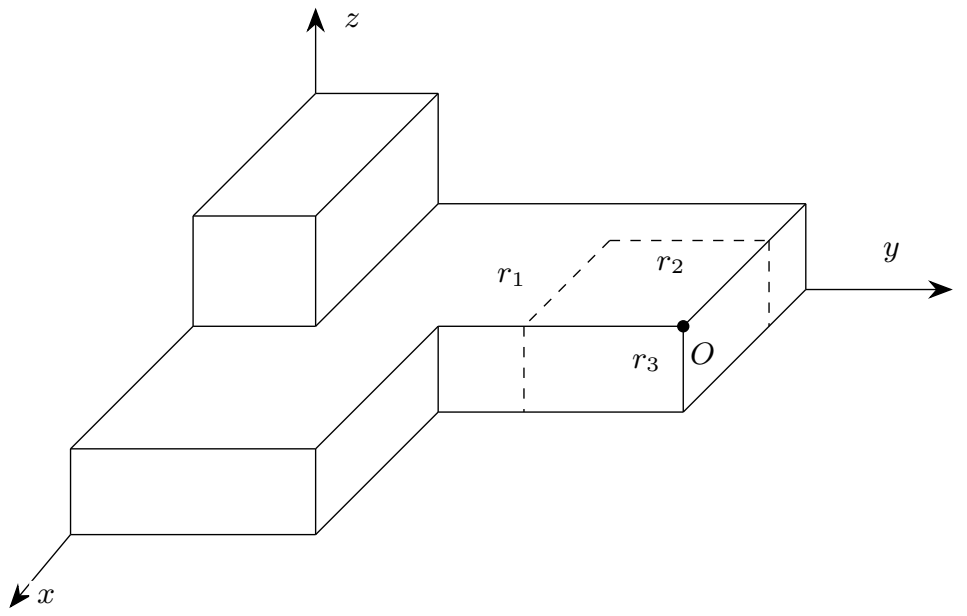
\resizebox{0.5\textwidth}{!}{%
\begin{circuitikz}
\tikzstyle{every node}=[font=\fontsize{9.1pt}{11.8pt}\selectfont]
\draw (8.75,9.875) to[short] (7.5,8.625);
\draw (7.5,8.625) to[short] (5,8.625);
\draw (5,8.625) to[short] (3.75,7.375);
\draw (8.75,9.875) to[short] (5,9.875);
\draw (2.5,8.625) to[short] (1.25,7.375);
\draw (1.25,7.375) to[short] (3.75,7.375);
\draw (1.25,7.375) to[short] (1.25,6.5);
\draw (1.25,6.5) to[short] (3.75,6.5);
\draw (3.75,7.375) to[short] (3.75,6.5);
\draw (3.75,6.5) to[short] (5,7.75);
\draw (5,8.625) to[short] (5,7.75);
\draw (5,7.75) to[short] (7.5,7.75);
\draw (7.5,8.625) to[short] (7.5,7.75);
\draw (7.5,7.75) to[short] (8.75,9);
\draw (8.75,9.875) to[short] (8.75,9);
\draw (2.5,8.625) to[short] (3.75,8.625);
\draw (5,9.875) to[short] (3.75,8.625);
\draw (2.5,8.625) to[short] (2.5,9.75);
\draw (2.5,9.75) to[short] (3.75,9.75);
\draw (3.75,9.75) to[short] (3.75,8.625);
\draw (3.75,9.75) to[short] (5,11);
\draw (5,11) to[short] (5,9.875);
\draw (2.5,9.75) to[short] (3.75,11);
\draw (3.75,11) to[short] (5,11);
\draw [dashed] (6.75,9.5) -- (8.375,9.5);
\draw [dashed] (6.75,9.5) -- (5.875,8.625);
\draw [dashed] (5.875,8.625) -- (5.875,7.75);
\draw [dashed] (8.375,9.5) -- (8.375,8.625);
\node [font=\fontsize{9.1pt}{11.8pt}\selectfont, inner xsep=0.080cm, inner ysep=0.085cm, rounded corners=0.020cm] at (5.75,9.125) {$r_1$};
\node [font=\fontsize{9.1pt}{11.8pt}\selectfont, inner xsep=0.080cm, inner ysep=0.085cm, rounded corners=0.020cm] at (7.375,9.25) {$r_2$};
\node [font=\fontsize{9.1pt}{11.8pt}\selectfont, inner xsep=0.080cm, inner ysep=0.085cm, rounded corners=0.020cm] at (7.125,8.25) {$r_3$};
\node [font=\fontsize{9.1pt}{11.8pt}\selectfont, inner xsep=0.080cm, inner ysep=0.085cm, rounded corners=0.020cm] at (7.7,8.35) {$O$};

\draw [-{Stealth[scale=1.5]}, ] (1.25,6.5) -- (0.625,5.75);
\draw [-{Stealth[scale=1.5]}, ] (8.75,9) -- (10.25,9);
\draw [-{Stealth[scale=1.5]}, ] (3.75,11) -- (3.75,11.875);
\node [font=\fontsize{9.1pt}{11.8pt}\selectfont, fill={rgb,255:red,255; green,255; blue,255}, fill opacity=1, text opacity=1, inner xsep=0.080cm, inner ysep=0.085cm, rounded corners=0.020cm] at (1,5.875) {$x$};
\node [font=\fontsize{9.1pt}{11.8pt}\selectfont, fill={rgb,255:red,255; green,255; blue,255}, fill opacity=1, text opacity=1, inner xsep=0.080cm, inner ysep=0.085cm, rounded corners=0.020cm] at (9.625,9.375) {$y$};
\node [font=\fontsize{9.1pt}{11.8pt}\selectfont, fill={rgb,255:red,255; green,255; blue,255}, fill opacity=1, text opacity=1, inner xsep=0.080cm, inner ysep=0.085cm, rounded corners=0.020cm] at (4.125,11.75) {$z$};
\node at (7.5,8.625) [circ, color={rgb,255:red,0; green,0; blue,0}] {};
\end{circuitikz}
}%
\end{definition}

\begin{definition}
The \textbf{interior faces} $F^-$ are $$F^-:=\{m \in R^+| m_1=o_1-r_1, \text{ or } m_2=o_2-r_2, \text{ or } m_3=o_3-r_3\}.$$ The 
\textbf{exterior faces} $F^+$ are $$F^+:=\{m \in (R^+ \setminus F^-)| m_1=o_1+1, \text{ or } m_2=o_2+1, \text{ or } m_3=o_3+1\}.$$

%$F^+:=R^+ \setminus R,$ all the monomials in $R^+$ that are not in $R$. 
\end{definition}

The above definitions are the following intuition written down rigorously using coordinates. We want an outer corner $O$ along with the rectangular cuboid of size $r_1 \times r_2 \times r_3$ adjacent to it (see picture above). And we would like this outer corner $O$ to be isolated from the other outer corners of the staircase $E$. That's equivalent to imposing conditions on a larger rectangular cuboid $R+$ of size $(r_1+1) \times (r_2+1) \times (r_3+1)$ around $O$, and require that the exterior face $F^+$, which are the boxes to the upper, right, and front direction of $O$, to be not contained in $E$. Similarly we also have the following definition.

\begin{definition}
For an inner corner $U=(u_1, u_2, u_3)$, define its \textbf{extended rectangular inner corner of size} $(t_1, t_2, t_3)$ at $U$ is a rectangular cuboid $R^+$ %$R: R^+ \cap I$, 
such that 
$$
R^+ = \left\{ 
m=(m_1,m_2,m_3):
\begin{array}{l}
u_1-1\leq m_1\leq u_1+t_1,\\
u_2-1\leq m_2\leq u_2+t_2,\\
u_3-1\leq m_3\leq u_3+t_3
\end{array} \right\},$$
such that every lattice point $m\in R^+$ satisfies
\begin{enumerate}
    \item $m\in E$ or $m_i <0$ for some $i$, if $m_1=u_1-1$, or $m_2=u_2-1$, or $m_3=u_3-1$;
    \item $m\in I$ otherwise.
\end{enumerate}
\end{definition}

The following Lemma gives a condition where there must be at least one bounded connected component. 

\begin{lemma}
\label{bigimportantlemma}
Let $R:= R^+ \setminus (F^+ \cup F^-)$. For a staircase $E$, suppose $O=(o_1,o_2,o_3)$ is an extended rectangular outer corner of size $s=(s_1, s_2, s_3)$, and $U=(u_1, u_2,u_3)$ is an extended inner corner of the staircase $E$. If there exists some vector $\alpha=(-n_1,-n_2, p)$ such that $(O-\alpha) \in I$ and $F^- -\alpha \subset (E \cup S_{-xyz})$ where $S_{-xyz}$ denotes the union of the lattice points where at least one of the $x,y,z$ coordinates is negative, then $I \setminus (I-\alpha)$ has at least one bounded connected component as a subset of $A:=(R-\alpha) \cap I$.

Consequently, if $U=(u_1, u_2,u_3)$ is an extended inner corner of the staircase $E$ of size $t=(t_1, t_2, t_3)$, and $O=(o_1,o_2,o_3)$ is an extended outer corner of size $s=(s_1, s_2, s_3)$, then the condition of having $O+(1,1,1)=s+\alpha+U$, and $s \leq t$ is a special case of the general statement above, which implies that $I \setminus (I-\alpha)$ has at least one bounded connected component.

The cases of other signatures $npn$ and $pnn$ holds by symmetry. 
\end{lemma}

\begin{proof}
By assumption, $(O-\alpha) \in I$, so the set $(R-\alpha) \cap I$ is not empty with the element $O-\alpha$. Because $R \subset E$, $(R-\alpha) \cap I$ is a non-empty subset of the finite component $I_{\textit{fin}}=(E-\alpha) \cap I$.

By Remark \ref{rmkpathconnect}, it's sufficient to show that $(R-\alpha) \cap I$ is not connected to the infinite component $I_\infty$ by any lattice path. 

Because $F^--\alpha \subset (E \cup S_{-xyz})$, $(F^--\alpha) \cap I = \emptyset$. So $(F^--\alpha) \cap I_\infty = \emptyset$ and $(F^--\alpha) \cap I_{\textit{fin}} = \emptyset$.

Because $F^+ \subset I$, $(F^+ - \alpha) \cap (E-\alpha)=\emptyset$, so $(F^+-\alpha) \cap I_{\textit{fin}} = \emptyset$. Because all the monomials in $F^+$ have positive $x,y,z$ coordinates, $F^+ \cap U = \emptyset$, so $(F^+-\alpha) \cap (U-\alpha)=\emptyset$, and $(F^+-\alpha) \cap I_\infty=\emptyset$.

Taking the union, $(F^--\alpha) \cup (F^+ - \alpha)$ is the boundary of $A$ that contains no lattice point in either $I_{\textit{fin}}$ or $I_\infty$. So $A$ cannot connect to $I_\infty$ by any lattice path.

\end{proof}

Project the top and bottom layers to the $XY$ plane, and denote the outer corners of the top layer as $(x_1,y_1),...,(x_k,y_k)$, and the outer corners of the bottom layer as $(g_1,l_1),...,(g_r,l_r)$, such that $x_1 >x_2 > ... > x_k$, $y_k > ... > y_1$, and that $g_1 > g_2 >...> g_r$, $l_r>...>l_1$. We draws the monomials of a staircase as lattice points.

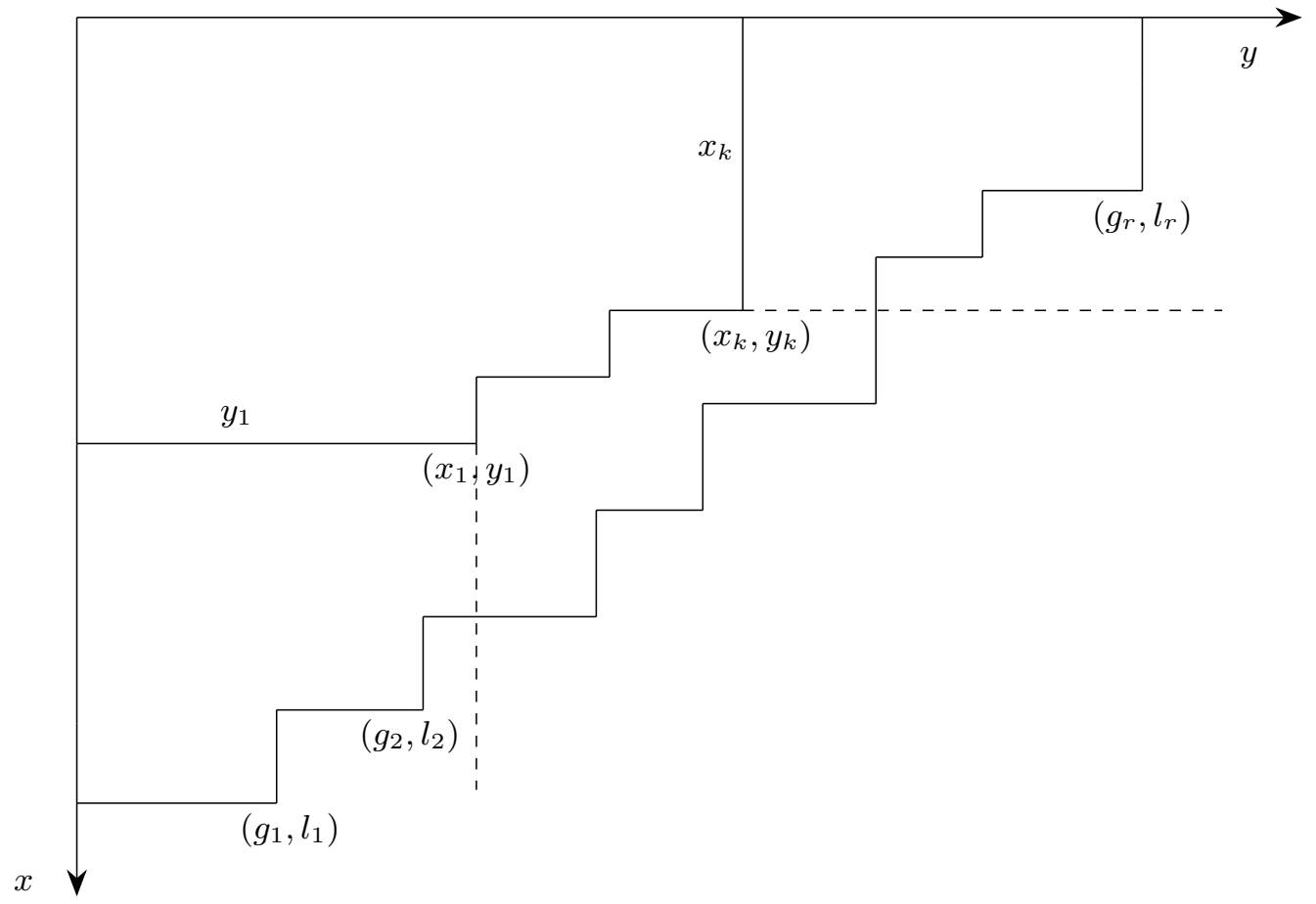
\resizebox{0.7\textwidth}{!}{%
\begin{circuitikz}
\tikzstyle{every node}=[font=\fontsize{9.1pt}{11.8pt}\selectfont]
\draw (0,11.125) to[short] (6.25,11.125);
\draw (0,11.125) to[short] (0,7.125);
\draw (6.25,11.125) to[short] (10,11.125);
\draw (10,11.125) to[short] (10,10.5);
\draw (6.25,11.125) to[short] (6.25,8.375);
\draw (6.25,8.375) to[short] (5,8.375);
\draw (5,8.375) to[short] (5,7.75);
\draw (5,7.75) to[short] (3.75,7.75);
\draw (3.75,7.75) to[short] (3.75,7.125);
\draw (3.75,7.125) to[short] (0,7.125);
\draw (0,7.125) to[short] (0,4.5);
\draw [dashed] (3.75,7.125) -- (3.75,3.875);
\draw [dashed] (6.25,8.375) -- (10.75,8.375);
\draw (10,10.5) to[short] (10,9.5);
\draw (10,9.5) to[short] (8.5,9.5);
\draw (8.5,9.5) to[short] (8.5,8.875);
\draw (8.5,8.875) to[short] (7.5,8.875);
\draw (7.5,8.875) to[short] (7.5,7.5);
\draw (7.5,7.5) to[short] (5.875,7.5);
\draw (5.875,7.5) to[short] (5.875,6.5);
\draw (5.875,6.5) to[short] (4.875,6.5);
\draw (4.875,6.5) to[short] (4.875,5.5);
\draw (4.875,5.5) to[short] (3.25,5.5);
\draw (3.25,5.5) to[short] (3.25,4.625);
\draw (3.25,4.625) to[short] (1.875,4.625);
\draw (1.875,4.625) to[short] (1.875,3.75);
\draw (1.875,3.75) to[short] (0,3.75);
\draw (0,4.5) to[short] (0,3.75);
\draw [-{Stealth[scale=1.5]}, ] (0,3.75) -- (0,2.875);
\draw [-{Stealth[scale=1.5]}, ] (10,11.125) -- (11.5,11.125);
\node [font=\fontsize{9.1pt}{11.8pt}\selectfont, inner xsep=0.080cm, inner ysep=0.085cm, rounded corners=0.020cm] at (-0.5,3) {$x$};
\node [font=\fontsize{9.1pt}{11.8pt}\selectfont, inner xsep=0.080cm, inner ysep=0.085cm, rounded corners=0.020cm] at (11,10.75) {$y$};
\node [font=\fontsize{9.1pt}{11.8pt}\selectfont, inner xsep=0.080cm, inner ysep=0.085cm, rounded corners=0.020cm] at (2,3.5) {$(g_1,l_1)$};
\node [font=\fontsize{9.1pt}{11.8pt}\selectfont, inner xsep=0.080cm, inner ysep=0.085cm, rounded corners=0.020cm] at (3.125,4.375) {$(g_2,l_2)$};
\node [font=\fontsize{9.1pt}{11.8pt}\selectfont, inner xsep=0.080cm, inner ysep=0.085cm, rounded corners=0.020cm] at (10,9.25) {$(g_r,l_r)$};
\node [font=\fontsize{9.1pt}{11.8pt}\selectfont, inner xsep=0.080cm, inner ysep=0.085cm, rounded corners=0.020cm] at (3.75,6.875) {$(x_1,y_1)$};
\node [font=\fontsize{9.1pt}{11.8pt}\selectfont, inner xsep=0.080cm, inner ysep=0.085cm, rounded corners=0.020cm] at (6.375,8.125) {$(x_k,y_k)$};
\node [font=\fontsize{9.1pt}{11.8pt}\selectfont, inner xsep=0.080cm, inner ysep=0.085cm, rounded corners=0.020cm] at (1.5,7.375) {$y_1$};
\node [font=\fontsize{9.1pt}{11.8pt}\selectfont, inner xsep=0.080cm, inner ysep=0.085cm, rounded corners=0.020cm] at (6,9.875) {$x_k$};
\end{circuitikz}
}%

\begin{lemma}
\label{sidecomp}
If the top and bottom layers of a staircase $E$ satisfy either $l_1< l_2 \leq y_1$, or $g_r < g_{r-1} \leq x_k$, %i.e. there are at least two outer corners of the bottom layer either both to the left or both to the right of the dashed lines, 
then there are at least two different vectors producing two bounded connected components of the staircase $E$.
\end{lemma}

\begin{proof}
Suppose without loss of generality that $l_1< l_2 \leq y_1$. There are two or more outer corners to the left of the dashed line. Consider the boxes $C_1=(g_1, l_1, v_1)$ and $C_2=(g_2, l_2, v_2)$ such that $C_1, C_2 \in E$, but $C_1+(0,0,1)=(g_1, l_1, v_1+1) \notin E$, and $C_2+(0,0,1)=(g_2, l_2, v_2+1) \notin E$. So $C_1, C_2$ are outer corners of $E$. By construction $v_1, v_2 \leq h_z$.

%Now we try to find a specific inner corner of the staircase $E$. The right-most box in the top layer next to the $y$-axis is the box $(1,y_k)$. The box to the right of it, $(1, y_k+1)$ is no longer in the top layer, but it is in the bottom layer with 3 dimensional coordinates $(1,y_k+1,1)$. Note that the box $(1,y_k+1)$ might also be in some layers higher than the bottom layer. 

Consider the box $N=(0,y_k+1,w)$, where $w$ is the smallest $z$ coordinate, such that the box $(0,y_k+1,w-1) \in E$ or $w=0$, but $(0,y_k+1,w) \notin E$. So $N=(0,y_k+1,w)$ is an inner corner of the staircase $E$. 

Define $s_1=-w+v_1$ and $j_1=1$ if $w>v_1$; and $s_1=-1$ and $j_1=v_1-w+1$ if $w \leq v_1$. 

Define $s_2=-w+v_2$ and $j_1=1$ if $w>v_2$; and $s_2=-1$ and $j_2=v_2-w+1$ if $w \leq v_2$. 

By Lemma \ref{bigimportantlemma}, the vector $\alpha=(g_1, -y_k+l_1-1, s_1)$, the outer corner $C_1$ of size $(1,1,l_1)$, and the inner corner $N$ produce a bounded connected component; the vector $\beta=(g_2, -y_k+l_2-1, s_2)$, the outer corner $C_2$ of size $(1,1,l_2)$, and the inner corner $N$ produce another bounded connected component. 

%We can shift corners $c_1$ by a vector $\alpha=(\alpha_1, \alpha_2, \alpha_3)$, such that $\alpha_1<0, \alpha_2>0, \alpha_3>0$ such that  $c_1+ \alpha =Q$. And similarly we can shift $c_2$ by $\alpha'=(\alpha'_1, \alpha'_2, \alpha'_3)$, such that $\alpha'_1<0, \alpha'_2>0, \alpha'_3>0$ such that  $c_2+ \alpha' =Q$. By Lemma \ref{bigimportantlemma}, there are two bounded connected components. Contradiction.
%The only situation to be careful is when $w$ is smaller than the $z$-coordinate of either one of $c_1$ or $c_2$, or both. In this case, we shift by the vector which is exactly $+1$ in the $z$ direction, and shift in the $x$ direction by another appropriate positive number of boxes, so that the corner $c_1$ or $c_2$ is exactly having the same $x$ coordinate as the inner corner $(x_1+1, 1,w)$. Then the corner $c_1$ or $c_2$ along with a strip of boxes are a bounded connected component.
\end{proof}

\section{When the top layer is a rectangle}

Consider a minimized staircase $E$. In this section, we assume that all of the top $x,y,z$ layers have one unique outer corner. We prove a few lemmas that list all possible shapes of the bottom layer.

\begin{lemma}
Consider a staircase $E$ whose corresponding ideal $I \in \Hilb^d(\mathbb{A}^3)$ has tangent space dimension $\dim_k T_{nnp}(I)=1$. 

If the top layer of $E$ is a unique $1$ by $1$ box, then the bottom layer of $E$ has one unique inner corner.
\end{lemma}

\begin{proof}
By Lemma \ref{innercorbottomlayer}, the bottom layer must have at least one inner corner.

Assume for contradiction that the bottom layer has $2$ or more inner corners, then either one of the two inner corners must have coordinates $(n_1, n_2, 0)$ such that $n_1, n_2 \geq 1$. For every inner corner $(n_1, n_2, 0)$, consider the unique outer corner $(0,0, h_z)$ of size $(1,1,1)$, and the $nnp$ vector $(-n_1, -n_2, h_z)$, by Lemma \ref{bigimportantlemma}, there is a bounded connected component for each of the two inner corners, contradiction.
\end{proof}

\begin{lemma}
\label{two_inner_corners}
Consider a staircase $E$ whose corresponding ideal $I \in \Hilb^d(\mathbb{A}^3)$ has tangent space dimension $\dim_k T_{nnp}(I)=1$. 

Let the top layer of $E$ be a rectangle of size $a$ by $b$. If either $a= 0$ or $b=0$, then the bottom layer of $E$ has one unique inner corner. If $a,b \geq 1$, then the bottom layer of $E$ either has one unique inner corner, or has exactly two inner corners $q=(a+1, w, 0)$ and $t=(r, b+1, 0)$ for some $w <b, r <a$.

\end{lemma}

\begin{proof}
By Lemma \ref{innercorbottomlayer}, the bottom layer must have at least one inner corner.

Now suppose for contradiction that the bottom layer has two or more inner corners. The pictures below uses the lattice points to draw the staircase's top layer and bottom layer, and the outer corner of the top layer has coordinates $(a, b, h_z)$. 

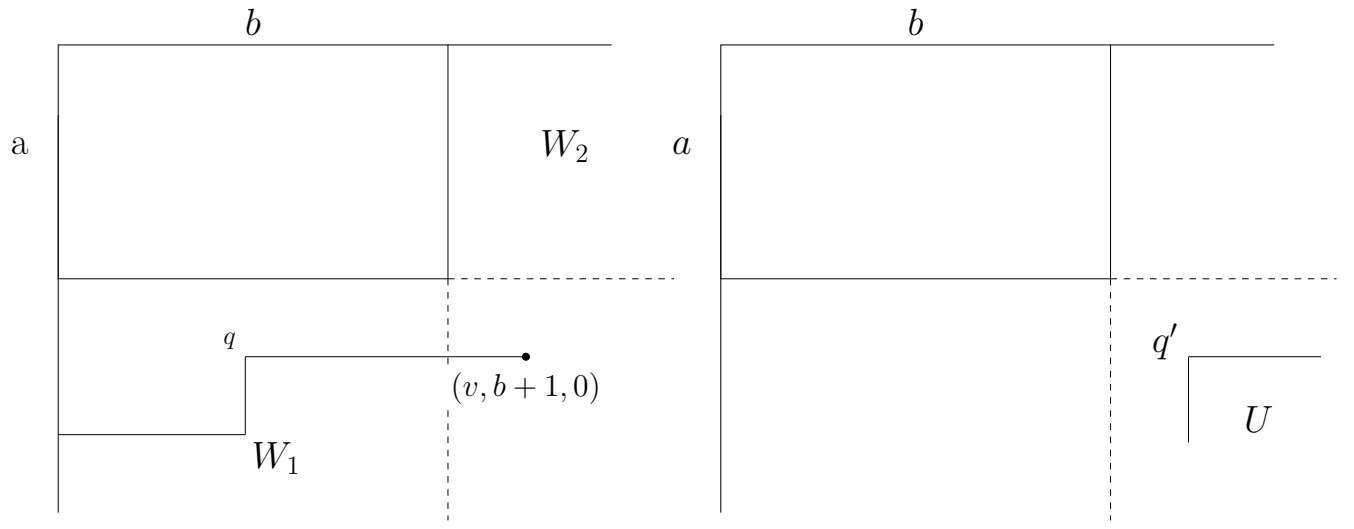
\resizebox{0.7\textwidth}{!}{%
\begin{circuitikz}
\tikzstyle{every node}=[font=\fontsize{18.2pt}{23.7pt}\selectfont]
\draw  (16.875,16) rectangle (23.125,12.25);
\draw [short] (23.125,16) -- (25.75,16);
\draw [short] (16.875,14.875) -- (16.875,8.5);
\draw [dashed] (23.125,12.25) -- (23.125,8.375);
\draw [dashed] (23.125,12.25) -- (26.75,12.25);
\draw [short] (16.875,9.75) -- (19.875,9.75);
\draw [short] (19.875,9.75) -- (19.875,11);
\draw [short] (19.875,11) -- (24.375,11);
\node [font=\fontsize{18.2pt}{23.7pt}\selectfont, fill={rgb,255:red,255; green,255; blue,255}, fill opacity=1, text opacity=1, inner xsep=0.080cm, inner ysep=0.085cm, rounded corners=0.020cm] at (20.375,9.375) {$W_1$};
\node [font=\fontsize{18.2pt}{23.7pt}\selectfont, fill={rgb,255:red,255; green,255; blue,255}, fill opacity=1, text opacity=1, inner xsep=0.080cm, inner ysep=0.085cm, rounded corners=0.020cm] at (25,14.375) {$W_2$};
\node [font=\fontsize{18.2pt}{23.7pt}\selectfont, fill={rgb,255:red,255; green,255; blue,255}, fill opacity=1, text opacity=1, inner xsep=0.080cm, inner ysep=0.085cm, rounded corners=0.020cm] at (16.25,14.375) {a};
\node [font=\fontsize{18.2pt}{23.7pt}\selectfont, fill={rgb,255:red,255; green,255; blue,255}, fill opacity=1, text opacity=1, inner xsep=0.080cm, inner ysep=0.085cm, rounded corners=0.020cm] at (20,16.375) {$b$};
\node [font=\fontsize{14.2pt}{18.5pt}\selectfont, fill={rgb,255:red,255; green,255; blue,255}, fill opacity=1, text opacity=1, inner xsep=0.080cm, inner ysep=0.085cm, rounded corners=0.020cm] at (24.375,10.5) {$(v, b+1,0)$};
\begin{scope}[transparency group, opacity=1]
\node at (24.375,11) [circ, color={rgb,255:red,0; green,0; blue,0}] {};
\end{scope}
\node [font=\fontsize{13.1pt}{17.0pt}\selectfont, fill={rgb,255:red,255; green,255; blue,255}, fill opacity=1, text opacity=1, inner xsep=0.080cm, inner ysep=0.085cm, rounded corners=0.000cm] at (19.625,11.25) {$q$};
\draw  (27.5,16) rectangle (33.75,12.25);
\draw [short] (33.75,16) -- (36.375,16);
\draw [short] (27.5,14.875) -- (27.5,8.5);
\draw [dashed] (33.75,12.25) -- (33.75,8.375);
\draw [dashed] (33.75,12.25) -- (37.375,12.25);
\node [font=\fontsize{18.2pt}{23.7pt}\selectfont, fill={rgb,255:red,255; green,255; blue,255}, fill opacity=1, text opacity=1, inner xsep=0.080cm, inner ysep=0.085cm, rounded corners=0.020cm] at (36.125,10) {$U$};
\node [font=\fontsize{18.2pt}{23.7pt}\selectfont, fill={rgb,255:red,255; green,255; blue,255}, fill opacity=1, text opacity=1, inner xsep=0.080cm, inner ysep=0.085cm, rounded corners=0.020cm] at (26.875,14.375) {$a$};
\node [font=\fontsize{18.2pt}{23.7pt}\selectfont, fill={rgb,255:red,255; green,255; blue,255}, fill opacity=1, text opacity=1, inner xsep=0.080cm, inner ysep=0.085cm, rounded corners=0.020cm] at (30.625,16.375) {$b$};
\draw (35,11) to[short] (35,9.625);
\draw (35,11) to[short] (37.125,11);
\node [font=\fontsize{18.2pt}{23.7pt}\selectfont, fill={rgb,255:red,255; green,255; blue,255}, fill opacity=1, text opacity=1, inner xsep=0.080cm, inner ysep=0.085cm, rounded corners=0.020cm] at (34.625,11.25) {$q'$};
\end{circuitikz}
}%

Assume that there is an inner corner $q=(v,w,0)$ in region $W_1$ so $w <b$, then by Lemma \ref{sidecomp}, there can only be at most one inner corner in region $W_1$, and the horizontal edge that emits from this inner corner must go straight to the right direction, until at least the point $(v, b+1,0)$. The inner corner $q$, outer corner $(a,b, h_z)$ of size $(1, b-w+2, 1)$ and the vector $(a-v, -1, h_z)$ satisfy Lemma \ref{bigimportantlemma}, so there is a bounded connected component. When the vertical edge that emits from $q$ has length equal to $2$ or more, then the inner corner $q$, outer corner $(a, b, h_z)$ of size $(2, b-w+2, 1)$ and the vector $(a-v-1, -1, h_z)$ satisfy Lemma \ref{bigimportantlemma}, so there is another bounded connected component, contradiction. So the length of the vertical edge that emits from $q$ must have length exactly $1$.

The same applies to region $W_2$ by symmetry, and there is at most one inner corner in $W_2$ with one bounded connected component. Denote again the inner corner in $W_2$ as $t=(r,s,0)$ where $r <a$, then the vector is $(-1, b-s, h_z)$

Assume that there is an inner corner $q'=(v',w',0)$ in region $U$ so $v' \geq a$ and $w' \geq b$, then the outer corner $(a, b, h_z)$ of size $(1,1,1)$, the inner corner $q'$ and the $nnp$ vector $(a-v', b'-w',h_z)$ satisfy Lemma \ref{bigimportantlemma}, so there is a bounded connected component. Note that if there is another inner corner in region $U$, then the same recipe of the inner corner, outer corner and $nnp$ vector again work to define another bounded connected component, so there can only be at most one inner corner in $U$. 

So we can conclude, if there are $2$ or more inner corners, and only one bounded connected component, then the two (or more) inner corners must produce exactly the same bounded connected component. In particular, the $nnp$ vectors must be the same. It is not possible for the two $nnp$ vectors for the inner corners in region $W_1$ and $U$ to be the same, because for $a-v-1=a-v'$, we must have that $v'=v+1$, but we must have $v'<v$ for $q$ and $q'$ to both be inner corners in their corresponding region of a partition. And neither can the two $nnp$ vectors for the inner corners in region $W_2$ and $U$ to be the same for the same reason. So we can only consider the possibility that the vectors for inner corners in region $W_1$ and $W_2$ to be the same, so we have $a-v=-1$ and $b-s=-1$. So the two inner corners in region $W_1$ and $W_2$ must satisfy that $q=(a+1, w, 0)$ and $t=(r, b+1, 0)$. Note that in order for these two inner corners to fit inside a partition, we must have $a \geq 1$ and $b \geq 1$. So in the case when either $a$ or $b$ is equal to $0$, where the top layer is either a single row or a single column of boxes, the bottom layer must have one unique inner corner.
\end{proof}

\begin{lemma}
\label{onerow}
Consider a staircase $E$ whose corresponding ideal $I \in \Hilb^d(\mathbb{A}^3)$ has tangent space dimension $\dim_k T_{nnp}(I)=1$.

If the top layer is a row of $b$ boxes, then the bottom layer can only be one of the three following types. 

\begin{enumerate}
\item An arbitrary partition with an unique size $1$ by $1$ inner corner.

\item A partition $(b+1,b+1,...,b+1,s,s,...,s)$ with one unique inner corner, and $s$ is any arbitrary number less than $b$. There can any arbitrary number of $s$'s and $b+1$'s.

\item A partition $(t,t,...,t,1,1,...,1)$ for any arbitrary number of $t$ and $1$'s, where $t > b$. 
\end{enumerate}
The picture shows the top layer and the bottom layer for each of the three types. The horizontal axis is the $y$ axis, and the vertical axis is the $x$ axis.

\

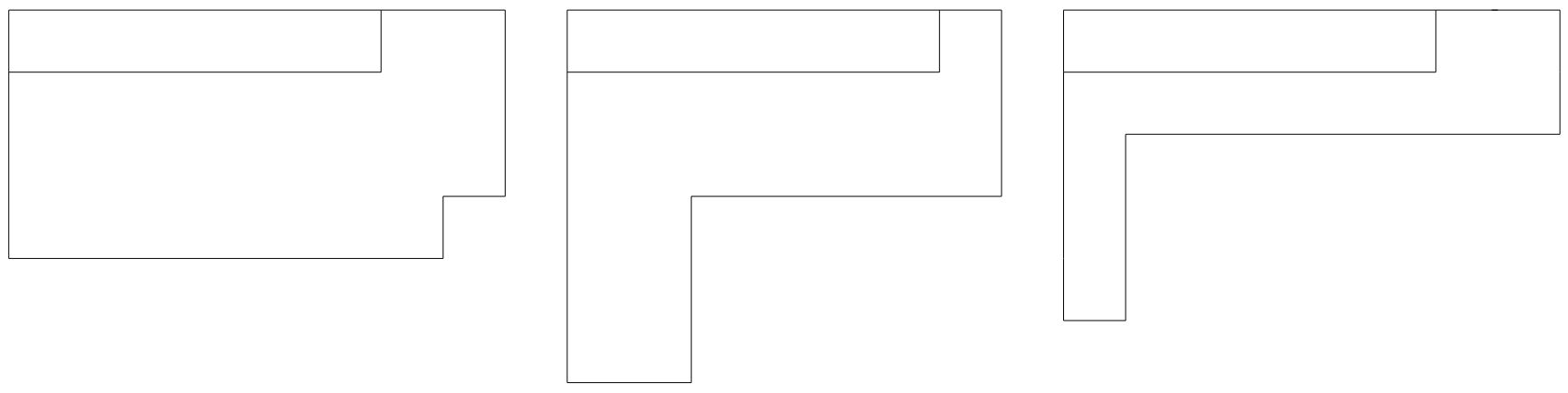
\resizebox{0.9\textwidth}{!}{%
\begin{circuitikz}
\tikzstyle{every node}=[font=\fontsize{18.2pt}{23.7pt}\selectfont]
\draw (-1.25,11.125) to[short] (6.25,11.125);
\draw (-1.25,11.125) to[short] (-1.25,9.875);
\draw (-1.25,9.875) to[short] (6.25,9.875);
\draw (6.25,11.125) to[short] (6.25,9.875);
\draw (-1.25,9.875) to[short] (-1.25,6.125);
\draw (-1.25,6.125) to[short] (7.5,6.125);
\draw (7.5,7.375) to[short] (7.5,6.125);
\draw (7.5,7.375) to[short] (8.75,7.375);
\draw (6.25,11.125) to[short] (8.75,11.125);
\draw (8.75,11.125) to[short] (8.75,7.375);
\draw (10,11.125) to[short] (17.5,11.125);
\draw (10,11.125) to[short] (10,9.875);
\draw (10,9.875) to[short] (17.5,9.875);
\draw (17.5,11.125) to[short] (17.5,9.875);
\draw (20,11.125) to[short] (27.5,11.125);
\draw (20,11.125) to[short] (20,9.875);
\draw (20,9.875) to[short] (27.5,9.875);
\draw (27.5,11.125) to[short] (27.5,9.875);
\draw (20,9.875) to[short] (20,6.125);
\draw (20,6.125) to[short] (20,4.875);
\draw (20,4.875) to[short] (21.25,4.875);
\draw (27.5,11.125) to[short] (28.75,11.125);
\draw (17.5,11.125) to[short] (18.75,11.125);
\draw (18.75,11.125) to[short] (18.75,7.375);
\draw (12.5,7.375) to[short] (12.5,3.625);
\draw (12.5,3.625) to[short] (10,3.625);
\draw (10,9.875) to[short] (10,3.625);
\draw (12.5,7.375) to[short] (18.75,7.375);
\draw (28.625,11.125) to[short] (30,11.125);
\draw (30,11.125) to[short] (30,9.875);
\draw (21.25,4.875) to[short] (21.25,8.625);
\draw (21.25,8.625) to[short] (30,8.625);
\draw (30,9.875) to[short] (30,8.625);
\end{circuitikz}
}%

\end{lemma}

\begin{proof}
By Lemma \ref{two_inner_corners}, there can only be one unique inner corner of the bottom layer. So the bottom layer is a partition that can be written as $(b',...,b',s,...,s)$, for some $b' > s$. We now discuss all the possibilities. 

When $b+2 \leq b'$ and $2 \leq s \leq b'-2$, then no matter how many $b'$ or $s$ we have, there are two or more bounded connected component as the following two cases. Denote the number of $b'$ as $n_{b'}$, and the inner corner of the bottom layer is denoted as $U=(n_{b'}, s, 0)$. 

\begin{enumerate}
\item If $s \leq b$, consider the outer corner $O=(0, b-1, h_z)$ of the top layer, the inner corner $U$ of size $(1, b'-s, 1)$, and the vectors $(-n_{b'}, -1, h_z)$ for the outer corner $O$ of size $(1,b-s+1,1)$, and the vector $(-n_{b'}, -2, h_z)$ for the outer corner $O$ of size $(1,b-s+2,1)$. Each of the two vectors induces a bounded connected component by Lemma \ref{bigimportantlemma}. Contradiction.

\item If $s > b$, consider the outer corner $O=(0, b-1, h_z)$ of the top layer, the inner corner $U$, and the vectors $(-n_{b'}, b-s-1, h_z)$ for the outer corner $O$ of size $(1,1,1)$, and the vector $(-n_{b'}, b-s-2, h_z)$ for the outer corner $O$ of size $(1,2,1)$. Each of the two vectors induces a bounded connected component by Lemma \ref{bigimportantlemma}. Contradiction.
\end{enumerate}

So for the bottom layer to have a unique inner corner, we negate the assumption: $b+2 \leq b'$ and $2 \leq s \leq b'-2$. Note that we always have $b' \geq b+1$ and $s < b'$ by construction. This gives the three possibilities. (1) $s+1 = b'$ and $b'$ is arbitrary. (2) $b'=b+1$ and $s$ is arbitrary, (3) $s=1$ and $b' \geq b+1$ is arbitrary. These are the three cases in the statement of the lemma.

\end{proof}

\begin{corollary}
\label{equiv_cond_onerow}
Denote the largest $x$ coordinate of the bottom layer as $h_x$, largest $y$ coordinate of the bottom layer as $h_y$. Denote the coordinate of the unique inner corner of the bottom layer as $(n_x,n_y,0)$. Denote the unique outer corner of the top layer as $(0, b, h_z)$, (note the shift of notation, for a simpler expression the corner is not $(0, b-1, h_z)$). The above Lemma \ref{onerow} is equivalent to the following statement. 

Consider a staircase $E$ whose corresponding ideal $I \in \Hilb^d(\mathbb{A}^3)$ has tangent space dimension $\dim_k T_{nnp}(I)=1$. If the top layer has a unique outer corner $(0, b, h_z)$, then the inner corner of the bottom layer is unique, and that $\min(n_y,h_y-b, h_y-n_y+1)=1$. 

By symmetry, if the top layer has a unique outer corner $(a, 0, h_z)$, we have $\min(n_x,h_x-a, h_x-n_x+1)=1$.   
\end{corollary}

\begin{lemma}
\label{rectangle}
Consider a staircase $E$ whose corresponding ideal $I \in \Hilb^d(\mathbb{A}^3)$ has tangent space dimension $\dim_k T_{nnp}(I)=1$. 

If the top layer is a rectangle of $a$ by $b$ boxes, with $a, b \geq 2$, then the bottom layer must be one of the four following cases, in addition to containing the $a$ by $b$ top layer.

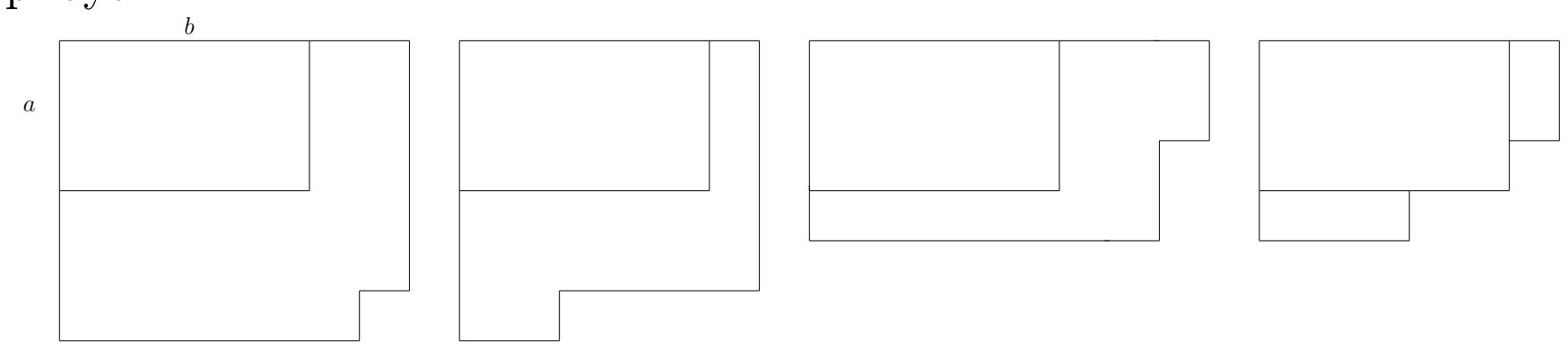
\resizebox{0.9\textwidth}{!}{%
\begin{circuitikz}
\tikzstyle{every node}=[font=\fontsize{18.2pt}{23.7pt}\selectfont]
\draw (1.25,8.625) to[short] (7.5,8.625);
\draw (7.5,8.625) to[short] (7.5,4.875);
\draw (7.5,4.875) to[short] (1.25,4.875);
\draw (1.25,8.625) to[short] (1.25,4.875);
\draw (7.5,8.625) to[short] (10,8.625);
\draw (10,8.625) to[short] (10,2.375);
\draw (10,2.375) to[short] (8.75,2.375);
\draw (8.75,2.375) to[short] (8.75,1.125);
\draw (8.75,1.125) to[short] (1.25,1.125);
\draw (1.25,4.875) to[short] (1.25,1.125);
\draw (31.25,8.625) to[short] (37.5,8.625);
\draw (37.5,8.625) to[short] (37.5,4.875);
\draw (37.5,4.875) to[short] (31.25,4.875);
\draw (31.25,8.625) to[short] (31.25,4.875);
\draw (37.5,8.625) to[short] (38.75,8.625);
\draw (38.75,8.625) to[short] (38.75,6.125);
\draw (38.75,6.125) to[short] (37.5,6.125);
\draw (31.25,4.875) to[short] (31.25,3.625);
\draw (31.25,3.625) to[short] (35,3.625);
\draw (35,4.875) to[short] (35,3.625);
\draw (11.25,8.625) to[short] (17.5,8.625);
\draw (17.5,8.625) to[short] (17.5,4.875);
\draw (17.5,4.875) to[short] (11.25,4.875);
\draw (11.25,8.625) to[short] (11.25,4.875);
\draw (20,8.625) to[short] (26.25,8.625);
\draw (26.25,8.625) to[short] (26.25,4.875);
\draw (26.25,4.875) to[short] (20,4.875);
\draw (20,8.625) to[short] (20,4.875);
\draw (26.25,8.625) to[short] (28.75,8.625);
\draw (17.5,8.625) to[short] (18.75,8.625);
\draw (18.75,8.625) to[short] (18.75,3.625);
\draw (18.75,3.625) to[short] (18.75,2.375);
\draw (18.75,2.375) to[short] (13.75,2.375);
\draw (13.75,2.375) to[short] (13.75,1.125);
\draw (13.75,1.125) to[short] (11.25,1.125);
\draw (11.25,4.875) to[short] (11.25,1.125);
\draw (20,5) to[short] (20,3.625);
\draw (20,3.625) to[short] (27.5,3.625);
\draw (28.625,8.625) to[short] (30,8.625);
\draw (27.375,3.625) to[short] (28.75,3.625);
\draw (28.75,6.125) to[short] (28.75,3.625);
\draw (28.75,6.125) to[short] (30,6.125);
\draw (30,8.625) to[short] (30,6.125);

\node [font=\fontsize{18.2pt}{23.7pt}\selectfont, fill={rgb,255:red,255; green,255; blue,255}, fill opacity=1, text opacity=1, inner xsep=0.080cm, inner ysep=0.085cm, rounded corners=0.020cm] at (4.5,9) {$b$};

\node [font=\fontsize{18.2pt}{23.7pt}\selectfont, fill={rgb,255:red,255; green,255; blue,255}, fill opacity=1, text opacity=1, inner xsep=0.080cm, inner ysep=0.085cm, rounded corners=0.020cm] at (0.5,7) {$a$};
\end{circuitikz}
}%

\begin{enumerate}

\item \label{rec:1} A partition with a unique size $1$ by $1$ inner corner.

\item \label{rec:2} A partition $(b+1, b+1, ..., b+1, s)$ for any arbitrary $ 1 \leq s \leq b$, and arbitrary number of $b+1$ parts.

\item \label{rec:3} A partition $(t, t,t,...,t, t-1, t-1, ..., t-1)$ for any arbitrary $t \geq b+1$ and the total number of parts of both $t$ and $t+1$ is $a+1$.

\item \label{rec:4} Any arbitrary partition with exactly $a+1$ parts and maximum length of the parts being $b+1$.

\end{enumerate}
\end{lemma}

\begin{proof}
Assume first that there is one unique inner corner, then the bottom partition is denoted as $(t,t,...,t,s,s,...,s)$ for an arbitrary number of $t$ and $s$ parts. By assumption $s < t$ and $b < t$. Denote the number of rows with a total of $t$ boxes as $n_t$, the number of rows with a total of $s$ boxes as $n_s$ so the inner corner of the bottom layer is $N:=(n_t, s, 0)$. Define numbers $m_a:=\max(1, n_t+1-a)$ and 
$m_b=\max(1,s+1-b)$.

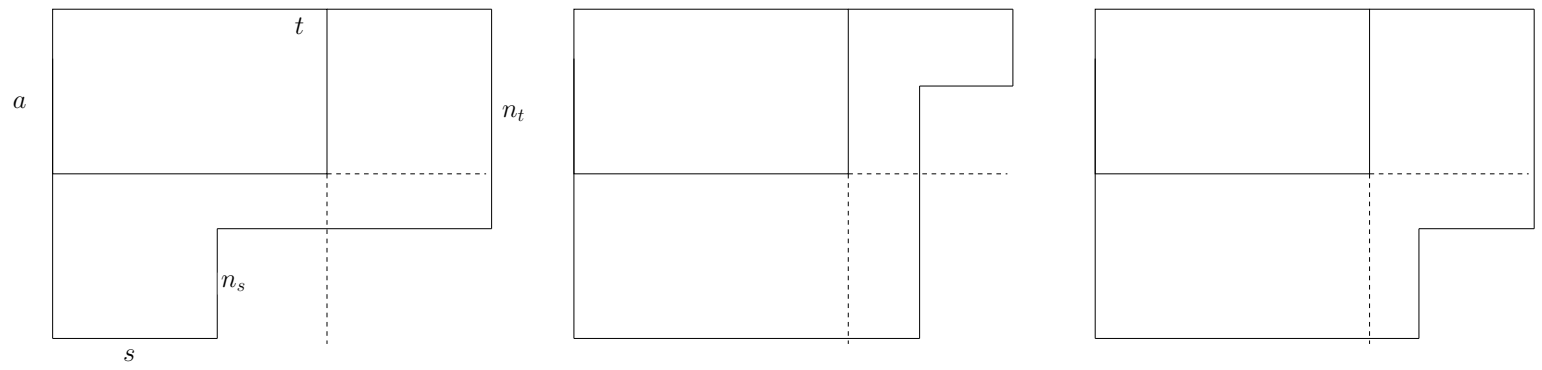
\resizebox{0.8\textwidth}{!}{%
\begin{circuitikz}
\tikzstyle{every node}=[font=\fontsize{18.2pt}{23.7pt}\selectfont]
\draw  (5,15.625) rectangle (11.25,11.875);
\draw [short] (11.25,15.625) -- (15,15.625);
\draw [short] (5,14.5) -- (5,8.125);
\draw [dashed] (11.25,11.875) -- (11.25,8);
\draw [dashed] (11.25,11.875) -- (14.875,11.875);
\draw [short] (15,15.625) -- (15,10.625);
\draw [short] (15,10.625) -- (8.75,10.625);
\draw [short] (8.75,10.625) -- (8.75,8.125);
\draw [short] (8.75,8.125) -- (5,8.125);
\draw  (28.75,15.625) rectangle (35,11.875);
\draw [short] (35,15.625) -- (38.75,15.625);
\draw [short] (28.75,14.5) -- (28.75,8.125);
\draw [dashed] (35,11.875) -- (35,8);
\draw [dashed] (35,11.875) -- (38.625,11.875);
\draw [short] (38.75,15.625) -- (38.75,10.625);
\draw [short] (38.75,10.625) -- (36.125,10.625);
\draw [short] (36.125,10.625) -- (36.125,8.125);
\draw [short] (36.125,8.125) -- (28.75,8.125);
\draw  (16.875,15.625) rectangle (23.125,11.875);
\draw [short] (23.125,15.625) -- (26.875,15.625);
\draw [short] (16.875,14.5) -- (16.875,8.125);
\draw [dashed] (23.125,11.875) -- (23.125,8);
\draw [dashed] (23.125,11.875) -- (26.75,11.875);
\draw [short] (26.875,15.625) -- (26.875,13.875);
\draw [short] (26.875,13.875) -- (24.75,13.875);
\draw [short] (24.75,13.875) -- (24.75,8.125);
\draw [short] (24.75,8.125) -- (16.875,8.125);

\node [font=\fontsize{18.2pt}{23.7pt}\selectfont, fill={rgb,255:red,255; green,255; blue,255}, fill opacity=1, text opacity=1, inner xsep=0.080cm, inner ysep=0.085cm, rounded corners=0.020cm] at (4.25,13.5) {$a$};
\node [font=\fontsize{18.2pt}{23.7pt}\selectfont, fill={rgb,255:red,255; green,255; blue,255}, fill opacity=1, text opacity=1, inner xsep=0.080cm, inner ysep=0.085cm, rounded corners=0.020cm] at (7.875,16.375) {$b$};
\node [font=\fontsize{18.2pt}{23.7pt}\selectfont, fill={rgb,255:red,255; green,255; blue,255}, fill opacity=1, text opacity=1, inner xsep=0.080cm, inner ysep=0.085cm, rounded corners=0.005cm] at (6.75,7.75) {$s$};
\node [font=\fontsize{18.2pt}{23.7pt}\selectfont, fill={rgb,255:red,255; green,255; blue,255}, fill opacity=1, text opacity=1, inner xsep=0.080cm, inner ysep=0.085cm, rounded corners=0.000cm] at (10.625,15.25) {$t$};
\node [font=\fontsize{18.2pt}{23.7pt}\selectfont, fill={rgb,255:red,255; green,255; blue,255}, fill opacity=1, text opacity=1, inner xsep=0.080cm, inner ysep=0.085cm, rounded corners=0.020cm] at (15.5,13.25) {$n_t$};
\node [font=\fontsize{18.2pt}{23.7pt}\selectfont, fill={rgb,255:red,255; green,255; blue,255}, fill opacity=1, text opacity=1, inner xsep=0.080cm, inner ysep=0.085cm, rounded corners=0.020cm] at (9.125,9.375) {$n_s$};
\end{circuitikz}
}%

Suppose for the purpose of contradiction that $t \geq b+2$ and $s \leq  t-2$. 
%Because $a \geq 2$,  there are $2$ or more vectors for $2$ or more bounded connected component, given as follows. 

Define the number $l_1:=1$ if $n_t \geq a$, and $l_1=a-n_t+1$ if $n_t <a$. Define the number $l_2:=1$ if $s \geq b$, and $l_2=b-s+1$ if $s <b$. 

By Lemma \ref{bigimportantlemma}, the outer corner the top layer $(a-1,b-1, h_z)$ of size $(l_1, l_2,1)$, the inner corner $N$, the vector $(-m_a, -m_b,h_z)$ create abounded connected component. Similarly, because $t-s \geq 2$ and $b \geq 2$, the outer corner $(a-1,b-1, h_z)$ of size $(l_1, l_2+1,1)$, the inner corner $N$, the vector $(-m_a,-m_b-1,h_z)$ create another bounded connected component. 

So we conclude that we must have the negation of the conditions, $t \geq b+2$ and $s \leq  t-2$. Because $t>b$, the negation of $t \geq b+2$ is $t=b+1$. Because and $s \leq t-1$, the negation of $s \leq  t-2$ is $s=t-1$. So we must have the following condition:
\begin{align}
\label{conditioni}
\text{Either } t=b+1 \text{ or } s=t-1. 
\tag{i}
\end{align}%
By symmetry of exchanging $x$ and $y$, we must also have the symmetry of the above condition:
\begin{align}
\label{conditionii}
\text{Either } n_t+n_s=a+1 \text{ or } n_s=1.
\tag{ii} 
\end{align}%
Condition (\ref{conditioni}) and condition (\ref{conditionii}) must hold simultaneously. 

Combining all the possible conditions, we have the following possibilities:

When $t=b+1$ and $n_t+n_s=a+1$, we obtain (\ref{rec:1}) with the size of the bottom layer being $a+1$ by $b+1$; when $s=t-1$ and $n_s=1$, we obtain every possible size of (\ref{rec:1}).

When $t=b+1$ and $n_s=1$, we obtain (\ref{rec:2}).

When $s=t-1$ and $n_t+n_s=a+1$, we obtain (\ref{rec:3}). 

If there are two or more inner corners, by Lemma \ref{two_inner_corners}, there can only be exactly two inner corners, $(a, w, 0)$ and $(r, b,0)$ for some $w<b$ and $r<a$ (note that the notation is shifted from $a+1$ to $a$ because the top layer in this lemma has coordinates $(a-1, b-1, h_z)$ instead of $(a, b, h_z)$ as in Lemma \ref{two_inner_corners}). This gives case (\ref{rec:4}). 

%there cannot be more than one inner corner to the left (or top) of each dashed line, otherwise there are two or more outer corners and there are more than two bounded connected component. So there can be maximum one inner corner to the left or top of the dashed line. If there's an inner corner to the left or top of the dashed lines, then each inner corner creates a bounded connected component unless we have (4), classified as follows. Denote the inner corner as $(o_1, o_2)$, and the outer corner $(a,b)$ satisfy the condition of 

%So the bottom layer can only have exactly two inner corners, with three outer corners, and that exactly the outer corner in the middle coincide with the outer corner of the top layer, so that the conditions of Lemma\ref{sidecomp} is not satisfied.

\end{proof}

For the case (\ref{rec:1})(\ref{rec:2})(\ref{rec:3}) in the above Lemma \ref{rectangle}, there is a unique inner corner in the bottom layer. Denote the coordinate of the unique inner corner as $(n_x, n_y,0)$. 
\begin{corollary}
\label{equiv_cond_rec}
Denote the unique outer corner of the top layer as $(a,b,h_z)$, (note the shift of notation, for a simpler expression the corner is not $(a-1, b-1, h_z)$). Denote the maximum $x$ and $y$ coordinates as $h_x, h_y$. The above Lemma \ref{rectangle} is equivalent the condition of both $\min(h_y-b, h_y-n_y+1)=1$ and that $\min(h_x-a, h_x-n_x+1)=1$.   
\end{corollary}

\section{When top layer has more than two outer corners}

In this section, we classify the staircases whose top layer has at least two
outer corners. We will show that, under the assumption that
$I \setminus (I - \alpha)$ has a unique bounded connected component for all $nnp$ vectors $\alpha$, the staircase $E$ must be what we call an opera house.

%The argument proceeds by analyzing the possible shapes of the bottom layer. We first show that the bottom layer must be exactly the region occupied by the top layer together with two additional rectangles each adjacent to the $x$ and $y$ coordinate axes. In the three dimensional staircase, the staircase is the top layer stacking on itself many times, forming the main part of the opera house, along with two wings of the opera house adjacent to the coordinate axis. Finally, we prove that the wings must be rectangular cuboids and provide the conditions on the numerical values of each part of the opera house.

\begin{definition}
\label{def_opera_house}
An opera house is the following staircase. It is a set-theoretic union of three sets of boxes: the main part and the two wings part. See the picture below.

The \textbf{main part} of the opera house is a pile of the same two dimensional partition $M$ stacking on top of itself $h_z+1$ layers in the $z$ direction, if the maximum $z$ coordinate is $h_z$. Denote $h_z$ as the \textbf{height} of the opera house. For the partition $M$, say its maximum $x$ and $y$ coordinates are $m_x, m_y$ accordingly. Denote the number of boxes in the $(m_x+1)$th row as \textbf{length of the  xth shoulder} $s_x$, and the number of boxes in the $(m_y+1)$th column as \textbf{length of the  yth shoulder} $s_y$.

There are also two \textbf{wings} of the opera house, each of the shape of a rectangular cuboid. Denote the wing near the $x$ axis as \textbf{the left wing}, and the wing near the $y$ axis as \textbf{the right wing}.

The unique outer corners of the left and right wings are denoted as $l=(l_1+ m_x ,l_2,l_3)$ of the left wing, and $r=(r_1,r_2+m_y,r_3)$ of the right wing. The wings satisfy $l_2 \leq s_x$, and $r_1 \leq s_y$. The height of the left wing is $l_3$, and the height of the right wing is $r_3$.

The numerical values of $l$ and $r$ are of the $3$ following types: (See Example \ref{wingstype} below.)
\begin{enumerate}
\item two spread-out rows: $l_2=r_1=1$, and $l_1,r_2$ arbitrary.
\item two folded rows: $l_1=r_2=1$, and $l_2,r_1$ arbitrary.
\item an unique box and an arbitrary rectangle: $l=(1,1,0)$ and $r$ arbitrary, or $r=(1,1,0)$ and $l$ arbitrary.
\end{enumerate}
At least one of the left and right wings must have height either $0$ or $h_z-1$. i.e. either $l_3=0, h_z-1$, or $r_3=0,h_z-1$. (Recall that height is the maximum $z$ coordinate, not the number of layers.)

\end{definition}

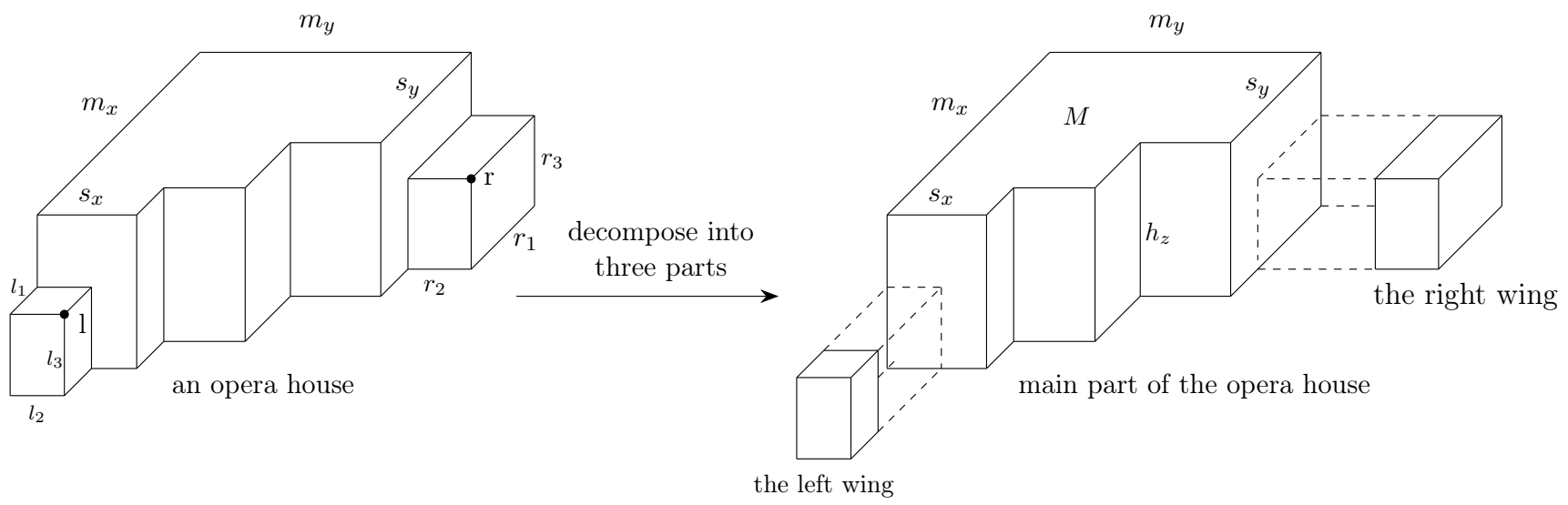
\resizebox{0.9\textwidth}{!}{%
\begin{circuitikz}
\tikzstyle{every node}=[font=\fontsize{10.2pt}{13.3pt}\selectfont]
\draw (3.75,11) to[short] (7.5,11);
\draw (7.5,11) to[short] (6.25,9.75);
\draw (6.25,9.75) to[short] (5,9.75);
\draw (5,9.75) to[short] (4.375,9.125);
\draw (4.375,9.125) to[short] (3.25,9.125);
\draw (3.75,11) to[short] (1.5,8.75);
\draw (1.5,8.75) to[short] (2.875,8.75);
\draw (4.625,9.375) to[short] (4.625,9.375);
\draw (3.25,9.125) to[short] (2.875,8.75);
\draw (3.25,9.125) to[short] (3.25,7);
\draw (4.375,9.125) to[short] (4.375,7);
\draw (5,9.75) to[short] (5,7.625);
\draw (1.5,8.75) to[short] (1.5,6.625);
\draw (2.875,8.75) to[short] (2.875,6.625);
\draw (6.25,9.75) to[short] (6.25,7.625);
\draw (7.5,11) to[short] (7.5,8.875);
\draw (7.5,8.875) to[short] (6.25,7.625);
\draw (6.25,7.625) to[short] (5,7.625);
\draw (5,7.625) to[short] (4.375,7);
\draw (4.375,7) to[short] (3.25,7);
\draw (3.25,7) to[short] (2.875,6.625);
\draw (1.5,6.625) to[short] (2.875,6.625);
\node [font=\fontsize{10.8pt}{14.1pt}\selectfont, inner xsep=0.080cm, inner ysep=0.085cm, rounded corners=0.020cm] at (5.75,6.375) {main part of the opera house};
\draw (0.25,6.5) to[short] (0.25,5.375);
\draw (1,6.5) to[short] (1,5.375);
\draw (0.25,6.5) to[short] (1,6.5);
\draw (0.25,5.375) to[short] (1,5.375);
\draw (0.625,6.875) to[short] (0.25,6.5);
\draw (1.375,6.875) to[short] (1,6.5);
\draw (1.375,6.875) to[short] (1.375,5.75);
\draw (1.375,5.75) to[short] (1,5.375);
\draw (0.625,6.875) to[short] (1.375,6.875);
\draw [dashed] (0.625,6.875) -- (1.5,7.75);
\draw [dashed] (1.375,5.75) -- (2.25,6.625);
\draw [dashed] (1.375,6.875) -- (2.25,7.75);
\draw [dashed] (2.25,7.75) -- (2.25,6.625);
\draw [dashed] (1.5,7.75) -- (2.25,7.75);
\draw [dashed] (1.375,6.5) -- (1.625,6.75);
\node [font=\fontsize{10.2pt}{13.3pt}\selectfont, inner xsep=0.080cm, inner ysep=0.085cm, rounded corners=0.020cm] at (0.625,5) {the left wing};
\draw (10,10.125) to[short] (9.125,9.25);
\draw (10,10.125) to[short] (10,8.875);
\draw (9.125,9.25) to[short] (9.125,8);
\draw (10,8.875) to[short] (9.125,8);
\draw (9.125,10.125) to[short] (10,10.125);
\draw (8.25,8) to[short] (9.125,8);
\draw (8.25,9.25) to[short] (9.125,9.25);
\draw (9.125,10.125) to[short] (8.25,9.25);
\draw (8.25,9.25) to[short] (8.25,8);
\draw [dashed] (6.625,8) -- (8.25,8);
\draw [dashed] (6.75,9.25) -- (8.375,9.25);
\draw [dashed] (7.5,10.125) -- (9.125,10.125);
\draw [dashed] (6.625,9.25) -- (6.625,8.125);
\draw [dashed] (6.625,9.25) -- (7.5,10.125);
\draw [dashed] (7.5,8.875) -- (8.25,8.875);
\node [font=\fontsize{11.9pt}{15.5pt}\selectfont, inner xsep=0.080cm, inner ysep=0.085cm, rounded corners=0.020cm] at (9.5,7.625) {the right wing};
\draw (-8,11) to[short] (-4.25,11);
\draw (-4.25,11) to[short] (-5.5,9.75);
\draw (-5.5,9.75) to[short] (-6.75,9.75);
\draw (-6.75,9.75) to[short] (-7.375,9.125);
\draw (-7.375,9.125) to[short] (-8.5,9.125);
\draw (-8,11) to[short] (-10.25,8.75);
\draw (-10.25,8.75) to[short] (-8.875,8.75);
\draw (-7.125,9.375) to[short] (-7.125,9.375);
\draw (-8.5,9.125) to[short] (-8.875,8.75);
\draw (-8.5,9.125) to[short] (-8.5,7);
\draw (-7.375,9.125) to[short] (-7.375,7);
\draw (-6.75,9.75) to[short] (-6.75,7.625);
\draw (-10.25,8.75) to[short] (-10.25,7.75);
\draw (-8.875,8.75) to[short] (-8.875,6.625);
\draw (-5.5,9.75) to[short] (-5.5,7.625);
\draw (-4.25,11) to[short] (-4.25,10.125);
\draw (-5.125,8) to[short] (-5.5,7.625);
\draw (-5.5,7.625) to[short] (-6.75,7.625);
\draw (-6.75,7.625) to[short] (-7.375,7);
\draw (-7.375,7) to[short] (-8.5,7);
\draw (-8.5,7) to[short] (-8.875,6.625);
\draw (-9.5,6.625) to[short] (-8.875,6.625);
\draw (-10.625,7.375) to[short] (-10.625,6.25);
\draw (-9.875,7.375) to[short] (-9.875,6.25);
\draw (-10.625,7.375) to[short] (-9.875,7.375);
\draw (-10.625,6.25) to[short] (-9.875,6.25);
\draw (-10.25,7.75) to[short] (-10.625,7.375);
\draw (-9.5,7.75) to[short] (-9.875,7.375);
\draw (-9.5,7.75) to[short] (-9.5,6.625);
\draw (-9.5,6.625) to[short] (-9.875,6.25);
\draw (-10.25,7.75) to[short] (-9.5,7.75);
\draw (-3.375,10.125) to[short] (-4.25,9.25);
\draw (-3.375,10.125) to[short] (-3.375,8.875);
\draw (-4.25,9.25) to[short] (-4.25,8);
\draw (-3.375,8.875) to[short] (-4.25,8);
\draw (-4.25,10.125) to[short] (-3.375,10.125);
\draw (-5.125,8) to[short] (-4.25,8);
\draw (-5.125,9.25) to[short] (-4.25,9.25);
\draw (-4.25,10.125) to[short] (-5.125,9.25);
\draw (-5.125,9.25) to[short] (-5.125,8);
\node [font=\fontsize{10.8pt}{14.1pt}\selectfont, inner xsep=0.080cm, inner ysep=0.085cm, rounded corners=0.020cm] at (-7.125,6.375) {an opera house};
\draw [-{Stealth[scale=1.5]}, ] (-3.625,7.625) -- (0,7.625);
\node [font=\fontsize{10.8pt}{14.1pt}\selectfont, inner xsep=0.080cm, inner ysep=0.085cm, rounded corners=0.020cm] at (-1.625,8.5) {decompose into};
\node [font=\fontsize{10.8pt}{14.1pt}\selectfont, inner xsep=0.080cm, inner ysep=0.085cm, rounded corners=0.020cm] at (-1.625,8) {three parts};
\node [font=\fontsize{10.8pt}{14.1pt}\selectfont, inner xsep=0.080cm, inner ysep=0.085cm, rounded corners=0.000cm] at (-9.625,7.25) {l};
\node [font=\fontsize{10.8pt}{14.1pt}\selectfont, inner xsep=0.080cm, inner ysep=0.085cm, rounded corners=0.000cm] at (-4,9.25) {r};
\begin{scope}[transparency group, opacity=1]
\node at (-9.875,7.375) [circ, color={rgb,255:red,0; green,0; blue,0}] {};
\end{scope}
\begin{scope}[transparency group, opacity=1]
\node at (-4.25,9.25) [circ, color={rgb,255:red,0; green,0; blue,0}] {};
\end{scope}
\node [font=\fontsize{10.8pt}{14.1pt}\selectfont, inner xsep=0.080cm, inner ysep=0.085cm, rounded corners=0.020cm] at (2.375,10.25) {$m_x$};
\node [font=\fontsize{10.8pt}{14.1pt}\selectfont, inner xsep=0.080cm, inner ysep=0.085cm, rounded corners=0.020cm] at (5.375,11.375) {$m_y$};
\node [font=\fontsize{10.8pt}{14.1pt}\selectfont, inner xsep=0.080cm, inner ysep=0.085cm, rounded corners=0.013cm] at (2.25,9) {$s_x$};
\node [font=\fontsize{10.8pt}{14.1pt}\selectfont, inner xsep=0.080cm, inner ysep=0.085cm, rounded corners=0.013cm] at (6.625,10.5) {$s_y$};
\node [font=\fontsize{10.8pt}{14.1pt}\selectfont, inner xsep=0.080cm, inner ysep=0.085cm, rounded corners=0.020cm] at (-9.375,10.25) {$m_x$};
\node [font=\fontsize{10.8pt}{14.1pt}\selectfont, inner xsep=0.080cm, inner ysep=0.085cm, rounded corners=0.020cm] at (-6.375,11.375) {$m_y$};
\node [font=\fontsize{10.8pt}{14.1pt}\selectfont, inner xsep=0.080cm, inner ysep=0.085cm, rounded corners=0.013cm] at (-5.125,10.5) {$s_y$};
\node [font=\fontsize{10.8pt}{14.1pt}\selectfont, inner xsep=0.080cm, inner ysep=0.085cm, rounded corners=0.013cm] at (-9.5,9) {$s_x$};
\node [font=\fontsize{8.0pt}{10.4pt}\selectfont, inner xsep=0.080cm, inner ysep=0.085cm, rounded corners=0.000cm] at (-10.5,7.75) {$l_1$};
\node [font=\fontsize{10.8pt}{14.1pt}\selectfont, inner xsep=0.080cm, inner ysep=0.085cm, rounded corners=0.010cm] at (-3.5,8.375) {$r_1$};
\node [font=\fontsize{8.0pt}{10.4pt}\selectfont, inner xsep=0.080cm, inner ysep=0.085cm, rounded corners=0.000cm] at (-10.25,6) {$l_2$};
\node [font=\fontsize{8.0pt}{10.4pt}\selectfont, inner xsep=0.080cm, inner ysep=0.085cm, rounded corners=0.000cm] at (-10,6.75) {$l_3$};
\node [font=\fontsize{10.2pt}{13.3pt}\selectfont, inner xsep=0.080cm, inner ysep=0.085cm, rounded corners=0.003cm] at (-4.75,7.75) {$r_2$};
\node [font=\fontsize{10.2pt}{13.3pt}\selectfont, inner xsep=0.080cm, inner ysep=0.085cm, rounded corners=0.003cm] at (-3.125,9.5) {$r_3$};
\node [font=\fontsize{10.2pt}{13.3pt}\selectfont, inner xsep=0.080cm, inner ysep=0.085cm, rounded corners=0.020cm] at (4.125,10.125) {$M$};
\node [font=\fontsize{10.2pt}{13.3pt}\selectfont, inner xsep=0.080cm, inner ysep=0.085cm, rounded corners=0.019cm] at (5.25,8.5) {$h_z$};
\end{circuitikz}
}%

%Denote the \textbf{inner corners} of the bottom layer of the opera house as $x^{d_1}y^{e_1}, ..., x^{d_k}y^{e_k}$, where $d_1<...<d_k$ and $e_1>...>e_k$ (See picture below). \Ytodo{add picture, add ref}

\begin{example}
\label{wingstype}
We show some examples of the three type wings:
\begin{enumerate}
\item two spread-out rows:

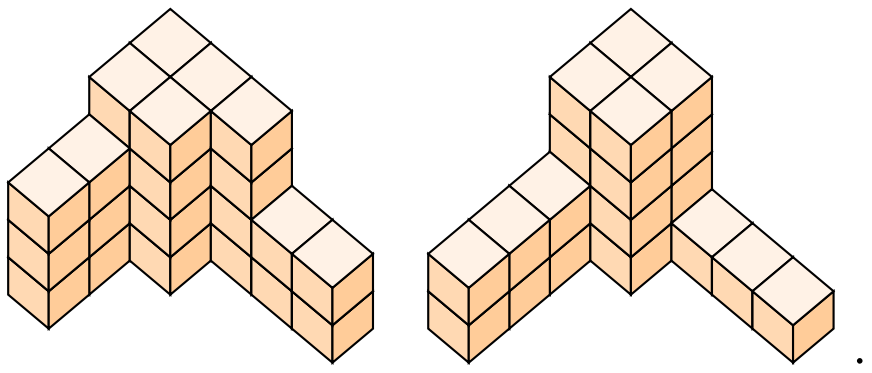
\monomialIdealPicture{{4,4,3,3},{4,4},{4},{2},{2}} \monomialIdealPicture{{4,4,2,2,2},{4,4},{1},{1},{1}}.
\item two folded rows:

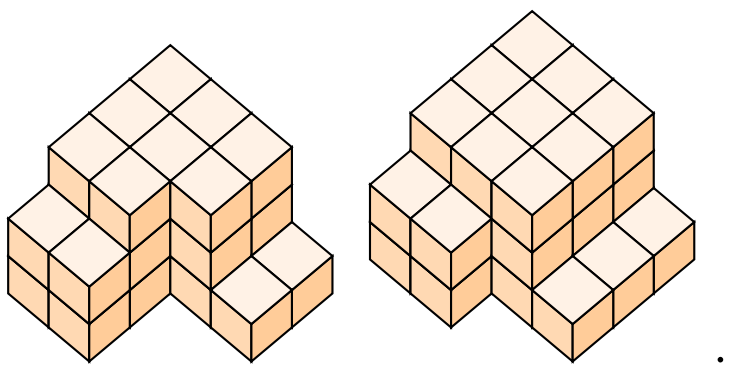
\monomialIdealPicture{{3,3,3,2},{3,3,3,2},{3,3},{1,1}}\monomialIdealPicture{{3,3,3,2},{3,3,3,2},{3,3,3},{1,1,1}}.
\item an unique box and an arbitrary rectangle:

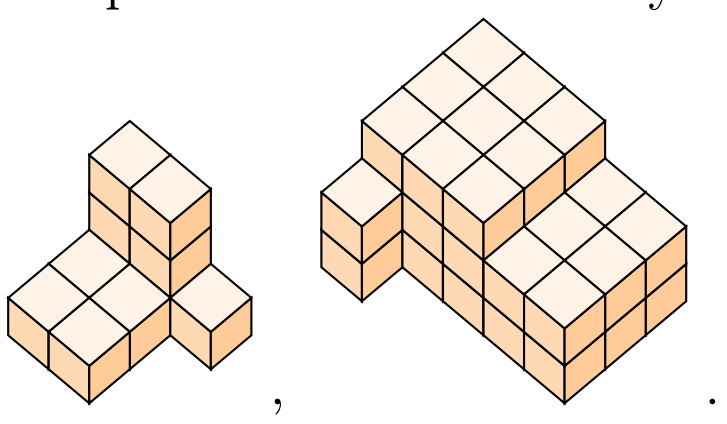
\monomialIdealPicture{{3,1,1},{3,1,1},{1}}, \monomialIdealPicture{{3,3,3,2},{3,3,3},{3,3,3},{2,2,2},{2,2,2}}.
\end{enumerate}

\end{example}

\begin{definition}
\label{big_pic}
Project the top layer of $E$ to the $XY$ plane. We use a box to denote an element of the staircases. Let
$(x_1,y_1),\ldots,(x_k,y_k)$ be its outer corners, ordered so that
$x_1>\cdots>x_k$ and $y_1<\cdots<y_k$. Divide the plane into the regions $A$, $B$, $C$, and $V$ as follows.

\begin{figure}[!ht]
\centering
\resizebox{0.6\textwidth}{!}{%
\begin{circuitikz}
\tikzstyle{every node}=[font=\fontsize{18.2pt}{23.7pt}\selectfont]
\draw (-0.125,13.625) to[short] (10.75,13.625);
\draw (10.75,13.625) to[short] (10.75,9.875);
\draw (10.75,9.875) to[short] (7.5,9.875);
\draw (7.5,9.875) to[short] (7.5,7.375);
\draw (7.5,7.375) to[short] (3.75,7.375);
\draw (3.75,7.375) to[short] (3.75,5.75);
\draw (3.75,5.75) to[short] (0,5.75);
\draw (0,13.625) to[short] (-1.375,13.625);
\draw (-1.375,13.625) to[short] (-1.375,4.875);
\draw (-1.375,4.875) to[short] (-1.375,3);
\draw (-1.375,3) to[short] (0,3);
\draw (0,5.75) to[short] (0,3);
\draw (11.25,13.125) to[short] (11.25,9.375);
\draw (11.25,9.375) to[short] (8,9.375);
\draw (8,9.375) to[short] (8,6.875);
\draw (8,6.875) to[short] (4.25,6.875);
\draw (4.25,6.875) to[short] (4.25,5.25);
\draw (4.25,5.25) to[short] (0.5,5.25);
\draw (-0.875,2.5) to[short] (0.5,2.5);
\draw (0.5,5.25) to[short] (0.5,2.5);
\draw [dashed] (0.5,2.5) -- (0.5,0);
\draw [dashed] (0.375,2.5) -- (7.625,2.5);
\draw [dashed] (4.25,5.25) -- (4.25,0);
\draw [dashed] (4.25,5.25) -- (9.75,5.25);
\draw [dashed] (8,6.875) -- (8,0);
\draw [dashed] (8,6.875) -- (14.75,6.875);
\draw [dashed] (11.25,9.375) -- (11.25,2.5);
\draw [dashed] (11.25,9.375) -- (16.625,9.375);
\draw [dashed] (11.25,13.125) -- (17.25,13.125);
\draw [dashed] (-0.875,2.5) -- (-0.875,0);
\draw [dashed] (7.375,2.5) -- (13,2.5);
\draw [dashed] (9.625,5.25) -- (15.25,5.25);
\draw [dashed] (11.25,2.5) -- (11.25,0);
\draw (-0.875,3) to[short] (-0.875,2.5);
\draw (11.25,13.125) to[short] (10.75,13.125);
\node [font=\fontsize{18.2pt}{23.7pt}\selectfont, inner xsep=0.080cm, inner ysep=0.085cm, rounded corners=0.020cm] at (2.625,10.5) {$A$};
\node [font=\fontsize{10.2pt}{13.3pt}\selectfont, inner xsep=0.080cm, inner ysep=0.085cm, rounded corners=0.000cm] at (4,7.125) {$B$};
\node [font=\fontsize{10.2pt}{13.3pt}\selectfont, inner xsep=0.080cm, inner ysep=0.085cm, rounded corners=0.020cm] at (2.5,3.625) {$C_1$};
\node [font=\fontsize{10.2pt}{13.3pt}\selectfont, inner xsep=0.080cm, inner ysep=0.085cm, rounded corners=0.020cm] at (6,6) {$C_1$};
\draw (-1.375,3) to[short] (-1.375,2.375);
\draw (-1.375,2.5) to[short] (-0.875,2.5);
\draw (10.75,13.625) to[short] (11.25,13.625);
\draw (11.25,13.625) to[short] (11.25,13.125);
\node [font=\fontsize{10.2pt}{13.3pt}\selectfont, inner xsep=0.080cm, inner ysep=0.085cm, rounded corners=0.020cm] at (9.625,8.125) {$C_1$};
\node [font=\fontsize{10.2pt}{13.3pt}\selectfont, inner xsep=0.080cm, inner ysep=0.085cm, rounded corners=0.020cm] at (6.25,3.75) {$C_2$};
\node [font=\fontsize{10.2pt}{13.3pt}\selectfont, inner xsep=0.080cm, inner ysep=0.085cm, rounded corners=0.020cm] at (2.375,1) {$C_2$};
\node [font=\fontsize{10.2pt}{13.3pt}\selectfont, inner xsep=0.080cm, inner ysep=0.085cm, rounded corners=0.020cm] at (9.625,6.125) {$C_2$};
\node [font=\fontsize{10.2pt}{13.3pt}\selectfont, inner xsep=0.080cm, inner ysep=0.085cm, rounded corners=0.020cm] at (13.125,8.125) {$C_2$};
\node [font=\fontsize{10.2pt}{13.3pt}\selectfont, inner xsep=0.080cm, inner ysep=0.085cm, rounded corners=0.020cm] at (13.125,6) {$C_3$};
\node [font=\fontsize{10.2pt}{13.3pt}\selectfont, inner xsep=0.080cm, inner ysep=0.085cm, rounded corners=0.020cm] at (9.625,3.875) {$C_3$};
\node [font=\fontsize{10.2pt}{13.3pt}\selectfont, inner xsep=0.080cm, inner ysep=0.085cm, rounded corners=0.020cm] at (6,1.125) {$C_3$};
\node [font=\fontsize{10.2pt}{13.3pt}\selectfont, inner xsep=0.080cm, inner ysep=0.085cm, rounded corners=0.020cm] at (9.625,1.125) {$C_3$};
\node [font=\fontsize{10.2pt}{13.3pt}\selectfont, inner xsep=0.080cm, inner ysep=0.085cm, rounded corners=0.020cm] at (13.25,3.75) {$C_3$};
\node [font=\fontsize{10.2pt}{13.3pt}\selectfont, inner xsep=0.080cm, inner ysep=0.085cm, rounded corners=0.020cm] at (-0.25,1.125) {$V_2$};
\node [font=\fontsize{10.2pt}{13.3pt}\selectfont, inner xsep=0.080cm, inner ysep=0.085cm, rounded corners=0.020cm] at (13.25,11.25) {$V_1$};
\node [font=\fontsize{10.2pt}{13.3pt}\selectfont, inner xsep=0.080cm, inner ysep=0.085cm, rounded corners=0.020cm] at (-1.125,2.75) {$\beta_2$};
\node [font=\fontsize{10.2pt}{13.3pt}\selectfont, inner xsep=0.080cm, inner ysep=0.085cm, rounded corners=0.020cm] at (11,13.375) {$\beta_1$};
\end{circuitikz}
}%
%\caption{Your Caption}
%\label{fig:my_label}
\end{figure}

The region $A$ is the part occupied by the top layer.

The region $B$ is the region immediately above and to the right of the
shifted top layer under translation by $(1,1)$. For each outer corner $(x_i,y_i)$ of the top layer, the corresponding outer corner of the translation is $(x_i+1,y_i+1)$.

Divide the region $C$ into three types of regions as follows. For each outer
corner of the region $B$, draw a vertical dashed line and a horizontal dashed
line, as shown in the figure. Let $C_2$ be the union of the rectangles that
share a corner with $B$, and let $C_1$ be the union of the rectangles
lying between $B$ and $C_2$. The region $C_3$ consists of the remaining parts of $C$.

Let $V=V_1\cup V_2$ denote the two unbounded rectangular regions
adjacent to the coordinate axes.

Denote the two boxes $\beta_1:=(x_1+1, 0, 0)$ and $\beta_2:=(0,y_k+1,0)$. 

\end{definition}

\begin{remark}
For the vector $\alpha=(-1,-1,h_z)$ and staircase $E$ with the top layer as drawn, the corresponding infinite component $I_\infty$ is the L shaped wall of width $1$. And the bottom layer of $I_\infty$ consists of one row of boxes adjacent to $\beta_1$ and one column of boxes adjacent to $\beta_2$. Note that in the picture above, we labeled every region except for the bottom layer of $I_\infty$.
\end{remark}

\begin{lemma}
\label{V1V2}
The bottom layer of a minimized staircase $E$ contains the boxes $\beta_1$ and $\beta_2$. 
\end{lemma}
\begin{proof}
Suppose that $\beta_1$ is not contained in the bottom layer. Since the bottom layer is a partition, no box to the right
of $\beta_1$ can be contained in the bottom layer. Similarly, if $\beta_2$ is not contained in the bottom layer,
then no box below $\beta_2$ can be contained in the bottom layer.

It follows that the two outer corners $(x_1,y_1,h_z)$ and
$(x_k,y_k,h_z)$ are dominating corners of the staircase $E$, contradiction.
\end{proof}

\begin{corollary}
\label{beta12}
After translating the staircase $E$ by $-\alpha=(1,1, -h_z)$, the finite component $I_{fin}$ cannot be connected to the infinite wall component of $I_\infty$ by any lattice path.
\end{corollary}

\begin{proof}
Consider the boxes in the translated staircase $E-\alpha$ that have non-negative $x,y,z$ coordinate, and are not in the original staircase $E$, which is exactly the boxes in region $A$, so the finite component $I_{fin}$ is a subset of the boxes in region $A$. Because the boxes $\beta_1, \beta_2 \in E$, they separate the translated top layer from the infinite wall component. 

Now we show that $S_B$ is not connected to $I_\infty$. Because $I_{fin}$ is a subset in the shifted top layer, there are no boxes in $I_{fin}$ to the top, right, or front of them. Because $\beta_1$ and $\beta_2$ are in the bottom layer, and the rest of the left and back side of the boundary of $B$ is shared with region $A$, the left and back side of $B$ is not in $I_{fin}$. The bottom side of $S_B$ has negative $z$ coordinates and cannot be in $I_{fin}$.
\end{proof}

\begin{lemma}
\label{reg_B}
If the edge of the bottom layer enters region $B$, then it must eventually leaves the region $B$ and cuts out some boxes that are in region $B$ but not in the bottom layer of the staircase. These boxes, denoted as $S_B$ is a bounded connected component.
\end{lemma}

\begin{proof}
First we show that $S_B$ is in $I_{fin}$. Shift the staircase $E$ by the vector $(-1,-1, h_z)$. The finite component $I_{fin}$ only consists of boxes in the top layer of $E$, because the rest of the staircase $E$ is shifted to have negative $z$ coordinates and therefore not contained in $I_{fin}$. We also know that the boxes in the shifted top layer that are also in region $A$ is not contained in $I_{fin}$ because they are contained in $E$. Finally, $S_B$ is contained in $I_{fin}$ because they are not in the bottom layer by construction and therefore in $I$. 

Now we show that $S_B$ is not connected to $I_\infty$. Because $S_B$ is in the shifted top layer, there are no boxes in $I_{fin}$ to the top, right, or front of them. Because $S_B$ is exactly "cut out" by the boundaries of the bottom layer: by the edge of an inner corner in $C_1$ that cuts horizontally into $B$ and then cutting out of $B$, and because $\beta_1$ and $\beta_2$ are in the bottom layer, the left, and back side of $S_B$ is in the bottom layer of $E$ and cannot be in $I_{fin}$. The bottom side of $S_B$ has negative $z$ coordinates and cannot be in $I_{fin}$. So $S_B$ is completely separated from $I_\infty$ because there's no box on any side of $S_B$ that connects it to $I_\infty$.
\end{proof}

\begin{lemma}
\label{Bcompconnected}
If there is an inner corner in region $B$, then it must follow the shared boundary between region $A$ and $B$ for some number of steps, and completely leaves the region $B$ and never intersect it again. Equivalently, the boxes in region $B$ outside of the bottom layer, if exists, must be path connected as a unique connected component.   
\end{lemma}

\begin{proof}
Again consider the vector $\alpha=(-1, -1, h_z)$, and the finite component $I_{\textit{fin}}$ consists exactly of the monomials in $B$ that is not in the bottom layer. By Corollary \ref{beta12}, the finite component is not connected to the infinite component $I_\infty$ by any lattice path. So region $B$ removing the boxes in the bottom layer, if exists, must be a unique connected component. 
\end{proof}

\begin{lemma}
\label{nocornerinV}
The bottom layer of a minimized staircase $E$ with minimal tangent space dimension cannot have any inner corners in region $V_1 \cup V_2$.
\end{lemma}

\begin{proof}
Without loss of generality, assume that there is an inner corner in the region $V_2$. The proof where there is an inner corner in region $V_1$ is similar. 

When $y_1 \geq 2$, subdivide the region $V_2$ into $V_2^1$ and $V_2^2$, where $V_2^2$ is a column of boxes of width $1$. When $y_1=1$, the entire region $V_2$ itself is a column of boxes of width $1$.

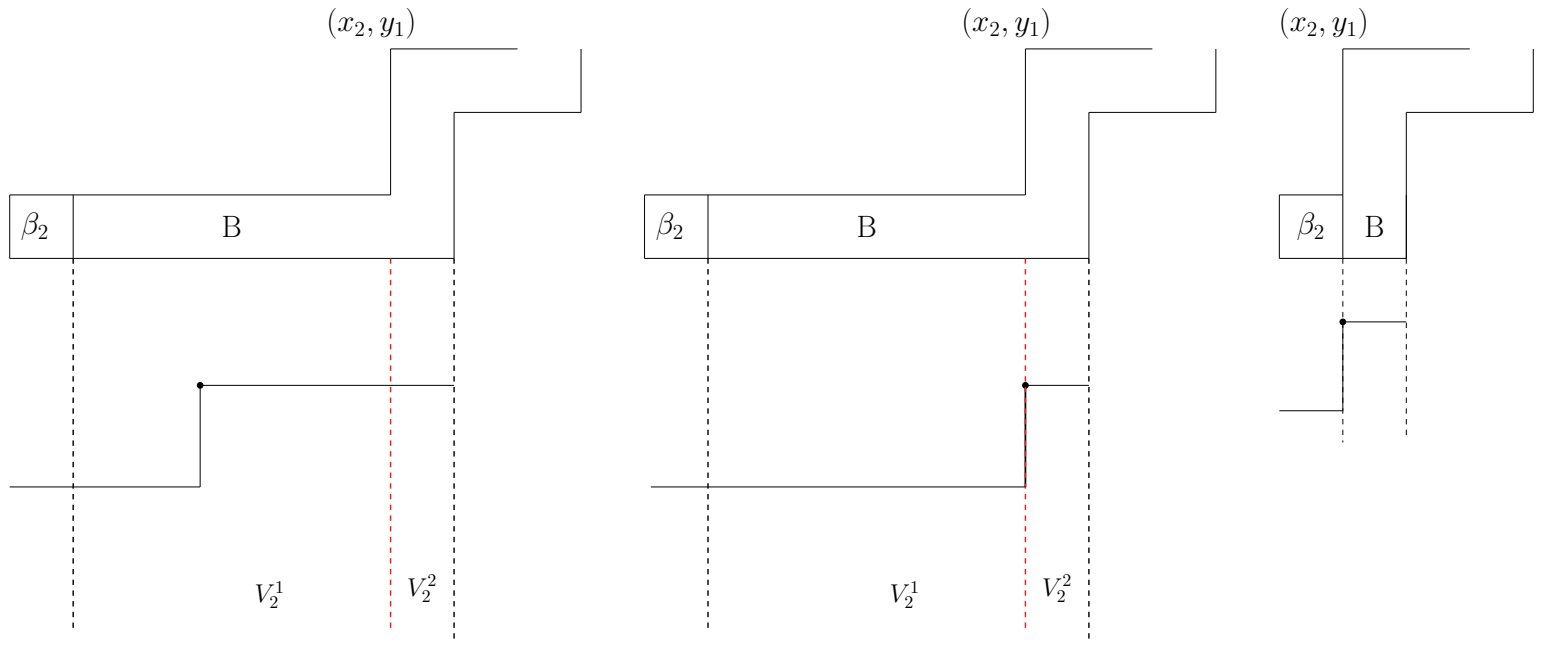
\resizebox{0.8\textwidth}{!}{%
\begin{circuitikz}
\tikzstyle{every node}=[font=\fontsize{18.2pt}{23.7pt}\selectfont]
\draw (2.5,11.125) to[short] (2.5,9.875);
\draw (2.5,9.875) to[short] (3.75,9.875);
\draw (3.75,9.875) to[short] (3.75,11.125);
\draw (2.5,11.125) to[short] (10,11.125);
\draw (3.75,9.875) to[short] (11.25,9.875);
\draw (10,11.125) to[short] (10,14);
\draw (11.25,9.875) to[short] (11.25,12.75);
\draw [line width=0.8pt, dashed] (3.75,9.875) -- (3.75,2.5);
\draw [line width=0.8pt, dashed] (11.25,9.875) -- (11.25,2.375);
\begin{scope}[transparency group, opacity=1]
\node at (6.25,7.375) [circ, color={rgb,255:red,0; green,0; blue,0}] {};
\end{scope}
\draw (6.25,7.375) to[short] (11.25,7.375);
\draw (6.25,7.375) to[short] (6.25,5.375);
\node [font=\fontsize{18.2pt}{23.7pt}\selectfont, fill={rgb,255:red,255; green,255; blue,255}, fill opacity=1, text opacity=1, inner xsep=0.080cm, inner ysep=0.085cm, rounded corners=0.020cm] at (3,10.5) {$\beta_2$};
\node [font=\fontsize{18.2pt}{23.7pt}\selectfont, fill={rgb,255:red,255; green,255; blue,255}, fill opacity=1, text opacity=1, inner xsep=0.080cm, inner ysep=0.085cm, rounded corners=0.020cm] at (6.875,10.5) {B};
\node [font=\fontsize{14.2pt}{18.5pt}\selectfont, fill={rgb,255:red,255; green,255; blue,255}, fill opacity=1, text opacity=1, inner xsep=0.080cm, inner ysep=0.085cm, rounded corners=0.020cm] at (7.625,3.25) {$V_2^1$};
\draw [ color={rgb,255:red,237; green,28; blue,28}, draw opacity=1, line width=0.8pt, dashed] (10,9.875) -- (10,2.5);
\draw (6.25,5.375) to[short] (2.5,5.375);
\node [font=\fontsize{14.2pt}{18.5pt}\selectfont, fill={rgb,255:red,255; green,255; blue,255}, fill opacity=1, text opacity=1, inner xsep=0.080cm, inner ysep=0.085cm, rounded corners=0.020cm] at (10.625,3.375) {$V_2^2$};
\draw (27.5,11.125) to[short] (28.75,11.125);
\draw (28.75,9.875) to[short] (30,9.875);
\draw (30,9.875) to[short] (30,11.125);
\draw (27.5,11.125) to[short] (27.5,9.875);
\draw (27.5,9.875) to[short] (28.75,9.875);
\draw (28.75,11.125) to[short] (28.75,9.875);
\draw [dashed] (28.75,9.875) -- (28.75,6.25);
\draw [dashed] (30,9.875) -- (30,6.375);
\draw (15,11.125) to[short] (15,9.875);
\draw (15,9.875) to[short] (16.25,9.875);
\draw (16.25,9.875) to[short] (16.25,11.125);
\draw (15,11.125) to[short] (22.5,11.125);
\draw (16.25,9.875) to[short] (23.75,9.875);
\draw [line width=0.8pt, dashed] (16.25,9.875) -- (16.25,2.5);
\draw [line width=0.8pt, dashed] (23.75,9.875) -- (23.75,2.375);
\begin{scope}[transparency group, opacity=1]
\node at (22.5,7.375) [circ, color={rgb,255:red,0; green,0; blue,0}] {};
\end{scope}
\draw [ line width=1pt](22.5,7.375) to[short] (22.5,5.375);
\node [font=\fontsize{18.2pt}{23.7pt}\selectfont, fill={rgb,255:red,255; green,255; blue,255}, fill opacity=1, text opacity=1, inner xsep=0.080cm, inner ysep=0.085cm, rounded corners=0.020cm] at (15.5,10.5) {$\beta_2$};
\node [font=\fontsize{18.2pt}{23.7pt}\selectfont, fill={rgb,255:red,255; green,255; blue,255}, fill opacity=1, text opacity=1, inner xsep=0.080cm, inner ysep=0.085cm, rounded corners=0.020cm] at (19.375,10.5) {B};
\node [font=\fontsize{14.2pt}{18.5pt}\selectfont, fill={rgb,255:red,255; green,255; blue,255}, fill opacity=1, text opacity=1, inner xsep=0.080cm, inner ysep=0.085cm, rounded corners=0.020cm] at (20.125,3.25) {$V_2^1$};
\draw [ color={rgb,255:red,237; green,28; blue,28}, draw opacity=1, line width=0.8pt, dashed] (22.5,9.875) -- (22.5,2.5);
\draw (22.5,5.375) to[short] (15.125,5.375);
\node [font=\fontsize{14.2pt}{18.5pt}\selectfont, fill={rgb,255:red,255; green,255; blue,255}, fill opacity=1, text opacity=1, inner xsep=0.080cm, inner ysep=0.085cm, rounded corners=0.020cm] at (23.125,3.375) {$V_2^2$};
\draw (22.5,7.375) to[short] (23.75,7.375);
\node [font=\fontsize{18.2pt}{23.7pt}\selectfont, fill={rgb,255:red,255; green,255; blue,255}, fill opacity=1, text opacity=1, inner xsep=0.080cm, inner ysep=0.085cm, rounded corners=0.020cm] at (28.125,10.5) {$\beta_2$};
\node [font=\fontsize{18.2pt}{23.7pt}\selectfont, fill={rgb,255:red,255; green,255; blue,255}, fill opacity=1, text opacity=1, inner xsep=0.080cm, inner ysep=0.085cm, rounded corners=0.020cm] at (29.375,10.5) {B};
\begin{scope}[transparency group, opacity=1]
\node at (28.75,8.625) [circ, color={rgb,255:red,0; green,0; blue,0}] {};
\end{scope}
\draw (28.75,8.625) to[short] (30,8.625);
\draw (28.75,8.625) to[short] (28.75,6.875);
\draw (28.75,6.875) to[short] (27.5,6.875);
\draw (11.25,12.75) to[short] (13.75,12.75);
\draw (10,14) to[short] (12.5,14);
\draw (13.75,12.75) to[short] (13.75,14);
\node [font=\fontsize{18.2pt}{23.7pt}\selectfont, inner xsep=0.080cm, inner ysep=0.085cm, rounded corners=0.020cm] at (9.625,14.5) {$(x_2,y_1)$};
\draw (22.5,11.125) to[short] (22.5,14);
\draw (23.75,9.875) to[short] (23.75,12.75);
\draw (23.75,12.75) to[short] (26.25,12.75);
\draw (22.5,14) to[short] (25,14);
\draw (26.25,12.75) to[short] (26.25,14);
\node [font=\fontsize{18.2pt}{23.7pt}\selectfont, inner xsep=0.080cm, inner ysep=0.085cm, rounded corners=0.020cm] at (22.125,14.5) {$(x_2,y_1)$};
\draw (28.75,11.125) to[short] (28.75,14);
\draw (30,9.875) to[short] (30,12.75);
\draw (30,12.75) to[short] (32.5,12.75);
\draw (28.75,14) to[short] (31.25,14);
\draw (32.5,12.75) to[short] (32.5,14);
\node [font=\fontsize{18.2pt}{23.7pt}\selectfont, inner xsep=0.080cm, inner ysep=0.085cm, rounded corners=0.020cm] at (28.375,14.5) {$(x_2,y_1)$};
\end{circuitikz}
}%

By Lemma \ref{sidecomp}, we discuss the restriction of the inner corners by their location, and list all the possibilities.

\begin{enumerate}
 \item Suppose there is an inner corner $q=(v,w,0)$ in the region $V_2^1$. The vertical edge that emanates from $q$ that initially goes downwards must eventually turn left, creating an outer corner. Note that this edge cannot create more outer corners along the way because of Lemma \ref{sidecomp}, so this edge must go straight left. Say it has the $x$ coordinate $c_1$. The horizontal edge that emanates from $q$ that initially goes to the right must not create any extra outer corner $O=(o_1,o_2,0)$ to the left of, or exactly at the location of, the red dashed line, by Lemma \ref{sidecomp}, where $o_2 \leq y_1$. So this edge that emanates from the inner corner $q$ must go to the right for at least $y_1-w+2$ of steps of boxes, from the corner $q$ with $y$ coordinate $w$ to the box with $y$ coordinate $y_1+1$.
 
\item Suppose there is an inner corner $q=(v,y_1+1,0)$ exactly at the red dashed line, the boundary between $V_2^1$ and $V_2^2$. Again, the vertical edge that emanates from $q$ that initially goes downwards must eventually turn left and go straight left, creating a unique outer corner with $y$ coordinate $y_1$. Denote its $x$ coordinate as $c_1$, then the outer corner is $(c_1, y_1,0)$. The horizontal edge that emanates from $q$ must go to the right for at least $1$ step.

\item Suppose the width of $V_2$ is $1$, then an inner corner $q$ in $V_2$ must have coordinates $q=(v, 1, 0)$ for some $v$. The horizontal edge that emanates from $q$ must go to the right for at least $1$ step. The vertical edge that emanates from $q$ that initially goes downwards must eventually turn left, creating an outer corner. Again, denote its $x$ coordinate as $c_1$, then the outer corner is $(c_1, y_1,0)$
\end{enumerate}

These are the three possible cases of the location of the inner corner $q$. In all of these three cases above, we find two triples of an inner corner, outer corner and a vector that satisfies Lemma \ref{bigimportantlemma} as follows.

\begin{enumerate}
\item Consider the inner corner $(x_2+1, y_1+1,o_3)$ of $E$ for some $o_3$ such that $(x_2+1, y_1+1,o_3+1) \in E$ and $(x_2+1, y_1+1,o_3+1) \in I$. We know that this is an inner corner, because the top layer does not contain the box $(x_2+1, y_1+1, h_z)$, but the bottom layer contains the box $(x_2+1, y_1+1, 0)$. And the top layer contains $(x_2, y_1+1, h_z)$, $(x_2+1, y_1, h_z)$ and $(x_2, y_1, h_z)$.

Outer corner $(c_1, w, 0)$ of size $(1,1,1)$. For all the three possible cases above, $c_1 > x_2+1$, and $w < y_1$.  

So the vector $(c_1-x_2-1,w-y_1-1, -o_3)$ of signature $pnn$, with the inner and outer corners above, form a bounded connected component.

\item Consider the inner corner $(0, y_k+1, o_4)$ for some $o_4$ such that $(0, y_k+1, o_4-1) \in E$ and $(0, y_k+1, o_4) \in I$. This box is not shown in the picture of this lemma, but it is the box labeled as $\beta_1$ in the bigger picture in Definition \ref{big_pic}. We know that this is an inner corner, because the top layer does not contain the box $(0, y_k+1, h_z)$, but the bottom layer contains the box $\beta_2=(0, y_k+1, 0)$. And the top layer contains $(0, y_k, h_z)$.

Consider the outer corner $(c_1, w, 0)$ of size $(1,1,1)$. For all the three possible cases above $w < y_k+1$. So the vector $(c_1, w-(y_k+1), -o_4)$ of signature $pnn$, the other corner and the inner corners listed above, form a bounded connected component.
\end{enumerate}

\end{proof}

\begin{lemma}
\label{nocornerinC}
Consider a minimized staircase $E$ with minimum tangent space dimension. If the top layer of $E$ has at least $2$ outer corners, then the bottom layer of $E$ cannot have any inner corners in region $C$.
\end{lemma}

\begin{proof}
Suppose, for contradiction, that the bottom layer has an inner corner $q=(v,w,0)$ in region $C$. We consider the three regions $C_1$, $C_2$, and $C_3$ separately.

\begin{enumerate}
\item 
First suppose that $q$ lies in $C_3$. Since the top layer has at least two outer corners, there are at least two indices $i$ of outer corners of the top layer $(x_i, y_i, h_z)$ such that $v > x_i$ and $w > y_i$.

For each such $i$, the outer corner $(x_i,y_i,h_z)$, the inner corner $q$, and the $nnp$ vector $\alpha=(v-x_i, w-y_i, h_z)$ satisfy the conditions of Lemma \ref{bigimportantlemma}. Hence each such pair gives a bounded connected component, contradiction. There should be no inner corners in $C_3$.

%If $(v,w)$ is far enough away from region $A,B$ and is in region $C_3$, then $v \geq x_i'$ and $w \geq y_i'$ for at least two $i \in \{1,...,k\}$. In this case, every single corner $(x_i, y_i)$ and $(v,w)$ satisfy the conditions of Lemma \ref{bigimportantlemma}, so there are at least $2$ connected components. So there cannot be any inner corner of the bottom partition in region $C_3$.

\item
Suppose that $q=(v,w,0)$ lies in $C_2$, including the upper and left boundaries of $C_2$.. Let $T_i$ be the rectangle in $C_2$ containing the outer corner $(x_i+1,y_i+1)$ of region $B$. So we have that $x_i +1\leq v$ and $y_i +1 \leq w$. 

The outer corner  $(x_i,y_i, h_z)$ of size $(1,1,1)$, the inner corner $q$, and the $nnp$ vector $\alpha=(x_i-v,y_i-w, h_z)$ satisfy the conditions of Lemma \ref{bigimportantlemma}, so there is a bounded connected component. 

We note that there cannot be two inner corners inner corner $q$ in $C_2$. If so, then each of them, together with the corresponding outer corner(s) and $nnp$ vectors, would satisfy the conditions of Lemma \ref{bigimportantlemma}, giving at least two bounded connected components.%Indeed, if there were two inner corners in the same rectangle $T_i$, then each of them, together with the outer corner $(x_i,y_i,h_z)$, would satisfy the conditions of Lemma \ref{bigimportantlemma}, giving at least two bounded connected components.

Let $U,D,L,R$ denote the four up, down, left, and right sides of $T_i$. Since $q$ is an inner corner, two edges emanate from $q$, one in the right direction and the other in the downward direction. We show that at least one of these two edges leave $T_i$ produces another bounded connected component.

Suppose first that the downward edge leaves $T_i$ through $D$. Since there can be no additional inner corner in $C_2$, this edge must continue straight downward until it leaves $T_i$. Let $(t_1^D,t_2^D,0)$ denote the box to the right of the edge where it leaves $T_i$. Then
$t_1^D=x_{i-1}+1$ and
$t_2^D=w$.

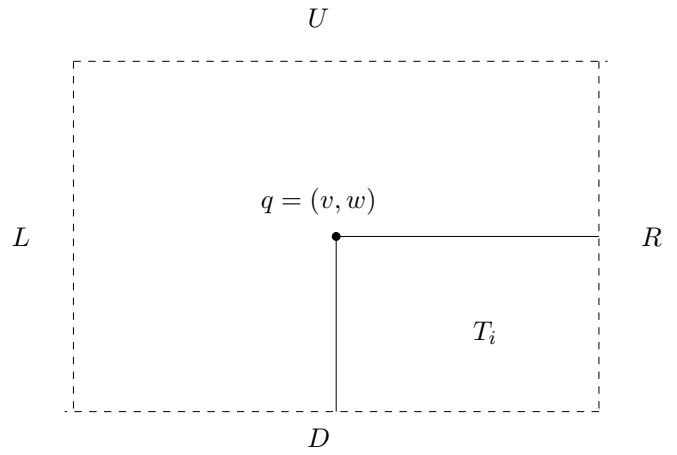
\resizebox{0.35\textwidth}{!}{%
\begin{circuitikz}
\tikzstyle{every node}=[font=\fontsize{10.8pt}{14.1pt}\selectfont]
\draw (6.25,9.875) to[short] (6.25,7.375);
\begin{scope}[transparency group, opacity=1]
\node at (6.25,9.875) [circ, color={rgb,255:red,0; green,0; blue,0}] {};
\end{scope}
\node [font=\fontsize{10.8pt}{14.1pt}\selectfont, inner xsep=0.080cm, inner ysep=0.085cm, rounded corners=0.020cm] at (6,10.375) {$q=(v,w)$};
\draw [dashed] (2.5,12.375) -- (10.125,12.375);
\draw [dashed] (10,12.375) -- (10,7.375);
\draw [dashed] (10,7.375) -- (2.375,7.375);
\draw [dashed] (2.5,12.375) -- (2.5,7.375);
\node [font=\fontsize{10.8pt}{14.1pt}\selectfont, inner xsep=0.080cm, inner ysep=0.085cm, rounded corners=0.020cm] at (6,13) {$U$};
\node [font=\fontsize{10.8pt}{14.1pt}\selectfont, inner xsep=0.080cm, inner ysep=0.085cm, rounded corners=0.020cm] at (6,7) {$D$};
\node [font=\fontsize{10.8pt}{14.1pt}\selectfont, inner xsep=0.080cm, inner ysep=0.085cm, rounded corners=0.020cm] at (1.75,9.875) {$L$};
\node [font=\fontsize{10.8pt}{14.1pt}\selectfont, inner xsep=0.080cm, inner ysep=0.085cm, rounded corners=0.020cm] at (10.75,9.875) {$R$};
\node [font=\fontsize{10.8pt}{14.1pt}\selectfont, inner xsep=0.080cm, inner ysep=0.085cm, rounded corners=0.020cm] at (8.375,8.5) {$T_i$};
\draw (6.25,9.875) to[short] (10,9.875);
\end{circuitikz}
}%

Because $q \in C_2$, we have that $w < y_{i-1}$. Consequently, the inner corner $q=(v,w,0)$ of size $(x_{i-1}+2-v, 1, 1)$, the outer corner $(x_{i-1},y_{i-1},h_z)$ of size $(x_{i-1}+2-v,1,1)$, and the $nnp$ vector $(-1, y_{i-1}-w, h_z)$ satisfy the conditions of Lemma \ref{bigimportantlemma}. This gives a second bounded connected component, contradiction.

Similarly, if the rightward edge leaves $T_i$ through $R$, then inner corner $q=(v,w,0)$ of size $(1, y_{i+1}+2-w,1)$, the outer corner $(x_{i+1},y_{i+1},h_z)$ of size $(1, y_{i+1}+2-w,1)$, and the $nnp$ vector $(x_{i+1}-v,-1, h_z)$ satisfy the conditions of Lemma \ref{bigimportantlemma}, again giving a second bounded connected component, contradiction.

It remains to consider the possibility that one of the edges leaves $T_i$ through $L$ or $U$. Suppose, for example, that the edge leaves through $L$. It then either enters the adjacent rectangle $T_i'$ in region $C_1$, or it enters $V_2$. 

We show that if this edge enters $C_1$, then there is a second bounded connected component. Because there are at least two outer corners in the top layer, every rectangle $T_i$ in region $C_2$ is adjacent to at least one rectangle in region $C_1$. So it will be enough to conclude that any inner corner in region $C_2$ induces two bounded connected componen.

If the edge in region $T_i$ of $C_2$ enters $T_i'$ in $C_1$ and subsequently leaves $T_i'$ through its bottom side $D'$, then there must be an inner corner in $T_i'$. Choosing  such inner corner with largest $x$-coordinate, and denote it as $q'=(v',w',0)$. The outer corner $(x_{i-1},y_{i-1},h_z)$ of size $(x_{i-1}+2-v',1, 1)$, the inner corner $q'$, and the $nnp$ vector $(-1,y_{i-1}-w', h_z)$ satisfy the conditions of Lemma \ref{bigimportantlemma}, giving a second bounded connected component.
%and the shifted outer corner is $(x_{i-1}+1,w',1)$

\hspace*{-1cm} 
\resizebox{0.5\textwidth}{!}
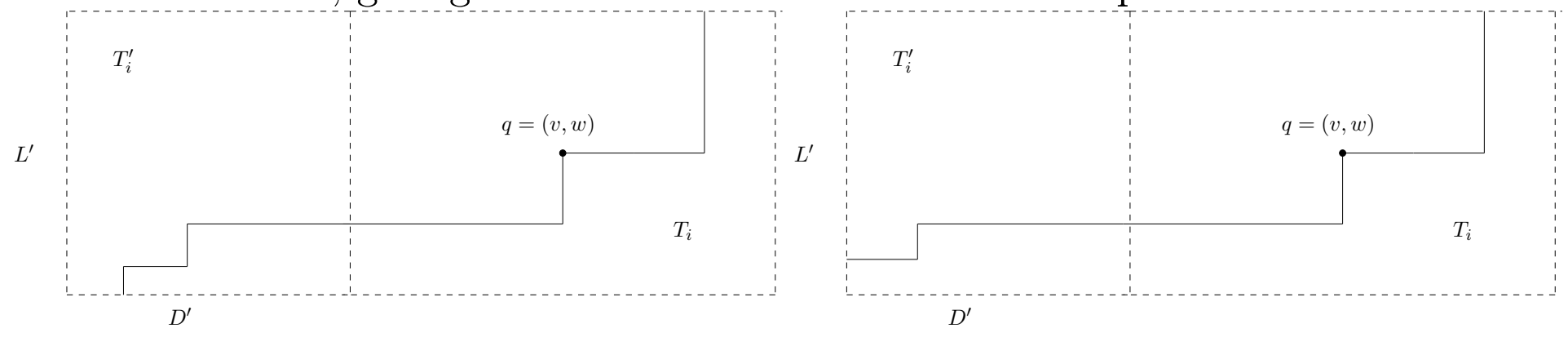
{%
\begin{circuitikz}
\tikzstyle{every node}=[font=\fontsize{10.8pt}{14.1pt}\selectfont]
\draw (6.25,9.875) to[short] (6.25,8.625);
\begin{scope}[transparency group, opacity=1]
\node at (6.25,9.875) [circ, color={rgb,255:red,0; green,0; blue,0}] {};
\end{scope}
\node [font=\fontsize{10.8pt}{14.1pt}\selectfont, inner xsep=0.080cm, inner ysep=0.085cm, rounded corners=0.020cm] at (6,10.375) {$q=(v,w)$};
\draw [dashed] (2.5,12.375) -- (10.125,12.375);
\draw [dashed] (10,12.375) -- (10,7.375);
\draw [dashed] (10,7.375) -- (2.375,7.375);
\draw [dashed] (2.5,12.375) -- (2.5,7.375);
\node [font=\fontsize{10.8pt}{14.1pt}\selectfont, inner xsep=0.080cm, inner ysep=0.085cm, rounded corners=0.020cm] at (8.375,8.5) {$T_i$};
\draw (6.25,9.875) to[short] (7.5,9.875);
\draw (7.5,9.875) to[short] (8.75,9.875);
\draw (8.75,9.875) to[short] (8.75,12.375);
\draw (6.25,8.625) to[short] (2.375,8.625);
\draw [dashed] (2.5,12.375) -- (-2.5,12.375);
\draw [dashed] (-2.5,12.375) -- (-2.5,7.375);
\draw [dashed] (-2.5,7.375) -- (2.5,7.375);
\node [font=\fontsize{10.8pt}{14.1pt}\selectfont, inner xsep=0.080cm, inner ysep=0.085cm, rounded corners=0.020cm] at (-3.25,9.875) {$L'$};
\draw (2.5,8.625) to[short] (-0.375,8.625);
\node [font=\fontsize{10.8pt}{14.1pt}\selectfont, inner xsep=0.080cm, inner ysep=0.085cm, rounded corners=0.020cm] at (-0.5,7) {$D'$};
\node [font=\fontsize{10.8pt}{14.1pt}\selectfont, inner xsep=0.080cm, inner ysep=0.085cm, rounded corners=0.020cm] at (-1.5,11.5) {$T'_i$};
\draw (-0.375,8.625) to[short] (-0.375,7.875);
\draw (-0.375,7.875) to[short] (-1.5,7.875);
\draw (-1.5,7.875) to[short] (-1.5,7.375);
\end{circuitikz}
}%
\resizebox{0.5\textwidth}{!}{%
\begin{circuitikz}
\tikzstyle{every node}=[font=\fontsize{10.8pt}{14.1pt}\selectfont]
\draw (6.25,9.875) to[short] (6.25,8.625);
\begin{scope}[transparency group, opacity=1]
\node at (6.25,9.875) [circ, color={rgb,255:red,0; green,0; blue,0}] {};
\end{scope}
\node [font=\fontsize{10.8pt}{14.1pt}\selectfont, inner xsep=0.080cm, inner ysep=0.085cm, rounded corners=0.020cm] at (6,10.375) {$q=(v,w)$};
\draw [dashed] (2.5,12.375) -- (10.125,12.375);
\draw [dashed] (10,12.375) -- (10,7.375);
\draw [dashed] (10,7.375) -- (2.375,7.375);
\draw [dashed] (2.5,12.375) -- (2.5,7.375);
\node [font=\fontsize{10.8pt}{14.1pt}\selectfont, inner xsep=0.080cm, inner ysep=0.085cm, rounded corners=0.020cm] at (8.375,8.5) {$T_i$};
\draw (6.25,9.875) to[short] (7.5,9.875);
\draw (7.5,9.875) to[short] (8.75,9.875);
\draw (8.75,9.875) to[short] (8.75,12.375);
\draw (6.25,8.625) to[short] (2.375,8.625);
\draw [dashed] (2.5,12.375) -- (-2.5,12.375);
\draw [dashed] (-2.5,12.375) -- (-2.5,7.375);
\draw [dashed] (-2.5,7.375) -- (2.5,7.375);
\node [font=\fontsize{10.8pt}{14.1pt}\selectfont, inner xsep=0.080cm, inner ysep=0.085cm, rounded corners=0.020cm] at (-3.25,9.875) {$L'$};
\node [font=\fontsize{10.8pt}{14.1pt}\selectfont, inner xsep=0.080cm, inner ysep=0.085cm, rounded corners=0.020cm] at (-0.5,7) {$D'$};
\node [font=\fontsize{10.8pt}{14.1pt}\selectfont, inner xsep=0.080cm, inner ysep=0.085cm, rounded corners=0.020cm] at (-1.5,11.5) {$T'_i$};
\draw (2.375,8.625) to[short] (-1.25,8.625);
\draw (-1.25,8.625) to[short] (-1.25,8);
\draw (-1.25,8) to[short] (-2.5,8);
\end{circuitikz}
}%

On the other hand, if the edge leaves $T_i'$ through its left side $L'$, then either it continues to the left or goes downwards. Eventually, this edge can either (1) enters and cuts through the region $B$, following the boundary between $A$ and $B$ until it reaches the box $(x_{i-1}, y_{i-1}+1,0)$ in region $B$, or (2) it creates an inner corner when it turns downward along the side $L'$ and eventually reaches the box $(x_{i-1}+1, y_{i-1}+2,0)$ in region $C_1$. 

If the edge reaches the box $(x_{i-1}+1, y_{i-1}+2,0)$, denote the inner corner on the side $L'$ as $q'=(v', y_{i-1}+2,0)$. The outer corner $(x_{i-1},y_{i-1},h_z)$ of size $(x_{i-1}+2-v',1, 1)$, inner corner $q'$, and vector $\alpha=(-1,-2,h_z)$ by Lemma \ref{bigimportantlemma} creates a bounded connected component. 

If the edge continues to the left and enters region $B$ (and eventually must leave the region $B$ and never enter region $A$, because monomials in region $A$ are the top layer, which must be contained in the bottom layer), then the boxes it cuts out in the region $B$ is a bounded connected component by Lemma \ref{reg_B}, obtained with the $nnp$ vector $(-1,-1,h_z)$.

The cases where the edges leave $T_i$ through $U$ or $R$ are symmetric. Hence an inner corner in $C_2$ always gives two bounded connected components, contradiction. Therefore there are no inner corners in $C_2$.

\item
Assume that an inner corner $q=(v,w,0)$ is in the rectangle $T_i'$ in the region $C_1$, and recall that $T_i'$ is the rectangle between the corners $(x_i+1,y_i+1)$ and $(x_{i+1}+1,y_{i+1}+1)$. Denote the four sides of the rectangle $T_i'$ as $U',D',L',R'$ correspondingly.

For each edge of an inner corner $q$, assume that it leaves the rectangle $T_i'$ at the side $R'$ with the box below the edge having $x$ coordinate $v'$ (See picture below). This box has $y$ coordinate $y_{i+1}+1$. Either we have $v'=v$, or $v'<v$ and there's at least another inner corner on this edge, and the inner corner with the smallest $x$ coordinate is exactly $v'$. Denote the inner corner on this edge with $x$ coordinate $v'$ as $q'=(v', w', 0)$. Then the outer corner $(x_{i+1}, y_{i+1},h_z)$ of size $(1, y_{i+1}+2-w',1)$, inner corner $q'$, and vector $(x_{i+1}-v',-1,h_z)$ satisfy the conditions of Lemma \ref{bigimportantlemma}, so there is a bounded connected component. The same argument applies when the edge leaves the rectangle $T_i'$ at the top-right vertex on the edge $R'$, at the box $(x_{i+1}+1,y_{i+1}+1,0)$. 

A similar argument applies when the edge of $q$ leaves at the side $D'$ with $x$ coordinate $x_i+1$, in which case we use the inner corner $p=(v'', w'', 0)$ where $p$ is the inner corner on this edge with the smallest $y$ coordinate, and the outer corner $(x_{i}, y_{i},h_z)$ of size $(x_i+2-v'',1,1)$, and $nnp$ vector $(-1, y_i-w'', h_z)$, creating a bounded connected component. This also applies when the edge leaves exactly at the bottom-left vertex $(x_{i},y_{i},0)$ of the rectangle. 

If there $q$ is the unique inner corner in $T_i'$, then $q=q'=q''$. The vectors we used for each case, $(x_{i+1}-v',-1,h_z)$ and $(-1, y_i-w'', h_z)$ are not equal, because it implies $x_{i+1}+1=v$ and $y_i+1=w$. And the box $q=(x_{i+1}+1, y_i+1,0)$ is in region $B$. (The top-left corner of $T_i$ is $(x_{i+1}+2, y_i+2,0)$, and the presumed inner corner $q$ is to the top-left of the corner.) So even when the inner corner $q$ is unique, the two different edges each produces a different bounded connected component.

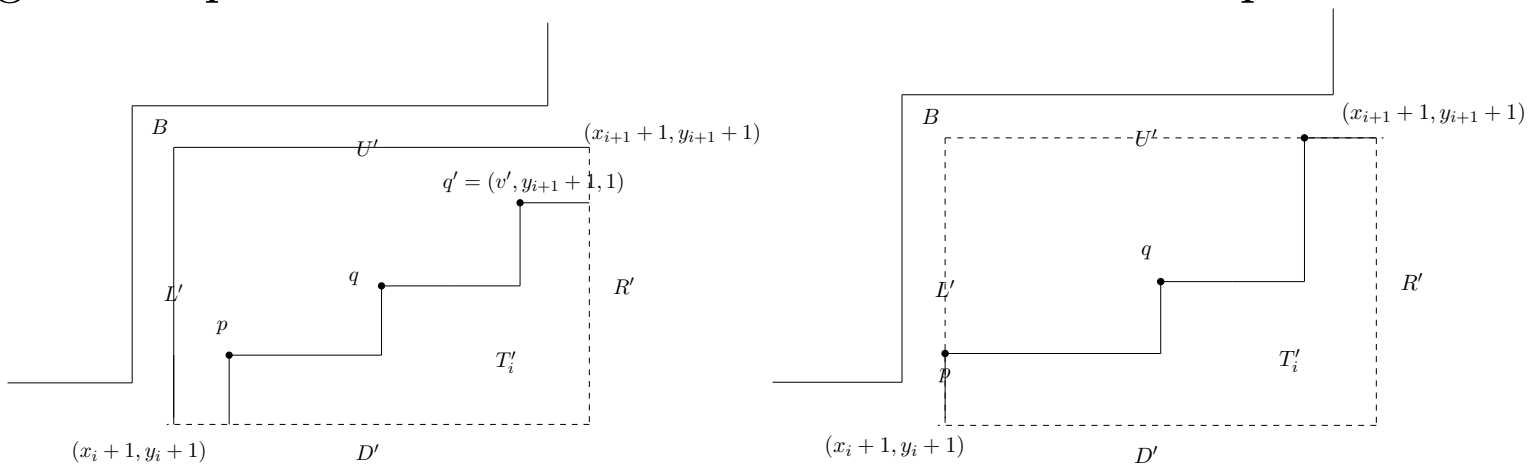
\resizebox{0.4\textwidth}{!}{%
\begin{circuitikz}
\tikzstyle{every node}=[font=\fontsize{10.8pt}{14.1pt}\selectfont]
\draw (6.25,9.875) to[short] (6.25,8.625);
\begin{scope}[transparency group, opacity=1]
\node at (6.25,9.875) [circ, color={rgb,255:red,0; green,0; blue,0}] {};
\end{scope}
\node [font=\fontsize{10.8pt}{14.1pt}\selectfont, inner xsep=0.080cm, inner ysep=0.085cm, rounded corners=0.020cm] at (9,11.75) {$q'=(v',y_{i+1}+1,1)$};
\draw [short] (2.5,12.375) -- (10,12.375);
\draw [dashed] (10,12.375) -- (10,7.375);
\draw [dashed] (10,7.375) -- (2.375,7.375);
\draw [short] (2.5,12.375) -- (2.5,7.375);
\draw (6.25,9.875) to[short] (7.5,9.875);
\draw (7.5,9.875) to[short] (8.75,9.875);
\draw (8.75,9.875) to[short] (8.75,11.375);
\node [font=\fontsize{10.8pt}{14.1pt}\selectfont, inner xsep=0.080cm, inner ysep=0.085cm, rounded corners=0.020cm] at (8.5,8.5) {$T_i'$};
\node [font=\fontsize{10.8pt}{14.1pt}\selectfont, inner xsep=0.080cm, inner ysep=0.085cm, rounded corners=0.020cm] at (6,6.875) {$D'$};
\node [font=\fontsize{10.8pt}{14.1pt}\selectfont, inner xsep=0.080cm, inner ysep=0.085cm, rounded corners=0.020cm] at (6,12.375) {$U'$};
\node [font=\fontsize{10.8pt}{14.1pt}\selectfont, inner xsep=0.080cm, inner ysep=0.085cm, rounded corners=0.020cm] at (2.5,9.75) {$L'$};
\node [font=\fontsize{10.8pt}{14.1pt}\selectfont, inner xsep=0.080cm, inner ysep=0.085cm, rounded corners=0.020cm] at (10.625,9.875) {$R'$};
\node [font=\fontsize{10.8pt}{14.1pt}\selectfont, inner xsep=0.080cm, inner ysep=0.085cm, rounded corners=0.020cm] at (11.5,12.625) {$(x_{i+1}+1,y_{i+1}+1)$};
\node [font=\fontsize{10.8pt}{14.1pt}\selectfont, inner xsep=0.080cm, inner ysep=0.085cm, rounded corners=0.020cm] at (1.875,6.875) {$(x_i+1,y_i+1)$};
\draw (8.75,11.375) to[short] (10,11.375);
\draw (6.25,8.625) to[short] (3.5,8.625);
\draw (2.5,8.625) to[short] (2.5,7.5);
\begin{scope}[transparency group, opacity=1]
\node at (8.75,11.375) [circ, color={rgb,255:red,0; green,0; blue,0}] {};
\end{scope}
\begin{scope}[transparency group, opacity=1]
\node at (3.5,8.625) [circ, color={rgb,255:red,0; green,0; blue,0}] {};
\end{scope}
\node [font=\fontsize{10.8pt}{14.1pt}\selectfont, inner xsep=0.080cm, inner ysep=0.085cm, rounded corners=0.020cm] at (3.375,9.125) {$p$};

\draw (1.75,8.125) to[short] (1.75,13.125);
\draw (1.75,13.125) to[short] (9.25,13.125);
\draw (1.75,8.125) to[short] (-0.5,8.125);
\draw (9.25,13.125) to[short] (9.25,14.625);
\node [font=\fontsize{10.8pt}{14.1pt}\selectfont, inner xsep=0.080cm, inner ysep=0.085cm, rounded corners=0.020cm] at (2.25,12.75) {$B$};
\node [font=\fontsize{10.8pt}{14.1pt}\selectfont, inner xsep=0.080cm, inner ysep=0.085cm, rounded corners=0.000cm] at (5.75,10) {$q$};
\draw (3.5,8.625) to[short] (3.5,7.375);
\end{circuitikz}
}%
\resizebox{0.4\textwidth}{!}{%
\begin{circuitikz}
\tikzstyle{every node}=[font=\fontsize{10.8pt}{14.1pt}\selectfont]
\draw (6.25,9.875) to[short] (6.25,8.625);
\begin{scope}[transparency group, opacity=1]
\node at (6.25,9.875) [circ, color={rgb,255:red,0; green,0; blue,0}] {};
\end{scope}
\node [font=\fontsize{10.8pt}{14.1pt}\selectfont, inner xsep=0.080cm, inner ysep=0.085cm, rounded corners=0.020cm] at (6,10.375) {$q$};
\draw [dashed] (2.5,12.375) -- (10.125,12.375);
\draw [dashed] (10,12.375) -- (10,7.375);
\draw [dashed] (10,7.375) -- (2.375,7.375);
\draw [dashed] (2.5,12.375) -- (2.5,7.375);
\draw (6.25,9.875) to[short] (7.5,9.875);
\draw (7.5,9.875) to[short] (8.75,9.875);
\draw (8.75,9.875) to[short] (8.75,12.375);
\node [font=\fontsize{10.8pt}{14.1pt}\selectfont, inner xsep=0.080cm, inner ysep=0.085cm, rounded corners=0.020cm] at (8.5,8.5) {$T_i'$};
\node [font=\fontsize{10.8pt}{14.1pt}\selectfont, inner xsep=0.080cm, inner ysep=0.085cm, rounded corners=0.020cm] at (6,6.875) {$D'$};
\node [font=\fontsize{10.8pt}{14.1pt}\selectfont, inner xsep=0.080cm, inner ysep=0.085cm, rounded corners=0.020cm] at (6,12.375) {$U'$};
\node [font=\fontsize{10.8pt}{14.1pt}\selectfont, inner xsep=0.080cm, inner ysep=0.085cm, rounded corners=0.020cm] at (2.5,9.75) {$L'$};
\node [font=\fontsize{10.8pt}{14.1pt}\selectfont, inner xsep=0.080cm, inner ysep=0.085cm, rounded corners=0.020cm] at (10.625,9.875) {$R'$};
\node [font=\fontsize{10.8pt}{14.1pt}\selectfont, inner xsep=0.080cm, inner ysep=0.085cm, rounded corners=0.020cm] at (1.625,7) {$(x_{i}+1,y_{i}+1)$};
\node [font=\fontsize{10.8pt}{14.1pt}\selectfont, inner xsep=0.080cm, inner ysep=0.085cm, rounded corners=0.020cm] at (11,12.8) {$(x_{i+1}+1,y_{i+1}+1)$};
\draw (8.75,12.375) to[short] (10,12.375);
\draw (6.25,8.625) to[short] (2.5,8.625);
\draw (2.5,8.625) to[short] (2.5,7.5);
\begin{scope}[transparency group, opacity=1]
\node at (8.75,12.375) [circ, color={rgb,255:red,0; green,0; blue,0}] {};
\end{scope}
\begin{scope}[transparency group, opacity=1]
\node at (2.5,8.625) [circ, color={rgb,255:red,0; green,0; blue,0}] {};
\end{scope}
\node [font=\fontsize{10.8pt}{14.1pt}\selectfont, inner xsep=0.080cm, inner ysep=0.085cm, rounded corners=0.020cm] at (2.5,8.25) {$p$};
\draw (1.75,8.125) to[short] (1.75,13.125);
\draw (1.75,13.125) to[short] (9.25,13.125);
\draw (1.75,8.125) to[short] (-0.5,8.125);
\draw (9.25,13.125) to[short] (9.25,14.625);
\node [font=\fontsize{10.8pt}{14.1pt}\selectfont, inner xsep=0.080cm, inner ysep=0.085cm, rounded corners=0.020cm] at (2.25,12.75) {$B$};
\end{circuitikz}
}%

%shifting the corner $(x_i', y_i')$ up by an appropriate vector $(0,+b)$ until $(x_i', y_i'+b)$ hits the box right below where the edge leaves the side $R'$. This creates a bounded connected component. Note that this argument also works when the edge leaves the rectangle $T_i'$ at the top-right vertex. A similar argument applies when the edge of this inner corner in this rectangle leaves at the side $D$, in which case we shift the corner $(x_{i+1}', y_{i+1}')$ to the right to create a bounded connected component. This also applies when the edge leaves at the bottom-left vertex of the rectangle.

If the edges of $q$ leaves the rectangle $T_i'$ at the sides $L'$ or $U'$, it either enters region $B$, or it follows the sides of the rectangle and leaves the rectangle at the bottom-left or top-right vertex which we discussed already. If the edge enters region $B$, then by Lemma \ref{reg_B}, there's a bounded connected component created on this side.

In any of the cases above, each edge of the inner corner $q$ creates $2$ bounded connected components for the $nnp$ signature. This concludes the contradiction for having an inner corner in region $C_1$. 

%and we use the same argument for the sides $R,D$, and claim that each edge creates a bounded connected component.
\end{enumerate}

\end{proof}

\begin{corollary}
\label{shapeofbottomlayer}
Assume that the minimized staircase $E$ has minimal dimension for tangent space $T_{nnp}$, then there exist a set of path connected boxes in region $B$ that are not contained in the bottom layer, and these boxes form a bounded connected component.    
\end{corollary}

\begin{proof}
By Lemma \ref{nocornerinV} and Lemma \ref{nocornerinC}, there are no inner corners of the bottom layer of $E$ in region $C$ or $V$. By Lemma \ref{innercorbottomlayer}, there is at least one inner corner in the bottom layer. Region $A$ as the projection of the top layer must be contained completely in the bottom layer. So the bottom layer must have at least one inner corner in region $B$. 

Consider the inner corner in the bottom layer of $E$ in $B$, then by Lemma \ref{Bcompconnected}, if there are more than boxes in region $B$ but not in the bottom layer, then they must be path connected and form a bounded connected component. 
\end{proof}

\begin{lemma}
Assume that $E$ is a minimized staircase with minimal dimension for tangent space $T_{nnp}$ and the height of $E$ is $h_z$, then $E$ must be a union of two parts. One part is the top layer stacking on itself $h_z+1$ times, and the other part are two staircases with smaller heights each near the $x$ and $y$ axis.
\end{lemma}

\begin{proof}
By Corollary \ref{shapeofbottomlayer}, the boundary of the bottom layer of $E$ is as follows. It must start somewhere on the $x$ axis with coordinate larger than or equal to $x_1+1$ to contain the box $\beta_2$, and then enters and leaves the region $B$, following the outline between region $A$ and $B$, and then it ends somewhere on the $y$ axis with coordinate larger than or equal to $y_k+1$ to contain the box $\beta_1$.

When the bottom layer follows the outline of $A$, the boundaries of the top and bottom layer are the same, so the top and bottom layer coincide. This only happens somewhere in region $B$, not at the starting point on the $x$ axis or at the ending point near the $y$ axis.

For the boundary of the bottom layer near the $x$ and $y$ axis, the bottom layer must have more extra boxes than the top layer (at least, $\beta_1$ and $\beta_2$ are in the bottom layer but not in the top layer). So the middle layers have flexibility as long as they both contain the top layer and are contained in the bottom layer. These two parts are the wings part of the staircase $E$, near the $x$ and $y$ axis.
\end{proof}

Denote the numbers of rows and columns of the left wing as $p_1, q_1$, and the numbers of rows and columns of the right wing as $p_2,q_2$ accordingly. Note that these numbers are the lengths and widths of the bottom layers of the wings, because the bottom layer is the largest partition among all layers. We now prove some lemmas restricting these numbers, and the shapes of the wings.

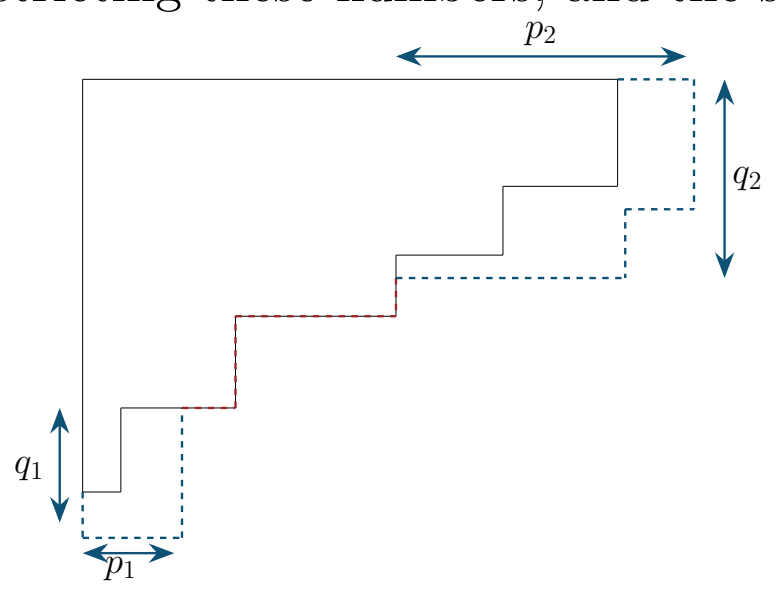
\resizebox{0.4\textwidth}{!}{%
\begin{circuitikz}
\tikzstyle{every node}=[font=\fontsize{18.2pt}{23.7pt}\selectfont]
\draw (6.25,9.75) to[short] (15,9.75);
\draw (15,9.75) to[short] (15,8);
\draw (15,8) to[short] (13.125,8);
\draw (13.125,8) to[short] (13.125,6.875);
\draw (13.125,6.875) to[short] (11.375,6.875);
\draw (6.25,9.75) to[short] (6.25,3);
\draw (6.25,3) to[short] (6.875,3);
\draw (6.875,3) to[short] (6.875,4.375);
\draw (6.875,4.375) to[short] (8.75,4.375);
\draw (8.75,4.375) to[short] (8.75,5.875);
\draw (8.75,5.875) to[short] (11.375,5.875);
\draw (11.375,6.875) to[short] (11.375,5.875);
\draw [ color={rgb,255:red,15; green,81; blue,116}, draw opacity=1, line width=1.1pt, dashed] (6.25,3) -- (6.25,2.25);
\draw [ color={rgb,255:red,13; green,82; blue,117}, draw opacity=1, line width=1.1pt, dashed] (6.25,2.25) -- (7.875,2.25);
\draw [ color={rgb,255:red,11; green,83; blue,118}, draw opacity=1, line width=1.1pt, dashed] (7.875,2.25) -- (7.875,4.375);
\draw [ color={rgb,255:red,155; green,38; blue,38}, draw opacity=1, line width=1.1pt, dashed] (7.875,4.375) -- (8.75,4.375);
\draw [ color={rgb,255:red,155; green,38; blue,38}, draw opacity=1, line width=1.1pt, dashed] (8.75,4.375) -- (8.75,5.875);
\draw [ color={rgb,255:red,155; green,38; blue,38}, draw opacity=1, line width=1.1pt, dashed] (8.75,5.875) -- (11.375,5.875);
\draw [ color={rgb,255:red,155; green,38; blue,38}, draw opacity=1, line width=1.1pt, dashed] (11.375,6.5) -- (11.375,5.875);
\draw [ color={rgb,255:red,15; green,81; blue,116}, draw opacity=1, line width=1.1pt, dashed] (11.375,6.5) -- (15.125,6.5);
\draw [ color={rgb,255:red,13; green,82; blue,117}, draw opacity=1, line width=1.1pt, dashed] (15.125,7.625) -- (15.125,6.5);
\draw [ color={rgb,255:red,13; green,82; blue,117}, draw opacity=1, line width=1.1pt, dashed] (15,9.75) -- (16.25,9.75);
\draw [ color={rgb,255:red,13; green,82; blue,117}, draw opacity=1, line width=1.1pt, dashed] (16.25,9.75) -- (16.25,7.625);
\draw [ color={rgb,255:red,13; green,82; blue,117}, draw opacity=1, line width=1.1pt, dashed] (15.125,7.625) -- (16.25,7.625);
\node [font=\fontsize{18.2pt}{23.7pt}\selectfont, inner xsep=0.080cm, inner ysep=0.085cm, rounded corners=0.020cm] at (6.875,1.75) {$p_1$};
\node [font=\fontsize{18.2pt}{23.7pt}\selectfont, inner xsep=0.080cm, inner ysep=0.085cm, rounded corners=0.020cm] at (5.375,3.375) {$q_1$};
\draw [ color={rgb,255:red,13; green,82; blue,117}, draw opacity=1, line width=1.1pt, {Stealth[scale=1.5]}-{Stealth[scale=1.5]}, ] (6.25,2) -- (7.75,2);
\draw [ color={rgb,255:red,13; green,82; blue,117}, draw opacity=1, line width=1.1pt, {Stealth[scale=1.5]}-{Stealth[scale=1.5]}, ] (5.875,4.375) -- (5.875,2.5);
\draw [ color={rgb,255:red,13; green,82; blue,117}, draw opacity=1, line width=1.1pt, {Stealth[scale=1.5]}-{Stealth[scale=1.5]}, ] (16.75,9.75) -- (16.75,6.5);
\draw [ color={rgb,255:red,13; green,82; blue,117}, draw opacity=1, line width=1.1pt, {Stealth[scale=1.5]}-{Stealth[scale=1.5]}, ] (11.375,10.125) -- (16.125,10.125);
\node [font=\fontsize{18.2pt}{23.7pt}\selectfont, fill={rgb,255:red,255; green,255; blue,255}, fill opacity=1, text opacity=1, inner xsep=0.080cm, inner ysep=0.085cm, rounded corners=0.020cm] at (13.75,10.5) {$p_2$};
\node [font=\fontsize{18.2pt}{23.7pt}\selectfont, fill={rgb,255:red,255; green,255; blue,255}, fill opacity=1, text opacity=1, inner xsep=0.080cm, inner ysep=0.085cm, rounded corners=0.020cm] at (17.125,8.125) {$q_2$};
\end{circuitikz}
}%
\begin{lemma}
\label{size_of_side_comp}
For a minimized staircase $E$ with minimal dimension for the tangent space $T_{nnp}$, the inner corners of the top layer cannot be located above the bottom layers of the wings. In other words, $p_1 \leq y_1$ and $q_2 \leq x_k$.
\end{lemma}

\begin{proof}
Assume for contradiction that $q_1 > x_k$, and the other case, $p_2 > y_1$ can be proved similarly. 

Consider the image of the projection of the inner corner of the top layer $(x_2+1, y_1+1, h_z)$ in the $z$ direction to the bottom layer $(x_2+1, y_1+1,0)$.

When $p_1 > y_1$, this box $(x_2+1, y_1+1,0)$ is contained in the bottom layer of the left wing.  Denote $N=(x_2+1, y_1+1,u)$ for some $u \leq h_z$ such that $(x_2+1, y_1+1,u) \in I$ and $(x_2+1, y_1+1,u-1) \in E$. Because $(x_2, y_1+1,h_z), (x_2+1, y_1,h_z), (x_2, y_1,h_z)$ and all the boxes below them are in $E$, $N$ is an inner corner of the staircase $E$.

Consider the outer corners with the largest $y$ coordinate of the right wing. If we slice the layers in the $y$ direction instead of the usual 
$z$ direction, these outer corners are on the right-most slice, and there can be more than one outer corners with the largest $y$ coordinates. Among these, we pick the outer corner with the largest $x$ coordinate, and denote it as $C=(v, w, t)$, where $v < x_2$, $w>y_1$ and $t <h_z$. Note that we must have $(v+1, w, i) \in I$ for any $i \geq 0$, otherwise there is at least one more outer corner on the right-most slice with largest $x$ coordinate. So we have an outer corner $C$ with size $(1, 1, t+1)$. 

If $t < u$, then the outer corner $C$ of size $(1,1,1)$, the inner corner $N=(x_2+1, y_1+1,u)$, and the $npn$ vector $\alpha_1=(v-x_2, w-y_1, t-u)$ satisfy the condition of Lemma \ref{bigimportantlemma}, so there is a bounded connected component. If $t \geq u$, then the outer corner $C$ of size $(1,1,t-u+2)$, the inner corner $N=(x_2+1, y_1+1,u)$, and the $npn$ vector $\alpha_1=(v-x_2, w-y_1, -1)$ satisfy the condition of Lemma \ref{bigimportantlemma}, so there is a bounded connected component. 

Consider another the inner corner $N'=(x_1+1, 0, u')$ where $(x_1+1, 0, u') \in I$ and $(x_1+1, 0, u'-1) \in E$. This is the inner corner whose projection in the $z$ direction lands on the box $\beta_2$. We know that this must be an inner corner because $\beta_2$ is not in the top layer, so the maximum $z$-height at this box is strictly less than $h_z$, and the box $(x_1, 0, h_z)$ right next to it has height $h_z$.

If $t < u'$, then the outer corner $C=(v, w, t)$ of size $(1,1,1)$, the inner corner $N'$, and the vector $\alpha_1=(v-x_1, w, t-u')$ of signature $npn$ also satisfy the condition of Lemma \ref{bigimportantlemma}, so there is another bounded connected component. If $t \geq u'$, then the outer corner $C=(v, w, t)$ of size $(1,1,t-u'+2)$, the inner corner $N'$, and the vector $\alpha_1=(v-x_1, w, -1)$ of signature $npn$ also satisfy the condition of Lemma \ref{bigimportantlemma}, so there is another bounded connected component. 

In all of the cases above, we show that there are at least two bounded connected components, contradiction. So $p_1 \leq y_1$ and $q_2 \leq x_k$.
\end{proof}

\begin{corollary}
\label{bot_rec_wings}
For a minimized staircase $E$ with minimal tangent space dimension for $T_{nnp}$, the bottom layer of the two wings must be two rectangles and cannot have any inner corners.
\end{corollary}

\begin{proof}
This directly follow from Lemma \ref{size_of_side_comp} and Lemma \ref{sidecomp}. The picture on the left is not possible. The bottom layer of the wings must be the two rectangles, as shown on the right picture.%also rules out the possibility of having more than $1$ outer corner in the bottom layer each of the wing area (see the picture on the left). So the bottom layer of the wings must be exactly two rectangles. We further refine the numerical condition of the rectangles. 
\end{proof}

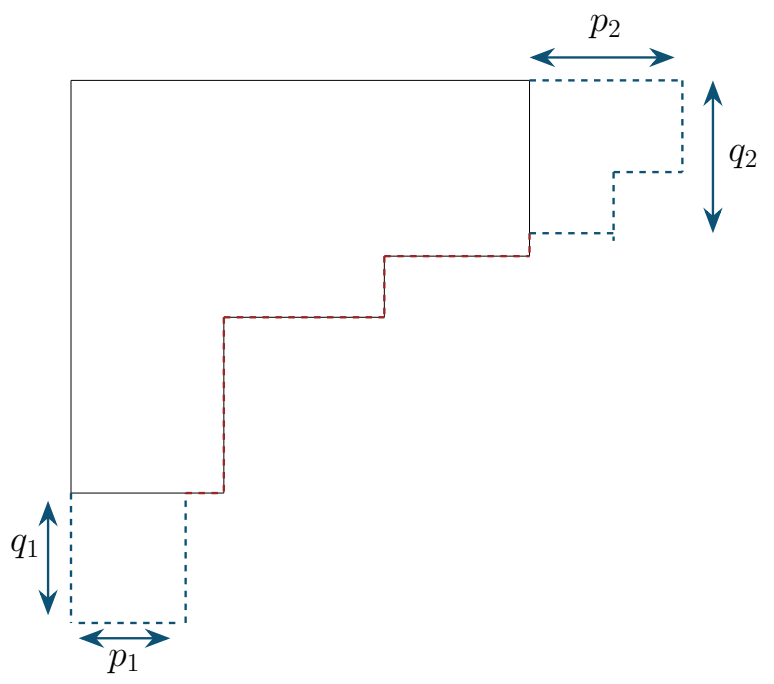
\resizebox{0.4\textwidth}{!}{%
\begin{circuitikz}
\tikzstyle{every node}=[font=\fontsize{18.2pt}{23.7pt}\selectfont]
\draw (6.25,9.75) to[short] (13.75,9.75);
\draw (13.75,9.75) to[short] (13.75,6.875);
\draw (13.75,6.875) to[short] (11.375,6.875);
\draw (6.25,9.75) to[short] (6.25,3);
\draw (6.25,3) to[short] (8.75,3);
\draw (8.75,3) to[short] (8.75,5.875);
\draw (8.75,5.875) to[short] (11.375,5.875);
\draw (11.375,6.875) to[short] (11.375,5.875);
\draw [ color={rgb,255:red,15; green,81; blue,116}, draw opacity=1, line width=1.1pt, dashed] (6.25,3) -- (6.25,0.875);
\draw [ color={rgb,255:red,13; green,82; blue,117}, draw opacity=1, line width=1.1pt, dashed] (6.375,0.875) -- (8,0.875);
\draw [ color={rgb,255:red,11; green,83; blue,118}, draw opacity=1, line width=1.1pt, dashed] (8.125,0.875) -- (8.125,3);
\draw [ color={rgb,255:red,155; green,38; blue,38}, draw opacity=1, line width=1.1pt, dashed] (8.75,5.875) -- (11.375,5.875);
\draw [ color={rgb,255:red,155; green,38; blue,38}, draw opacity=1, line width=1.1pt, dashed] (11.375,6.875) -- (11.375,5.875);
\draw [ color={rgb,255:red,15; green,81; blue,116}, draw opacity=1, line width=1.1pt, dashed] (13.75,7.25) -- (15.125,7.25);
\draw [ color={rgb,255:red,13; green,82; blue,117}, draw opacity=1, line width=1.1pt, dashed] (15.125,8.25) .. controls (15.125,7.75) and (15.125,8.625) .. (15.125,7.125);
\draw [ color={rgb,255:red,13; green,82; blue,117}, draw opacity=1, line width=1.1pt, dashed] (13.75,9.75) -- (16.25,9.75);
\draw [ color={rgb,255:red,13; green,82; blue,117}, draw opacity=1, line width=1.1pt, dashed] (16.25,9.75) -- (16.25,8.25);
\draw [ color={rgb,255:red,13; green,82; blue,117}, draw opacity=1, line width=1.1pt, dashed] (15.125,8.25) -- (16.25,8.25);
\node [font=\fontsize{18.2pt}{23.7pt}\selectfont, inner xsep=0.080cm, inner ysep=0.085cm, rounded corners=0.020cm] at (7.125,0.25) {$p_1$};
\node [font=\fontsize{18.2pt}{23.7pt}\selectfont, inner xsep=0.080cm, inner ysep=0.085cm, rounded corners=0.020cm] at (5.5,2.125) {$q_1$};
\draw [ color={rgb,255:red,13; green,82; blue,117}, draw opacity=1, line width=1.1pt, {Stealth[scale=1.5]}-{Stealth[scale=1.5]}, ] (6.375,0.625) -- (7.875,0.625);
\draw [ color={rgb,255:red,13; green,82; blue,117}, draw opacity=1, line width=1.1pt, {Stealth[scale=1.5]}-{Stealth[scale=1.5]}, ] (5.875,2.875) -- (5.875,1);
\draw [ color={rgb,255:red,13; green,82; blue,117}, draw opacity=1, line width=1.1pt, {Stealth[scale=1.5]}-{Stealth[scale=1.5]}, ] (16.75,9.75) -- (16.75,7.25);
\draw [ color={rgb,255:red,13; green,82; blue,117}, draw opacity=1, line width=1.1pt, {Stealth[scale=1.5]}-{Stealth[scale=1.5]}, ] (13.75,10.125) -- (16.125,10.125);
\node [font=\fontsize{18.2pt}{23.7pt}\selectfont, fill={rgb,255:red,255; green,255; blue,255}, fill opacity=1, text opacity=1, inner xsep=0.080cm, inner ysep=0.085cm, rounded corners=0.020cm] at (15,10.625) {$p_2$};
\node [font=\fontsize{18.2pt}{23.7pt}\selectfont, fill={rgb,255:red,255; green,255; blue,255}, fill opacity=1, text opacity=1, inner xsep=0.080cm, inner ysep=0.085cm, rounded corners=0.020cm] at (17.25,8.5) {$q_2$};
\draw [ color={rgb,255:red,155; green,38; blue,38}, draw opacity=1, line width=1.1pt, dashed] (8.125,3) -- (8.75,3);
\draw [ color={rgb,255:red,155; green,38; blue,38}, draw opacity=1, line width=1.1pt, dashed] (8.75,5.875) -- (8.75,3);
\draw [ color={rgb,255:red,155; green,38; blue,38}, draw opacity=1, line width=1.1pt, dashed] (11.375,6.875) -- (13.75,6.875);
\draw [ color={rgb,255:red,155; green,38; blue,38}, draw opacity=1, line width=1.1pt, dashed] (13.75,7.25) -- (13.75,6.875);
\end{circuitikz}
}%
\resizebox{0.4\textwidth}{!}
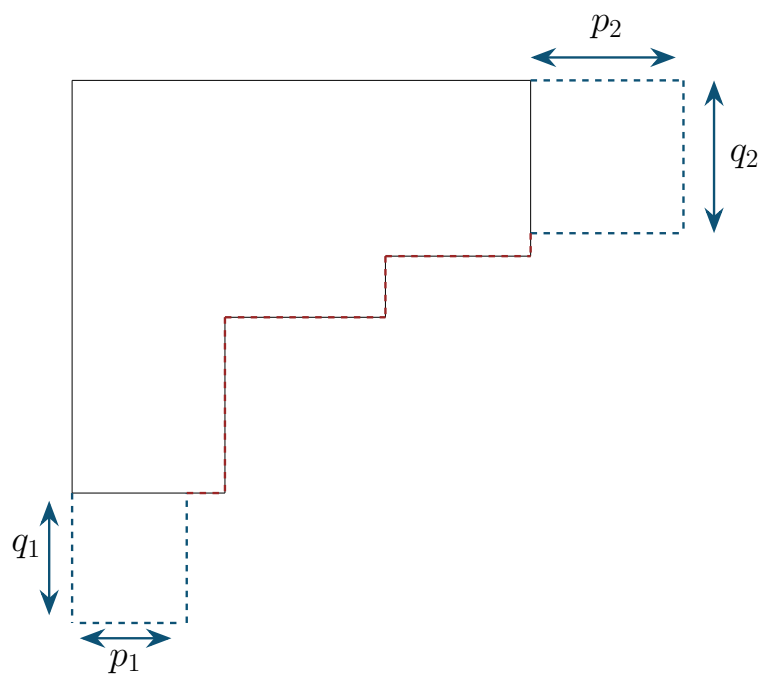
{%
\begin{circuitikz}
\tikzstyle{every node}=[font=\fontsize{18.2pt}{23.7pt}\selectfont]
\draw (6.25,9.75) to[short] (13.75,9.75);
\draw (13.75,9.75) to[short] (13.75,6.875);
\draw (13.75,6.875) to[short] (11.375,6.875);
\draw (6.25,9.75) to[short] (6.25,3);
\draw (6.25,3) to[short] (8.75,3);
\draw (8.75,3) to[short] (8.75,5.875);
\draw (8.75,5.875) to[short] (11.375,5.875);
\draw (11.375,6.875) to[short] (11.375,5.875);
\draw [ color={rgb,255:red,15; green,81; blue,116}, draw opacity=1, line width=1.1pt, dashed] (6.25,3) -- (6.25,0.875);
\draw [ color={rgb,255:red,13; green,82; blue,117}, draw opacity=1, line width=1.1pt, dashed] (6.375,0.875) -- (8,0.875);
\draw [ color={rgb,255:red,11; green,83; blue,118}, draw opacity=1, line width=1.1pt, dashed] (8.125,0.875) -- (8.125,3);
\draw [ color={rgb,255:red,155; green,38; blue,38}, draw opacity=1, line width=1.1pt, dashed] (8.75,5.875) -- (11.375,5.875);
\draw [ color={rgb,255:red,155; green,38; blue,38}, draw opacity=1, line width=1.1pt, dashed] (11.375,6.875) -- (11.375,5.875);
\draw [ color={rgb,255:red,15; green,81; blue,116}, draw opacity=1, line width=1.1pt, dashed] (13.75,7.25) -- (16.25,7.25);
\draw [ color={rgb,255:red,13; green,82; blue,117}, draw opacity=1, line width=1.1pt, dashed] (13.75,9.75) -- (16.25,9.75);
\draw [ color={rgb,255:red,13; green,82; blue,117}, draw opacity=1, line width=1.1pt, dashed] (16.25,9.75) -- (16.25,7.25);
\node [font=\fontsize{18.2pt}{23.7pt}\selectfont, inner xsep=0.080cm, inner ysep=0.085cm, rounded corners=0.020cm] at (7.125,0.25) {$p_1$};
\node [font=\fontsize{18.2pt}{23.7pt}\selectfont, inner xsep=0.080cm, inner ysep=0.085cm, rounded corners=0.020cm] at (5.5,2.125) {$q_1$};
\draw [ color={rgb,255:red,13; green,82; blue,117}, draw opacity=1, line width=1.1pt, {Stealth[scale=1.5]}-{Stealth[scale=1.5]}, ] (6.375,0.625) -- (7.875,0.625);
\draw [ color={rgb,255:red,13; green,82; blue,117}, draw opacity=1, line width=1.1pt, {Stealth[scale=1.5]}-{Stealth[scale=1.5]}, ] (5.875,2.875) -- (5.875,1);
\draw [ color={rgb,255:red,13; green,82; blue,117}, draw opacity=1, line width=1.1pt, {Stealth[scale=1.5]}-{Stealth[scale=1.5]}, ] (16.75,9.75) -- (16.75,7.25);
\draw [ color={rgb,255:red,13; green,82; blue,117}, draw opacity=1, line width=1.1pt, {Stealth[scale=1.5]}-{Stealth[scale=1.5]}, ] (13.75,10.125) -- (16.125,10.125);
\node [font=\fontsize{18.2pt}{23.7pt}\selectfont, fill={rgb,255:red,255; green,255; blue,255}, fill opacity=1, text opacity=1, inner xsep=0.080cm, inner ysep=0.085cm, rounded corners=0.020cm] at (15,10.625) {$p_2$};
\node [font=\fontsize{18.2pt}{23.7pt}\selectfont, fill={rgb,255:red,255; green,255; blue,255}, fill opacity=1, text opacity=1, inner xsep=0.080cm, inner ysep=0.085cm, rounded corners=0.020cm] at (17.25,8.5) {$q_2$};
\draw [ color={rgb,255:red,155; green,38; blue,38}, draw opacity=1, line width=1.1pt, dashed] (8.125,3) -- (8.75,3);
\draw [ color={rgb,255:red,155; green,38; blue,38}, draw opacity=1, line width=1.1pt, dashed] (8.75,5.875) -- (8.75,3);
\draw [ color={rgb,255:red,155; green,38; blue,38}, draw opacity=1, line width=1.1pt, dashed] (11.375,6.875) -- (13.75,6.875);
\draw [ color={rgb,255:red,155; green,38; blue,38}, draw opacity=1, line width=1.1pt, dashed] (13.75,7.25) -- (13.75,6.875);
\end{circuitikz}
}%

We further discuss the numerical conditions on the bottom layers of the wings.

\begin{lemma}
For a minimized staircase $E$ whose tangent space $T_{nnp}$ has minimal dimension $1$: 

$\min(p_1, p_2)=1$, and $\min(q_1, q_2)=1$.
\end{lemma}

\begin{proof}
As we discussed multiple times before, the top layer of the staircase $B$, shifted by the $nnp$ vector $(-1,-1,h_z)$, becomes a bounded connected component.

When $p_1 \geq 2$ and $p_2 \geq 2$, shifting the top layer of $E$ by the vector $(-1, -2, h_z)$, the bottom layer must cut the shifted top layer into a bounded connected component, contradiction. In order for there to be one unique bounded connected component, we must have that either $p_1=1$ or $p_2=1$. 

By similar reasoning with the $nnp$ vector $(-2, -1, h_z)$, we also have the other two conditions for $q_1$ and $q_2$.
\end{proof}

\begin{corollary}
\label{bottomlayerofwings}
Assume $E$ is a minimized staircase with tangent space $T_{nnp}$ having minimal dimension $1$, then the bottom layer of the wings can only have the following forms.

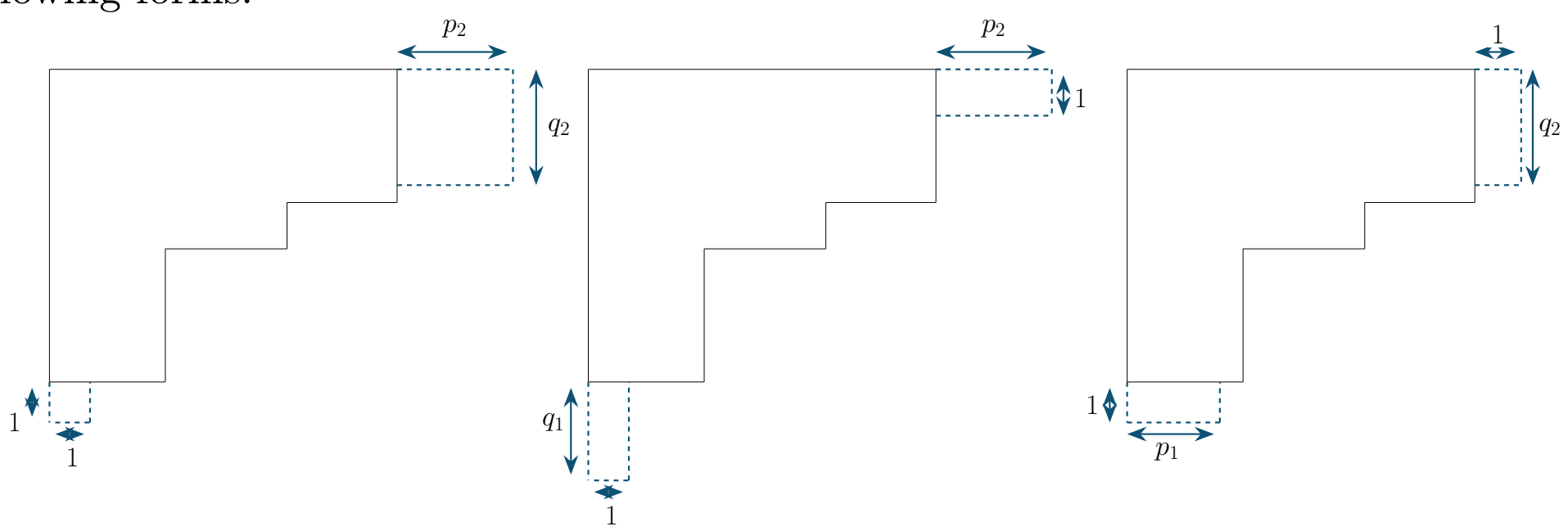
\resizebox{1\textwidth}{!}{%
\begin{circuitikz}
\tikzstyle{every node}=[font=\fontsize{18.2pt}{23.7pt}\selectfont]
\draw (6.25,9.75) to[short] (13.75,9.75);
\draw (13.75,9.75) to[short] (13.75,6.875);
\draw (13.75,6.875) to[short] (11.375,6.875);
\draw (6.25,9.75) to[short] (6.25,3);
\draw (6.25,3) to[short] (8.75,3);
\draw (8.75,3) to[short] (8.75,5.875);
\draw (8.75,5.875) to[short] (11.375,5.875);
\draw (11.375,6.875) to[short] (11.375,5.875);
\draw [ color={rgb,255:red,15; green,81; blue,116}, draw opacity=1, line width=1.1pt, dashed] (6.25,3) -- (6.25,0.875);
\draw [ color={rgb,255:red,13; green,82; blue,117}, draw opacity=1, line width=1.1pt, dashed] (6.375,0.875) -- (7.125,0.875);
\draw [ color={rgb,255:red,11; green,83; blue,118}, draw opacity=1, line width=1.1pt, dashed] (7.125,0.875) -- (7.125,3);
\draw [ color={rgb,255:red,15; green,81; blue,116}, draw opacity=1, line width=1.1pt, dashed] (13.75,8.75) -- (16.25,8.75);
\draw [ color={rgb,255:red,13; green,82; blue,117}, draw opacity=1, line width=1.1pt, dashed] (13.75,9.75) -- (16.25,9.75);
\draw [ color={rgb,255:red,13; green,82; blue,117}, draw opacity=1, line width=1.1pt, dashed] (16.25,9.75) -- (16.25,8.75);
\node [font=\fontsize{18.2pt}{23.7pt}\selectfont, inner xsep=0.080cm, inner ysep=0.085cm, rounded corners=0.020cm] at (6.75,0.125) {1};
\node [font=\fontsize{18.2pt}{23.7pt}\selectfont, inner xsep=0.080cm, inner ysep=0.085cm, rounded corners=0.020cm] at (5.5,2.125) {$q_1$};
\draw [ color={rgb,255:red,13; green,82; blue,117}, draw opacity=1, line width=1.1pt, {Stealth[scale=1.5]}-{Stealth[scale=1.5]}, ] (6.375,0.625) -- (7,0.625);
\draw [ color={rgb,255:red,13; green,82; blue,117}, draw opacity=1, line width=1.1pt, {Stealth[scale=1.5]}-{Stealth[scale=1.5]}, ] (5.875,2.875) -- (5.875,1);
\draw [ color={rgb,255:red,13; green,82; blue,117}, draw opacity=1, line width=1.1pt, {Stealth[scale=1.5]}-{Stealth[scale=1.5]}, ] (16.5,9.625) -- (16.5,8.75);
\draw [ color={rgb,255:red,13; green,82; blue,117}, draw opacity=1, line width=1.1pt, {Stealth[scale=1.5]}-{Stealth[scale=1.5]}, ] (13.75,10.125) -- (16.125,10.125);
\node [font=\fontsize{18.2pt}{23.7pt}\selectfont, fill={rgb,255:red,255; green,255; blue,255}, fill opacity=1, text opacity=1, inner xsep=0.080cm, inner ysep=0.085cm, rounded corners=0.020cm] at (15,10.625) {$p_2$};
\node [font=\fontsize{18.2pt}{23.7pt}\selectfont, fill={rgb,255:red,255; green,255; blue,255}, fill opacity=1, text opacity=1, inner xsep=0.080cm, inner ysep=0.085cm, rounded corners=0.020cm] at (16.875,9.125) {1};
\draw (17.875,9.75) to[short] (25.375,9.75);
\draw (25.375,9.75) to[short] (25.375,6.875);
\draw (25.375,6.875) to[short] (23,6.875);
\draw (17.875,9.75) to[short] (17.875,3);
\draw (17.875,3) to[short] (20.375,3);
\draw (20.375,3) to[short] (20.375,5.875);
\draw (20.375,5.875) to[short] (23,5.875);
\draw (23,6.875) to[short] (23,5.875);
\draw [ color={rgb,255:red,15; green,81; blue,116}, draw opacity=1, line width=1.1pt, dashed] (17.875,3) -- (17.875,2.125);
\draw [ color={rgb,255:red,13; green,82; blue,117}, draw opacity=1, line width=1.1pt, dashed] (17.875,2.125) -- (19.875,2.125);
\draw [ color={rgb,255:red,11; green,83; blue,118}, draw opacity=1, line width=1.1pt, dashed] (19.875,2.125) -- (19.875,3);
\draw [ color={rgb,255:red,15; green,81; blue,116}, draw opacity=1, line width=1.1pt, dashed] (25.375,7.25) -- (26.375,7.25);
\draw [ color={rgb,255:red,13; green,82; blue,117}, draw opacity=1, line width=1.1pt, dashed] (25.375,9.75) -- (26.375,9.75);
\draw [ color={rgb,255:red,13; green,82; blue,117}, draw opacity=1, line width=1.1pt, dashed] (26.375,9.75) -- (26.375,7.25);
\node [font=\fontsize{18.2pt}{23.7pt}\selectfont, inner xsep=0.080cm, inner ysep=0.085cm, rounded corners=0.020cm] at (18.75,1.5) {$p_1$};
\node [font=\fontsize{18.2pt}{23.7pt}\selectfont, inner xsep=0.080cm, inner ysep=0.085cm, rounded corners=0.020cm] at (17.125,2.5) {1};
\draw [ color={rgb,255:red,13; green,82; blue,117}, draw opacity=1, line width=1.1pt, {Stealth[scale=1.5]}-{Stealth[scale=1.5]}, ] (17.875,1.875) .. controls (18.875,1.875) and (18.75,1.875) .. (19.75,1.875);
\draw [ color={rgb,255:red,13; green,82; blue,117}, draw opacity=1, line width=1.1pt, {Stealth[scale=1.5]}-{Stealth[scale=1.5]}, ] (17.5,2.875) -- (17.5,2.125);
\draw [ color={rgb,255:red,13; green,82; blue,117}, draw opacity=1, line width=1.1pt, {Stealth[scale=1.5]}-{Stealth[scale=1.5]}, ] (26.625,9.75) -- (26.625,7.25);
\draw [ color={rgb,255:red,13; green,82; blue,117}, draw opacity=1, line width=1.1pt, {Stealth[scale=1.5]}-{Stealth[scale=1.5]}, ] (25.375,10.125) -- (26.25,10.125);
\node [font=\fontsize{18.2pt}{23.7pt}\selectfont, fill={rgb,255:red,255; green,255; blue,255}, fill opacity=1, text opacity=1, inner xsep=0.080cm, inner ysep=0.085cm, rounded corners=0.020cm] at (25.875,10.5) {1};
\node [font=\fontsize{18.2pt}{23.7pt}\selectfont, fill={rgb,255:red,255; green,255; blue,255}, fill opacity=1, text opacity=1, inner xsep=0.080cm, inner ysep=0.085cm, rounded corners=0.020cm] at (27,8.5) {$q_2$};
\draw (-5.375,9.75) to[short] (2.125,9.75);
\draw (2.125,9.75) to[short] (2.125,6.875);
\draw (2.125,6.875) to[short] (-0.25,6.875);
\draw (-5.375,9.75) to[short] (-5.375,3);
\draw (-5.375,3) to[short] (-2.875,3);
\draw (-2.875,3) to[short] (-2.875,5.875);
\draw (-2.875,5.875) to[short] (-0.25,5.875);
\draw (-0.25,6.875) to[short] (-0.25,5.875);
\draw [ color={rgb,255:red,15; green,81; blue,116}, draw opacity=1, line width=1.1pt, dashed] (-5.375,3) -- (-5.375,2.125);
\draw [ color={rgb,255:red,13; green,82; blue,117}, draw opacity=1, line width=1.1pt, dashed] (-5.25,2.125) -- (-4.5,2.125);
\draw [ color={rgb,255:red,11; green,83; blue,118}, draw opacity=1, line width=1.1pt, dashed] (-4.5,2.125) -- (-4.5,3);
\draw [ color={rgb,255:red,15; green,81; blue,116}, draw opacity=1, line width=1.1pt, dashed] (2.125,7.25) -- (4.625,7.25);
\draw [ color={rgb,255:red,13; green,82; blue,117}, draw opacity=1, line width=1.1pt, dashed] (2.125,9.75) -- (4.625,9.75);
\draw [ color={rgb,255:red,13; green,82; blue,117}, draw opacity=1, line width=1.1pt, dashed] (4.625,9.75) -- (4.625,7.25);
\node [font=\fontsize{18.2pt}{23.7pt}\selectfont, inner xsep=0.080cm, inner ysep=0.085cm, rounded corners=0.020cm] at (-4.875,1.375) {1};
\node [font=\fontsize{18.2pt}{23.7pt}\selectfont, inner xsep=0.080cm, inner ysep=0.085cm, rounded corners=0.020cm] at (-6.125,2.125) {1};
\draw [ color={rgb,255:red,13; green,82; blue,117}, draw opacity=1, line width=1.1pt, {Stealth[scale=1.5]}-{Stealth[scale=1.5]}, ] (-5.25,1.875) -- (-4.625,1.875);
\draw [ color={rgb,255:red,13; green,82; blue,117}, draw opacity=1, line width=1.1pt, {Stealth[scale=1.5]}-{Stealth[scale=1.5]}, ] (-5.75,2.875) -- (-5.75,2.25);
\draw [ color={rgb,255:red,13; green,82; blue,117}, draw opacity=1, line width=1.1pt, {Stealth[scale=1.5]}-{Stealth[scale=1.5]}, ] (5.125,9.75) -- (5.125,7.25);
\draw [ color={rgb,255:red,13; green,82; blue,117}, draw opacity=1, line width=1.1pt, {Stealth[scale=1.5]}-{Stealth[scale=1.5]}, ] (2.125,10.125) -- (4.5,10.125);
\node [font=\fontsize{18.2pt}{23.7pt}\selectfont, fill={rgb,255:red,255; green,255; blue,255}, fill opacity=1, text opacity=1, inner xsep=0.080cm, inner ysep=0.085cm, rounded corners=0.020cm] at (3.375,10.625) {$p_2$};
\node [font=\fontsize{18.2pt}{23.7pt}\selectfont, fill={rgb,255:red,255; green,255; blue,255}, fill opacity=1, text opacity=1, inner xsep=0.080cm, inner ysep=0.085cm, rounded corners=0.020cm] at (5.625,8.5) {$q_2$};
\end{circuitikz}
}%
\begin{enumerate}
\item If $p_1=q_1=1$, then $p_2$ and $q_2$ can be any arbitrary number. If $p_1 \geq 2$, and $p_2 \geq 2$, then $p_1=q_1=1$. (One box of size $1$ by $1$, and the other arbitrary.)

\item If $p_1=1$, and $p_2 \geq 2$, then $q_1$ can be any arbitrary number, and $q_2=1$. (Spread out wings.)

\item If $p_1 \geq 2$, and $p_2=1$, then $q_2$ can be any arbitrary number, and $q_1=1$. (Folded wings.)
 
\end{enumerate}
\end{corollary}

%So far, we have proved that if the top layer has $2$ or more corners, then the staircase must have two parts. The main part is some partition stacking on itself $l$ times, and the wings part of the partition is some chunks of staircase attached to the main part next to the $x$ and $y$ axis, where the bottom layer of each wing part must be a rectangular shape on the $XY$ plane, and in the 2d bottom layer, the number of rows of the top-right rectangular wing must be smaller than or equal to $x_1$, and the number of columns of the bottom-left rectangular wing must be smaller than or equal to $y_k$. We now discuss if the 3d wings could have any corners.

\begin{lemma}
\label{cuboid}
Assume that the minimized staircase $E$ has minimal tangent space  $T_{nnp}$ dimension equal to $1$. 

The wings of $E$ as a three dimensional staircase cannot have any inner corners, and each must be a solid rectangular cuboid.   
\end{lemma}

\begin{proof}
We discuss case by case as listed in Corollary \ref{bottomlayerofwings}.

\begin{enumerate}
\item Suppose that the bottom layer of the wings is a single box and another rectangle.

Assume without loss of generality that the left wing's bottom layer is a single box, so it must be a vertical strip of boxes. The vertical strip of boxes must be a cuboid cannot have any inner corners. 

The right wing's bottom layer is a rectangle, so we assume for contradiction that the right wing has an inner corner and is not a rectangular cuboid. 
%Pick an inner corner of the top-right wing, and denote it as $N=(v,w,t)$.

%There are $3$ cases to consider. Suppose this vertical strip of boxes, the bottom-left wing, has height $h$.

%\begin{enumerate}
%\item If $h=1$, then we have one unique box as the bottom-left wing. This box is an outer corner with coordinates $O=(x_1+1, 1, 1)$. The outer corner $O$ of size $(1,1,1)$, the inner corner $N=(v,w,t)$, and the $pnn$ vector $\alpha_1=(x_1+1-v, 1-w,1-t)$ gives a bounded connected component by Lemma \ref{bigimportantlemma}.

%Additionally, assume that the biggest $z$ coordinate of the top-right wing is $s$. The box $N'=(1, y_k+1, z+1)$ is an inner corner. The outer corner $O=(x_1+1, 1, 1)$ of size $(1,1,1)$, the inner corner $N'$, and the $pnn$ vector $\alpha_2=(x_1,-y_k, -z)$ gives another bounded connected component by Lemma \ref{bigimportantlemma}, contradiction.

%\item If $h=h_z-1$, then the height of bottom-left wing is just one box shorter than the main part of the opera house. This box $M=(x_1+1, 1, h_z)$ is an inner corner.

%Now we take the $x$ axis as the direction of height, and the bottom-left wing becomes the top layer of the staircase. The bottom layer is the layer of the staircase above the $YZ$ plane.

Now consider the staircase $E$ sliced into layers in the $x$ direction, so the vertical strip of boxes as the left wing is now the top $x$ layer.  By Lemma \ref{onerow}, the bottom layer on the $YZ$ plane must be a partition with exactly one inner corner, in one of the three possible shapes classified. The wings by construction is shorter than the main part of the partition, so the location where the right wing connects to the main part of the opera house forms an inner corner of the bottom layer on the $YZ$ plane. Therefore, the bottom layer on the $YZ$ plane of the right wing must be a rectangle without any inner corner. 

We assume for contradiction that the right wing as a staircase has at least one inner corner, but because the $YZ$ bottom layer of the right wing is a rectangle, the inner corner must be strictly above the $YZ$ plane. In other words, any of the inner corners' $x$ coordinate must be strictly greater than $0$.

We list out each case as in Lemma \ref{onerow} by considering the top and bottom layers in the $x$ direction. By Corollary \ref{bot_rec_wings}, the $XY$ bottom layer of the right wing is also a rectangle, and the right wing in the third picture only has one layer, so we rule out the case of the inner corner in the third picture. 

Because both the $XY$ and $YZ$ bottom layers are rectangles, and we assumed for contradiction that there must be a three dimensional inner corner in the right wing, the top $y$ layer must have an inner corner, and the reason is as follows. If the top $y$ layer is a rectangle with no inner corner, and either the length or width is smaller than the bottom $y$ layer of the right wing (not the bottom $y$ layer of the entire staircase), then either the $XY$ and $YZ$ bottom layers of the right wing has an inner corner, contradiction. If the top $y$ layer is a rectangle with no inner corner, and has the same size as the bottom $y$ layer of the right wing, then the right wing is a rectangular cuboid and has no inner corner.

\hspace*{-0.5cm}
\resizebox{0.9\textwidth}{!}
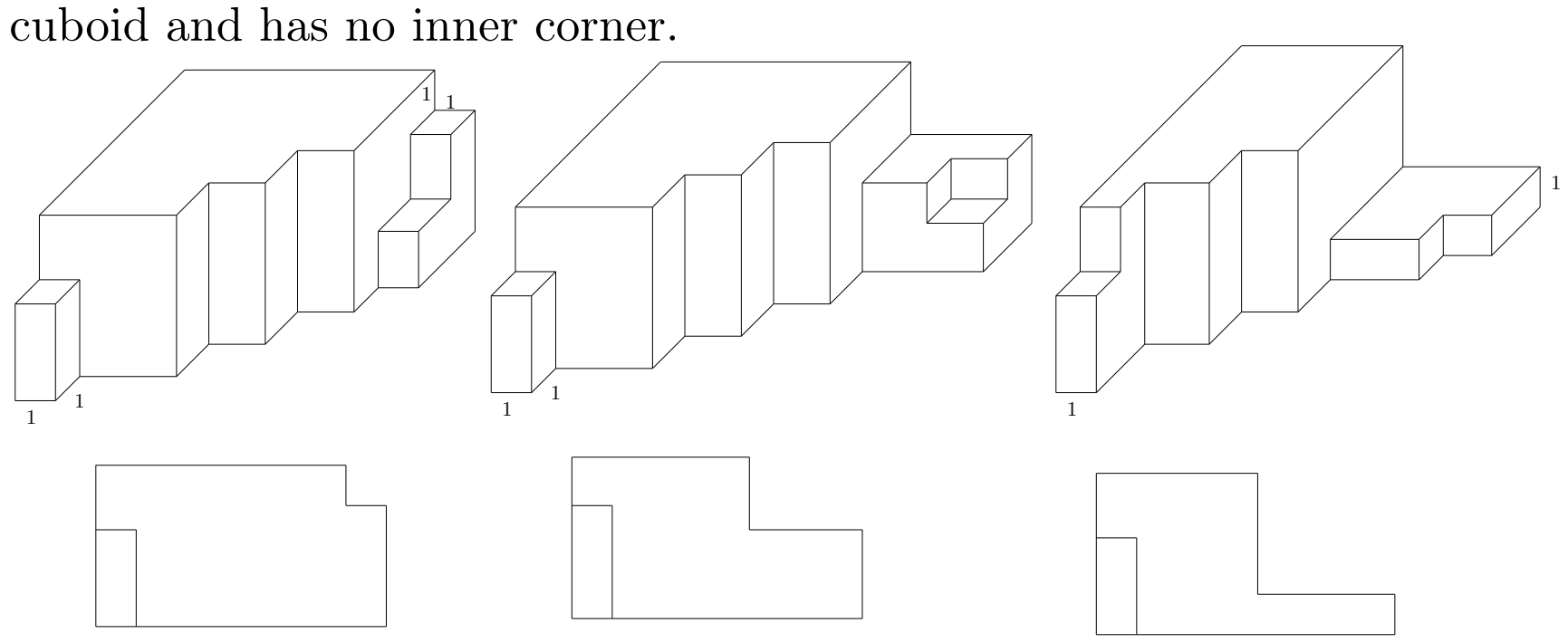
{%
\begin{circuitikz}
\tikzstyle{every node}=[font=\fontsize{18.2pt}{23.7pt}\selectfont]
\draw (14.125,10.875) to[short] (18,10.875);
\draw (11.875,8.625) to[short] (14.125,10.875);
\draw (11.875,8.625) to[short] (14,8.625);
\draw (11.875,7.625) to[short] (12.5,7.625);
\draw (11.5,7.25) to[short] (11.875,7.625);
\draw (15.875,9.625) to[short] (16.75,9.625);
\draw (16.75,9.625) to[short] (18,10.875);
\draw (17.625,9.875) to[short] (18,10.25);
\draw (17.625,9.875) to[short] (18.25,9.875);
\draw (18,10.25) to[short] (18.625,10.25);
\draw (17.125,8.375) to[short] (17.125,7.5);
\draw (11.5,7.25) to[short] (11.5,5.75);
\draw (11.875,8.625) to[short] (11.875,7.625);
\draw (15.875,9.625) to[short] (15.875,7.125);
\draw (14,8.625) to[short] (14,6.125);
\draw (14.5,9.125) to[short] (14.5,6.625);
\draw (14.5,6.625) to[short] (15.375,6.625);
\draw (15.375,6.625) to[short] (15.875,7.125);
\draw (15.875,7.125) to[short] (16.75,7.125);
\draw (16.75,9.625) to[short] (16.75,7.125);
\draw (18,10.875) to[short] (18,10.25);
\draw (16.75,7.125) to[short] (17.125,7.5);
\draw (18.25,9.875) to[short] (18.625,10.25);
\draw (17.75,8.375) to[short] (17.75,7.5);
\draw (17.125,7.5) to[short] (17.75,7.5);
\draw (18.625,10.25) to[short] (18.625,8.375);
\draw (17.75,7.5) to[short] (18.625,8.375);
\draw (11.5,5.75) to[short] (12.125,5.75);
\draw (12.5,6.125) to[short] (14,6.125);
\draw (14,6.125) to[short] (14.5,6.625);
\node [font=\fontsize{9.1pt}{11.8pt}\selectfont, inner xsep=0.080cm, inner ysep=0.085cm, rounded corners=0.000cm] at (11.75,5.5) {$1$};
\draw (12.125,7.25) to[short] (12.125,5.75);
\draw (12.125,7.25) to[short] (12.5,7.625);
\draw (11.5,7.25) to[short] (12.125,7.25);
\draw (12.5,7.625) to[short] (12.5,6.125);
\draw (12.125,5.75) to[short] (12.5,6.125);
%\node [font=\fontsize{9.1pt}{11.8pt}\selectfont, inner xsep=0.080cm, inner ysep=0.085cm, rounded corners=0.000cm] at (17.875,-1.25) {$a$};
\draw (14,8.625) to[short] (14.5,9.125);
\draw (15.375,9.125) to[short] (15.375,6.625);
\draw (14.5,9.125) to[short] (15.375,9.125);
\draw (15.375,9.125) to[short] (15.875,9.625);
\draw (30.5,11.25) to[short] (33,11.25);
\draw (28,8.75) to[short] (30.5,11.25);
\draw (28,8.75) to[short] (28.625,8.75);
\draw (28,7.75) to[short] (28.625,7.75);
\draw (27.625,7.375) to[short] (28,7.75);
\draw (31.375,9.625) to[short] (33,11.25);
\draw (31.875,8.25) to[short] (33,9.375);
\draw (31.875,8.25) to[short] (33.25,8.25);
\draw (31.875,7.625) to[short] (33.25,7.625);
\draw (33.25,8.25) to[short] (33.25,7.625);
\draw (27.625,7.375) to[short] (27.625,5.875);
\draw (28,8.75) to[short] (28,7.75);
\draw (31.375,9.625) to[short] (31.375,7.125);
\draw (28.625,8.75) to[short] (28.625,7.75);
\draw (33,11.25) to[short] (33,9.375);
\draw (31.875,8.25) to[short] (31.875,7.625);
\draw (31.375,7.125) to[short] (31.875,7.625);
\draw (34.375,8) to[short] (35.125,8.75);
\draw (35.125,9.375) to[short] (35.125,8.75);
\draw (33,9.375) to[short] (35.125,9.375);
\draw (34.375,8.625) to[short] (35.125,9.375);
\draw (27.625,5.875) to[short] (28.25,5.875);
\draw (28.625,8.75) to[short] (29,9.125);
\draw (28.25,7.375) to[short] (28.25,5.875);
\draw (28.25,7.375) to[short] (28.625,7.75);
\draw (27.625,7.375) to[short] (28.25,7.375);
\draw (28.25,5.875) to[short] (29,6.625);
\draw (29,9.125) to[short] (29,6.625);
\draw (29,9.125) to[short] (30,9.125);
\draw (29,6.625) to[short] (30,6.625);
\draw (12.75,4.75) to[short] (16.625,4.75);
\draw (12.75,4.75) to[short] (12.75,2.25);
\draw (13.375,3.75) to[short] (13.375,2.25);
\draw (12.75,3.75) to[short] (13.375,3.75);
\draw (12.75,2.25) to[short] (13.375,2.25);
\draw (17.25,4.125) to[short] (17.25,2.25);
\draw (16.625,4.75) to[short] (16.625,4.125);
\draw (16.625,4.125) to[short] (17.25,4.125);
\draw (13.375,2.25) to[short] (17.25,2.25);
\draw (28.25,4.625) to[short] (28.25,2.125);
\draw (28.875,3.625) to[short] (28.875,2.125);
\draw (28.25,3.625) to[short] (28.875,3.625);
\draw (28.25,2.125) to[short] (28.875,2.125);
\draw (30.75,4.625) to[short] (30.75,2.75);
\draw (28.875,2.125) to[short] (32.875,2.125);
\draw (30.75,2.75) to[short] (32.875,2.75);
\draw (32.875,2.75) to[short] (32.875,2.125);
\draw (18.25,9.875) to[short] (18.25,8.875);
\draw (17.75,8.375) to[short] (18.25,8.875);
\draw (17.625,8.875) to[short] (18.25,8.875);
\draw (17.125,8.375) to[short] (17.75,8.375);
\draw (17.625,9.875) to[short] (17.625,8.875);
\draw (17.125,8.375) to[short] (17.625,8.875);
\draw (34.375,8.625) to[short] (34.375,8);
\draw (33.625,8.625) to[short] (34.375,8.625);
\draw (33.25,8.25) to[short] (33.625,8.625);
\draw (33.25,7.625) to[short] (33.625,8);
\draw (33.625,8) to[short] (34.375,8);
\draw (33.625,8.625) to[short] (33.625,8);
\node [font=\fontsize{9.1pt}{11.8pt}\selectfont, inner xsep=0.080cm, inner ysep=0.085cm, rounded corners=0.000cm] at (12.5,5.75) {$1$};
\node [font=\fontsize{9.1pt}{11.8pt}\selectfont, inner xsep=0.080cm, inner ysep=0.085cm, rounded corners=0.000cm] at (18.25,10.375) {$1$};
\node [font=\fontsize{9.1pt}{11.8pt}\selectfont, inner xsep=0.080cm, inner ysep=0.085cm, rounded corners=0.000cm] at (27.875,5.625) {$1$};
\node [font=\fontsize{9.1pt}{11.8pt}\selectfont, inner xsep=0.080cm, inner ysep=0.085cm, rounded corners=0.000cm] at (35.375,9.125) {$1$};
\node [font=\fontsize{9.1pt}{11.8pt}\selectfont, inner xsep=0.080cm, inner ysep=0.085cm, rounded corners=0.000cm] at (17.875,10.5) {$1$};
\draw (30,9.125) to[short] (30.5,9.625);
\draw (30,9.125) to[short] (30,6.625);
\draw (30.5,9.625) to[short] (31.375,9.625);
\draw (30.5,9.625) to[short] (30.5,7.125);
\draw (30.5,7.125) to[short] (31.375,7.125);
\draw (30,6.625) to[short] (30.5,7.125);
\draw (28.25,4.625) to[short] (30.75,4.625);
\draw (21.5,11) to[short] (25.375,11);
\draw (19.25,8.75) to[short] (21.5,11);
\draw (19.25,8.75) to[short] (21.375,8.75);
\draw (19.25,7.75) to[short] (19.875,7.75);
\draw (18.875,7.375) to[short] (19.25,7.75);
\draw (23.25,9.75) to[short] (24.125,9.75);
\draw (24.125,9.75) to[short] (25.375,11);
\draw (24.625,9.125) to[short] (25.375,9.875);
\draw (24.625,9.125) to[short] (25.625,9.125);
\draw (25.375,9.875) to[short] (27.25,9.875);
\draw (24.625,9.125) to[short] (24.625,7.75);
\draw (18.875,7.375) to[short] (18.875,5.875);
\draw (19.25,8.75) to[short] (19.25,7.75);
\draw (23.25,9.75) to[short] (23.25,7.25);
\draw (21.375,8.75) to[short] (21.375,6.25);
\draw (21.875,9.25) to[short] (21.875,6.75);
\draw (21.875,6.75) to[short] (22.75,6.75);
\draw (22.75,6.75) to[short] (23.25,7.25);
\draw (23.25,7.25) to[short] (24.125,7.25);
\draw (24.125,9.75) to[short] (24.125,7.25);
\draw (25.375,11) to[short] (25.375,9.875);
\draw (24.125,7.25) to[short] (24.625,7.75);
\draw (26.875,9.5) to[short] (27.25,9.875);
\draw (25.625,9.125) to[short] (25.625,8.5);
\draw (24.625,7.75) to[short] (26.5,7.75);
\draw (27.25,9.875) to[short] (27.25,8.5);
\draw (26.5,7.75) to[short] (27.25,8.5);
\draw (18.875,5.875) to[short] (19.5,5.875);
\draw (19.875,6.25) to[short] (21.375,6.25);
\draw (21.375,6.25) to[short] (21.875,6.75);
\node [font=\fontsize{9.1pt}{11.8pt}\selectfont, inner xsep=0.080cm, inner ysep=0.085cm, rounded corners=0.000cm] at (19.125,5.625) {$1$};
\draw (19.5,7.375) to[short] (19.5,5.875);
\draw (19.5,7.375) to[short] (19.875,7.75);
\draw (18.875,7.375) to[short] (19.5,7.375);
\draw (19.875,7.75) to[short] (19.875,6.25);
\draw (19.5,5.875) to[short] (19.875,6.25);
\draw (21.375,8.75) to[short] (21.875,9.25);
\draw (22.75,9.25) to[short] (22.75,6.75);
\draw (21.875,9.25) to[short] (22.75,9.25);
\draw (22.75,9.25) to[short] (23.25,9.75);
\draw (20.125,4.875) to[short] (22.875,4.875);
\draw (20.125,4.875) to[short] (20.125,2.375);
\draw (20.75,4.125) to[short] (20.75,2.375);
\draw (20.125,4.125) to[short] (20.75,4.125);
\draw (20.125,2.375) to[short] (20.75,2.375);
\draw (24.625,3.75) to[short] (24.625,2.375);
\draw (22.875,4.875) to[short] (22.875,3.75);
\draw (22.875,3.75) to[short] (24.625,3.75);
\draw (20.75,2.375) to[short] (24.625,2.375);
\draw (26,9.5) to[short] (26.875,9.5);
\draw (26.875,9.5) to[short] (26.875,8.875);
\draw (25.625,9.125) to[short] (26,9.5);
\node [font=\fontsize{9.1pt}{11.8pt}\selectfont, inner xsep=0.080cm, inner ysep=0.085cm, rounded corners=0.000cm] at (19.875,5.875) {$1$};
\draw (26,9.5) to[short] (26,8.875);
\draw (26,8.875) to[short] (26.875,8.875);
\draw (25.625,8.5) to[short] (26,8.875);
\draw (25.625,8.5) to[short] (26.5,8.5);
\draw (26.5,8.5) to[short] (26.875,8.875);
\draw (26.5,8.5) to[short] (26.5,7.75);
\end{circuitikz}
}%

Consider the $y$ axis as the direction to slice the staircase into layers. In the first two pictures, the top $y$ layer have at least $2$ outer corners, and the bottom layer has an inner corner in region $C$. By Lemma \ref{nocornerinC}, we have a contradiction. So we can also rule out the first two pictures. So there's no three dimensional inner corner in the right wing.

\item 
If the bottom layer of the wings are a row of boxes and a column of boxes, as classified in Corollary \ref{bottomlayerofwings} case (2), (3), then we have two wings that look like "walls" with thickness $1$, either folded or spread out. Suppose for contradiction that there are some inner corners on either one of wings. 

In Corollary \ref{bottomlayerofwings} case (3), the wings are "folded" as in the picture on the left. Suppose without loss of generality that the right wing has an inner corner, then the top $XZ$ layer has at least two outer corners, and the bottom $XZ$ layer has an inner corner in region $C$, as in Lemma \ref{nocornerinC}, contradiction.

In Corollary \ref{bottomlayerofwings} case (2), the wings are "spread out" as in the picture on the right. Suppose without loss of generality that the right wing has an inner corner. Consider slicing the staircases in the $x$ direction, so the bottom layer on the $YZ$ plane has at least $2$ outer corners. The top $x$ layer of the staircase left wing is a vertical column of boxes, with one unique outer corner. So this contradicts Lemma \ref{onerow}.

In both cases, we have contradiction.

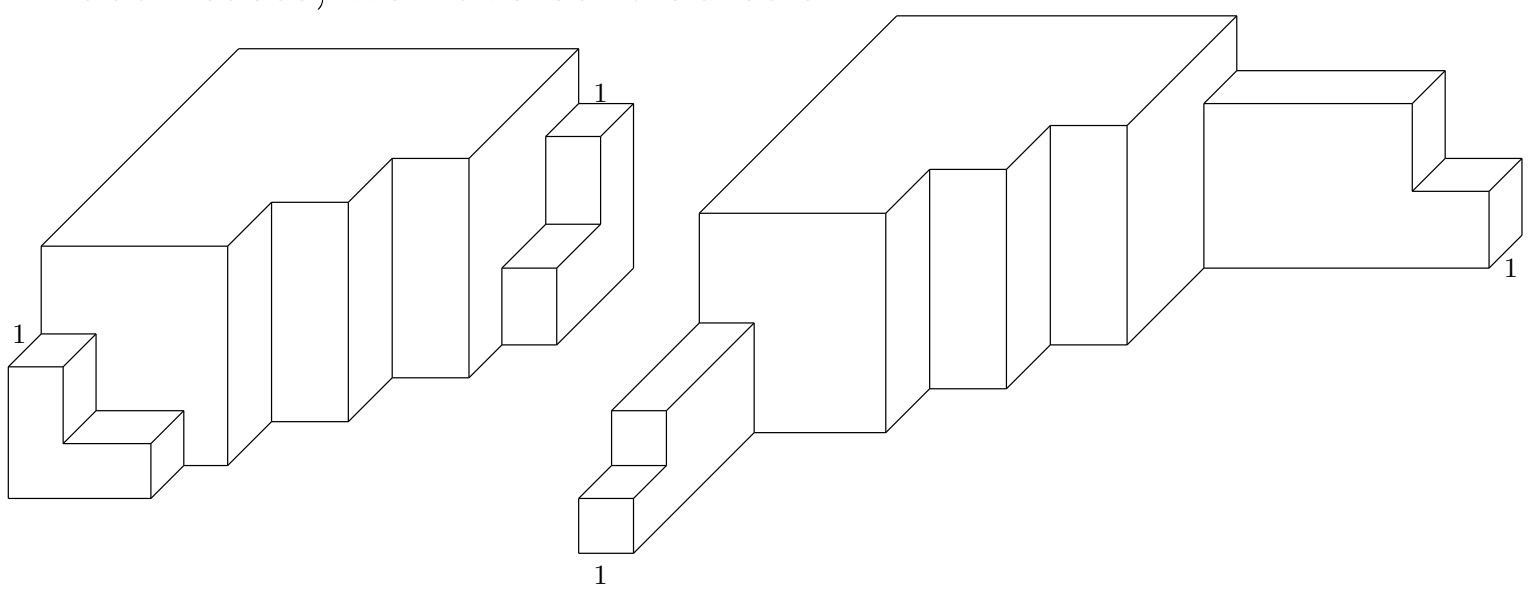
\resizebox{0.8\textwidth}{!}{%
\begin{circuitikz}
\tikzstyle{every node}=[font=\fontsize{9.1pt}{11.8pt}\selectfont]
\draw (14.125,10.875) to[short] (18,10.875);
\draw (11.875,8.625) to[short] (14.125,10.875);
\draw (11.875,8.625) to[short] (14,8.625);
\draw (11.875,7.625) to[short] (12.5,7.625);
\draw (11.5,7.25) to[short] (11.875,7.625);
\draw (15.875,9.625) to[short] (16.75,9.625);
\draw (16.75,9.625) to[short] (18,10.875);
\draw (17.625,9.875) to[short] (18,10.25);
\draw (17.625,9.875) to[short] (18.25,9.875);
\draw (18,10.25) to[short] (18.625,10.25);
\draw (17.125,8.375) to[short] (17.125,7.5);
\draw (11.5,7.25) to[short] (11.5,5.75);
\draw (11.875,8.625) to[short] (11.875,7.625);
\draw (15.875,9.625) to[short] (15.875,7.125);
\draw (14,8.625) to[short] (14,6.125);
\draw (14.5,9.125) to[short] (14.5,6.625);
\draw (14.5,6.625) to[short] (15.375,6.625);
\draw (15.375,6.625) to[short] (15.875,7.125);
\draw (15.875,7.125) to[short] (16.75,7.125);
\draw (16.75,9.625) to[short] (16.75,7.125);
\draw (18,10.875) to[short] (18,10.25);
\draw (16.75,7.125) to[short] (17.125,7.5);
\draw (18.25,9.875) to[short] (18.625,10.25);
\draw (17.75,8.375) to[short] (17.75,7.5);
\draw (17.125,7.5) to[short] (17.75,7.5);
\draw (18.625,10.25) to[short] (18.625,8.375);
\draw (17.75,7.5) to[short] (18.625,8.375);
\draw (11.5,5.75) to[short] (13.125,5.75);
\draw (13.5,6.125) to[short] (14,6.125);
\draw (14,6.125) to[short] (14.5,6.625);
\draw (13.125,6.375) to[short] (13.125,5.75);
\draw (12.125,7.25) to[short] (12.5,7.625);
\draw (11.5,7.25) to[short] (12.125,7.25);
\draw (12.5,7.625) to[short] (12.5,6.75);
\draw (13.125,5.75) to[short] (13.5,6.125);
%\node [font=\fontsize{9.1pt}{11.8pt}\selectfont, inner xsep=0.080cm, inner ysep=0.085cm, rounded corners=0.000cm] at (17.875,-1.25) {$a$};
\draw (14,8.625) to[short] (14.5,9.125);
\draw (15.375,9.125) to[short] (15.375,6.625);
\draw (14.5,9.125) to[short] (15.375,9.125);
\draw (15.375,9.125) to[short] (15.875,9.625);
\draw (18.25,9.875) to[short] (18.25,8.875);
\draw (17.75,8.375) to[short] (18.25,8.875);
\draw (17.625,8.875) to[short] (18.25,8.875);
\draw (17.125,8.375) to[short] (17.75,8.375);
\draw (17.625,9.875) to[short] (17.625,8.875);
\draw (17.125,8.375) to[short] (17.625,8.875);
\node [font=\fontsize{9.1pt}{11.8pt}\selectfont, inner xsep=0.080cm, inner ysep=0.085cm, rounded corners=0.000cm] at (18.25,10.375) {$1$};

\draw (12.125,7.25) to[short] (12.125,6.375);
\draw (12.125,6.375) to[short] (12.5,6.75);
\draw (12.125,6.375) to[short] (13.125,6.375);
\draw (13.125,6.375) to[short] (13.5,6.75);
\draw (13.5,6.75) to[short] (13.5,6.125);
\draw (12.5,6.75) to[short] (13.5,6.75);
\node [font=\fontsize{9.1pt}{11.8pt}\selectfont, inner xsep=0.080cm, inner ysep=0.085cm, rounded corners=0.000cm] at (11.625,7.625) {$1$};
\draw (21.625,11.25) to[short] (25.5,11.25);
\draw (19.375,9) to[short] (21.625,11.25);
\draw (19.375,9) to[short] (21.5,9);
\draw (19.375,7.75) to[short] (20,7.75);
\draw (18.375,6.75) to[short] (19.375,7.75);
\draw (23.375,10) to[short] (24.25,10);
\draw (24.25,10) to[short] (25.5,11.25);
\draw (25.125,10.25) to[short] (25.5,10.625);
\draw (25.125,10.25) to[short] (27.5,10.25);
\draw (25.5,10.625) to[short] (27.875,10.625);
\draw (19,6.75) to[short] (19,6.125);
\draw (19.375,9) to[short] (19.375,7.75);
\draw (23.375,10) to[short] (23.375,7.5);
\draw (21.5,9) to[short] (21.5,6.5);
\draw (22,9.5) to[short] (22,7);
\draw (22,7) to[short] (22.875,7);
\draw (22.875,7) to[short] (23.375,7.5);
\draw (23.375,7.5) to[short] (24.25,7.5);
\draw (24.25,10) to[short] (24.25,7.5);
\draw (25.5,11.25) to[short] (25.5,10.625);
\draw (24.25,7.5) to[short] (25.125,8.375);
\draw (27.5,10.25) to[short] (27.875,10.625);
\draw (27.875,10.625) to[short] (27.875,9.625);
\draw (18,5.75) to[short] (18.625,5.75);
\draw (20,6.5) to[short] (21.5,6.5);
\draw (21.5,6.5) to[short] (22,7);
\draw (18.625,5.75) to[short] (18.625,5.125);
\draw (19,6.75) to[short] (20,7.75);
\draw (18.375,6.75) to[short] (19,6.75);
\draw (20,7.75) to[short] (20,6.5);
\draw (18.625,5.125) to[short] (20,6.5);
\draw (21.5,9) to[short] (22,9.5);
\draw (22.875,9.5) to[short] (22.875,7);
\draw (22,9.5) to[short] (22.875,9.5);
\draw (22.875,9.5) to[short] (23.375,10);
\draw (27.5,10.25) to[short] (27.5,9.25);
\draw (27.5,9.25) to[short] (28.375,9.25);
\draw (25.125,10.25) to[short] (25.125,8.375);

\draw (18.375,6.75) to[short] (18.375,6.125);
\draw (18.625,5.75) to[short] (19,6.125);
\draw (18.375,6.125) to[short] (19,6.125);
\draw (18,5.75) to[short] (18,5.125);
\draw (18,5.125) to[short] (18.625,5.125);

\draw (27.5,9.25) to[short] (27.875,9.625);
\draw (28.375,9.25) to[short] (28.75,9.625);
\draw (27.875,9.625) to[short] (28.75,9.625);
\draw (28.75,9.625) to[short] (28.75,8.75);
\draw (28.375,9.25) to[short] (28.375,8.375);
\draw (28.375,8.375) to[short] (28.75,8.75);
\draw (25.125,8.375) to[short] (28.375,8.375);
\draw (18,5.75) to[short] (18.375,6.125);
\node [font=\fontsize{9.1pt}{11.8pt}\selectfont, inner xsep=0.080cm, inner ysep=0.085cm, rounded corners=0.000cm] at (18.25,4.875) {$1$};
\node [font=\fontsize{9.1pt}{11.8pt}\selectfont, inner xsep=0.080cm, inner ysep=0.085cm, rounded corners=0.000cm] at (28.625,8.375) {$1$};
\end{circuitikz}
}%

%then the top $XZ$ layer has at least two outer corners, and the bottom $XZ$ layer has an inner corner in region $C$, as in Lemma \ref{nocornerinC}, contradiction.

%Assume without loss of generality that the top-right wing has some inner corners. Then there are at least $2$ outer corners on the top-ring wing, and we can shift each of them to the inner corner right above the $y-z$ plane where the bottom-left wing connects to the main part of the opera house. So we have at least $2$ bounded connected component, contradiction.

%\item If the bottom layer of the bottom-right wing is a row of boxes and the top-left wing is a column of boxes. This situation is exactly the same as above, but the directions of the wings are in the perpendicular direction, and the wall-like wings are closely attached to the main part of the opera house. But we have the exact same argument as above. Assume without loss of generality that the top-right wing has two outer corners, and shift them to the inner corner right above the $y-z$ plane where the bottom-left wing connects to the main part of the opera house. So we have at least $2$ bounded connected component, contradiction.
\end{enumerate}

\end{proof}

\begin{lemma}
\label{height_of_wings}
For a minimized staircase $E$ with minimum tangent space $T_{nnp}$ dimension $1$, then one of the two wings must have heights either $0$ or $h_z-1$. (Recall that height is the maximum $z$ coordinate, not the number of layers.)
\end{lemma}

\begin{proof}
Assume for contradiction that both of the heights of the wings are not $0$ or $h_z-1$. Denote the height of the left wing as $m_x$, and the right wing as $m_y$. So $1 \leq m_x,m_y \leq l_z-2$. In this case, the vector $(-1,-1,h_z-1)$ also produces a bounded connected component. The bounded connected component is a subset of the region $A$ stacking on top of itself two layers. This collection of boxes, is not path connected to the infinite component $I_{\infty}$ because the boxes $\beta_1$ and $\beta_2$ both have height at least $1$.\end{proof}

At the end of the section, we conclude that all the Lemmas in this section together formulate exactly the definition of the opera houses. 

\begin{lemma}
An opera house has tangent space $T_{nnp}$ dimension exactly $1$.
\end{lemma}

\begin{proof}
The vector $(-1, -1, h_z)$ gives a bounded connected component in region $A$. For any other $nnp$ vector $(l_1, l_2, l_3)$, we want to show that there's no bounded connected component created.

For any $l_3 \geq h_z+1$, the entire shifted staircase $E -\alpha$ has negative $z$ coordinates so the bounded connected component is empty.

If either $l_1$ or $l_2$ is not $-1$, then whichever of the three shapes of the wings are, the finite component is path connected to the finite component. 

If $l_1=l_2=-1$ and $0 \leq l_3 \leq h_z-1$, then because one of the wings has height $0$ or $h_z-1$, the shifted staircase $E -\alpha$ is connected to the infinite component.
\end{proof}

\section{Classification when the top layers all have unique outer corners}

In this section, we consider staircases such that all the top $x,y,z$ layers have unique outer corners. If any of the top $x,y,z$ layers have two outer corners, they must be an opera house as classified in the last section.

\begin{lemma}
In the case of Lemma \ref{rectangle}(\ref{rec:4}), where the top layer has a unique outer corner $(a, b, h_z)$ and the bottom layer has largest $x$ coordinate $a+1$ and the largest $y$ coordinate $b+1$, the staircase is an opera house.   
\end{lemma}

\begin{proof}
There is not much variation for the shape of the middle layers, given that the top layer is only one column and one row different from the bottom layer. The staircase must consist of a rectangular cuboid which is the top layer stacking on itself $h_z+1$ layers, with the unique outer corner $(a, b, h_z)$. And the only part not completely determined by the top layers are the top $x$ and $y$ layers. Recall that the top $x$ and $y$ layers must both be rectangles, the staircase is an opera house with the top $x$ and $y$ layers being two wings.
\end{proof}

\begin{example}
An example of the above Lemma \monomialIdealPicture{{4,4,4,2},{4,4,4,2},{4,4,4},{1,1,1}}.   
\end{example}

\begin{corollary}
If a minimized staircase $E$ with minimal tangent space $T_{nnp}$ dimension $1$, and if $E$ is not an opera house, we can assume that the top and bottom layers of $E$ satisfy Lemma \ref{onerow} and Lemma \ref{rectangle} (\ref{rec:1})(\ref{rec:2})(\ref{rec:3}) but not (\ref{rec:4}), so from now on we can assume that the bottom layer has one unique inner corner. And by symmetry, all the bottom $x$ and $y$ layers also have one unique inner corner.
\end{corollary}

We now put some notations on the staircase $E$. As we will show in the process of giving notations, once we decide that the top $x,y,z$ layers each have one unique outer corner, and the bottom $x,y,z$ layers each have one inner corner, the coordinates of the inner corners of the bottom layers and the coordinates of the outer corners of the top layers are related. 

\begin{lemma}
\label{all_lengths_top_bottom}
For a minimized staircase $E$ with one unique outer corner for each of the top layers, and one unique inner corner for each of the bottom layers, there are only six variables $a,b,c,d,e,f$ together with three maximum lengths in each $x,y,z$ direction $h_z, h_y,h_z$ to completely characterize the top and bottom layers, as follows.  

The outer corner of the top $x$ layer is $(h_x, d, e)$; the outer corner of the top $y$ layer is $(c, h_y, f)$; the outer corner of the top $z$ layer is $(a, b, h_z)$. The inner corner of the bottom $x$ layer is $(0, b+1, e+1)$; the inner corner of the bottom $y$ layer is $(a+1, 0, f+1)$; the inner corner of the bottom $z$ layer is $(c+1,d+1, 0)$.
\end{lemma}

\begin{proof}
Denote the coordinates of the inner corner of the bottom $z$ layer as $(c+1,d+1,0)$. 

Denote the largest $x,y,z$ coordinates of the staircase as $h_x, h_y, h_z$. And denote the outer corners of the staircase as $(a, b, h_z), (c, h_y, f), (h_x, d, e)$ accordingly. Now we show that these completely determine the coordinates of the inner corners of the bottom $x,y,z$ layers. 

%Because we want the top $x,y,z$ layers to have one unique outer corner, the $h_y+1$th column and the $h_x+1$th row of the bottom layer must stack itself repeatedly in the $z$ direction, for the to $y$ and $x$ layers to be a rectangle. Denote the outer corners  $e+1$, and denote the number of times the $h_x+1$th row stack on itself as $f+1$. 

Consider the bottom $y$ layer on $XZ$ plane.

We know that the bottom $y$ layer must have the maximum $z$ coordinate $h_z$ and maximum $x$ coordinate $h_x$. We also know that the top $x$ layer's intersection with the $XZ$ plane, which is a vertical column of $f+1$ boxes, must also be in the bottom $y$ layer. Similarly, the top $z$ layer's intersection with the $XZ$ plane is a horizontal row of $a+1$ boxes, must also be in the bottom $y$ layer. Because the bottom $y$ layer must have one unique inner corner, and $f<h_z, $ by assumption, the inner corner must have coordinates $(a+1, 0, f+1)$.

The other inner corners are determined similarly.
\end{proof}

\begin{remark}
We draw some examples of the interaction of the coordinates of the inner and outer corners.

In the above Lemma \ref{all_lengths_top_bottom}, the lengths in the $x$ directions are $h_x, a, c$. The lengths in the $y$ directions are $h_y, b, d$. The lengths in the $z$ directions are $h_z, e, f$.
\end{remark}

\resizebox{0.8\textwidth}{!}{%
\begin{circuitikz}
\tikzstyle{every node}=[font=\fontsize{8.0pt}{10.4pt}\selectfont]
\draw (6.75,9.875) to[short] (6.375,9.5);
\draw (6.75,9.875) to[short] (4.625,9.875);
\draw (1.25,6) to[short] (3.375,6);
\draw (3.375,6.625) to[short] (3.375,6);
\draw (3.375,6) to[short] (5.125,7.75);
\draw (6.75,9.875) to[short] (6.75,8.125);
\draw (4.25,9.5) to[short] (6.375,9.5);
\draw (6.375,9.5) to[short] (6.375,7.75);
\draw (5.125,7.75) to[short] (6.375,7.75);
\draw (4.625,9.875) to[short] (4.25,9.5);
\draw (2.5,7.875) to[short] (2.5,9.75);
\draw (2.5,9.75) to[short] (3.75,9.75);
\draw (3.75,9.75) to[short] (3.75,7.875);
\draw (3.75,9.75) to[short] (4.625,10.625);
\draw (4.625,10.625) to[short] (4.625,9.875);
\draw (2.5,9.75) to[short] (3.375,10.625);
\draw (3.375,10.625) to[short] (4.625,10.625);
\draw [-{Stealth[scale=1.5]}, ] (1.25,6) -- (0.625,5.25);
\draw [-{Stealth[scale=1.5]}, ] (6.75,8.125) .. controls (7.5,8.125) and (7.375,8.125) .. (8,8.125) ;
\draw [-{Stealth[scale=1.5]}, ] (3.375,10.625) -- (3.375,11.5);
\draw (4.25,8.375) to[short] (3.75,7.875);
\draw (6.75,8.125) to[short] (6.375,7.75);
\draw (4.25,8.375) to[short] (5.125,8.375);
\draw (1.25,6.625) to[short] (3.375,6.625);
\draw (1.25,6.625) to[short] (1.25,6);
\draw (1.25,6.625) to[short] (2.5,7.875);
\draw (3.375,6.625) to[short] (5.125,8.375);
\draw (5.125,8.375) to[short] (5.125,7.75);
\draw (2.5,7.875) to[short] (3.75,7.875);
\draw (4.25,9.5) to[short] (4.25,8.375);
\draw (13.375,9) to[short] (11.875,7.5);
\draw (13.375,9) to[short] (11.875,9);
\draw (7.875,6.25) to[short] (9.125,6.25);
\draw (13.375,9) to[short] (13.375,8.375);
\draw (11.875,7.5) to[short] (11.875,6.875);
\draw (9.75,6.875) to[short] (11.875,6.875);
\draw (9.125,9) to[short] (9.125,10);
\draw (9.125,10) to[short] (11,10);
\draw (11,10) to[short] (11,8.125);
\draw (11,10) to[short] (11.875,10.875);
\draw (11.875,10.875) to[short] (11.875,9);
\draw (9.125,10) to[short] (10,10.875);
\draw (10,10.875) to[short] (11.875,10.875);
\draw [-{Stealth[scale=1.5]}, ] (7.875,6.25) -- (7.25,5.5);
\draw [-{Stealth[scale=1.5]}, ] (13.375,8.375) -- (14.875,8.375);
\draw [-{Stealth[scale=1.5]}, ] (10,10.875) -- (10,11.75);
\draw (10.375,9) to[short] (9.125,7.75);
\draw (13.375,8.375) to[short] (11.875,6.875);
\draw (9.75,7.5) to[short] (11.875,7.5);
\draw (7.875,7.75) to[short] (7.875,6.25);
\draw (7.875,7.75) to[short] (9.125,9);
\draw (9.75,7.5) to[short] (9.75,6.875);
\draw (9.125,7.75) to[short] (9.125,6.25);
\draw (7.875,7.75) to[short] (9.125,7.75);
\draw (9.125,9) to[short] (10.375,9);
\draw (9.75,6.875) to[short] (9.125,6.25);
\draw (11,8.125) to[short] (11.875,9);
\draw (10.375,8.125) to[short] (11,8.125);
\draw (10.375,8.125) to[short] (9.75,7.5);
\node [font=\fontsize{8.0pt}{10.4pt}\selectfont, fill={rgb,255:red,255; green,255; blue,255}, fill opacity=1, text opacity=1, inner xsep=0.080cm, inner ysep=0.085cm, rounded corners=0.020cm] at (2.875,10.625) {$a$};
\node [font=\fontsize{8.0pt}{10.4pt}\selectfont, fill={rgb,255:red,255; green,255; blue,255}, fill opacity=1, text opacity=1, inner xsep=0.080cm, inner ysep=0.085cm, rounded corners=0.020cm] at (4,10.875) {$b$};
\node [font=\fontsize{8.0pt}{10.4pt}\selectfont, fill={rgb,255:red,255; green,255; blue,255}, fill opacity=1, text opacity=1, inner xsep=0.080cm, inner ysep=0.085cm, rounded corners=0.020cm] at (6.875,7.75) {$c$};
\node [font=\fontsize{8.0pt}{10.4pt}\selectfont, fill={rgb,255:red,255; green,255; blue,255}, fill opacity=1, text opacity=1, inner xsep=0.080cm, inner ysep=0.085cm, rounded corners=0.020cm] at (2.125,5.75) {$d$};
\node [font=\fontsize{8.0pt}{10.4pt}\selectfont, fill={rgb,255:red,255; green,255; blue,255}, fill opacity=1, text opacity=1, inner xsep=0.080cm, inner ysep=0.085cm, rounded corners=0.020cm] at (1.625,6.25) {$f$};
\node [font=\fontsize{8.0pt}{10.4pt}\selectfont, fill={rgb,255:red,255; green,255; blue,255}, fill opacity=1, text opacity=1, inner xsep=0.080cm, inner ysep=0.085cm, rounded corners=0.020cm] at (7.125,8.875) {$e$};
\node [font=\fontsize{8.0pt}{10.4pt}\selectfont, fill={rgb,255:red,255; green,255; blue,255}, fill opacity=1, text opacity=1, inner xsep=0.080cm, inner ysep=0.085cm, rounded corners=0.020cm] at (10.75,11.125) {$b$};
\node [font=\fontsize{8.0pt}{10.4pt}\selectfont, fill={rgb,255:red,255; green,255; blue,255}, fill opacity=1, text opacity=1, inner xsep=0.080cm, inner ysep=0.085cm, rounded corners=0.020cm] at (9.375,10.75) {$a$};
\node [font=\fontsize{8.0pt}{10.4pt}\selectfont, fill={rgb,255:red,255; green,255; blue,255}, fill opacity=1, text opacity=1, inner xsep=0.080cm, inner ysep=0.085cm, rounded corners=0.020cm] at (7.5,6.75) {$f$};
\node [font=\fontsize{8.0pt}{10.4pt}\selectfont, fill={rgb,255:red,255; green,255; blue,255}, fill opacity=1, text opacity=1, inner xsep=0.080cm, inner ysep=0.085cm, rounded corners=0.020cm] at (8.5,6) {$d$};
\node [font=\fontsize{8.0pt}{10.4pt}\selectfont, fill={rgb,255:red,255; green,255; blue,255}, fill opacity=1, text opacity=1, inner xsep=0.080cm, inner ysep=0.085cm, rounded corners=0.020cm] at (13.125,7.5) {$c$};
\node [font=\fontsize{8.0pt}{10.4pt}\selectfont, fill={rgb,255:red,255; green,255; blue,255}, fill opacity=1, text opacity=1, inner xsep=0.080cm, inner ysep=0.085cm, rounded corners=0.020cm] at (13.75,8.75) {$e$};
\draw (10.375,9) to[short] (10.375,8.125);
\end{circuitikz}
}%

\begin{lemma}
\label{f=e}
If $f=e$ for a staircase $E$ described in Lemma \ref{all_lengths_top_bottom}, then $E$ is an opera house.
\end{lemma}

\begin{proof}
Consider the top and bottom $z$ layers, satisfying Lemma \ref{rectangle}(1)(2)(3), so either $a \leq c$ or $b \leq d$, or otherwise we must have (4). 

Suppose without loss of generality that $a \leq c$, and the proof the other case is exactly symmetric. We now show that the staircase is an opera house in the $y$ orientation by identifying the main part and the wings part of the opera house. We do so by dividing the boxes in the staircase into $3$ disjoint subsets, depending on $x$ and $z$ coordinate of the boxes.

Now we consider the top and bottom $y$ layers. The top layer has outer corner $(c, h_y, f)$, and the bottom layer has one unique inner corner $(a+1, 0,f+1)$. Because $e=f$ and $a \leq c$, the bottom $y$ layer in the $XZ$ plane and the top $y$ layer has the following shape as in the picture. 

Define the set $B_+$ as the boxes $b=(b_1, b_2,b_3)$ such that $0 \leq b_1 < a, e < b_3 \leq h_z$. Define the set $B_-$ as the boxes $b=(b_1, b_2,b_3)$ such that $c < b_1 \leq h_x, 0 \leq b_3 < e$. Define the set $B_0$ as the boxes $b=(b_1, b_2,b_3)$ such that $0 \leq b_1 \leq c, 0 \leq b_3 \leq e$, as shown in the picture. 

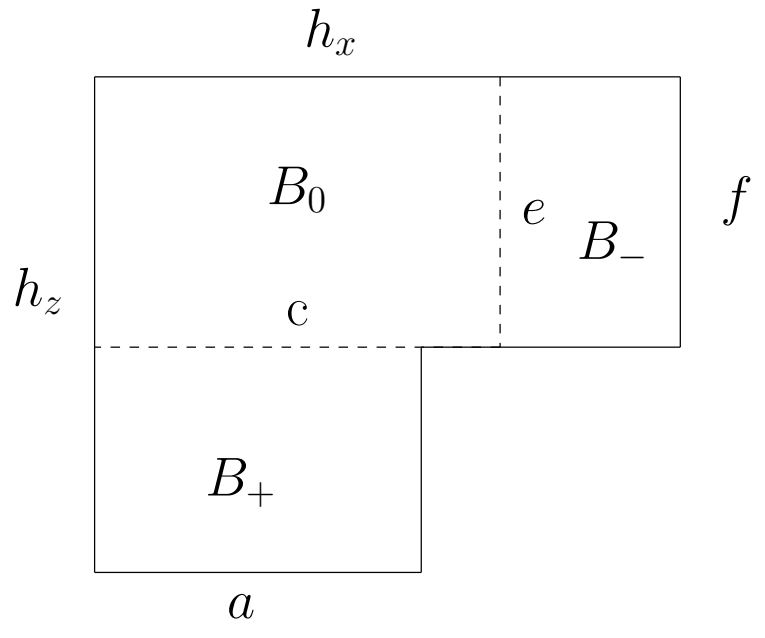
\resizebox{0.4\textwidth}{!}{%
\begin{circuitikz}
\tikzstyle{every node}=[font=\fontsize{18.2pt}{23.7pt}\selectfont]
\draw (23.75,16) to[short] (30.25,16);
\draw (23.75,16) to[short] (23.75,10.5);
\draw (23.75,10.5) to[short] (27.375,10.5);
\draw (27.375,10.5) to[short] (27.375,13);
\draw (27.375,13) to[short] (30.25,13);
\draw (30.25,16) to[short] (30.25,13);
\node [font=\fontsize{18.2pt}{23.7pt}\selectfont, fill={rgb,255:red,255; green,255; blue,255}, fill opacity=1, text opacity=1, inner xsep=0.080cm, inner ysep=0.085cm, rounded corners=0.020cm] at (25.375,10.125) {$a$};
\node [font=\fontsize{18.2pt}{23.7pt}\selectfont, fill={rgb,255:red,255; green,255; blue,255}, fill opacity=1, text opacity=1, inner xsep=0.080cm, inner ysep=0.085cm, rounded corners=0.000cm] at (30.875,14.625) {$f$};
\node [font=\fontsize{18.2pt}{23.7pt}\selectfont, fill={rgb,255:red,255; green,255; blue,255}, fill opacity=1, text opacity=1, inner xsep=0.080cm, inner ysep=0.085cm, rounded corners=0.020cm] at (26.375,16.5) {$h_x$};
\node [font=\fontsize{18.2pt}{23.7pt}\selectfont, fill={rgb,255:red,255; green,255; blue,255}, fill opacity=1, text opacity=1, inner xsep=0.080cm, inner ysep=0.085cm, rounded corners=0.020cm] at (23.125,13.625) {$h_z$};
\draw [dashed] (28.25,16) -- (28.25,13);
\draw [dashed] (28.25,13) -- (23.75,13);
\node [font=\fontsize{18.2pt}{23.7pt}\selectfont, fill={rgb,255:red,255; green,255; blue,255}, fill opacity=1, text opacity=1, inner xsep=0.080cm, inner ysep=0.085cm, rounded corners=0.020cm] at (26,13.375) {c};
\node [font=\fontsize{18.2pt}{23.7pt}\selectfont, fill={rgb,255:red,255; green,255; blue,255}, fill opacity=1, text opacity=1, inner xsep=0.080cm, inner ysep=0.085cm, rounded corners=0.020cm] at (28.625,14.5) {$e$};
\node [font=\fontsize{18.2pt}{23.7pt}\selectfont, fill={rgb,255:red,255; green,255; blue,255}, fill opacity=1, text opacity=1, inner xsep=0.080cm, inner ysep=0.085cm, rounded corners=0.020cm] at (25.375,11.5) {$B_+$};
\node [font=\fontsize{18.2pt}{23.7pt}\selectfont, fill={rgb,255:red,255; green,255; blue,255}, fill opacity=1, text opacity=1, inner xsep=0.080cm, inner ysep=0.085cm, rounded corners=0.020cm] at (29.5,14.125) {$B_-$};
\node [font=\fontsize{18.2pt}{23.7pt}\selectfont, fill={rgb,255:red,255; green,255; blue,255}, fill opacity=1, text opacity=1, inner xsep=0.080cm, inner ysep=0.085cm, rounded corners=0.020cm] at (26,14.75) {$B_0$};
\end{circuitikz}
}%

The main part of the opera house are all the boxes in $B_0$. Because top $y$ layer is a rectangle of size $c$ by $e$, this rectangle repeating itself $h_y+1$ times is a rectangular cuboid as the main part of the opera house.

The boxes in $B_+$ must be a rectangular cubiod. This is because the top $z$ layer of $B_+$ is also the top $z$ layer of the entire staircase. Because the top $z$ layer is the smallest possible partition among all layers in the $z$ direction, and the top $z$ layer is already reaching maximum coordinates $a, b$. %The $B_+$ wing of the opera house has size $a \times b \times (h_z-f)$.

By similar argument, the boxes in $B_-$ must also be a rectangular cubiod, as the top $x$ layer of $B_-$ is of maximum coordinates $f, d$.%, which has size $(h_x-c)\times d \times f$.

%If $e=f \geq 1$, by Corollary \ref{equiv_cond_rec}, we must have $\min(h_y-b, h_y-d)=1$, $\min(h_x-a, h_x-c)=1$, and $\min(h_z-f, h_z-e)=1$. Because we assumed that $a \leq c$, we must have $h_x-c=1$. Because $f=e$, we also must have that $h_z-f=1$. So the staircase is an opera house whose wings are folded. The condition $\min(h_y-b, h_y-d)=1$ ensures that at least one of the $y$ heights, either $b$ or $d$, of the wings are $h_y-1$, as classified in Lemma \ref{height_of_wings}. 

If $e=f=0$, by Corollary \ref{equiv_cond_onerow}, we must have that $\min(a+1, h_x-a, h_x-c)=1$, so either $a=0$ or $h_x-c=1$, which means that the bottom layer of the wings are either spread out ($f=a=0$) or a box and an arbitrary rectangle ($h_x-c=1, f=0$).

If $e=f \geq 1$ and $c \geq 1$, by Corollary \ref{equiv_cond_rec}, we must have that $\min(h_x-a, h_x-c)=1$, and $\min(h_z-f, h_z-e)=1$. Because we assumed that $a \leq c$, we must have $h_x-c=1$. Because $f=e$, we also must have that $h_z-f=1$. So the bottom layer of the wings are folded.

If $e=f \geq 1$ and $c = 0$, then $a=0$, by Corollary \ref{equiv_cond_onerow}, we must have that $\min(f+1, h_z-e, h_z-f)=1$, so $h_z-e=h_z-f=1$, so the bottom layer of the wings are of the type a box and a rectangle.

Lemma \ref{height_of_wings} does not assume that the main part of the opera house has at least two outer corners, so we can use it here, which ensures that the height of the wings also must be that either one of them is $0$ or $h_y-1$. 

%because the inner corner of the bottom $z$ layer is $(c, d, 1)$, and the top $z$ layer has size $a \times b$, we must have both $\min(h_y-b, h_y-d)=1$ and that $\min(h_x-a, h_x-c)=1$. Because the inner corner of the bottom $x$ layer is $(1, b, e)$, and the top $x$ layer has size $f\times d$, we must have both $\min(h_y-b, h_y-d)=1$ and that $\min(h_z-f, h_z-e)=1$. 

%Because the inner corner of the bottom $x$ layer is $(1, b, e)$, and the top $x$ layer has size $f\times d$, we must have both $\min(h_y-b, h_y-d)=1$ and that $\min(h_x-a, h_x-c)=1$.
\end{proof}
\begin{example}
The following are opera houses, where the $y$ direction (instead of the $z$ direction) is the vertical direction.

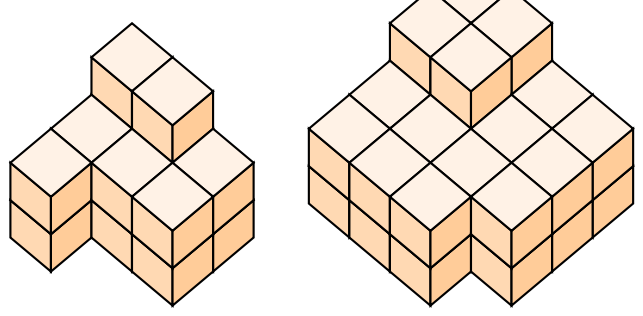
\monomialIdealPicture{{3,2,2},{3,2},{2,2}} \monomialIdealPicture{{3,3,2,2},{3,3,2,2},{2,2,2,2},{2,2,2}}    
\end{example}

\begin{definition}
Assume $f \neq e$ for a minimized staircase whose top layers all have one unique outer corner and the bottom layer all have one unique inner corner. If $f < e$, denote the staircase as \textbf{right-oriented}. If $f>e$, denote the staircase as \textbf{left-oriented}.
\end{definition}
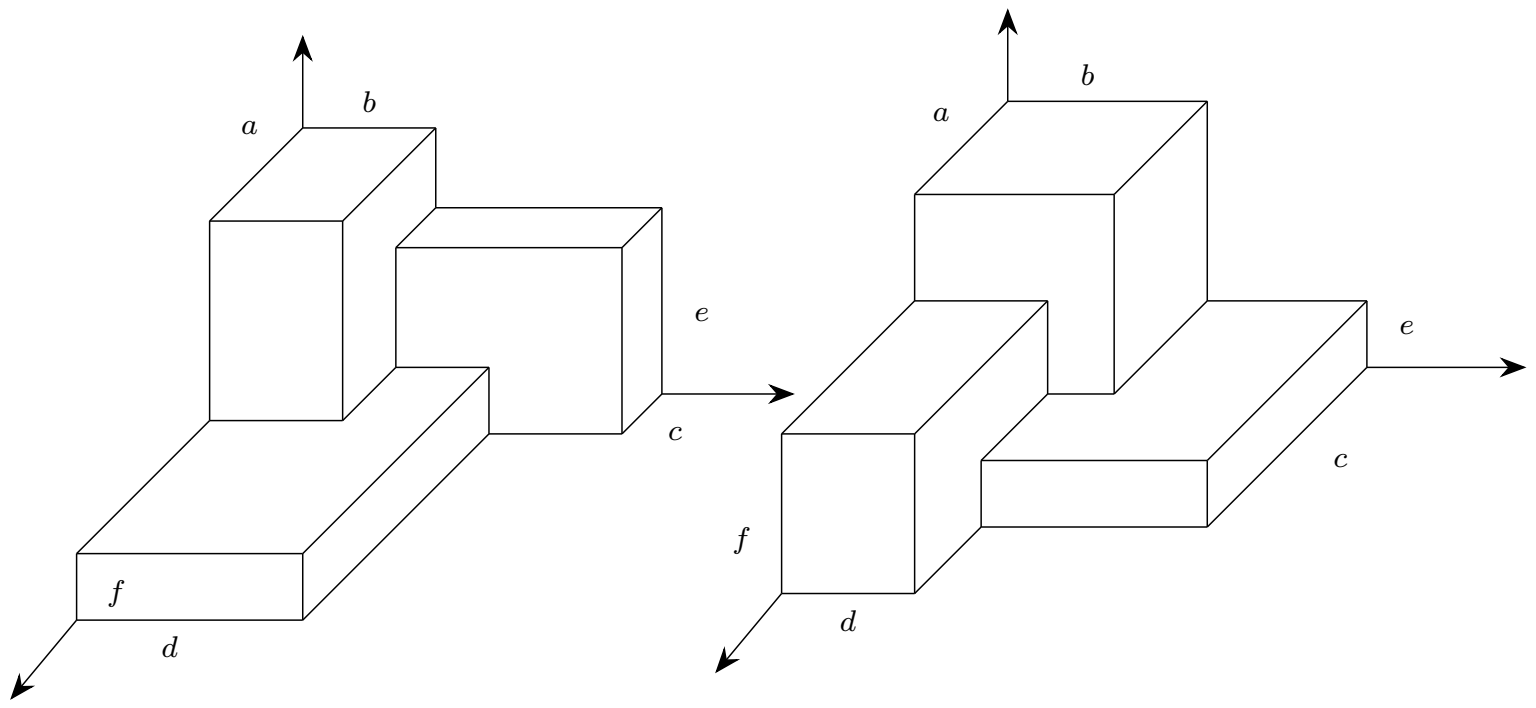
\resizebox{0.8\textwidth}{!}{%
\begin{circuitikz}
\tikzstyle{every node}=[font=\fontsize{8.0pt}{10.4pt}\selectfont]
\draw (6.75,9.875) to[short] (6.375,9.5);
\draw (6.75,9.875) to[short] (4.625,9.875);
\draw (1.25,6) to[short] (3.375,6);
\draw (3.375,6.625) to[short] (3.375,6);
\draw (3.375,6) to[short] (5.125,7.75);
\draw (6.75,9.875) to[short] (6.75,8.125);
\draw (4.25,9.5) to[short] (6.375,9.5);
\draw (6.375,9.5) to[short] (6.375,7.75);
\draw (5.125,7.75) to[short] (6.375,7.75);
\draw (4.625,9.875) to[short] (4.25,9.5);
\draw (2.5,7.875) to[short] (2.5,9.75);
\draw (2.5,9.75) to[short] (3.75,9.75);
\draw (3.75,9.75) to[short] (3.75,7.875);
\draw (3.75,9.75) to[short] (4.625,10.625);
\draw (4.625,10.625) to[short] (4.625,9.875);
\draw (2.5,9.75) to[short] (3.375,10.625);
\draw (3.375,10.625) to[short] (4.625,10.625);
\draw [-{Stealth[scale=1.5]}, ] (1.25,6) -- (0.625,5.25);
\draw [-{Stealth[scale=1.5]}, ] (6.75,8.125) .. controls (7.5,8.125) and (7.375,8.125) .. (8,8.125) ;
\draw [-{Stealth[scale=1.5]}, ] (3.375,10.625) -- (3.375,11.5);
\draw (4.25,8.375) to[short] (3.75,7.875);
\draw (6.75,8.125) to[short] (6.375,7.75);
\draw (4.25,8.375) to[short] (5.125,8.375);
\draw (1.25,6.625) to[short] (3.375,6.625);
\draw (1.25,6.625) to[short] (1.25,6);
\draw (1.25,6.625) to[short] (2.5,7.875);
\draw (3.375,6.625) to[short] (5.125,8.375);
\draw (5.125,8.375) to[short] (5.125,7.75);
\draw (2.5,7.875) to[short] (3.75,7.875);
\draw (4.25,9.5) to[short] (4.25,8.375);
\draw (13.375,9) to[short] (11.875,7.5);
\draw (13.375,9) to[short] (11.875,9);
\draw (7.875,6.25) to[short] (9.125,6.25);
\draw (13.375,9) to[short] (13.375,8.375);
\draw (11.875,7.5) to[short] (11.875,6.875);
\draw (9.75,6.875) to[short] (11.875,6.875);
\draw (9.125,9) to[short] (9.125,10);
\draw (9.125,10) to[short] (11,10);
\draw (11,10) to[short] (11,8.125);
\draw (11,10) to[short] (11.875,10.875);
\draw (11.875,10.875) to[short] (11.875,9);
\draw (9.125,10) to[short] (10,10.875);
\draw (10,10.875) to[short] (11.875,10.875);
\draw [-{Stealth[scale=1.5]}, ] (7.875,6.25) -- (7.25,5.5);
\draw [-{Stealth[scale=1.5]}, ] (13.375,8.375) -- (14.875,8.375);
\draw [-{Stealth[scale=1.5]}, ] (10,10.875) -- (10,11.75);
\draw (10.375,9) to[short] (9.125,7.75);
\draw (13.375,8.375) to[short] (11.875,6.875);
\draw (9.75,7.5) to[short] (11.875,7.5);
\draw (7.875,7.75) to[short] (7.875,6.25);
\draw (7.875,7.75) to[short] (9.125,9);
\draw (9.75,7.5) to[short] (9.75,6.875);
\draw (9.125,7.75) to[short] (9.125,6.25);
\draw (7.875,7.75) to[short] (9.125,7.75);
\draw (9.125,9) to[short] (10.375,9);
\draw (9.75,6.875) to[short] (9.125,6.25);
\draw (11,8.125) to[short] (11.875,9);
\draw (10.375,8.125) to[short] (11,8.125);
\draw (10.375,8.125) to[short] (9.75,7.5);
\node [font=\fontsize{8.0pt}{10.4pt}\selectfont, fill={rgb,255:red,255; green,255; blue,255}, fill opacity=1, text opacity=1, inner xsep=0.080cm, inner ysep=0.085cm, rounded corners=0.020cm] at (2.875,10.625) {$a$};
\node [font=\fontsize{8.0pt}{10.4pt}\selectfont, fill={rgb,255:red,255; green,255; blue,255}, fill opacity=1, text opacity=1, inner xsep=0.080cm, inner ysep=0.085cm, rounded corners=0.020cm] at (4,10.875) {$b$};
\node [font=\fontsize{8.0pt}{10.4pt}\selectfont, fill={rgb,255:red,255; green,255; blue,255}, fill opacity=1, text opacity=1, inner xsep=0.080cm, inner ysep=0.085cm, rounded corners=0.020cm] at (6.875,7.75) {$c$};
\node [font=\fontsize{8.0pt}{10.4pt}\selectfont, fill={rgb,255:red,255; green,255; blue,255}, fill opacity=1, text opacity=1, inner xsep=0.080cm, inner ysep=0.085cm, rounded corners=0.020cm] at (2.125,5.75) {$d$};
\node [font=\fontsize{8.0pt}{10.4pt}\selectfont, fill={rgb,255:red,255; green,255; blue,255}, fill opacity=1, text opacity=1, inner xsep=0.080cm, inner ysep=0.085cm, rounded corners=0.020cm] at (1.625,6.25) {$f$};
\node [font=\fontsize{8.0pt}{10.4pt}\selectfont, fill={rgb,255:red,255; green,255; blue,255}, fill opacity=1, text opacity=1, inner xsep=0.080cm, inner ysep=0.085cm, rounded corners=0.020cm] at (7.125,8.875) {$e$};
\node [font=\fontsize{8.0pt}{10.4pt}\selectfont, fill={rgb,255:red,255; green,255; blue,255}, fill opacity=1, text opacity=1, inner xsep=0.080cm, inner ysep=0.085cm, rounded corners=0.020cm] at (10.75,11.125) {$b$};
\node [font=\fontsize{8.0pt}{10.4pt}\selectfont, fill={rgb,255:red,255; green,255; blue,255}, fill opacity=1, text opacity=1, inner xsep=0.080cm, inner ysep=0.085cm, rounded corners=0.020cm] at (9.375,10.75) {$a$};
\node [font=\fontsize{8.0pt}{10.4pt}\selectfont, fill={rgb,255:red,255; green,255; blue,255}, fill opacity=1, text opacity=1, inner xsep=0.080cm, inner ysep=0.085cm, rounded corners=0.020cm] at (7.5,6.75) {$f$};
\node [font=\fontsize{8.0pt}{10.4pt}\selectfont, fill={rgb,255:red,255; green,255; blue,255}, fill opacity=1, text opacity=1, inner xsep=0.080cm, inner ysep=0.085cm, rounded corners=0.020cm] at (8.5,6) {$d$};
\node [font=\fontsize{8.0pt}{10.4pt}\selectfont, fill={rgb,255:red,255; green,255; blue,255}, fill opacity=1, text opacity=1, inner xsep=0.080cm, inner ysep=0.085cm, rounded corners=0.020cm] at (13.125,7.5) {$c$};
\node [font=\fontsize{8.0pt}{10.4pt}\selectfont, fill={rgb,255:red,255; green,255; blue,255}, fill opacity=1, text opacity=1, inner xsep=0.080cm, inner ysep=0.085cm, rounded corners=0.020cm] at (13.75,8.75) {$e$};
\draw (10.375,9) to[short] (10.375,8.125);
\end{circuitikz}
}%

\begin{example}
The first staircase is left-oriented, and the second one is right-oriented.

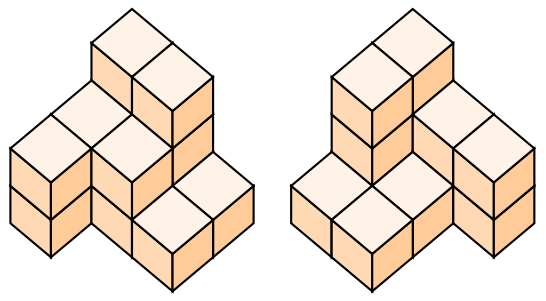
\monomialIdealPicture{{3,2,2},{3,2},{1,1}}\monomialIdealPicture{{3,3,1},{2,1,1},{2}}
\end{example}

\begin{lemma}
\label{left_right_orientation_inequality}
Assume that $E$ is a minimized staircase whose top layers all have one unique outer corner, and bottom layers all have one unique inner corner. Assume that $E$ is not an opera house and left-oriented, then $E$ must also satisfy that $c >a$ and $b>d$, other than $f>e$. Similarly, a right-oriented staircase $E$ also satisfy that $c<a$ and $b<d$, other than $f<e$.
\end{lemma}

\begin{proof}
By Lemma \ref{f=e}, we can also assume that $a \neq c$ and $b \neq d$, otherwise the staircase is an opera house.

Assume that the staircase is left-oriented so $f >e$. Assume for contradiction that we also add the condition that $a>c$, then we want to show that in either the case of $b >d$ or $b <d$, we have an opera house. As shown in the following pictures, if $b >d$, then we draw the top and bottom $z$ layers; if $b <d$, then we draw the top and bottom $x$ layers. And again we divide all the boxes in the staircase into three sets, $B_0, B_+, B_-$, and $B_0$ is the main part of the opera house in both cases. %The $a \times b$ rectangle as the top $z$ layer repeating itself $h_z$ times is a rectangular cuboid as the main part of the opera house in the case when $b > d$; the $d \times f$ rectangle as the top $x$ layer repeating itself $h_x$ times is a rectangular cuboid as the main part of the opera house in the case when $b < d$. 
Now we prove that $B_+,B_-$ are the wings.

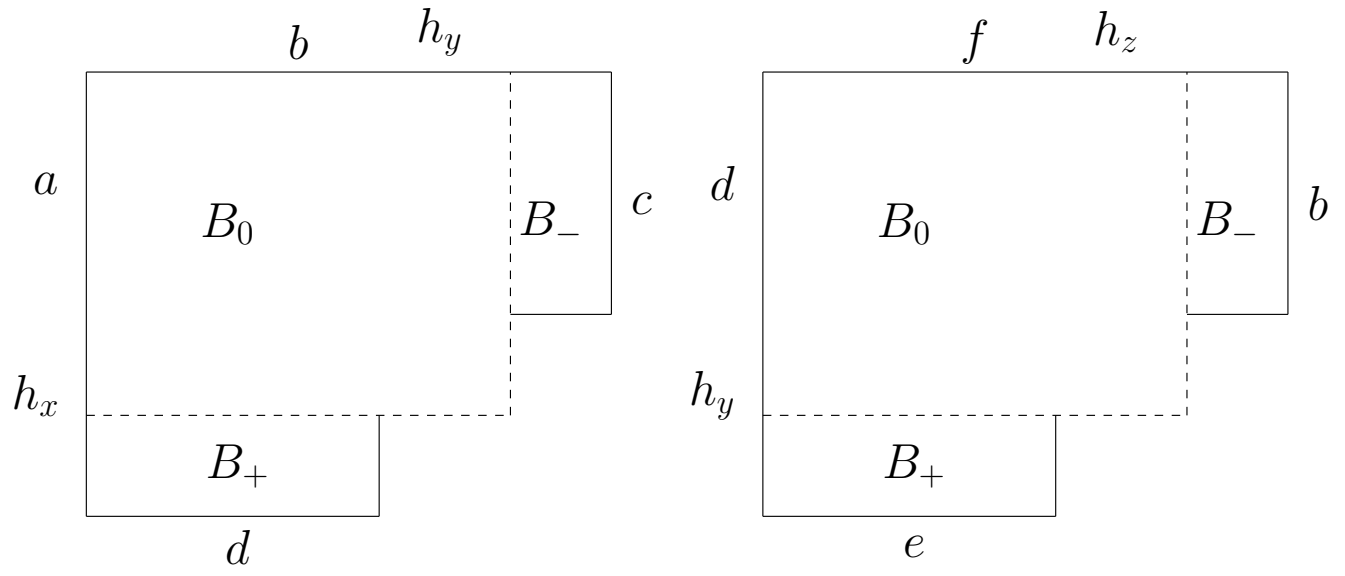
\resizebox{0.7\textwidth}{!}{%
\begin{circuitikz}
\tikzstyle{every node}=[font=\fontsize{18.2pt}{23.7pt}\selectfont]
\draw (23.75,16) to[short] (30.25,16);
\draw (23.75,16) to[short] (23.75,10.5);
\draw (23.75,10.5) to[short] (27.375,10.5);
\draw (27.375,10.5) to[short] (27.375,11.75);
\draw (29,13) to[short] (30.25,13);
\draw (30.25,16) to[short] (30.25,13);
\node [font=\fontsize{18.2pt}{23.7pt}\selectfont, fill={rgb,255:red,255; green,255; blue,255}, fill opacity=1, text opacity=1, inner xsep=0.080cm, inner ysep=0.085cm, rounded corners=0.020cm] at (23.125,12) {$h_x$};
\node [font=\fontsize{18.2pt}{23.7pt}\selectfont, fill={rgb,255:red,255; green,255; blue,255}, fill opacity=1, text opacity=1, inner xsep=0.080cm, inner ysep=0.085cm, rounded corners=0.020cm] at (28.125,16.5) {$h_y$};
\node [font=\fontsize{18.2pt}{23.7pt}\selectfont, fill={rgb,255:red,255; green,255; blue,255}, fill opacity=1, text opacity=1, inner xsep=0.080cm, inner ysep=0.085cm, rounded corners=0.020cm] at (25.625,11.125) {$B_+$};
\node [font=\fontsize{18.2pt}{23.7pt}\selectfont, fill={rgb,255:red,255; green,255; blue,255}, fill opacity=1, text opacity=1, inner xsep=0.080cm, inner ysep=0.085cm, rounded corners=0.020cm] at (29.5,14.125) {$B_-$};
\draw [dashed] (23.75,11.75) -- (29,11.75);
\draw [dashed] (29,11.75) -- (29,16);
\node [font=\fontsize{18.2pt}{23.7pt}\selectfont, fill={rgb,255:red,255; green,255; blue,255}, fill opacity=1, text opacity=1, inner xsep=0.080cm, inner ysep=0.085cm, rounded corners=0.020cm] at (23.25,14.625) {$a$};
\node [font=\fontsize{18.2pt}{23.7pt}\selectfont, fill={rgb,255:red,255; green,255; blue,255}, fill opacity=1, text opacity=1, inner xsep=0.080cm, inner ysep=0.085cm, rounded corners=0.020cm] at (26.375,16.375) {$b$};
\node [font=\fontsize{18.2pt}{23.7pt}\selectfont, fill={rgb,255:red,255; green,255; blue,255}, fill opacity=1, text opacity=1, inner xsep=0.080cm, inner ysep=0.085cm, rounded corners=0.020cm] at (30.625,14.375) {$c$};
\node [font=\fontsize{18.2pt}{23.7pt}\selectfont, fill={rgb,255:red,255; green,255; blue,255}, fill opacity=1, text opacity=1, inner xsep=0.080cm, inner ysep=0.085cm, rounded corners=0.020cm] at (25.625,10.125) {$d$};
\node [font=\fontsize{18.2pt}{23.7pt}\selectfont, fill={rgb,255:red,255; green,255; blue,255}, fill opacity=1, text opacity=1, inner xsep=0.080cm, inner ysep=0.085cm, rounded corners=0.020cm] at (25.5,14.125) {$B_0$};
\draw (32.125,16) to[short] (38.625,16);
\draw (32.125,16) to[short] (32.125,10.5);
\draw (32.125,10.5) to[short] (35.75,10.5);
\draw (35.75,10.5) to[short] (35.75,11.75);
\draw (37.375,13) to[short] (38.625,13);
\draw (38.625,16) to[short] (38.625,13);
\node [font=\fontsize{18.2pt}{23.7pt}\selectfont, fill={rgb,255:red,255; green,255; blue,255}, fill opacity=1, text opacity=1, inner xsep=0.080cm, inner ysep=0.085cm, rounded corners=0.020cm] at (31.5,12) {$h_y$};
\node [font=\fontsize{18.2pt}{23.7pt}\selectfont, fill={rgb,255:red,255; green,255; blue,255}, fill opacity=1, text opacity=1, inner xsep=0.080cm, inner ysep=0.085cm, rounded corners=0.020cm] at (36.5,16.5) {$h_z$};
\node [font=\fontsize{18.2pt}{23.7pt}\selectfont, fill={rgb,255:red,255; green,255; blue,255}, fill opacity=1, text opacity=1, inner xsep=0.080cm, inner ysep=0.085cm, rounded corners=0.020cm] at (34,11.125) {$B_+$};
\node [font=\fontsize{18.2pt}{23.7pt}\selectfont, fill={rgb,255:red,255; green,255; blue,255}, fill opacity=1, text opacity=1, inner xsep=0.080cm, inner ysep=0.085cm, rounded corners=0.020cm] at (37.875,14.125) {$B_-$};
\draw [dashed] (32.125,11.75) -- (37.375,11.75);
\draw [dashed] (37.375,11.75) -- (37.375,16);
\node [font=\fontsize{18.2pt}{23.7pt}\selectfont, fill={rgb,255:red,255; green,255; blue,255}, fill opacity=1, text opacity=1, inner xsep=0.080cm, inner ysep=0.085cm, rounded corners=0.020cm] at (31.625,14.625) {$d$};
\node [font=\fontsize{18.2pt}{23.7pt}\selectfont, fill={rgb,255:red,255; green,255; blue,255}, fill opacity=1, text opacity=1, inner xsep=0.080cm, inner ysep=0.085cm, rounded corners=0.000cm] at (34.75,16.375) {$f$};
\node [font=\fontsize{18.2pt}{23.7pt}\selectfont, fill={rgb,255:red,255; green,255; blue,255}, fill opacity=1, text opacity=1, inner xsep=0.080cm, inner ysep=0.085cm, rounded corners=0.020cm] at (39,14.375) {$b$};
\node [font=\fontsize{18.2pt}{23.7pt}\selectfont, fill={rgb,255:red,255; green,255; blue,255}, fill opacity=1, text opacity=1, inner xsep=0.080cm, inner ysep=0.085cm, rounded corners=0.020cm] at (34,10.125) {$e$};
\node [font=\fontsize{18.2pt}{23.7pt}\selectfont, fill={rgb,255:red,255; green,255; blue,255}, fill opacity=1, text opacity=1, inner xsep=0.080cm, inner ysep=0.085cm, rounded corners=0.020cm] at (33.875,14.125) {$B_0$};
\end{circuitikz}
}%

%In the case of $b>d$, the $B_-$ wing has size $(h_y-b) \times c \times e$, and the $B_+$ wing has size $(h_x-a) \times d \times f$. 

%In the case of $b<d$, the $B_-$ wing has size $(h_z-f) \times b \times a$, and the $B_+$ wing has size $(h_y-d) \times e \times c$.

When $b>d$, by Corollary \ref{equiv_cond_rec}, we must have $\min(h_y-b, h_y-d)=1$, and $\min(h_x-a, h_x-c)=1$.%, and $\min(h_z-f, h_z-e)=1$. 
Because we assumed that $a > c$, we must have $h_x-a=1$. %Because $f>e$, we also must have that $h_z-f=1$.
Because $b>d$, we also must have that $h_y-b=1$. 
%When $b<d$, we also must have that $h_y-d=1$.
Lemma \ref{height_of_wings} ensures that the height of the wings also must be that either one of them is $0$ or $h_z-1$. 

Similarly for $b<d$, we must have that $h_y-d=1$ and $h_z-f=1$, and the height of the wings also must be that either one of them is $0$ or $h_x-1$.

So in both cases $b>d$ or $b<d$, the staircase is an opera house whose wings are folded.

%Finally, when $b>d$, the heights of the wings are $h_z-e$ and $h_z-f=1$. When $b<d$, the heights of the wings are $h_x-a$ and $h_x-c=1$. So the height conditions of the wings, as stated in Lemma \ref{height_of_wings} are also satisfyed.

%The boxes in $B_+$ must also be a rectangular cuboid. We know that the top $x$ layer of $B_+$ is also the top $z$ layer of the entire staircase. Because the inner corner of the bottom $y$ layer is assumed to be $(a+1, 1, f+1)$, so the longest row of the top $x$ layer is $f$, and because the top $x$ layer must be a rectangle with one unique outer corner, the top $x$ layer of $B_+$ is a $f \times d$ rectangle. The $B_+$ wing of the opera house has size $f \times d \times (h_x-a)$.

%By similar argument, the boxes in $B_-$ must also be a rectangular cuboid, which has size $c \times e \times (h_y-b)$.

A similar argument applies for the right-oriented staircases.
\end{proof}

%Once we know that the bottom $z$ layer has an unique inner corner $(n_1, n_2)$, then we know that the top $x$ layer has unique outer corner with $y$ coordinate $n_2$, and that the top $y$ layer has unique outer corner with $x$ coordinate $n_1$.

\begin{definition}
For a minimized staircase $E$ whose top $x,y,z$ layers has one unique outer corner, and the numerics of the outer corners are $a,b,c,d,e,f$ as in Lemma \ref{all_lengths_top_bottom}, consider all the boxes $e=(x,y,z)$ in the staircase $E$ that satisfy any one of the three conditions:
\begin{enumerate}
\item $x \leq a$ and  $y \leq b$ and  $z \leq h_z$.

\item $x \leq c$ and  $y \leq h_y$ and  $z \leq e$.

\item $x \leq h_x$ and  $y \leq d$ and  $z \leq f$.
\end{enumerate}

Then denote these boxes as the \textbf{rotating wall} part of the staircase $E$. Denote the rest of the boxes of the staircase $E$ as the \textbf{interior boxes}. Denote the smaller set of values, $b,c,f$ if right oriented, and $a,d,e$ if left oriented, as the thickness of the three walls.
\end{definition}

\begin{example}
In the example below, the rotating walls of 
the ideal $(x^5, y^5, z^5, x^4y, y^4z,x^2z^4,x^2y^2z,x^3yz,x^2yz^2)$ is the following. The picture on the left is the original ideal, and the picture in the middle is the rotating wall part. The interior boxes part is on the right.

\monomialIdealPicture  {{5,5,4,4,4},{5,5,3,2},{5,5,2,2},{5,5,1,1},{1,1,1,1}} \monomialIdealPicture{{5,5,4,4,4},{5,5,1,1},{5,5,1,1},{5,5,1,1},{1,1,1,1}} \monomialIdealPicture{{2,1},{1,1}}
  
\end{example}

Lemma \ref{left_right_orientation_inequality},  Corollaries \ref{equiv_cond_onerow} and \ref{equiv_cond_rec} are equivalent to the following.

\begin{lemma}
\label{rotatingwallscondition}
When a minimized staircase $E$ has one unique outer corner for each of the top $x,y,z$ layers, and one unique inner corner for each of the bottom layers, and $E$ is not an opera house, we have the following numerical conditions on the rotating walls part.

\begin{enumerate}
\item If the rotating walls all have thickness $0$, then they can have any arbitrary size.
\item If at least two of the three the rotating walls have thickness greater than or equal to $1$, then $a=h_x-1$,$d=h_y-1$ and $e=h_z-1$ if the staircase is right oriented, and $c=h_x-1$,$b=h_y-1$ and $f=h_z-1$ if the staircase is left oriented.

\item If one of the three rotating walls have thickness greater than or equal to $1$, then the two numerics adjacent to the thickness and in parallel with the thickness is fixed as the maximum coordinate of the stairacse minus 1. For example, setting $b>0$ fixes $a=h_x-1$,$d=h_y-1$.
\end{enumerate} 

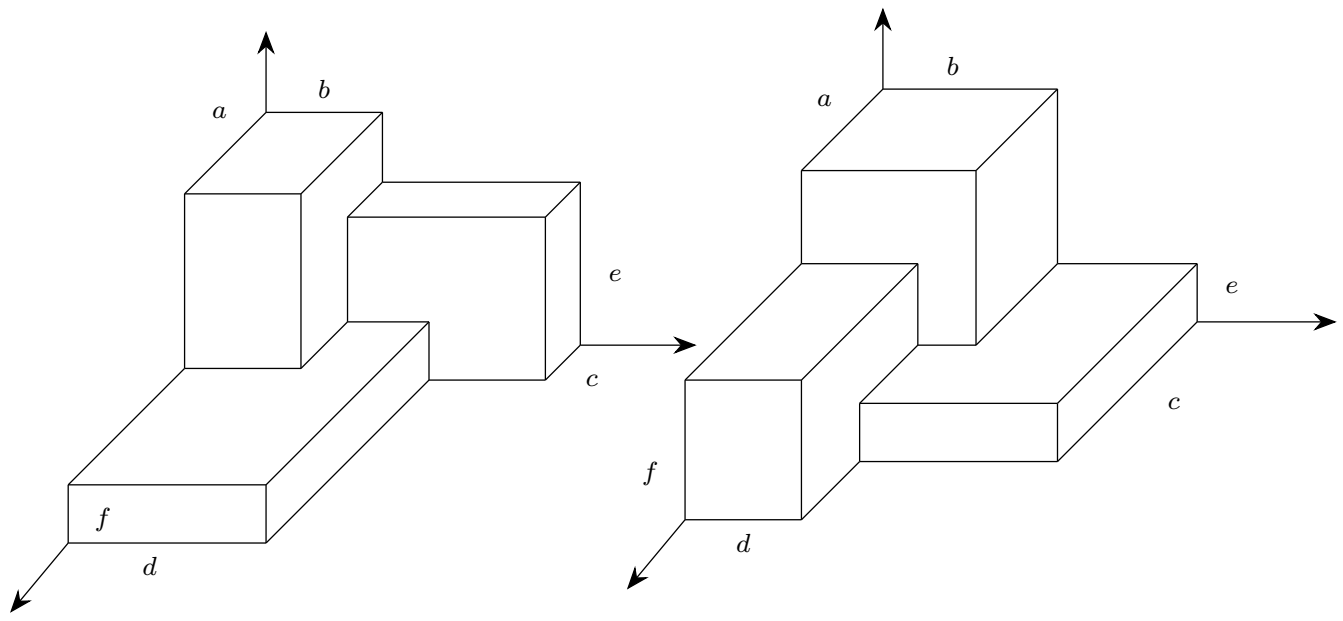
\resizebox{0.7\textwidth}{!}{%
\begin{circuitikz}
\tikzstyle{every node}=[font=\fontsize{8.0pt}{10.4pt}\selectfont]
\draw (6.75,9.875) to[short] (6.375,9.5);
\draw (6.75,9.875) to[short] (4.625,9.875);
\draw (1.25,6) to[short] (3.375,6);
\draw (3.375,6.625) to[short] (3.375,6);
\draw (3.375,6) to[short] (5.125,7.75);
\draw (6.75,9.875) to[short] (6.75,8.125);
\draw (4.25,9.5) to[short] (6.375,9.5);
\draw (6.375,9.5) to[short] (6.375,7.75);
\draw (5.125,7.75) to[short] (6.375,7.75);
\draw (4.625,9.875) to[short] (4.25,9.5);
\draw (2.5,7.875) to[short] (2.5,9.75);
\draw (2.5,9.75) to[short] (3.75,9.75);
\draw (3.75,9.75) to[short] (3.75,7.875);
\draw (3.75,9.75) to[short] (4.625,10.625);
\draw (4.625,10.625) to[short] (4.625,9.875);
\draw (2.5,9.75) to[short] (3.375,10.625);
\draw (3.375,10.625) to[short] (4.625,10.625);
\draw [-{Stealth[scale=1.5]}, ] (1.25,6) -- (0.625,5.25);
\draw [-{Stealth[scale=1.5]}, ] (6.75,8.125) .. controls (7.5,8.125) and (7.375,8.125) .. (8,8.125) ;
\draw [-{Stealth[scale=1.5]}, ] (3.375,10.625) -- (3.375,11.5);
\draw (4.25,8.375) to[short] (3.75,7.875);
\draw (6.75,8.125) to[short] (6.375,7.75);
\draw (4.25,8.375) to[short] (5.125,8.375);
\draw (1.25,6.625) to[short] (3.375,6.625);
\draw (1.25,6.625) to[short] (1.25,6);
\draw (1.25,6.625) to[short] (2.5,7.875);
\draw (3.375,6.625) to[short] (5.125,8.375);
\draw (5.125,8.375) to[short] (5.125,7.75);
\draw (2.5,7.875) to[short] (3.75,7.875);
\draw (4.25,9.5) to[short] (4.25,8.375);
\draw (13.375,9) to[short] (11.875,7.5);
\draw (13.375,9) to[short] (11.875,9);
\draw (7.875,6.25) to[short] (9.125,6.25);
\draw (13.375,9) to[short] (13.375,8.375);
\draw (11.875,7.5) to[short] (11.875,6.875);
\draw (9.75,6.875) to[short] (11.875,6.875);
\draw (9.125,9) to[short] (9.125,10);
\draw (9.125,10) to[short] (11,10);
\draw (11,10) to[short] (11,8.125);
\draw (11,10) to[short] (11.875,10.875);
\draw (11.875,10.875) to[short] (11.875,9);
\draw (9.125,10) to[short] (10,10.875);
\draw (10,10.875) to[short] (11.875,10.875);
\draw [-{Stealth[scale=1.5]}, ] (7.875,6.25) -- (7.25,5.5);
\draw [-{Stealth[scale=1.5]}, ] (13.375,8.375) -- (14.875,8.375);
\draw [-{Stealth[scale=1.5]}, ] (10,10.875) -- (10,11.75);
\draw (10.375,9) to[short] (9.125,7.75);
\draw (13.375,8.375) to[short] (11.875,6.875);
\draw (9.75,7.5) to[short] (11.875,7.5);
\draw (7.875,7.75) to[short] (7.875,6.25);
\draw (7.875,7.75) to[short] (9.125,9);
\draw (9.75,7.5) to[short] (9.75,6.875);
\draw (9.125,7.75) to[short] (9.125,6.25);
\draw (7.875,7.75) to[short] (9.125,7.75);
\draw (9.125,9) to[short] (10.375,9);
\draw (9.75,6.875) to[short] (9.125,6.25);
\draw (11,8.125) to[short] (11.875,9);
\draw (10.375,8.125) to[short] (11,8.125);
\draw (10.375,8.125) to[short] (9.75,7.5);
\node [font=\fontsize{8.0pt}{10.4pt}\selectfont, fill={rgb,255:red,255; green,255; blue,255}, fill opacity=1, text opacity=1, inner xsep=0.080cm, inner ysep=0.085cm, rounded corners=0.020cm] at (2.875,10.625) {$a$};
\node [font=\fontsize{8.0pt}{10.4pt}\selectfont, fill={rgb,255:red,255; green,255; blue,255}, fill opacity=1, text opacity=1, inner xsep=0.080cm, inner ysep=0.085cm, rounded corners=0.020cm] at (4,10.875) {$b$};
\node [font=\fontsize{8.0pt}{10.4pt}\selectfont, fill={rgb,255:red,255; green,255; blue,255}, fill opacity=1, text opacity=1, inner xsep=0.080cm, inner ysep=0.085cm, rounded corners=0.020cm] at (6.875,7.75) {$c$};
\node [font=\fontsize{8.0pt}{10.4pt}\selectfont, fill={rgb,255:red,255; green,255; blue,255}, fill opacity=1, text opacity=1, inner xsep=0.080cm, inner ysep=0.085cm, rounded corners=0.020cm] at (2.125,5.75) {$d$};
\node [font=\fontsize{8.0pt}{10.4pt}\selectfont, fill={rgb,255:red,255; green,255; blue,255}, fill opacity=1, text opacity=1, inner xsep=0.080cm, inner ysep=0.085cm, rounded corners=0.020cm] at (1.625,6.25) {$f$};
\node [font=\fontsize{8.0pt}{10.4pt}\selectfont, fill={rgb,255:red,255; green,255; blue,255}, fill opacity=1, text opacity=1, inner xsep=0.080cm, inner ysep=0.085cm, rounded corners=0.020cm] at (7.125,8.875) {$e$};
\node [font=\fontsize{8.0pt}{10.4pt}\selectfont, fill={rgb,255:red,255; green,255; blue,255}, fill opacity=1, text opacity=1, inner xsep=0.080cm, inner ysep=0.085cm, rounded corners=0.020cm] at (10.75,11.125) {$b$};
\node [font=\fontsize{8.0pt}{10.4pt}\selectfont, fill={rgb,255:red,255; green,255; blue,255}, fill opacity=1, text opacity=1, inner xsep=0.080cm, inner ysep=0.085cm, rounded corners=0.020cm] at (9.375,10.75) {$a$};
\node [font=\fontsize{8.0pt}{10.4pt}\selectfont, fill={rgb,255:red,255; green,255; blue,255}, fill opacity=1, text opacity=1, inner xsep=0.080cm, inner ysep=0.085cm, rounded corners=0.020cm] at (7.5,6.75) {$f$};
\node [font=\fontsize{8.0pt}{10.4pt}\selectfont, fill={rgb,255:red,255; green,255; blue,255}, fill opacity=1, text opacity=1, inner xsep=0.080cm, inner ysep=0.085cm, rounded corners=0.020cm] at (8.5,6) {$d$};
\node [font=\fontsize{8.0pt}{10.4pt}\selectfont, fill={rgb,255:red,255; green,255; blue,255}, fill opacity=1, text opacity=1, inner xsep=0.080cm, inner ysep=0.085cm, rounded corners=0.020cm] at (13.125,7.5) {$c$};
\node [font=\fontsize{8.0pt}{10.4pt}\selectfont, fill={rgb,255:red,255; green,255; blue,255}, fill opacity=1, text opacity=1, inner xsep=0.080cm, inner ysep=0.085cm, rounded corners=0.020cm] at (13.75,8.75) {$e$};
\draw (10.375,9) to[short] (10.375,8.125);
\end{circuitikz}
}%

\end{lemma}

\begin{proof}
The proof is basically a case by case discussion that combines all the conditions in the lemma and corollaries.

If the staircase is right-oriented, then $c<a, b<d, f<e$:
 
\begin{enumerate}
\item Project $E$ to the $XY$ plane and look at conditions on the top and bottom $z$ layers given in the Corollaries. 
We know $a>0$, and if $b>0$, then $\min(h_x-a, h_x-c)=h_x-a=1$ and $\min(h_y-b, h_y-d)=h_y-d=1$; and if $b=0$, then $\min(c+1, h_x-c, h_x-a)=\min(c+1, h_x-a)=1$.
So either $b>0$ and $a=h_x-1$,$d=h_y-1$; or $b=0$ and $\min(c+1, h_x-a)=1$.

\item Project $E$ to the $XZ$ plane and look at conditions on the top and bottom $y$ layers. Similarly, either $c>0$ and $a=h_x-1$, $e=h_z-1$; or $c=0$ and $\min(f+1, h_z-e)=1$.

\item Project $E$ to the $YZ$ plane and look at conditions on the top and bottom $x$ layers. Either $f>0$ and $d=h_y-1$, $e=h_z-1$; or $f=0$ and $\min(b+1, h_y-d)=1$.

%For the $y$ direction, $\min(1+b, h_y-d)=1$. For the $z$ direction, $\min(1+f, h_z-e)=1$.
\end{enumerate}

So we can make the conclusion:

\begin{enumerate}
\item If $b=c=f=0$, then the choices of $a,d,e$ are unrestricted.

\item If exactly one of the $b,c,f$ are greater than $0$, then two of the $a,d,f$ are fixed as follows. Setting $b>0$ fixes $a=h_x-1$,$d=h_y-1$. Setting $c>0$ fixes $a=h_x-1$, $e=h_z-1$. Setting $f>0$ fixes $d=h_y-1$, $e=h_z-1$.

\item If at least two of the $b, c, f$ are greater than $0$, then $a=h_x-1$,$d=h_y-1$ and $e=h_z-1$ are determined.

\end{enumerate}

If the staircase is left-oriented, recall that $f>e, c>a, b>d$, we can make a similar computation and draw a similar conclusion.

\end{proof}

\begin{definition} 
\cite{AJR1}
Let $E$ be a staircase. A \textbf{singularizing triple} for $I$ is a triple of outer corners $\{a,b,c\}=\{(a_1, a_2, a_3), (b_1, b_2, b_3), (c_1, c_2, c_3)\}$ of $E$ such that $a_1 > b_1,c_1; b_2 > a_2,c_2; c_3 > a_3,b_3.$
\end{definition}

\begin{definition}
\label{rot_stair_def}
A rotating staircase $E$ is a minimized staircase with one unique outer corner for each of the top $x,y,z$ layers, and one unique inner corner for each of the bottom layers. The rotating staircase $E$ has two parts, the rotating wall part and the interior boxes part, where each part satisfy the numerical conditions below.

\begin{enumerate}
\item The rotating walls part is the union of three rectangular cuboids, where the unique outer corners of the rectangular cuboids are $(a,b, h_z),(c,h_y, e), (h_x,d, f)$ satisfying $a,c < h_x$, and $b,d < h_y$, and $e,f < h_z$. 

\item \label{rot_stair_2} The numerical values of the rectangular cuboids in the rotating walls satisfy those in Lemma \ref{rotatingwallscondition}.

\item The interior boxes part can be any arbitrary partition as long as their projections in all of the $x,y,z$ directions are fully contained in the corresponding bottom layers. In other words, the interior boxes are confined in the region marked by the dashed lines in the picture above.

\item \label{rot_stair_4} The interior boxes have no singularizing triples.
\end{enumerate}

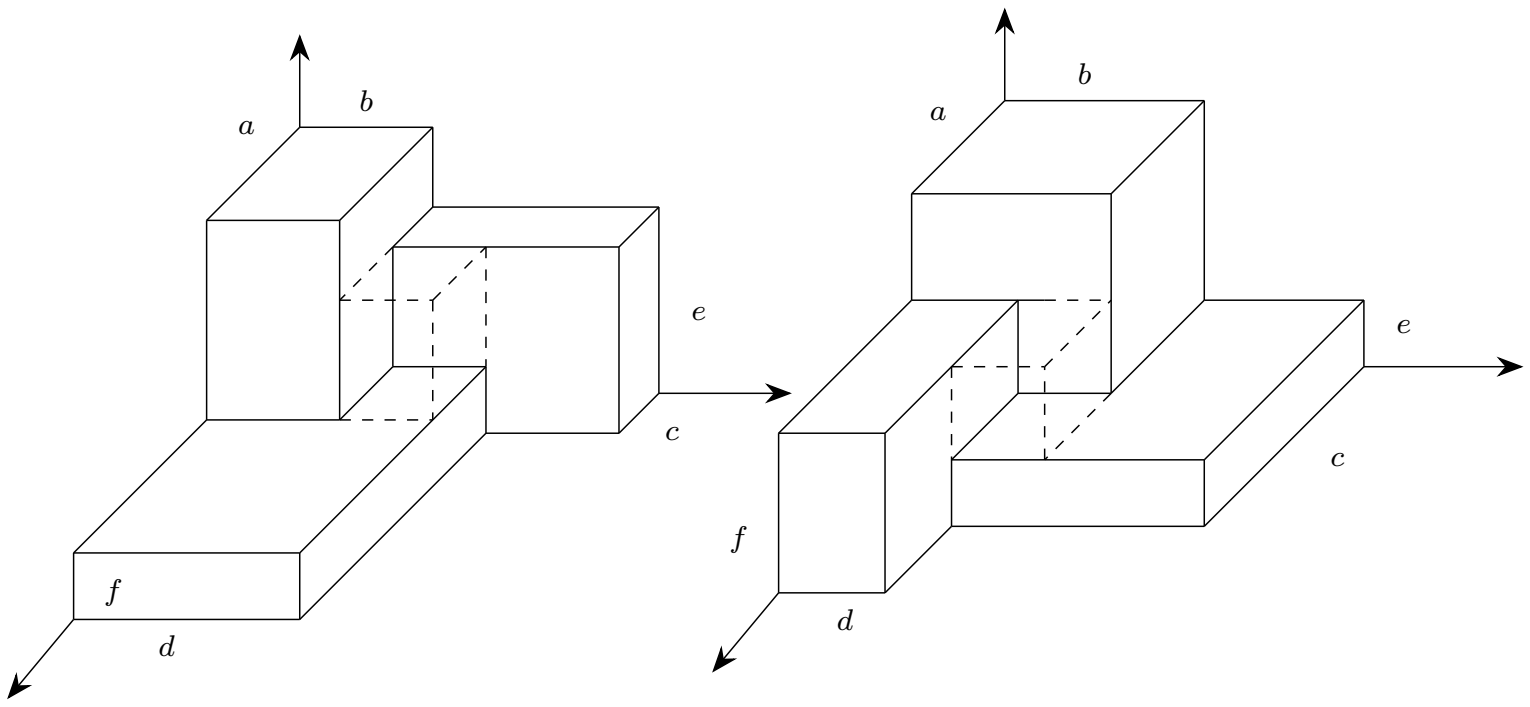
\resizebox{0.8\textwidth}{!}{%
\begin{circuitikz}
\tikzstyle{every node}=[font=\fontsize{18.2pt}{23.7pt}\selectfont]
\draw (6.75,9.875) to[short] (6.375,9.5);
\draw (6.75,9.875) to[short] (4.625,9.875);
\draw (1.25,6) to[short] (3.375,6);
\draw (3.375,6.625) to[short] (3.375,6);
\draw (3.375,6) to[short] (5.125,7.75);
\draw (6.75,9.875) to[short] (6.75,8.125);
\draw (4.25,9.5) to[short] (6.375,9.5);
\draw (6.375,9.5) to[short] (6.375,7.75);
\draw (5.125,7.75) to[short] (6.375,7.75);
\draw (4.625,9.875) to[short] (4.25,9.5);
\draw (2.5,7.875) to[short] (2.5,9.75);
\draw (2.5,9.75) to[short] (3.75,9.75);
\draw (3.75,9.75) to[short] (3.75,7.875);
\draw (3.75,9.75) to[short] (4.625,10.625);
\draw (4.625,10.625) to[short] (4.625,9.875);
\draw (2.5,9.75) to[short] (3.375,10.625);
\draw (3.375,10.625) to[short] (4.625,10.625);
\draw [-{Stealth[scale=1.5]}, ] (1.25,6) -- (0.625,5.25);
\draw [-{Stealth[scale=1.5]}, ] (6.75,8.125) .. controls (7.5,8.125) and (7.375,8.125) .. (8,8.125) ;
\draw [-{Stealth[scale=1.5]}, ] (3.375,10.625) -- (3.375,11.5);
\draw (4.25,8.375) to[short] (3.75,7.875);
\draw (6.75,8.125) to[short] (6.375,7.75);
\draw (4.25,8.375) to[short] (5.125,8.375);
\draw (1.25,6.625) to[short] (3.375,6.625);
\draw (1.25,6.625) to[short] (1.25,6);
\draw (1.25,6.625) to[short] (2.5,7.875);
\draw (3.375,6.625) to[short] (5.125,8.375);
\draw (5.125,8.375) to[short] (5.125,7.75);
\draw (2.5,7.875) to[short] (3.75,7.875);
\draw (4.25,9.5) to[short] (4.25,8.375);
\draw (13.375,9) to[short] (11.875,7.5);
\draw (13.375,9) to[short] (11.875,9);
\draw (7.875,6.25) to[short] (8.875,6.25);
\draw (13.375,9) to[short] (13.375,8.375);
\draw (11.875,7.5) to[short] (11.875,6.875);
\draw (9.5,6.875) to[short] (11.875,6.875);
\draw (9.125,9) to[short] (9.125,10);
\draw (9.125,10) to[short] (11,10);
\draw (11,10) to[short] (11,8.125);
\draw (11,10) to[short] (11.875,10.875);
\draw (11.875,10.875) to[short] (11.875,9);
\draw (9.125,10) to[short] (10,10.875);
\draw (10,10.875) to[short] (11.875,10.875);
\draw [-{Stealth[scale=1.5]}, ] (7.875,6.25) -- (7.25,5.5);
\draw [-{Stealth[scale=1.5]}, ] (13.375,8.375) -- (14.875,8.375);
\draw [-{Stealth[scale=1.5]}, ] (10,10.875) -- (10,11.75);
\draw (10.125,9) to[short] (8.875,7.75);
\draw (13.375,8.375) to[short] (11.875,6.875);
\draw (9.5,7.5) to[short] (11.875,7.5);
\draw (7.875,7.75) to[short] (7.875,6.25);
\draw (7.875,7.75) to[short] (9.125,9);
\draw (9.5,7.5) to[short] (9.5,6.875);
\draw (8.875,7.75) to[short] (8.875,6.25);
\draw (7.875,7.75) to[short] (8.875,7.75);
\draw (9.125,9) to[short] (10.375,9);
\draw (9.5,6.875) to[short] (8.875,6.25);
\draw (11,8.125) to[short] (11.875,9);
\draw (10.125,8.125) to[short] (11,8.125);
\draw (10.125,8.125) to[short] (9.5,7.5);
\node [font=\fontsize{8.0pt}{10.4pt}\selectfont, fill={rgb,255:red,255; green,255; blue,255}, fill opacity=1, text opacity=1, inner xsep=0.080cm, inner ysep=0.085cm, rounded corners=0.020cm] at (2.875,10.625) {$a$};
\node [font=\fontsize{8.0pt}{10.4pt}\selectfont, fill={rgb,255:red,255; green,255; blue,255}, fill opacity=1, text opacity=1, inner xsep=0.080cm, inner ysep=0.085cm, rounded corners=0.020cm] at (4,10.875) {$b$};
\node [font=\fontsize{8.0pt}{10.4pt}\selectfont, fill={rgb,255:red,255; green,255; blue,255}, fill opacity=1, text opacity=1, inner xsep=0.080cm, inner ysep=0.085cm, rounded corners=0.020cm] at (6.875,7.75) {$c$};
\node [font=\fontsize{8.0pt}{10.4pt}\selectfont, fill={rgb,255:red,255; green,255; blue,255}, fill opacity=1, text opacity=1, inner xsep=0.080cm, inner ysep=0.085cm, rounded corners=0.020cm] at (2.125,5.75) {$d$};
\node [font=\fontsize{8.0pt}{10.4pt}\selectfont, fill={rgb,255:red,255; green,255; blue,255}, fill opacity=1, text opacity=1, inner xsep=0.080cm, inner ysep=0.085cm, rounded corners=0.020cm] at (1.625,6.25) {$f$};
\node [font=\fontsize{8.0pt}{10.4pt}\selectfont, fill={rgb,255:red,255; green,255; blue,255}, fill opacity=1, text opacity=1, inner xsep=0.080cm, inner ysep=0.085cm, rounded corners=0.020cm] at (7.125,8.875) {$e$};
\node [font=\fontsize{8.0pt}{10.4pt}\selectfont, fill={rgb,255:red,255; green,255; blue,255}, fill opacity=1, text opacity=1, inner xsep=0.080cm, inner ysep=0.085cm, rounded corners=0.020cm] at (10.75,11.125) {$b$};
\node [font=\fontsize{8.0pt}{10.4pt}\selectfont, fill={rgb,255:red,255; green,255; blue,255}, fill opacity=1, text opacity=1, inner xsep=0.080cm, inner ysep=0.085cm, rounded corners=0.020cm] at (9.375,10.75) {$a$};
\node [font=\fontsize{8.0pt}{10.4pt}\selectfont, fill={rgb,255:red,255; green,255; blue,255}, fill opacity=1, text opacity=1, inner xsep=0.080cm, inner ysep=0.085cm, rounded corners=0.020cm] at (7.5,6.75) {$f$};
\node [font=\fontsize{8.0pt}{10.4pt}\selectfont, fill={rgb,255:red,255; green,255; blue,255}, fill opacity=1, text opacity=1, inner xsep=0.080cm, inner ysep=0.085cm, rounded corners=0.020cm] at (8.5,6) {$d$};
\node [font=\fontsize{8.0pt}{10.4pt}\selectfont, fill={rgb,255:red,255; green,255; blue,255}, fill opacity=1, text opacity=1, inner xsep=0.080cm, inner ysep=0.085cm, rounded corners=0.020cm] at (13.125,7.5) {$c$};
\node [font=\fontsize{8.0pt}{10.4pt}\selectfont, fill={rgb,255:red,255; green,255; blue,255}, fill opacity=1, text opacity=1, inner xsep=0.080cm, inner ysep=0.085cm, rounded corners=0.020cm] at (13.75,8.75) {$e$};
\draw (10.125,9) to[short] (10.125,8.125);
\draw [dashed] (4.25,9.5) -- (3.75,9);
\draw [dashed] (3.75,7.875) -- (4.625,7.875);
\draw [dashed] (3.75,9) -- (4.625,9);
\draw [dashed] (4.625,9) -- (4.625,7.875);
\draw [dashed] (4.625,9) -- (5.125,9.5);
\draw [dashed] (5.125,9.5) -- (5.125,8.375);
\draw [dashed] (10.375,9) -- (11,9);
\draw [dashed] (9.5,8.375) -- (9.5,7.5);
\draw [dashed] (11,8.125) -- (10.375,7.5);
\draw [dashed] (9.5,8.375) -- (10.375,8.375);
\draw [dashed] (10.375,7.5) -- (10.375,8.375);
\draw [dashed] (10.375,8.375) -- (11,9);
\end{circuitikz}
}%
\end{definition}

\begin{example}
\label{rot_star_examples}
The following are all rotating staircases.

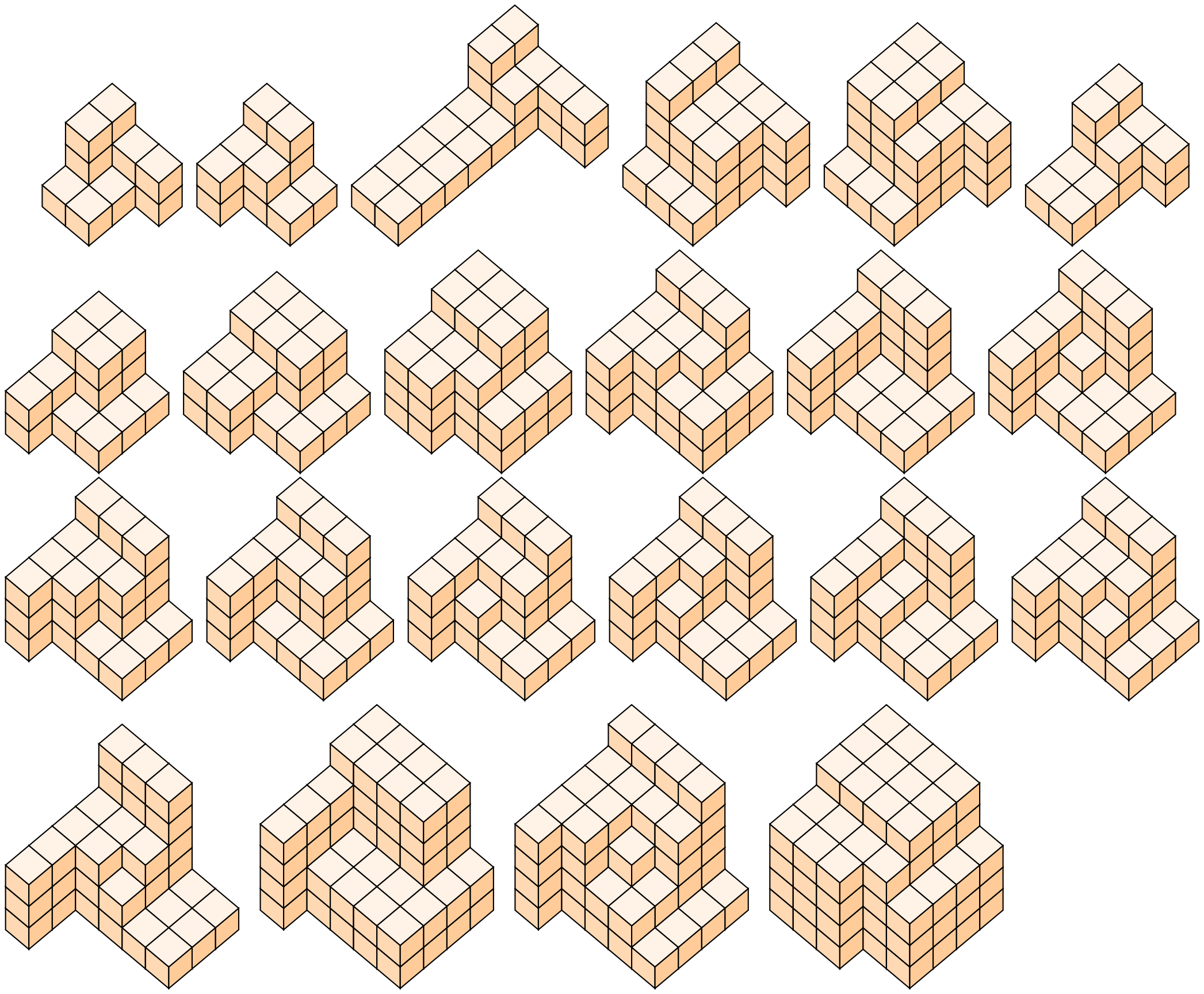
\monomialIdealPicture{{3,3,1},{2,1,1},{2}}\monomialIdealPicture{{3,2,2},{3,2},{1,1}}\monomialIdealPicture{{3,3,1,1,1,1,1},{2,2,1,1,1,1,1},{2},{2}}\monomialIdealPicture{{4,4,4,1},{3,3,3,1},{3,3,3,1},{3}}\monomialIdealPicture{{4,4,4,1},{4,4,4,1},{3,3,3,1},{3}}\monomialIdealPicture{{3,3,1,1},{2,2,1,1},{2}}\monomialIdealPicture{{3,3,2,2},{3,3,1},{1,1,1}}\monomialIdealPicture{{3,3,2,2},{3,3,2,2},{3,3,1},{1,1,1}}\monomialIdealPicture{{4,4,3,3},{4,4,3,3},{4,4,3},{2,2,2}}\monomialIdealPicture{{4,3,3,3},{4,3,3},{4,3,2},{2,2,2}}\monomialIdealPicture{{4,3,3,3},{4,1,1},{4,1,1},{1,1,1}}\monomialIdealPicture{{4,3,3,3},{4,2,1},{4,1,1},{1,1,1}}\monomialIdealPicture{{4,3,3,3},{4,3,3},{4,3,1},{1,1,1}}\monomialIdealPicture{{4,3,3,3},{4,3,1},{4,3,1},{1,1,1}}\monomialIdealPicture{{4,3,3,3},{4,3,2},{4,3,1},{1,1,1}}\monomialIdealPicture{{4,3,3,3},{4,3,2},{4,1,1},{1,1,1}}\monomialIdealPicture{{4,3,3,3},{4,2,2},{4,1,1},{1,1,1}}\monomialIdealPicture{{4,3,3,3},{4,3,3},{4,3,2},{1,1,1}}\monomialIdealPicture{{5,3,3,3,3},{5,3,3},{5,3,2},{1,1,1},{1,1,1}}\monomialIdealPicture{{5,5,4,4,4},{5,5,2,2},{5,5,2,2},{5,5,2,2},{2,2,2,2}}\monomialIdealPicture{{5,4,4,4,4},{5,4,4,4},{5,4,3,2},{5,4,2,2},{1,1,1,1}}\monomialIdealPicture{{5,5,5,4,4},{5,5,5,4,4},{5,5,5,4,4},{5,5,5,4},{3,3,3,3}} 
\end{example}

\begin{example}
We give a few counterexamples of staircases with $d$ boxes that look like rotating staircases, but one of the four conditions are not satisfied in Definition \ref{rot_stair_def}, and the tangent dimension is bigger than $3d+6$.

The following staircases do not satisfy the numerical conditions of the rotating walls (\ref{rot_stair_2}):

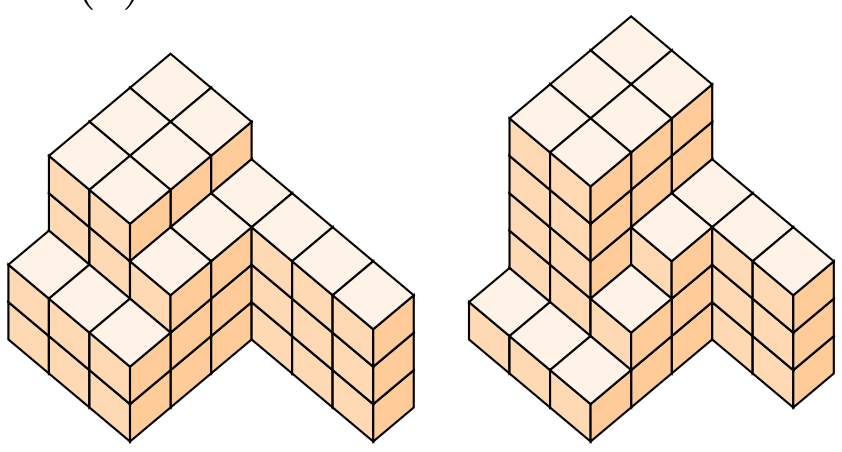
\monomialIdealPicture{{4,4,4,2},{4,4,4,2},{3,3,3,2},{3},{3},{3}} \monomialIdealPicture{{5,5,5,1},{5,5,5,1},{3,3,2,1},{3},{3}}.

The following staircases do not satisfy the condition that the interior boxes have a singularizing triple, the condition (\ref{rot_stair_4}):

\monomialIdealPicture{{4,3,3,3},{4,3,2},{4,2,1},{1,1,1}} \monomialIdealPicture{{5,4,4,4},{5,4,3},{5,2,1},{5,1,1},{1,1,1}}.

\end{example}

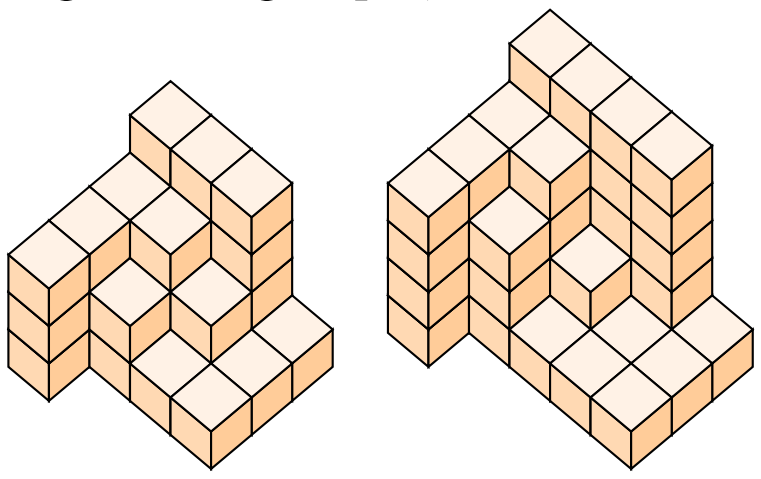
\begin{lemma}
\label{rot_stair_one_comp}
A rotating staircase $E$ as defined in Definition \ref{rot_stair_def} whose corresponding ideal is $I \in \Hilb^d(\mathbb{A}^3)$ has tangent dimension exactly equal to $3d+6$.
\end{lemma}

\begin{proof}
The outline of the proof is that there is exactly one bounded connected component for each $nnp, npn, pnn$ tangent spaces. Without loss of generality, we assume that we have a left-oriented, minimized, rotating staircase $E$, and prove the case for the $nnp$ signature, and the other orientation and other signatures follow from a similar argument. By assumption, $c>a, b>d$, and $f>e$.

For the $nnp$ vector, take the outer corner $(a, b,h_z)$ of size $(1, b-d+1, 1)$ of the top $z$ layer, the inner corner $(c+1,d+1,0)$ of the bottom layer, and the $nnp$ vector $\beta=(-c-1+a, -1,h_z)$ satisfy the condition of Lemma \ref{bigimportantlemma}, so there is a bounded connected component.

Now we show that there are no other bounded connected component for vectors of type $nnp$. 

Suppose there is another bounded connected component, then it must contain a shifted outer corner $O-\alpha$ of the shifted staircase $E-\alpha$. The outer corner could come from the shifted rotating walls part of $E-\alpha$, or the shifted interior boxes part of $E-\alpha$. We discuss each case.

We first consider the outer corners of the rotating walls.

Suppose we consider another $nnp$ vector, denoted by $\alpha=(-n_1, -n_2, p) \neq \beta$. We want to show that any box in $(E- \alpha) \cap I$ is connected to the infinite component.

First, note that the two outer corners $(h_x, d,f)$ and $(c, h_y, e)$, if in $(E- \alpha) \cap I$, when shifted by some $nnp$ vectors, must be connected to the infinite component via the edge $((h_x, d,f)-\alpha) \rightarrow ((h_x, 0,f)-\alpha)$, and $((c, h_y, e)-\alpha) \rightarrow ((0, h_y, e)-\alpha)$. So the only outer corner of the rotating walls part that could be contained in a bounded connected component, not connected to the infinite component, is $O=(a, b, h_z)$. 

When $h_z-e \leq p \leq h_z$, the shifter outer corner $(a, b, h_z)-\alpha$ has $z$ coordinate between $0$ and $e$, so we must have that $n_1 > c+1-a$ and $n_2 > 1$, otherwise the shifted outer corner $O=(a, b, h_z) -\alpha$ is in $E$ and not in $(E- \alpha) \cap I$. %By Lemma \ref{rotatingwallscondition}, we have that $\min(a+1, h_x-c)=1$, and $\min(e+1, h_z-f)=1$. 
So the outer corner $O$ shifted by $\alpha$ is $(a, b, h_z)-(-n_1, -n_2, p)>(c+1, b+1,h_z-p)$. 

Now we use Lemma \ref{rotatingwallscondition} to construct a path connecting the shifted outer corner $O-\alpha$ to the infinite component. If $d>0$, then $f=h_x-1$ and $b=h_y-1$, so $b+1=h_y$ and the shifted outer corner $O-\alpha$ has $y$ coordinate greater than $h_y$. So the shifted outer corner $O-\alpha$ is connected to the infinite component by the edge $((a, b, h_z)-\alpha) \rightarrow ((0, b, h_z)-\alpha)$.

If $d=0$, the shifted outer corner $O-\alpha$ has $y$ coordinate greater than $b+d+1$, so the shifted outer corner $O-\alpha$ is connected to the infinite component by the edge $((a, b, h_z)-\alpha) \rightarrow (a, 0, h_z)-\alpha)$. 

%If $h_x-c=1$,  then $c+1 >h_x$; if $a=0$, then $n_1 > c$; and in both cases 

When $h_z-e > p$, the shifter outer corner $(a, b, h_z)-\alpha$ has $z$ coordinate between $e$ and $h_z$, for any $nnp$ vector, the outer corner $O-\alpha$ is connected to the infinite component by the edge $((a,b, h_z)-\alpha) \rightarrow ((0,b,h_z)-\alpha)$.

%Suppose $p'< h_z-\min(e,f)$. The outer corner $(a, b, h_z)$, if is an element of $(E-\alpha) \cap I$, in other words not removed because it is not in $I$, then it is path connected to infinite component by going to the negative $y$ direction if it is right oriented, and it is path connected to infinite component by going to the negative $x$ direction if it is left oriented. Suppose $h_z-\min(e,f) \leq p' < h_z$. The outer corner $(a, b, h_z)$,

%So we conclude that $p'=h_z$ for there to be a bounded connected component. But if the rotating staircase satisfy its defining numerical properties Lemma \ref{rotatingwallscondition}, then there is no bounded connected component 

%infinite component by .  so there is no extra bounded connected component induced by the outer corner of the rotating walls part. If there is a bounded connected component, then it must contain at least one of the three corners. So if every corner must be connected to the infinite component, there's no bounded connected component.

Suppose for the purpose of contradiction that there exists an outer corner $O'$ in interior boxes part $T$, that after shifting by an $nnp$ vector $\alpha$, is also contained in another bounded connected component $B$.

Consider the interior boxes $T$ as a staircase itself. 

If the bounded connected component $B$ of $E-\alpha$, containing $O'-\alpha$, is fully contained in the dashed box, then the bounded connected component of $E$ is also a bounded connected component of $T$. By assumption, the interior boxes $T$ has no singularizing triple, so by Theorem \ref{smoothconditionsAJR}, $T$ shouldn't have any bounded connected components, contradiction. 

If the bounded connected component $B$ of $E-\alpha$ is not fully contained in the dashed box, we discuss all the possibilities of its location. The dashed box has three sides, the front side, the up side, and the right side, where $B$ could sit. 

\begin{enumerate}
\item Because the vector $\alpha$ has signature $nnp$,  the bounded connected component $B$ of $E-\alpha$ cannot be higher than the up side of the dashed box. 

\item If some box $b$ of bounded connected component $B$ of $E-\alpha$ to on the front side of the dashed box, then by construction, $B$ is connected to the infinite wall part by any path that goes from $b$ straight down to the the negative $z$ direction.

\item If some box $b$ of bounded connected component $B$ of $E-\alpha$ is on the right side of the dashed box, then by construction, $B$ is connected to the infinite wall part by any path that goes from $b$ to the the negative $z$ direction.

\end{enumerate}

In all cases above, we have a contradiction. There are no other bounded connected components of signature $nnp$.

\end{proof}

\begin{definition}
\label{rot_stair_definition_comb}
A \textbf{rotating singular staircase} is a minimized staircase satisfying the conditions (1)(2) and (3) in Definition \ref{rot_stair_def} , but has at least one singularizing triple in the interior boxes part.
\end{definition}

\begin{lemma}
\label{rot_sing_stair}
An rotating singular staircase has tangent dimension strictly greater than $3d+6$.
\end{lemma}

\begin{proof}
As we showed in the proof Lemma \ref{rot_stair_one_comp}, if a minimized staircase $E$ satisfies the conditions (1)(2) and (3) in Definition \ref{rot_stair_def}, then there is at least one bounded connected component for each $nnp$, $npn$, $pnn$ signature, coming from the top and bottom layers of the rotating wall part. 

Consider the interior boxes part itself as a staircase $T$ with corresponding ideal $I_T$. Because $T$ has a singularizing triple, by Theorem \ref{smoothconditionsAJR}, there is at least one bounded connected component of $(T-\alpha) \cap I_T$. The bounded connected component of $((T-\alpha) \cap I_T) +(a,d,e)$ is also a bounded connected component of $(E-\alpha) \cap I$. So the tangent space dimension is strictly greater than $3d+6$.
%and $(T-\alpha) \cap I_T$ shifted by the vector $(a,d,e)$ is fully contained in the dashed box.  Because $(T-\alpha) \cap I_T$ shifted by the vector $(a,d,e)$ is a subset of $(E-\alpha) \cap I$, the bounded connected component of $((T-\alpha) \cap I_T) +(a,d,e)$ is also a bounded connected component of $(E-\alpha) \cap I$. This is the extra bounded connected component we find.
\end{proof}

Now, we have shown that all of the four conditions in the Definition \ref{rot_stair_def} of the rotating walls part are necessary for the tangent dimension to be exactly $3d+6$. If the staircase does not satisfy the numerical conditions (1)(2), we have two or more bounded connected components. (3) is tautological as a partition. The necessity of (4) is given as Lemma \ref{rot_sing_stair}.

%If the staircase does not satisfy the numerical conditions (1)(2), we have two or more bounded connected components. \Ytodo{add some explanation, refer to the previous lemmas } 

%\begin{definition}
%Take the bottom layer in $x,y$ or $z$ directions, and find how many times it repeats itself, without losing any boxes, when it goes up to the top. Define these boxes as the "wall" part of the staircase. Define the rest of the boxes as the "interior" part. In other words, the wall parts are $\{(x,y,z)| x \leq x_0, y \leq y_0, z \leq z_0\}$ where $x_0,y_0,z_0$ are the times that the bottom layer repeats itself, in the $x,y,z$ directions. The interior part is $\{(x,y,z)| x \geq x_0+1, y \geq y_0+1, z \geq z_0+1\}$.
%\end{definition}

\section{linkage}

\subsection{Pairs of linking staircases}

\begin{remark}
\label{linkedincube}
Consider taking linkage of a monomial ideal $I$ of $\mathbb{C}[x,y,z]$ using the regular sequence $(x^a,y^b,z^c) \subset I$. Geometrically, we are taking complement of the staircase $E$ of the ideal $I$ inside a rectangular cuboid of size $a$ by $b$ by $c$. %\rtodo{cite somewhere as a lemma/remark or prove this?}
\end{remark}

\begin{corollary}
\label{cor_linkage_pairs}
The complement of an opera house $E$, where the maximum $x,y,z$ axis are $h_x=a, h_y=b, h_z=c$, inside a minimum fitting cube of size $a \times b \times c$ is another opera house $F$. The complement of a rotating staircase $R$, where the maximum $x,y,z$ axis are $h_x=a, h_y=b, h_z=c$, inside a minimum fitting cube of size $a \times b \times c$ is another rotating staircase $G$.

Therefore, both the opera houses and rotating staircases come in pairs that are linked to each other, and the union of two staircases in a pair is a cube.    
\end{corollary}

\begin{proof}
This boils down to going through the conditions defining the opera houses and rotating boxes in the last section, and confirming that the complemental shape in a minimum fitting cube again satisfies each of the conditions. It is an easy check but long to write down all the conditions again, so we omit it here.
\end{proof}

We now give some examples of the above Corollary.

\begin{example}
\begin{enumerate}
\item Pairs of opera houses 

(a). \monomialIdealPicture{{2,2,2,2,1},{2,2,2},{1,1}}
\monomialIdealPicture{{2,2,2,1,1},{2,2},{1}} 

(b). \monomialIdealPicture{{3,3,3,1},{3,3,3},{3,3,3},{2}}
\monomialIdealPicture{{3,3,3,2},{3},{3},{1}}

\item Pairs of rotating staircases

(a).\monomialIdealPicture{{3,3,1},{2,1,1},{2}}\monomialIdealPicture{{3,2,2},{3,2},{1,1}} 

(b).\monomialIdealPicture{{3,3,1,1},{2,2,1,1},{2}}\monomialIdealPicture{{3,3,2,2},{3,3,1},{1,1,1}}

\end{enumerate}
\end{example}

Now we will link any opera house to a staircase of the shape of a tripod. 

\begin{definition}
A tripod is a staircase $E$ whose corresponding ideal has form $I_{tri(t_1, t_2, t_3)} := (x^{t_1},y^{t_2},z^{t_3},xy,xz,yz)$ for some $t_1, t_2, t_3  \geq  2$.
\end{definition}

\begin{example}
The following are tripods.
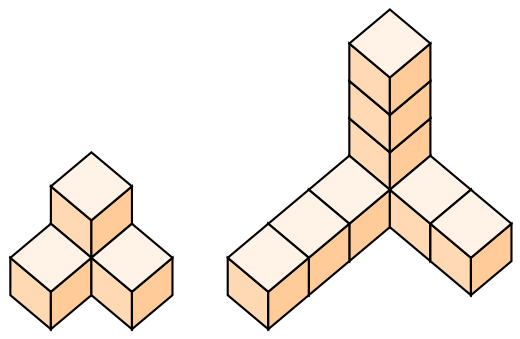
\monomialIdealPicture{{2,1},{1}} \monomialIdealPicture{{4,1,1,1},{1},{1}}
\end{example}

\begin{lemma}\cite{AJR1}[Proposition 5.12.]
All tripods are in the linkage class of the ideal $(x,y,z)^2$.
\end{lemma}

The strategy of the proof boils down to a few steps, where each step links an opera house to another simpler opera house. More specifically, for any opera house that is not already a tripod, we link it to another opera house with one fewer inner corners in the bottom layer. So inductively, we can link a complicated opera house to one with a unique inner corner in the bottom layer. Then we reduce the height of both of the wings of the opera house to $1$. Finally, such opera houses are directly linked to a tripod. In the following subsections, we explicitly write and prove the above steps.

\subsection{Notations of the opera house.\\}
Denote the number of boxes on the bottom layer of each wing that connects to the main part as the \textbf{width of the wing}. For example, for all opera houses of two spread-out rows, both of the wings have width $1$. 

For an opera house with the inner corners in the bottom layers denoted as $(d_1, e_1, 0), ..., (d_k, e_k, 0)$, where $d_1 > d_2 > ... >d_k$ and $e_1 < e_2 < ... <e_k$. The generators of the ideal $I$ of the opera house $E$, other than $x^{d_1}y^{e_1}, ..., x^{d_k}y^{e_k}$,  are denoted as $x^a, y^b, z^c$ and $x^{p_1}z^{q_1}$, $y^{p_2}z^{q_2}$, where $p_1 \geq d_1$ and $p_2 \geq e_k$. 

The following picture shows as an example, the case when $p_1=d_1, p_2=d_2$ and the case when $p_1>d_1, p_2>d_2$.  

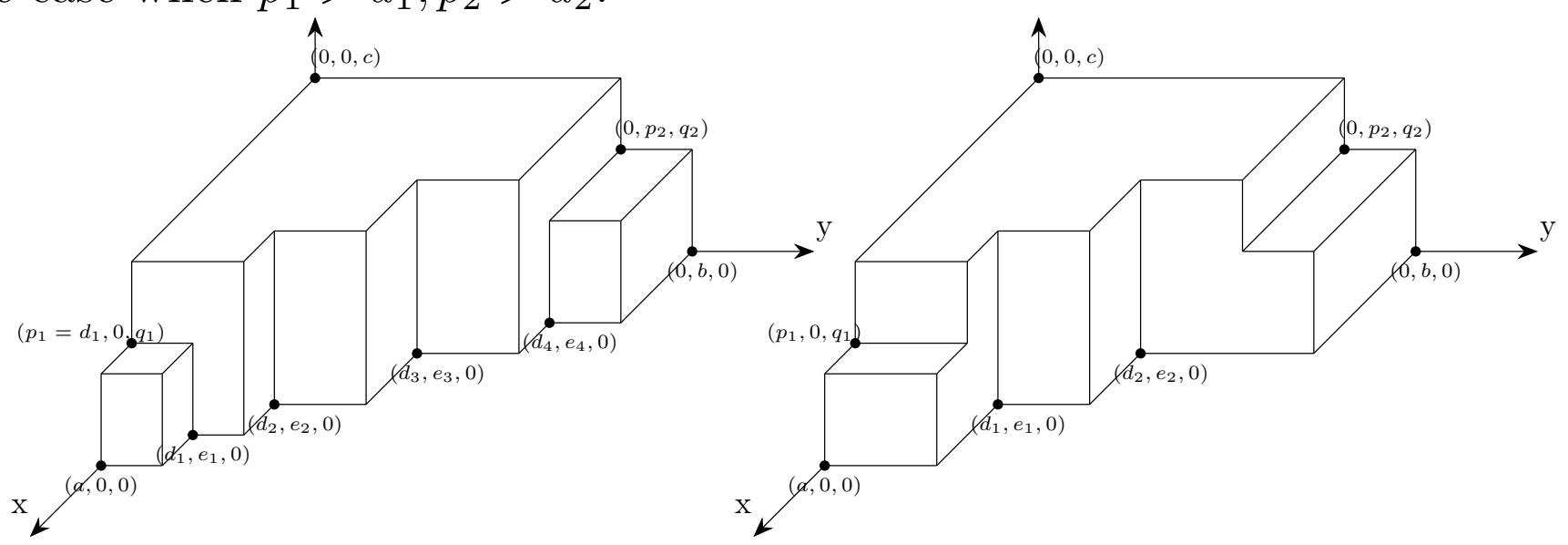
\resizebox{0.9\textwidth}{!}{%
\begin{circuitikz}
\tikzstyle{every node}=[font=\fontsize{18.2pt}{23.7pt}\selectfont]
\draw (-8,11) to[short] (-4.25,11);
\draw (-4.25,11) to[short] (-5.5,9.75);
\draw (-5.5,9.75) to[short] (-6.75,9.75);
\draw (-6.75,9.75) to[short] (-7.375,9.125);
\draw (-7.375,9.125) to[short] (-8.5,9.125);
\draw (-8,11) to[short] (-10.25,8.75);
\draw (-10.25,8.75) to[short] (-8.875,8.75);
\draw (-7.125,9.375) to[short] (-7.125,9.375);
\draw (-8.5,9.125) to[short] (-8.875,8.75);
\draw (-8.5,9.125) to[short] (-8.5,7);
\draw (-7.375,9.125) to[short] (-7.375,7);
\draw (-6.75,9.75) to[short] (-6.75,7.625);
\draw (-10.25,8.75) to[short] (-10.25,7.75);
\draw (-8.875,8.75) to[short] (-8.875,6.625);
\draw (-5.5,9.75) to[short] (-5.5,7.625);
\draw (-4.25,11) to[short] (-4.25,10.125);
\draw (-5.125,8) to[short] (-5.5,7.625);
\draw (-5.5,7.625) to[short] (-6.75,7.625);
\draw (-6.75,7.625) to[short] (-7.375,7);
\draw (-7.375,7) to[short] (-8.5,7);
\draw (-8.5,7) to[short] (-8.875,6.625);
\draw (-9.5,6.625) to[short] (-8.875,6.625);
\draw (-10.625,7.375) to[short] (-10.625,6.25);
\draw (-9.875,7.375) to[short] (-9.875,6.25);
\draw (-10.625,7.375) to[short] (-9.875,7.375);
\draw (-10.625,6.25) to[short] (-9.875,6.25);
\draw (-10.25,7.75) to[short] (-10.625,7.375);
\draw (-9.5,7.75) to[short] (-9.875,7.375);
\draw (-9.5,7.75) to[short] (-9.5,6.625);
\draw (-9.5,6.625) to[short] (-9.875,6.25);
\draw (-10.25,7.75) to[short] (-9.5,7.75);
\draw (-3.375,10.125) to[short] (-4.25,9.25);
\draw (-3.375,10.125) to[short] (-3.375,8.875);
\draw (-4.25,9.25) to[short] (-4.25,8);
\draw (-3.375,8.875) to[short] (-4.25,8);
\draw (-4.25,10.125) to[short] (-3.375,10.125);
\draw (-5.125,8) to[short] (-4.25,8);
\draw (-5.125,9.25) to[short] (-4.25,9.25);
\draw (-4.25,10.125) to[short] (-5.125,9.25);
\draw (-5.125,9.25) to[short] (-5.125,8);
\begin{scope}[transparency group, opacity=1]
\node at (-9.5,6.625) [circ, color={rgb,255:red,0; green,0; blue,0}] {};
\end{scope}
\begin{scope}[transparency group, opacity=1]
\node at (-8.5,7) [circ, color={rgb,255:red,0; green,0; blue,0}] {};
\end{scope}
\begin{scope}[transparency group, opacity=1]
\node at (-6.75,7.625) [circ, color={rgb,255:red,0; green,0; blue,0}] {};
\end{scope}
\begin{scope}[transparency group, opacity=1]
\node at (-5.125,8) [circ, color={rgb,255:red,0; green,0; blue,0}] {};
\end{scope}
\draw [-{Stealth[scale=1.5]}, ] (-10.625,6.25) -- (-11.5,5.375);
\draw [-{Stealth[scale=1.5]}, ] (-3.375,8.875) -- (-1.875,8.875);
\draw [-{Stealth[scale=1.5]}, ] (-8,11) -- (-8,11.75);
\node [font=\fontsize{10.8pt}{14.1pt}\selectfont, inner xsep=0.080cm, inner ysep=0.085cm, rounded corners=0.000cm] at (-11.625,5.75) {x};
\node [font=\fontsize{10.8pt}{14.1pt}\selectfont, inner xsep=0.080cm, inner ysep=0.085cm, rounded corners=0.000cm] at (-1.75,9.125) {y};
\node [font=\fontsize{6.8pt}{8.9pt}\selectfont, inner xsep=0.080cm, inner ysep=0.085cm, rounded corners=0.020cm] at (-9.375,6.375) {$(d_1,e_1,0)$};
\node [font=\fontsize{6.8pt}{8.9pt}\selectfont, inner xsep=0.080cm, inner ysep=0.085cm, rounded corners=0.020cm] at (-4.875,7.75) {$(d_4,e_4,0)$};
\node [font=\fontsize{6.8pt}{8.9pt}\selectfont, inner xsep=0.080cm, inner ysep=0.085cm, rounded corners=0.020cm] at (-6.5,7.375) {$(d_3,e_3,0)$};
\node [font=\fontsize{6.8pt}{8.9pt}\selectfont, inner xsep=0.080cm, inner ysep=0.085cm, rounded corners=0.020cm] at (-8.25,6.75) {$(d_2,e_2,0)$};
\begin{scope}[transparency group, opacity=1]
\node at (-10.625,6.25) [circ, color={rgb,255:red,0; green,0; blue,0}] {};
\end{scope}
\begin{scope}[transparency group, opacity=1]
\node at (-3.375,8.875) [circ, color={rgb,255:red,0; green,0; blue,0}] {};
\end{scope}
\begin{scope}[transparency group, opacity=1]
\node at (-10.25,7.75) [circ, color={rgb,255:red,0; green,0; blue,0}] {};
\end{scope}
\begin{scope}[transparency group, opacity=1]
\node at (-4.25,10.125) [circ, color={rgb,255:red,0; green,0; blue,0}] {};
\end{scope}
\begin{scope}[transparency group, opacity=1]
\node at (-8,11) [circ, color={rgb,255:red,0; green,0; blue,0}] {};
\end{scope}
\node [font=\fontsize{6.8pt}{8.9pt}\selectfont, inner xsep=0.080cm, inner ysep=0.085cm, rounded corners=0.020cm] at (-10.625,6) {$(a,0,0)$};
\node [font=\fontsize{6.8pt}{8.9pt}\selectfont, inner xsep=0.080cm, inner ysep=0.085cm, rounded corners=0.020cm] at (-3.25,8.625) {$(0,b,0)$};
\node [font=\fontsize{6.8pt}{8.9pt}\selectfont, inner xsep=0.080cm, inner ysep=0.085cm, rounded corners=0.020cm] at (-7.625,11.25) {$(0,0,c)$};
\node [font=\fontsize{6.8pt}{8.9pt}\selectfont, inner xsep=0.080cm, inner ysep=0.085cm, rounded corners=0.020cm] at (-10.75,7.875) {$(p_1=d_1, 0, q_1)$};
\node [font=\fontsize{6.8pt}{8.9pt}\selectfont, inner xsep=0.080cm, inner ysep=0.085cm, rounded corners=0.020cm] at (-3.75,10.375) {$(0, p_2, q_2)$};
\draw (0.875,11) to[short] (4.625,11);
\draw (4.625,11) to[short] (3.375,9.75);
\draw (3.375,9.75) to[short] (2.125,9.75);
\draw (2.125,9.75) to[short] (1.5,9.125);
\draw (1.5,9.125) to[short] (0.375,9.125);
\draw (0.875,11) to[short] (-1.375,8.75);
\draw (-1.375,8.75) to[short] (0,8.75);
\draw (1.75,9.375) to[short] (1.75,9.375);
\draw (0.375,9.125) to[short] (0,8.75);
\draw (0.375,9.125) to[short] (0.375,7);
\draw (1.5,9.125) to[short] (1.5,7);
\draw (2.125,9.75) to[short] (2.125,7.625);
\draw (-1.375,8.75) to[short] (-1.375,7.75);
\draw (0,8.75) to[short] (0,7.75);
\draw (3.375,9.75) to[short] (3.375,8.875);
\draw (4.625,11) to[short] (4.625,10.125);
\draw (4.25,7.625) to[short] (2.125,7.625);
\draw (2.125,7.625) to[short] (1.5,7);
\draw (1.5,7) to[short] (0.375,7);
\draw (0.375,7) to[short] (-0.375,6.25);
\draw (-1.75,7.375) to[short] (-1.75,6.25);
\draw (-0.375,7.375) to[short] (-0.375,6.25);
\draw (-1.75,7.375) to[short] (-0.375,7.375);
\draw (-1.75,6.25) to[short] (-0.375,6.25);
\draw (-1.375,7.75) to[short] (-1.75,7.375);
\draw (5.5,10.125) to[short] (4.25,8.875);
\draw (5.5,10.125) to[short] (5.5,8.875);
\draw (4.25,8.875) to[short] (4.25,7.625);
\draw (5.5,8.875) to[short] (4.25,7.625);
\draw (4.625,10.125) to[short] (5.5,10.125);
\draw (3.375,8.875) to[short] (4.25,8.875);
\draw (4.625,10.125) to[short] (3.375,8.875);
\begin{scope}[transparency group, opacity=1]
\node at (0.375,7) [circ, color={rgb,255:red,0; green,0; blue,0}] {};
\end{scope}
\begin{scope}[transparency group, opacity=1]
\node at (2.125,7.625) [circ, color={rgb,255:red,0; green,0; blue,0}] {};
\end{scope}
\draw [-{Stealth[scale=1.5]}, ] (-1.75,6.25) -- (-2.625,5.375);
\draw [-{Stealth[scale=1.5]}, ] (5.5,8.875) -- (7,8.875);
\draw [-{Stealth[scale=1.5]}, ] (0.875,11) -- (0.875,11.75);
\node [font=\fontsize{10.8pt}{14.1pt}\selectfont, inner xsep=0.080cm, inner ysep=0.085cm, rounded corners=0.000cm] at (-2.75,5.75) {x};
\node [font=\fontsize{10.8pt}{14.1pt}\selectfont, inner xsep=0.080cm, inner ysep=0.085cm, rounded corners=0.000cm] at (7.125,9.125) {y};
\node [font=\fontsize{6.8pt}{8.9pt}\selectfont, inner xsep=0.080cm, inner ysep=0.085cm, rounded corners=0.020cm] at (2.375,7.375) {$(d_2,e_2,0)$};
\node [font=\fontsize{6.8pt}{8.9pt}\selectfont, inner xsep=0.080cm, inner ysep=0.085cm, rounded corners=0.020cm] at (0.625,6.75) {$(d_1,e_1,0)$};
\begin{scope}[transparency group, opacity=1]
\node at (-1.75,6.25) [circ, color={rgb,255:red,0; green,0; blue,0}] {};
\end{scope}
\begin{scope}[transparency group, opacity=1]
\node at (5.5,8.875) [circ, color={rgb,255:red,0; green,0; blue,0}] {};
\end{scope}
\begin{scope}[transparency group, opacity=1]
\node at (-1.375,7.75) [circ, color={rgb,255:red,0; green,0; blue,0}] {};
\end{scope}
\begin{scope}[transparency group, opacity=1]
\node at (4.625,10.125) [circ, color={rgb,255:red,0; green,0; blue,0}] {};
\end{scope}
\begin{scope}[transparency group, opacity=1]
\node at (0.875,11) [circ, color={rgb,255:red,0; green,0; blue,0}] {};
\end{scope}
\node [font=\fontsize{6.8pt}{8.9pt}\selectfont, inner xsep=0.080cm, inner ysep=0.085cm, rounded corners=0.020cm] at (-1.75,6) {$(a,0,0)$};
\node [font=\fontsize{6.8pt}{8.9pt}\selectfont, inner xsep=0.080cm, inner ysep=0.085cm, rounded corners=0.020cm] at (5.625,8.625) {$(0,b,0)$};
\node [font=\fontsize{6.8pt}{8.9pt}\selectfont, inner xsep=0.080cm, inner ysep=0.085cm, rounded corners=0.020cm] at (1.25,11.25) {$(0,0,c)$};
\node [font=\fontsize{6.8pt}{8.9pt}\selectfont, inner xsep=0.080cm, inner ysep=0.085cm, rounded corners=0.020cm] at (-1.875,7.875) {$(p_1, 0, q_1)$};
\node [font=\fontsize{6.8pt}{8.9pt}\selectfont, inner xsep=0.080cm, inner ysep=0.085cm, rounded corners=0.020cm] at (5.125,10.375) {$(0, p_2, q_2)$};
\draw (-1.375,7.75) to[short] (0,7.75);
\draw (0,7.75) to[short] (-0.375,7.375);
\end{circuitikz}
}%

\begin{definition}
\label{alg_opera_house}
An opera house is generated as $$I=(x^a, y^b, z^c, x^{d_1}y^{e_1}, ..., x^{d_k}y^{e_k}, x^{p_1}z^{q_1}, y^{p_2}z^{q_2}),$$
where $a > p_1 \geq d_1 > d_2 > ... >d_k$, and $e_1 < e_2 < ... <e_k \leq p_2 < b$, and $c>q_1, q_2$. Either one of $q_1, q_2$ is equal to $1$ or $c-1$.
\end{definition}

The height of the opera house is denoted as $h_z=c$. At least one of the two wings must have height $q_1, q_2$ to be equal to either $1$ or $c-1$. For all opera houses in which one of the two wings has height $h_z-1$, and the other wing does not have height $1$, there is an opera house linked to it in pairs such that one of the wings have height $1$ by Remark \ref{linkedincube}. So from now on, we only need to consider opera houses where one wing has exactly height $1$. Without loss of generality, set $q_1=1$, and the case of $q_2=1$ is similar by symmetry.

We conclude in end of this section that the ideal $I$ of an opera house $E$, at the start of the linkage process, is $$I=(x^a, y^b, z^c, x^{d_1}y^{e_1}, ..., x^{d_k}y^{e_k}, x^{p_1}z, y^{p_2}z^{q_2}),$$ where all these numerics satisfy all the conditions to be an opera house.

\subsection{Reducing the number of corners in the bottom layer.\\}

\textbf{Algorithm:}
Consider the complete intersection ideals $$G_1=(x^a+y^b+z^c, x^{d_1}y^{e_1},x^{p_1}z+y^{p_2}z^{q_2}),$$ 
$$G_2=(x^a+y^b+z^c, x^{d_1}y^{e_1},x^{d_1}z+y^{p_2-(e_2-e_1)}z^{q_2}).$$

Take a 2 step linkage:
Compute $I_2=G_1:I$;  $I_3=G_2:I_2$.

Before we compute the linkage, we prove the following useful Lemma.

\begin{lemma}
\label{lem_single_gen}
In a polynomial quotient ring $R=\frac{\mathbb{C}[x,y,z]}{(f_1,f_2)}$ for some polynomials $f_1, f_2$, where $a,b$ are not zero divisors, and $f_1, f_2, a$ forms a regular sequences, and $J$ an ideal in $R$ containing $a$, we have that $(b):((a):J))=(bJ):(a)$.
\end{lemma}

\begin{proof}
Denote $E:=(a):J$ and an element in $E$ as $e$. Denote $F:=(b):E$ and an element in $F$ as $f$. Denote $K:=(bJ):(a)$ and an element in $K$ as $k$. 

Suppose an element $k \in (bJ):(a)$, then $ak \subset bJ$ by definition. Similarly, $eJ \subset (a)$ and $fE \subset (b)$. In the following proof, we implicitly use that $R$ is commutative.

We first show one side of the inclusion: $(bJ):(a) \subset (b):((a):J))$. For any element $k \in K$, we want to show that $k \in F$. So we need to show that $kE \subset (b)$. Because we know that $ak=bj$ for some $j \in J$, and that $je = ar$ for some $r \in R$, we can compute $ake=bje= abr$ for some $r \in R$. Because $a$ is not a zero divisor, we cancel $a$ and have that $ke=br$, so $ke \in (b)$ and indeed $k \in F$.

Now we show the other side of the inclusion:  $(b):((a):J)) \subset (bJ):(a)$. For any element $f \in F$, we want to show that $f \in K$. Because $a \in E$ and $fE \subset (b)$, we have $fa=br_1$ for some $r_1 \in R$. Because $f \in F$, we have that $fe=br_2$ for any $e \in E$ and some $r_2 \in R$. Combining two equations, we have that $fae=br_1e=br_2a$. Because $b$ is not a zero divisor, $r_1e=r_2a$ for any $e$ and some $r_1, r_2 \in R$. So $r_1\in (a):E = (a):((a):J)=J$. (The last equality is the reflexivity property under linkage in $\mathbb{C}[x,y,z]$, that $(a,f_1,f_2):((a,f_1,f_2) : Q)=Q$ given that $f_1, f_2,a$ is a regular sequence, $a \in J$ and $Q=J+(f_1,f_2)$.) So $r_1 \in J$ and $fa = br_1 \in bJ$, so $f \in bJ :(a)$.
\end{proof}

\begin{lemma} 
The above algorithm links the ideal $I$ of the opera house to 

$$I_3 =
( x^{a-(p_1-d_1)}, y^{b-\delta}, z^c, x^{d_2}y^{e_2-\delta}, x^{d_3}y^{e_3-\delta}, \ldots, x^{d_k}y^{e_k-\delta},
x^{d_1}z,
y^{p_2-\delta}z^{q_2}), $$ where $\delta=e_2-e_1.$

\end{lemma}

\newpage
\resizebox{0.7\linewidth}{!}{%
    \begin{circuitikz}
      \tikzstyle{every node}=[font=\fontsize{6.8pt}{8.9pt}\selectfont]
      \draw (-8,-2) to[short] (-3.125,-2);
      \draw (-5.375,-3.25) to[short] (-5.875,-3.25);
      \draw (-5.875,-3.25) to[short] (-6.5,-3.875);
      \draw (-6.5,-3.875) to[short] (-8.5,-3.875);
      \draw (-8,-2) to[short] (-10.25,-4.25);
      \draw (-10.25,-4.25) to[short] (-8.875,-4.25);
      \draw (-7.125,-3.625) to[short] (-7.125,-3.625);
      \draw (-8.5,-3.875) to[short] (-8.875,-4.25);
      \draw (-8.5,-3.875) to[short] (-8.5,-6);
      \draw (-6.5,-3.875) to[short] (-6.5,-6);
      \draw (-5.375,-3.25) to[short] (-5.375,-5.375);
      \draw (-10.25,-4.25) to[short] (-10.25,-5.25);
      \draw (-8.875,-4.25) to[short] (-8.875,-6.375);
      \draw (-5.875,-3.25) to[short] (-5.875,-5.375);
      \draw (-3.625,-4.625) to[short] (-4,-5);
      \draw (-5.375,-5.375) to[short] (-5.875,-5.375);
      \draw (-5.875,-5.375) to[short] (-6.5,-6);
      \draw (-6.5,-6) to[short] (-8.5,-6);
      \draw (-8.5,-6) to[short] (-8.875,-6.375);
      \draw (-9.5,-6.375) to[short] (-8.875,-6.375);
      \draw (-10.625,-5.625) to[short] (-10.625,-6.75);
      \draw (-9.875,-5.625) to[short] (-9.875,-6.75);
      \draw (-10.625,-5.625) to[short] (-9.875,-5.625);
      \draw (-10.625,-6.75) to[short] (-9.875,-6.75);
      \draw (-10.25,-5.25) to[short] (-10.625,-5.625);
      \draw (-9.5,-5.25) to[short] (-9.875,-5.625);
      \draw (-9.5,-5.25) to[short] (-9.5,-6.375);
      \draw (-9.5,-6.375) to[short] (-9.875,-6.75);
      \draw (-10.25,-5.25) to[short] (-9.5,-5.25);
      \draw (-2.375,-3.25) to[short] (-2.875,-3.75);
      \draw (-2.875,-3.75) to[short] (-2.875,-4.625);
      \draw (-3.125,-3.25) to[short] (-2.375,-3.25);
      \draw (-3.625,-4.625) to[short] (-2.875,-4.625);
      \draw (-3.125,-3.25) to[short] (-3.625,-3.75);
      \draw (-3.625,-3.75) to[short] (-3.625,-4.625);
      \node at (-9.5,-6.375) [circ] {};
      \node at (-8.5,-6) [circ] {};
      \node at (-5.875,-5.375) [circ] {};
      \node at (-5,-5) [circ] {};
      \draw [-{Stealth[scale=1.5]}] (-10.625,-6.75) -- (-11.5,-7.625);
      \draw [-{Stealth[scale=1.5]}] (-2.375,-4.125) -- (-0.875,-4.125);
      \draw [-{Stealth[scale=1.5]}] (-8,-2) -- (-8,-1.25);
      \node [font=\fontsize{10.8pt}{14.1pt}\selectfont] at (-11.625,-7.25) {x};
      \node [font=\fontsize{10.8pt}{14.1pt}\selectfont] at (-1.25,-4.375) {y};
      \node at (-9.375,-6.625) {$(d_1,e_1,0)$};
      \node at (-4.75,-5.25) {$(d_4,e_4,0)$};
      \node at (-5.875,-5.625) {$(d_3,e_3,0)$};
      \node at (-8.25,-6.25) {$(d_2,e_2,0)$};
      \node at (-10.625,-6.75) [circ] {};
      \node at (-3.625,-4.625) [circ] {};
      \node at (-10.25,-5.25) [circ] {};
      \node at (-3.125,-3.25) [circ] {};
      \node at (-8,-2) [circ] {};
      \node at (-10.625,-7) {$(a,0,0)$};
      \node at (-2,-4.375) {$(0,b,0)$};
      \node at (-7.625,-1.75) {$(0,0,c)$};
      \node at (-10.75,-5.125) {$(p_1=d_1, 0, q_1=1)$};
      \node at (-2.5,-2.875) {$(0, p_2, q_2)$};
      \draw (-4,-5) to[short] (-5,-5);
      \draw (-5,-2.875) to[short] (-5.375,-3.25);
      \draw (-5,-5) to[short] (-5.375,-5.375);
      \draw (-5,-2.875) to[short] (-5,-5);
      \draw (-3.125,-2) to[short] (-3.125,-3.25);
      \draw (-4,-2.875) to[short] (-4,-5);
      \draw (-4,-2.875) to[short] (-5,-2.875);
      \draw (-3.125,-2) to[short] (-4,-2.875);
      \draw (-3.625,-3.75) to[short] (-2.875,-3.75);
      \draw (-2.375,-3.25) to[short] (-2.375,-4.125);
      \draw (-2.375,-4.125) to[short] (-2.875,-4.625);
      \node at (-3.5,-4.875) {$(d_5,e_5,0)$};
      \node at (-2.375,-4.125) [circ] {};
      \draw (-7.25,-10.75) to[short] (-3,-10.75);
      \draw (-5.375,-12.125) to[short] (-5.875,-12.125);
      \draw (-5.875,-12.125) to[short] (-6.5,-12.75);
      \draw (-6.5,-12.75) to[short] (-8.5,-12.75);
      \draw (-7.25,-10.75) to[short] (-9.625,-13.125);
      \draw (-9.625,-13.125) to[short] (-8.875,-13.125);
      \draw (-7.125,-12.5) to[short] (-7.125,-12.5);
      \draw (-8.5,-12.75) to[short] (-8.875,-13.125);
      \draw (-6.5,-12.75) to[short] (-6.5,-14.875);
      \draw (-5.375,-12.125) to[short] (-5.375,-14.25);
      \draw (-9.625,-13.125) to[short] (-9.625,-14.125);
      \draw (-8.5,-12.75) to[short] (-8.5,-14.875);
      \draw (-5.875,-12.125) to[short] (-5.875,-14.25);
      \draw (-3.5,-13.375) to[short] (-4,-13.875);
      \draw (-5.375,-14.25) to[short] (-5.875,-14.25);
      \draw (-5.875,-14.25) to[short] (-6.5,-14.875);
      \draw (-6.5,-14.875) to[short] (-8.5,-14.875);
      \draw (-8.5,-14.875) to[short] (-8.875,-15.25);
      \draw (-10,-14.5) to[short] (-10,-15.625);
      \draw (-9.25,-14.5) to[short] (-9.25,-15.625);
      \draw (-10,-14.5) to[short] (-9.25,-14.5);
      \draw (-10,-15.625) to[short] (-9.25,-15.625);
      \draw (-9.625,-14.125) to[short] (-10,-14.5);
      \draw (-8.875,-14.125) to[short] (-9.25,-14.5);
      \draw (-8.875,-13.125) to[short] (-8.875,-14.125);
      \draw (-8.875,-15.25) to[short] (-9.25,-15.625);
      \draw (-9.625,-14.125) to[short] (-8.875,-14.125);
      \draw (-2.25,-12) to[short] (-2.75,-12.5);
      \draw (-2.75,-12.5) to[short] (-2.75,-13.375);
      \draw (-3,-12) to[short] (-2.25,-12);
      \draw (-3.5,-13.375) to[short] (-2.75,-13.375);
      \draw (-3,-12) to[short] (-3.5,-12.5);
      \draw (-3.5,-12.5) to[short] (-3.5,-13.375);
      \node at (-8.5,-14.875) [circ] {};
      \node at (-5.875,-14.25) [circ] {};
      \node at (-5,-13.875) [circ] {};
      \draw [-{Stealth[scale=1.5]}] (-10,-15.625) -- (-10.875,-16.5);
      \draw [-{Stealth[scale=1.5]}] (-2.25,-12.875) -- (-0.75,-12.875);
      \draw [-{Stealth[scale=1.5]}] (-7.25,-10.75) -- (-7.25,-10);
      \node [font=\fontsize{10.8pt}{14.1pt}\selectfont] at (-11,-16.125) {x};
      \node [font=\fontsize{10.8pt}{14.1pt}\selectfont] at (-1.125,-13.125) {y};
      \node at (-4.75,-14.125) {$(d_4,e_4-\delta,0)$};
      \node at (-5.875,-14.5) {$(d_3,e_3-\delta,0)$};
      \node at (-8.25,-15.125) {$(d_2,e_2-\delta,0)$};
      \node at (-10,-15.625) [circ] {};
      \node at (-3.5,-13.375) [circ] {};
      \node at (-9.625,-14.125) [circ] {};
      \node at (-3,-12) [circ] {};
      \node at (-7.25,-10.75) [circ] {};
      \node at (-10,-15.875) {$(a,0,0)$};
      \node at (-1.875,-13.125) {$(0,b-\delta,0)$};
      \node at (-6.875,-10.5) {$(0,0,c)$};
      \node at (-10.125,-14) {$(p_1=d_1, 0, q_1=1)$};
      \node at (-2.5,-11.75) {$(0, p_2-\delta, q_2)$};
      \draw (-4,-13.875) to[short] (-5,-13.875);
      \draw (-5,-11.75) to[short] (-5.375,-12.125);
      \draw (-5,-13.875) to[short] (-5.375,-14.25);
      \draw (-5,-11.75) to[short] (-5,-13.875);
      \draw (-3,-10.75) to[short] (-3,-12);
      \draw (-4,-11.75) to[short] (-4,-13.875);
      \draw (-4,-11.75) to[short] (-5,-11.75);
      \draw (-3,-10.75) to[short] (-4,-11.75);
      \draw (-3.5,-12.5) to[short] (-2.75,-12.5);
      \draw (-2.25,-12) to[short] (-2.25,-12.875);
      \draw (-2.25,-12.875) to[short] (-2.75,-13.375);
      \node at (-3.375,-13.625) {$(d_5,e_5-\delta,0)$};
      \node at (-2.25,-12.875) [circ] {};
      \draw [dashed] (-6.5,-6) -- (-6.5,-12.875);
      \draw [dashed] (-8.5,-6) -- (-8.5,-12.75);
      \draw [color=staircasegray, draw opacity=1, dashed] (-9.5,-6.375) -- (-9.5,-9.25);
      \draw [color=staircasegray, draw opacity=1, dashed] (-8.875,-6.375) -- (-8.875,-9.25);
      \draw [color=staircasegray, draw opacity=1, dashed] (-7.25,-2) -- (-7.25,-7);
      \draw [color=staircasegray, draw opacity=1, dashed] (-6.625,-2) -- (-6.625,-7.125);
      \draw [color=staircasedarkgray, draw opacity=0.46, line width=0.2pt] (-9.5,-9.25) rectangle (-8.875,-11.375);
      \begin{scope}
        \clip (-9.5,-9.25) rectangle (-8.875,-11.375);
        \foreach \x in {-2.13,-1.97,...,2.13} {
          \draw[color=staircasedarkgray, draw opacity=0.46, line width=0.2pt, rotate around={45:(-9.1875,-10.3125)}] ([xshift=\x cm]-9.1875,-12.4375) -- ([xshift=\x cm]-9.1875,-8.1875);
        }
      \end{scope}
      \draw [color=staircasedarkgray, draw opacity=0.62, line width=0.2pt] (-9.5,-9.25) -- (-7.25,-7) -- (-6.625,-7) -- (-8.875,-9.25) -- cycle;
      \begin{scope}
        \clip (-9.5,-9.25) -- (-7.25,-7) -- (-6.625,-7) -- (-8.875,-9.25) -- cycle;
        \foreach \x in {-2.83,-2.67,...,2.83} {
          \draw[color=staircasedarkgray, draw opacity=0.62, line width=0.2pt, rotate around={62:(-8,-8.125)}] ([xshift=\x cm]-8,-11.95038523331378) -- ([xshift=\x cm]-8,-4.29961476668622);
        }
      \end{scope}
      \draw [color=staircasedarkgray, draw opacity=0.62, line width=0.2pt] (-8.875,-9.25) -- (-8.875,-11.375) -- (-6.625,-9.125) -- (-6.625,-7) -- cycle;
      \begin{scope}
        \clip (-8.875,-9.25) -- (-8.875,-11.375) -- (-6.625,-9.125) -- (-6.625,-7) -- cycle;
        \foreach \x in {-3.46,-3.26,...,3.46} {
          \draw[color=staircasedarkgray, draw opacity=0.62, line width=0.2pt, rotate around={45:(-7.75,-9.125)}] ([xshift=\x cm]-7.75,-13.584833581769303) -- ([xshift=\x cm]-7.75,-4.665166418230697);
        }
      \end{scope}
      \node [font=\fontsize{11.9pt}{15.5pt}\selectfont, color=black, text opacity=0.54] at (-10.375,-10.375) {removed};
      \draw (-5.25,-19.625) to[short] (-3,-19.625);
      \draw (-5.375,-21) to[short] (-5.875,-21);
      \draw (-5.875,-21) to[short] (-6.5,-21.625);
      \draw (-6.5,-21.625) to[short] (-7.25,-21.625);
      \draw (-5.25,-19.625) to[short] (-7.25,-21.625);
      \draw (-7.125,-21.375) to[short] (-7.125,-21.375);
      \draw (-8.5,-21.625) to[short] (-8.5,-21.625);
      \draw (-6.5,-21.625) to[short] (-6.5,-22.625);
      \draw (-5.375,-21) to[short] (-5.375,-23.125);
      \draw (-7.25,-21.625) to[short] (-7.25,-22.625);
      \draw (-5.875,-21) to[short] (-5.875,-23.125);
      \draw (-3.5,-22.25) to[short] (-4,-22.75);
      \draw (-5.375,-23.125) to[short] (-5.875,-23.125);
      \draw (-5.875,-23.125) to[short] (-6.5,-23.75);
      \draw (-7.625,-23) to[short] (-7.625,-24.125);
      \draw (-6.875,-23) to[short] (-6.875,-24.125);
      \draw (-7.625,-23) to[short] (-6.875,-23);
      \draw (-7.625,-24.125) to[short] (-6.875,-24.125);
      \draw (-7.25,-22.625) to[short] (-7.625,-23);
      \draw (-6.5,-22.625) to[short] (-6.875,-23);
      \draw (-6.5,-23.75) to[short] (-6.875,-24.125);
      \draw (-7.25,-22.625) to[short] (-6.5,-22.625);
      \draw (-2.25,-20.875) to[short] (-2.75,-21.375);
      \draw (-2.75,-21.375) to[short] (-2.75,-22.25);
      \draw (-3,-20.875) to[short] (-2.25,-20.875);
      \draw (-3.5,-22.25) to[short] (-2.75,-22.25);
      \draw (-3,-20.875) to[short] (-3.5,-21.375);
      \draw (-3.5,-21.375) to[short] (-3.5,-22.25);
      \node at (-5.875,-23.125) [circ] {};
      \node at (-5,-22.75) [circ] {};
      \draw [-{Stealth[scale=1.5]}] (-7.625,-24.125) -- (-8.5,-25);
      \draw [-{Stealth[scale=1.5]}] (-2.25,-21.75) -- (-0.75,-21.75);
      \draw [-{Stealth[scale=1.5]}] (-5.25,-19.625) -- (-5.25,-18.875);
      \node [font=\fontsize{10.8pt}{14.1pt}\selectfont] at (-11,-25.75) {x};
      \node [font=\fontsize{10.8pt}{14.1pt}\selectfont] at (-1.125,-22) {y};
      \node at (-4.75,-23) {$(d_4',e_4'-\delta',0)$};
      \node at (-5.75,-23.375) {$(d_3',e_3'-\delta',0)$};
      \node at (-7.625,-24.125) [circ] {};
      \node at (-3.5,-22.25) [circ] {};
      \node at (-7.25,-22.625) [circ] {};
      \node at (-3,-20.875) [circ] {};
      \node at (-5.25,-19.625) [circ] {};
      \node at (-7.625,-24.375) {$(a-p_1+d_2,0,0)$};
      \node at (-1.875,-22) {$(0,b'-\delta',0)$};
      \node at (-4.875,-19.375) {$(0,0,c)$};
      \node at (-7.75,-22.5) {$(d_2, 0, q_1=1)$};
      \node at (-2.5,-20.625) {$(0, p_2'-\delta', q_2)$};
      \draw (-4,-22.75) to[short] (-5,-22.75);
      \draw (-5,-20.625) to[short] (-5.375,-21);
      \draw (-5,-22.75) to[short] (-5.375,-23.125);
      \draw (-5,-20.625) to[short] (-5,-22.75);
      \draw (-3,-19.625) to[short] (-3,-20.875);
      \draw (-4,-20.625) to[short] (-4,-22.75);
      \draw (-4,-20.625) to[short] (-5,-20.625);
      \draw (-3,-19.625) to[short] (-4,-20.625);
      \draw (-3.5,-21.375) to[short] (-2.75,-21.375);
      \draw (-2.25,-20.875) to[short] (-2.25,-21.75);
      \draw (-2.25,-21.75) to[short] (-2.75,-22.25);
      \node at (-3.375,-22.5) {$(d_5',e_5'-\delta',0)$};
      \node at (-2.25,-21.75) [circ] {};
      \draw [dashed] (-6.5,-14.875) -- (-6.5,-21.625);
      \draw [color=black, draw opacity=0.41, dashed] (-8.5,-14.875) -- (-8.5,-17.875);
      \draw [color=staircasedarkgray, draw opacity=0.62, line width=0.2pt] (-8.875,-18.25) -- (-8.875,-20.375) -- (-8.5,-20) -- (-8.5,-17.875) -- cycle;
      \begin{scope}
        \clip (-8.875,-18.25) -- (-8.875,-20.375) -- (-8.5,-20) -- (-8.5,-17.875) -- cycle;
        \foreach \x in {-2.26,-2.06,...,2.26} {
          \draw[color=staircasedarkgray, draw opacity=0.62, line width=0.2pt, rotate around={12:(-8.625,-19.125)}] ([xshift=\x cm]-8.625,-22.388984276009793) -- ([xshift=\x cm]-8.625,-15.861015723990207);
        }
      \end{scope}
      \draw [color=staircasedarkgray, draw opacity=0.46, line width=0.2pt] (-9.625,-18.25) rectangle (-8.875,-20.375);
      \begin{scope}
        \clip (-9.625,-18.25) rectangle (-8.875,-20.375);
        \foreach \x in {-2.13,-1.97,...,2.13} {
          \draw[color=staircasedarkgray, draw opacity=0.46, line width=0.2pt, rotate around={45:(-9.25,-19.3125)}] ([xshift=\x cm]-9.25,-21.4375) -- ([xshift=\x cm]-9.25,-17.1875);
        }
      \end{scope}
      \draw [color=staircasegray, draw opacity=0.81, dashed] (-9.625,-14.125) -- (-9.625,-18.25);
      \draw [color=staircasedarkgray, draw opacity=0.62, line width=0.2pt] (-9.625,-18.25) -- (-7.25,-15.875) -- (-5.375,-15.875) -- (-7.25,-17.875) -- (-8.5,-17.875) -- (-8.875,-18.25) -- cycle;
      \begin{scope}
        \clip (-9.625,-18.25) -- (-7.25,-15.875) -- (-5.375,-15.875) -- (-7.25,-17.875) -- (-8.5,-17.875) -- (-8.875,-18.25) -- cycle;
        \foreach \x in {-3.84,-3.68,...,3.84} {
          \draw[color=staircasedarkgray, draw opacity=0.62, line width=0.2pt, rotate around={62:(-7.75,-17.375)}] ([xshift=\x cm]-7.75,-22.215447563520776) -- ([xshift=\x cm]-7.75,-12.534552436479224);
        }
      \end{scope}
      \draw [color=staircasegray, draw opacity=0.81, dashed] (-8.875,-14.125) -- (-8.875,-18.25);
      \node [font=\fontsize{11.9pt}{15.5pt}\selectfont, color=black, text opacity=0.54] at (-10.5,-19.25) {removed};
      \draw [color=black, draw opacity=0.48, line width=0.2pt] (-8.5,-17.875) -- (-8.5,-20) -- (-7.25,-20) -- (-7.25,-17.875) -- cycle;
      \begin{scope}
        \clip (-8.5,-17.875) -- (-8.5,-20) -- (-7.25,-20) -- (-7.25,-17.875) -- cycle;
        \foreach \x in {-2.23,-2.03,...,2.23} {
          \draw[color=black, draw opacity=0.48, line width=0.2pt, rotate around={45:(-7.875,-18.875)}] ([xshift=\x cm]-7.875,-22.107692682707253) -- ([xshift=\x cm]-7.875,-15.642307317292747);
        }
      \end{scope}
      \draw [color=black, draw opacity=0.48, line width=0.2pt] (-7.25,-17.875) -- (-7.25,-20) -- (-5.375,-17.875) -- (-5.375,-15.875) -- cycle;
      \begin{scope}
        \clip (-7.25,-17.875) -- (-7.25,-20) -- (-5.375,-17.875) -- (-5.375,-15.875) -- cycle;
        \foreach \x in {-3.29,-3.09,...,3.29} {
          \draw[color=black, draw opacity=0.48, line width=0.2pt, rotate around={12:(-6.25,-17.875)}] ([xshift=\x cm]-6.25,-22.16905651902912) -- ([xshift=\x cm]-6.25,-13.58094348097088);
        }
      \end{scope}
      \draw [color=staircasegray, draw opacity=1, dashed] (-7.25,-14.875) -- (-7.25,-17.875);
    \end{circuitikz}
  }%

\begin{proof}
We work in the quotient ring $R:=\frac{\mathbb{C}[x,y,z]}{(x^a+y^b+z^c, x^{d_1}y^{e_1})}$. So the ideal $G_1'=(x^{p_1}z+y^{p_2}z^{q_2})$ and $G_2'=(x^{d_1}z+y^{p_2-(e_2-e_1)}z^{q_2})$. The ideal $I$ becomes $$I=(x^a,y^b, x^{d_2}y^{e_2},...,x^{d_k}y^{e_k}, x^{p_1}z,y^{p_2}z^{q_2}).$$

Because the generators of $G_1'$ and $G_2'$ are not zero divisors, and in particular $G_1$ and $G_2$ are both generated by regular sequences, we use the above Lemma \ref{lem_single_gen} and compute

$$K:=G_2:(G_1:I)=G_2' \cdot I:G_1'.$$

In particular, a direct computation of generators shows that $I_3 \cdot G_1' \subset G_2' \cdot I$, so $I_3 \subset K$. 

The following equation of dimensions of vector spaces also holds by a direct computation of Grobner basis:

$$\dim_{\mathbb{C}}\frac{\mathbb{C}[x,y,z]}{I_3}+\dim_{\mathbb{C}}\frac{\mathbb{C}[x,y,z]}{G_1}=\dim_{\mathbb{C}}\frac{\mathbb{C}[x,y,z]}{I}+\dim_{\mathbb{C}}\frac{\mathbb{C}[x,y,z]}{G_2}.$$

So $I_3=K$.
\end{proof}

\begin{example}
\label{example_red_opera_house}
We show the consecutive reduction of an opera house that can be computed in Macaulay2, by running the above algorithm recursively.

We start from the following ideal $I$, with the corresponding $G_1$ and $G_2$. Note that $p_1=d_1$, $\delta=e_2-e_1= 1$, $a=12, b=9,c=6$.

\begin{verbatim}
I = ideal(x^12,y^9, z^6, x^7*y, x^5*y^2,x^3*y^4,x^2*y^5,x*y^7, x^7*z,y^8*z^4);
G1= ideal(x^12+y^9+z^6, x^7*y, x^7*z+y^8*z^4);
I2=G1:I
G2=ideal(x^12+z^6+y^9, x^7*y, x^7*z+y^7*z^4);
I3=G2:I2
I3= ideal(z^6,x^2*y^4,x^3*y^3,x^5*y,x*y^6,x^7*z,y^8,y^7*z^4,x^12)

\end{verbatim}

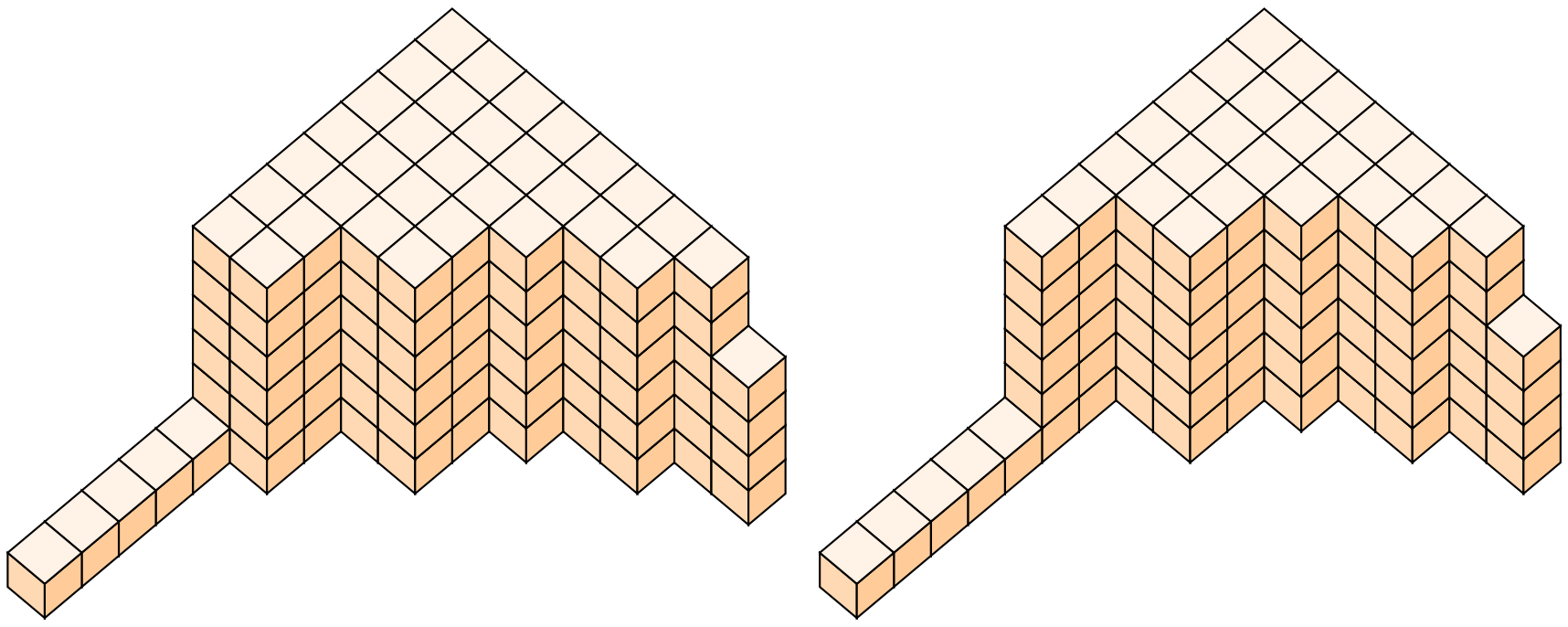
\monomialIdealPicture{{6,6,6,6,6,6,6,1,1,1,1,1},{6,6,6,6,6,6,6},{6,6,6,6,6},{6,6,6,6,6},{6,6,6},{6,6},{6,6},{6},{4}}\monomialIdealPicture{{6,6,6,6,6,6,6,1,1,1,1,1},{6,6,6,6,6},{6,6,6,6,6},{6,6,6},{6,6},{6,6},{6},{4}} 

Now we reduce again, setting $I=I_3$ from the above computation. Now $p_1=7, d_1=5$, $\delta=e_2-e_1=2$, and 

\begin{verbatim}
I = ideal(x^12,y^8,z^6,x*y^6,x^2*y^4,x^3*y^3,x^5*y, x^7*z,y^7*z^4);
G1= ideal(x^12+y^8+z^6, x^5*y, x^7*z+y^7*z^4);
I2=G1:I
G2=ideal(x^12+y^8+z^6, x^5*y, x^5*z+y^5*z^4);
I3=G2:I2
I3=ideal(x^2*y^2,x^3*y,x*y^4,z^6,x^5*z,y^6,y^5*z^4,x^10)
\end{verbatim}

Now the opera house is reduced to the following staircase, with the inner corner $x^5y$ removed.

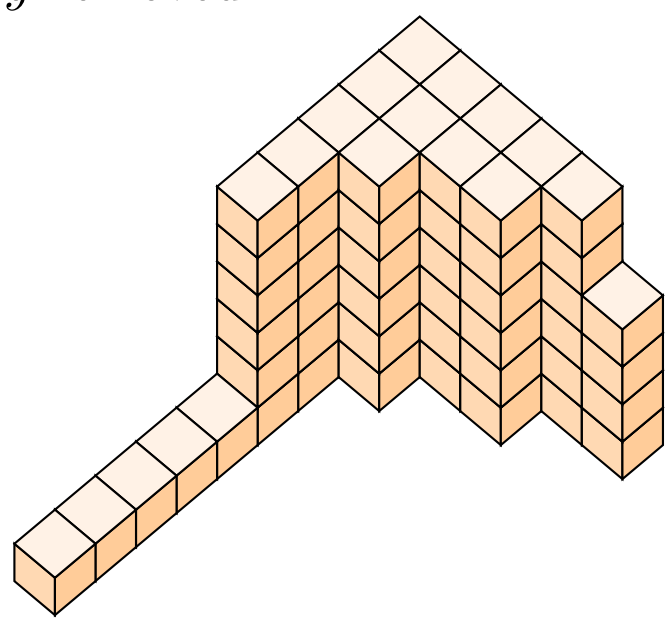
\monomialIdealPicture{{6,6,6,6,6,1,1,1,1,1},{6,6,6},{6,6},{6,6},{6},{4}}

Reduce two more times:
First using $e_2-e_1=1$, $p_1=5, d_1=3$, reduce the inner corner $x^3y$.
\begin{verbatim}
I = ideal(x^10,y^6,z^6, x*y^4,x^2*y^2,x^3*y, x^5*z,y^5*z^4);
G1= ideal(x^10+y^6+z^6, x^3*y, x^5*z+y^5*z^4);
I2=G1:I
G2=ideal(x^10+y^6+z^6, x^3*y, x^3*z+y^4*z^4);
I3=G2:I2
I3=ideal(x^2*y,x^3*z,x*y^3,y^5,z^6,y^4*z^4,x^8)
\end{verbatim}

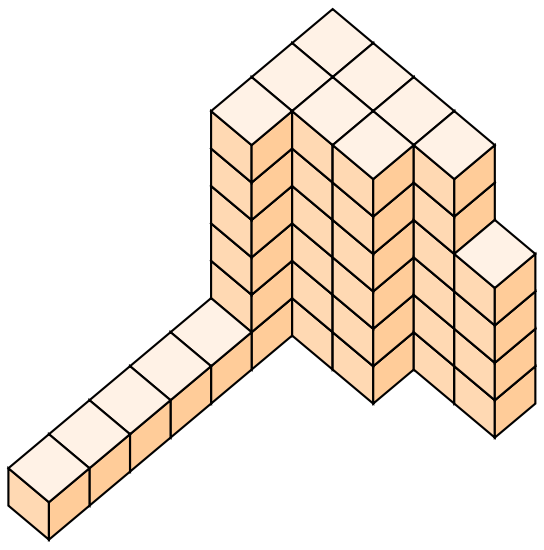
\monomialIdealPicture{{6,6,6,1,1,1,1,1},{6,6},{6,6},{6},{4}}

Then using $e_2-e_1=2$, $p_1=3, d_1=2$, reduce the inner corner $x^2y$.

\begin{verbatim}
I = ideal(x^8,y^5,z^6, x*y^3,x^2*y, x^3*z,y^4*z^4);
G1= ideal(x^8+y^5+z^6, x^2*y, x^3*z+y^4*z^4);
I2=G1:I
G2=ideal(x^8+y^5+z^6, x^2*y, x^2*z+y^2*z^4);
I3=G2:I2
I3=ideal(x*y,x^2*z,y^3,z^6,y^2*z^4,x^7)
\end{verbatim}
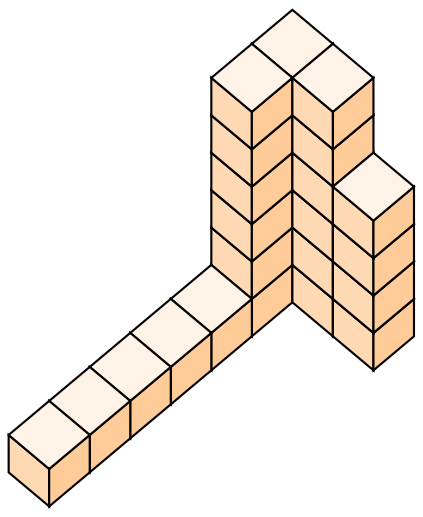
\monomialIdealPicture{{6,6,1,1,1,1,1},{6},{4}} 

Now we have one unique inner corner in the bottom layer. We move on to the reduction of the next algorithm.
\end{example}

\subsection{The last bit of the induction}
At this point, we reduce all opera houses of more than one inner corners in the bottom layer, to having one unique corner $x^{d}y^{e}$. The above algorithm does not work anymore, because the third generator of $G_2$ is no longer well defined. 

In particular, any such opera house is generated as 
$$I=(x^a, y^b,z^c, x^dy^e, x^{p_1}z^{q_1}, y^{p_2}z^{q_2}).$$

These opera houses have one unique pillar (possibly cut through by an inner corner of the bottom layer, if $p_1>d$ or $p_2>e$) in the middle, with two wings on the side. We draw a few examples to show that there are still many variations.

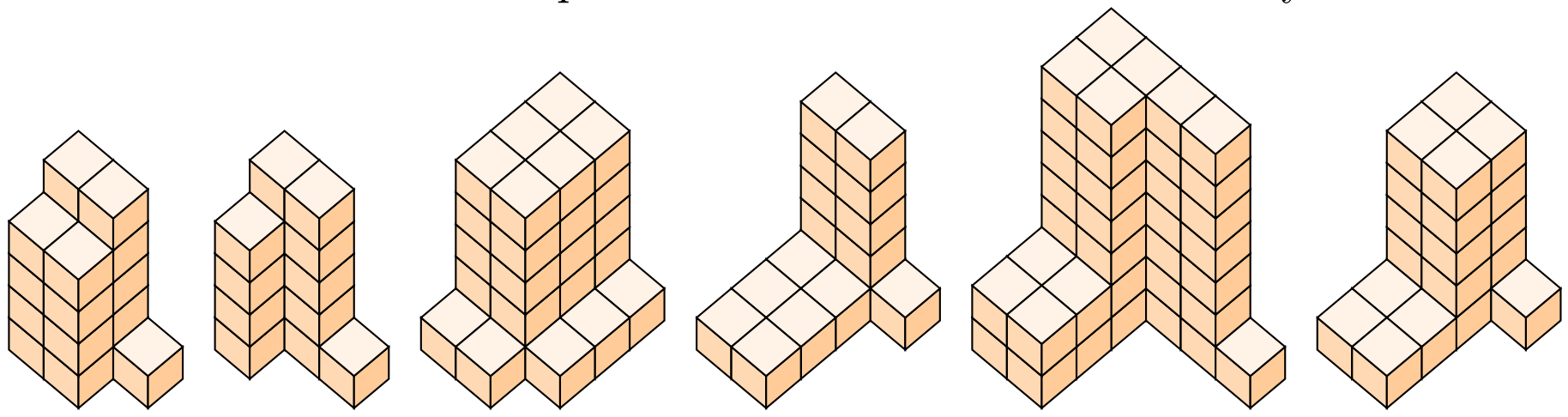
\monomialIdealPicture  {{5,4},{5,4},{1}}\monomialIdealPicture{{5,4},{5},{1}}\monomialIdealPicture{{5,5,5,1},{5,5,5,1},{1,1,1}}\monomialIdealPicture{{5,1,1,1},{5,1,1,1},{1}}\monomialIdealPicture{{7,7,2,2},{7,7,2,2},{7},{7},{1}}\monomialIdealPicture{{5,5,1,1},{5,5,1,1},{1}}

The algorithm in the last section, with a similar idea but some modification, still works when $p_1>d$ or $p_2>e$.

\textbf{Algorithm:}
When $p_1>d$ or $p_2>e$, consider the complete intersection ideals 

$$G_1=(x^a+y^b+z^c, x^dy^e,x^{p_1}z^{q_1}+y^{p_2}z^{q_2}).$$
$$G_2=(x^a+y^b+z^c, x^{d}y^{e},x^dz^{q_1}+y^ez^{q_2}).$$

Take a 2 step linkage:
Compute $I_2=G_1:I$;  $I_3=G_2:I_2$.

\begin{lemma}
 The above algorithm links the ideal $I$ of the opera house to 

$$I_3 =
( x^{a-(p_1-d)}, y^{b-(p_2-e)}, z^c, x^{d}y^{e}, x^{d}z^{q_1}, y^{e}z^{q_2}).$$

\end{lemma}

\begin{proof}
The proof strategy is also similar. We work in the quotient ring $R:=\frac{\mathbb{C}[x,y,z]}{(x^a+y^b+z^c, x^{d}y^{e})}$. So the ideal $G_1'=(x^{p_1}z^{q_1}+y^{p_2}z^{q_2})$ and $G_2'=(x^{d}z^{p_1}+y^{e}z^{q_2})$. The ideal $I$ becomes $$I=(x^a,y^b,x^{p_1}z^{q_1},y^{p_2}z^{q_2}).$$

Because the generators of $G_1'$ and $G_2'$ are not zero divisors, and in particular $G_1$ and $G_2$ are both generated by regular sequences, we use Lemma \ref{lem_single_gen} and compute

$$K:=G_2:(G_1:I)=G_2' \cdot I:G_1'.$$

And a direct computation of generators shows that $I_3 \cdot G_1' \subset G_2' \cdot I$, so $I_3 \subset K$. 

The following equation of dimensions of vector spaces also holds by a direct computation of Grobner basis:

$$\dim_{\mathbb{C}}\frac{\mathbb{C}[x,y,z]}{I_3}+\dim_{\mathbb{C}}\frac{\mathbb{C}[x,y,z]}{G_1}=\dim_{\mathbb{C}}\frac{\mathbb{C}[x,y,z]}{I}+\dim_{\mathbb{C}}\frac{\mathbb{C}[x,y,z]}{G_2}.$$

So $I_3=K$.

\end{proof}

\begin{example}
In the following example, $p_1=2>d=1$, and $p_2=3 > e=2$.
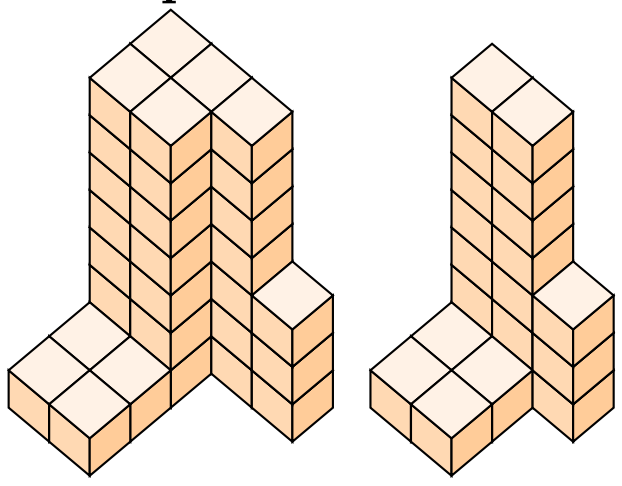
\monomialIdealPicture{{7,7,1,1},{7,7,1,1},{7},{3}}\monomialIdealPicture{{7,1,1},{7,1,1},{3}}

\begin{verbatim}
I = ideal(x^4,y^4,z^7,x*y^2,x^2*z,y^3*z^3)
G1= ideal(x*y^2,x^2*z+y^3*z^3,x^4+y^4+z^7);
I2=G1:I
G2=ideal(x^4+z^7+y^4, x*y^2, x*z+y^2*z^3);
I3=G2:I2
I3=ideal(x^3,y^3,z^7,x*z,x*y^2,y^2*z^3)
\end{verbatim}
\end{example}

\begin{example}
In the following example, $p_1=2>d=1$, and $p_2= e=1$.
\monomialIdealPicture{{5,5,3,3},{5,5,3,3},{1}} \monomialIdealPicture{{5,3,3},{5,3,3},{1}}
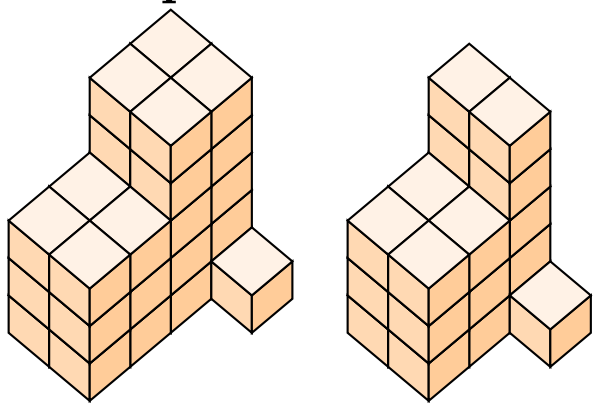
\begin{verbatim}
I = ideal(x*y^2,x^2*z^3,x^4,y^2*z,y^3,z^5)
G1= ideal(x*y^2,x^2*z^3+y^2*z,x^4+y^3+z^5);
I2=G1:I
G2=ideal(x*y^2,x*z^3+y^2*z,x^4+y^3+z^5);
I3=G2:I2
I3=ideal(y^2*z,y^3,x*y^2,x^3,x*z^3,z^5)
\end{verbatim}
\end{example}

The result of the above algorithm is an opera house with one unique inner corner $x^dy^e$ where $p_1=d$ and $p_2=e$.

However, when $p_1=d$ and $p_2=e$, the idea of the above algorithm doesn't work, even with modifications. For some unknown and interesting reason, we always get back the original ideal $I$ as the result of the 2 step linkage in the algorithm, so we develop another strategy.

\subsection{Reducing the height of the wings both to 1.}

Now we have reduced all opera houses to $I=(x^a,y^b, z^c, x^dy^e,x^dz^{q_1}, y^ez^{q_2})$, and either $q_1=1$ or $q_2=1$. We use the following algorithm to reduce both $q_1$ and $q_2$ to $1$.

Assume without loss of generality that $q_1=1$. 

\textbf{Algorithm:}

Consider the complete intersection ideal $G=(x^a+z^c, x^dz,y^b)$.
Take a 1 step linkage: $I_2=G:I$.

\begin{lemma}
 The above algorithm links the ideal $I$ of the opera house to 
 $$I_2= (x^a,y^b, z^{c-q_2+1}, x^dy^{b-e},x^dz, y^{b-e}z).$$
\end{lemma}

\begin{proof}
We can directly compute that $I_2 \cdot I \subset G$, and that the dimension satisfies $\dim_{\mathbb{C}}\frac{\mathbb{C}[x,y,z]}{(I_2)}+\dim_{\mathbb{C}}\frac{\mathbb{C}[x,y,z]}{(I)}=\dim_{\mathbb{C}}\frac{\mathbb{C}[x,y,z]}{(G)}$, so $I_2 \cdot I = G$.
\end{proof}

\begin{example}
In this example, $a=3, b=3, c=5,d=2, e=2, q_2=2$. The staircases $I$ and $I_2$ are shown below.

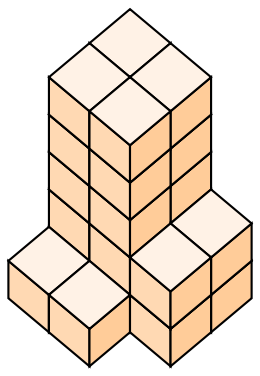
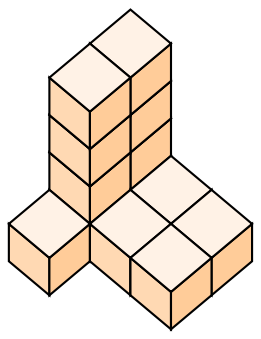
\monomialIdealPicture{{5,5,1},{5,5,1},{2,2}} 
\monomialIdealPicture{{4,4,1},{1,1},{1,1}}

\begin{verbatim}
I = ideal(x^2*y^2,x^2*z, x^3,y^2*z^2,y^3,z^5);
G=ideal(x^2*z,y^3,x^3+z^5);
I2=G:I
I2=ideal(y*z,x^2*z,y^3,x^2*y,x^3,z^4)
\end{verbatim}

\end{example}

\begin{example}
In this example, $a=3, b=3, c=5,d=2, e=1, q_2=3$.
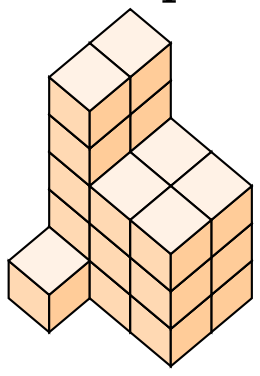
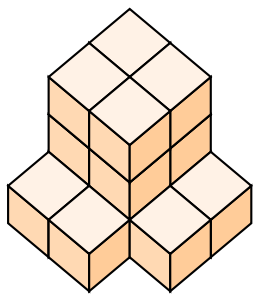
\monomialIdealPicture{{5,5,1},{3,3},{3,3}}
\monomialIdealPicture{{3,3,1},{3,3,1},{1,1}}
\begin{verbatim}
I=ideal(x^2*y,x^2*z,y*z^3,y^3,x^3,z^5);
G=ideal(x^2*z,x^3+z^5,y^3);
I2=G:I
I2=ideal(z^3,y^2*z,x^2*z,y^3,x^3,x^2*y^2)
\end{verbatim}
\end{example}

\subsection{Thin opera house.\\}

Now we have reduced everything to $I=(x^a,y^b, z^c, x^dy^e,x^dz, y^ez)$. There is an intermediate step that they are all linked to what we call thin opera house, before we reduce everything to a tripod. A \textbf{thin opera house} is an opera house $I$ such that at least one of the three monomials, $xy, xz, yz$, is contained in $I$. 

\begin{lemma}
Any opera house $I=(x^a,y^b, z^c, x^dy^e,x^dz, y^ez)$ where the height of both of the wings are $1$, is either already a tripod, or linked to a thin opera house generated as $$J=(x^a,y^b, z^c, xy,xz^{c-1}, yz^{c-1}).$$
\end{lemma}

\begin{proof}
We argue by discussing the possible types of wings that an opera house has.

If the wings of the staircase of $I$ are two folded rows with height both being $1$, then its complement in $(x^a, y^b, z^c)$ is a thin opera house containing the generator $xy$, and the height of the wings are $c-1$. In particular, the thin opera house is $(x^a,y^b, z^c, x^dy^e,x^dz, y^ez)$.

If the wings of the staircase of $I$ are two spread-out rows, then $d=e=1$ and $xy, yz, xz \in I$. The opera house is a tripod.  

If the wings' bottom layer is consisted of one box (of size $1$ by $1$) and a rectangle of arbitrary size, then either $d=1$ or $e=1$. 

Assume that $d=1$, and the case of $e=1$ follows similarly. Because $d=1$, the opera house is $(x^a, y^b,z^c, xy^e, xz,y^e)$. Because the right wing has size $1$ by $1$, $e=b-1$. Swap the variables $z$ with $y$, the opera house is now $(x^a, z^b,y^c, xz^e, xy, z^e)$, which is the proposed form.
\end{proof}

\begin{comment}
\begin{remark}
Then $I$ is generated as $I=(x^a,y^b, z^c, xy,xz^{q_1}, yz^{q_2})$. 
The possible height of the wings $q_1$ and $q_2$, is either $1$ or $c-1$. If the height is $1$ then we have a tripod already, so we assume that $q_1=q_2=c-1$.
\end{remark}
The possible length of the wings, is either $1$, $a-1$ (if it's an opera house in the $y$ direction) or $b-1$ (if it's an opera house in the $x$ direction).
Note that because we required $d=p_1$ and $e=p_2$, the height of the two wings (viewed as an opera house in either $x$ or $y$ directions) are the same. 
\begin{lemma}
\end{lemma}
\begin{proof}
We prove the case when $xz \in I$ and the opera house is in the $y$ direction, and the other case when $yz \in I$ follows similarly. Because we 
\end{proof}
\end{comment}

And a thin opera house is linked to a tripod, by applying to it the same algorithms we developed in the last sections.

\textbf{Algorithm:}
Consider the complete intersection ideal $G_1=(x^a+y^b+z^c, xy,xz^{c-1}+yz^{c-1})$.

Take a 1 step linkage: $I_2=G_1:J$.

Consider the complete intersection ideal that is a cube $G_2=(x^2,y^2,z^c)$.

Take another 1 step linkage: $I_3=G_2:I_2$.

\begin{lemma}
The above algorithm links any thin opera house $J$ to a tripod $$I_3=(x^2,y^2,z^c,xy,xz,yz).$$
\end{lemma}

\begin{proof}
We claim that $I_2=(x^2,y^2,z^c,xy,xz^{c-1},yz^{c-1})$ by computing $J \cdot I_2 \in G_1$, and $\dim_{\mathbb{C}}\frac{\mathbb{C}[x,y,z]}{(I_2)}+\dim_{\mathbb{C}}\frac{\mathbb{C}[x,y,z]}{(J)}=\dim_{\mathbb{C}}\frac{\mathbb{C}[x,y,z]}{(G_1)}$. Then $I_3$ is the complement of $I_2$ in a cube.
\end{proof}

\begin{example}
In this example, $a=4, b=6,c=4, q_1=q_2=3$.

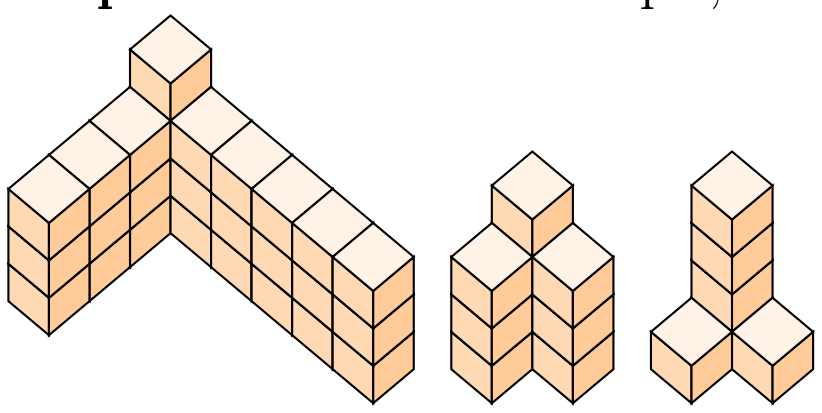
\monomialIdealPicture{{4,3,3,3},{3},{3},{3},{3},{3}}\monomialIdealPicture{{4,3},{3}}\monomialIdealPicture{{4,1},{1}}

\begin{verbatim}
I = ideal(x^4,y^6, z^4, x*y,x*z^3,y*z^3);
G1= ideal(x^4+y^6+z^4, x*y, x*z^3+y*z^3);
I2=G1:I
I2=ideal(y^2,x*y,x^2,z^4,y*z^3,x*z^3)
G2=ideal(x^2,y^2,z^4);
I3=G2:I2
I3=ideal(y*z,x*z,y^2,x*y,x^2,z^4)
\end{verbatim}
\end{example}

\subsection{Singularity of opera houses.}

So far, we have showed that all opera houses are in the same linkage class of a tripod. We further show in this section that all opera house have the same singularity type.

\begin{theorem}
For any monomial ideal $I \in \Hilb^d(\mathbb{A}^3)$, if its corresponding staircase $E$ is an opera house, then, the singularity at $[\mathrm{S}/\mathrm{I}] \in \Hilb^d(\mathbb{A}^3)$ is smoothly equivalent to the vertex of a cone over the Grassmannian $\mathrm{Gr}(2,6) \hookrightarrow \mathbb{P}^{14}$ in its Pl\''ucker embedding.
\end{theorem}

\begin{proof}
Because any ideal $I$ whose staircase $E$ is an opera house is in the same linkage class as a tripod, and by \cite{AJR1} \cite{Katz} the tripod has the singularity type of the vertex of a cone over the Grassmannian $\mathrm{Gr}(2,6) \hookrightarrow \mathbb{P}^{14}$, we apply Theorem \ref{AJR_singularity_type} of \cite{AJR1} and conclude the singularity type of all opera houses.
\end{proof}

\subsection{Linkage classes of rotating staircases.}

We have proved in Corollary \ref{cor_linkage_pairs}, that the rotating staircases are in pairs of linkage classes. We further prove in this section, that many pair are not linked to each other, and are not linked to a tripod. We use the following classification result of linkage classes.

\begin{proposition}
\cite{HU}
\label{prop_HU}
Let $R = k[X_1, \dots, X_n]$ be a graded polynomial ring with $\deg X_i = 1$ $(1 \le i \le n)$, let $I$ be a homogeneous $R$-ideal of grade $g > 0$ with a homogeneous minimal resolution as in (5.7), set $d = \min\{n_{1j}\}$, and assume that $\operatorname{grade}(I_d R) = g$. Write
$$R' = R_{(X_1, \dots, X_n)}, \quad I' = IR.$$
Then either
\begin{enumerate}
\item[a)] $\max\{n_{gj}\} > (g - 1)d$, and (if $k$ is infinite) there exists an $R'$-ideal $J'$, doubly linked to $I'$, such that $e(R'/J') < e(R'/I')$, and (for a suitable set of indeterminates $V$) $e(R'(V)/L^2(I')) < e(R'/I')$, or,
    
\item[b)] $\max\{n_{gj}\} \le (g - 1)d$, and for any $R'$-ideal $J'$ in the even linkage class of $I'$, $e(R'/J') \ge e(R'/I')$, and for any $m$, $e(R'(V)/L^{2m}(I')) = e(R'/I')$.
\end{enumerate}
\end{proposition}

\begin{definition}
Let $d$ be smallest total degree among the inner corners of the staircase $E$, and $m$ maximal total degree among the outer corners of the staircase $E$. If $m+3 \leq 2d$, then $E$ is called linkage-minimal.
\end{definition}

\begin{lemma}
Consider a rotating staircase $E$ in $\Hilb^{d_1}(\mathbb{A}^3)$ and its pair $F$ in $\Hilb^{d_2}(\mathbb{A}^3)$, that union to a minimal cube as in Corollary \ref{cor_linkage_pairs}. If $E$ and $F$ are both linkage-minimal, and there is another pair of ideals $E'$ in $\Hilb^{d_3}(\mathbb{A}^3)$ and $F'$ in $\Hilb^{d_4}(\mathbb{A}^3)$ that are also both linkage-minimal, and that the number of boxes $d_1$ and $d_2$ are distinct from both $d_3, d_4$ then $E,F$ are in different linkage classes of $E',F'$.
\end{lemma}

\begin{proof}
This lemma is a direct consequence of Proposition \ref{prop_HU}.

A polynomial ring is Cohen–Macaulay, so $g$ is equal to the height of $I$, which is $3$ for any ideal $I \in \Hilb^d(\mathbb{A}^3)$. 

Let $I$ denote the ideal of the rotating staircase $E$. The number $d$ is smallest total degree among the minimal generators of $I$, which are the inner corners of $E$. The number $n_{gj}=n_{3j} $ is equal to $3$ plus the maximal total degree among the outer corners of the staircase $E$. If $m+3 \leq 2d$ for $E$ and $F$, then the part b) of Proposition \ref{prop_HU} applies, and if the number of boxes of $E',F'$ are both not the same as either $E$ or $F$, then they are in different linkage classes.
\end{proof}

\begin{example}
\label{rot_star_link_class} 
All of the following pairs of rotating boxes are in different linkage classes.

\begin{enumerate}
\item Suppose we take $J_1=(z^4,xz^3,y^3z,y^4,x^3y,x^4,xyz^2,x^2yz,xy^2z)$, and its linkage
$I_1=(x^4,y^4,z^4):J$. The left picture is $J_1$ and the right picture is $I_1$. 
\monomialIdealPicture{{4,3,3,3},{4,2,1},{4,1,1},{1,1,1}} \monomialIdealPicture{{4,3,3,3},{4,3,3},{4,3,2},{1,1,1}} 

For $J_1$, $8 \leq  (3 - 1) \cdot 4$, and for $I_1$, $8 \leq (3- 1)\cdot 4$. So any ideal that is linked to both $I_1$ and $J_1$ must
have at least more boxes than $J_1$.

\item Similarly, we take $J_2=(z^4,xz^3,y^3z,y^4,x^3y,x^4,xyz^2,x^2yz);$, and its linkage
$I_2=(x^4,y^4,z^4):J_2$. The left picture is $J_2$ and the right picture is $I_2$. 
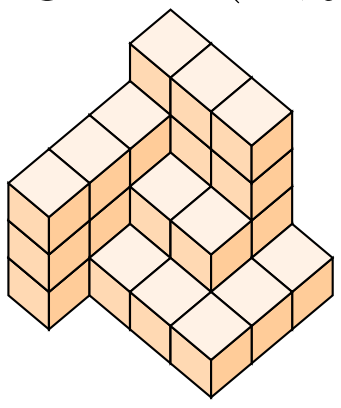
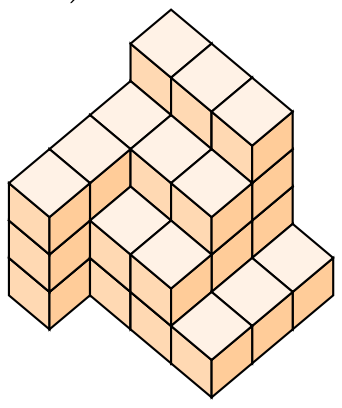
\monomialIdealPicture{{4,3,3,3},{4,2,1},{4,2,1},{1,1,1}} \monomialIdealPicture{{4,3,3,3},{4,3,2},{4,3,2},{1,1,1}} 

For $J_2$, $8 \leq  (3 - 1) \cdot 4$, and for $I_2$, $8 \leq (3- 1)\cdot 4$. So any ideal that is linked to both $I_2$ and $J_2$ must
have at least more boxes than $J_2$.

\item Similarly, we take $J_3=(z^4,xz^3,y^3z,y^4,x^3y,x^4,xyz^2);$, and its linkage
$I_3=(x^4,y^4,z^4):J_3$. The left picture is $J_3$ and the right picture is $I_3$. 
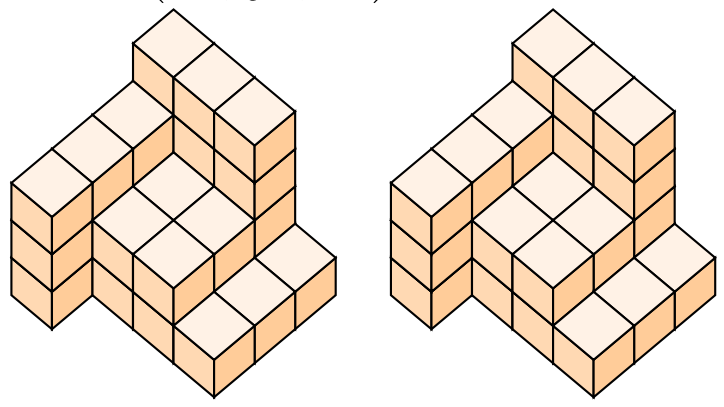
\monomialIdealPicture{{4,3,3,3},{4,2,2},{4,2,2},{1,1,1}} \monomialIdealPicture{{4,3,3,3},{4,2,2},{4,2,2},{1,1,1}} 

For $J_3$, $8 \leq  (3 - 1) \cdot 4$, and for $I_3$, $8 \leq (3- 1)\cdot 4$. So any ideal that is linked to both $I_3$ and $J_3$ must
have at least more boxes than $J_3$.
\end{enumerate}

\end{example}

\end{document}